\magnification=\magstep1
\input amstex
\documentstyle{amsppt}
\NoBlackBoxes
\vsize=7in
\hsize=5.5in
\def\bt{\mathop\boxtimes\limits}
\catcode`\@=11
\loadmathfont{rsfs}
\def\mycal{\mathfont@\rsfs}
\csname rsfs \endcsname
\topmatter
\title ON A CLASS OF TYPE II$_1$
FACTORS WITH BETTI NUMBERS INVARIANTS\endtitle
\vskip -.3in
\author SORIN POPA \endauthor
\rightheadtext{Betti numbers invariants}
\vskip -.3in
\affil University of California, Los Angeles\endaffil
\vskip -.3in
\address Math.Dept., UCLA, LA, CA 90095-155505\endaddress
\email popa@math.ucla.edu\endemail
\thanks Supported in part by a NSF Grant 0100883.\endthanks
\vskip -.3in
\abstract We prove that
a type II$_1$ factor $M$ can have at most
one Cartan subalgebra $A$ satisfying a combination of
rigidity and compact approximation properties. We use this result to
show that within the class $\Cal H \Cal T$ of
factors $M$ having such Cartan subalgebras $A \subset M$,
the Betti numbers of the standard equivalence
relation associated with $A \subset M$ ([G2]),
are in fact isomorphism
invariants for the factors $M$, $\beta^{^{HT}}_n(M),
n\geq 0$. The class $\Cal H\Cal T$ is closed
under amplifications and tensor
products, with the Betti numbers satisfying
$\beta^{^{HT}}_n(M^t)=
\beta^{^{HT}}_n(M)/t,
\forall t>0$, and a K{\"u}nneth type formula.
An example of a factor in the class $\Cal H\Cal T$ is given by
the group von Neumann factor
$M=L(\Bbb Z^2 \rtimes SL(2, \Bbb Z))$, for which
$\beta^{^{HT}}_1(M) = \beta_1(SL(2, \Bbb Z)) = 1/12$. Thus,
$M^t \not\simeq M, \forall t \neq 1$,
showing that the fundamental
group of $M$ is trivial.  This solves a long
standing problem of R.V. Kadison. Also, our results
bring some insight into a recent problem of A. Connes and answer a number
of open questions on von Neumann algebras.
\endabstract
\endtopmatter
\document
\vskip -.2in
\centerline{\bf Content}

\settabs\+\indent\indent & 10. \ & 1.1 \quad & More on rigidity
\cr \tenpoint{\+& 0. & Introduction \cr \+& 1. & Preliminaries \cr
\+& & 1.1. & Pointed correspondences \cr \+& & 1.2. & Completely
positive maps as Hilbert space operators \cr \+& & 1.3. & The
basic construction and its compact ideal space \cr \+& & 1.4. &
Discrete inclusions and bimodule decomposition \cr \+& 2. &
Relative Property H:  Definition and Examples \cr \+& 3. & More on
property H \cr \+& 4. & Rigid embeddings:  Definitions and
Properties \cr \+& 5. & More on rigid embeddings \cr \+& 6. & HT
subalgebras and the class $\Cal H \Cal T$ \cr \+& 7. & Subfactors
of an $\Cal H \Cal T$ factor \cr \+& 8. & Betti numbers for $\Cal
H \Cal T$ factors \cr \+& Appendix:  Some conjugacy results \cr}

\newpage
\heading 0. Introduction. \endheading

We consider in this paper the class of type II$_1$ factors
with maximal abelian $^*$-subalgebras
satisfying both a weak rigidity property,
in the spirit of Kazhdan
and Connes-Jones ([Kaz], [CJ]),
and a weak amenability property,
in the spirit of Haagerup's compact approximation property
([H]). Our main result shows that
a type II$_1$ factor $M$ can have at most
one such maximal abelian $^*$-subalgebra $A \subset M$,
up to unitary conjugacy.
Moreover, we prove that if $A\subset M$ satisfies these conditions then $A$ is
automatically a Cartan subalgebra of $M$, i.e., the normalizer of $A$ in $N$,
$\Cal N(A)=\{u \in M\mid uu^*=1, uAu^* = A\}$,
generates all the von Neumann algebra $M$.
In particular, $\Cal N(A)$
implements an ergodic measure preserving equivalence relation
on the standard probability space
$(X, \mu)$, with $A=L^\infty(X, \mu)$ ([FM]),
which up to orbit equivalence
only depends on the isomorphism class of $M$.

We call HT the Cartan subalgebras satisfying the combination
of the rigidity and compact approximation
properties and denote by $\Cal H\Cal T$ the class
of factors having HT Cartan subalgebras. Thus, our theorem
implies that if $M\in \Cal H \Cal T$,
then there exists a unique (up to isomorphism)
ergodic measure-preserving
equivalence relation $\Cal R^{^{HT}}_M$ on
$(X, \mu)$ associated with it, implemented by
the HT Cartan subalgebra of $M$. In particular, any
invariant for $\Cal R^{^{HT}}_M$
is an invariant for $M\in \Cal H\Cal T$.

In a recent paper ([G2]), D. Gaboriau introduced a notion
of $\ell^2$-Betti numbers for arbitrary countable
measure preserving equivalence relations
$\Cal R$, $\{\beta_n(\Cal R)\}_{n\geq 0}$,
starting from ideas of Atiyah ([A]) and Connes
([C4]), and generalizing the
notion of $L^2$-Betti numbers for measurable
foliations defined in ([C4]). His notion also generalizes the
$\ell^2$-Betti numbers for discrete
groups $\Gamma_0$ of Cheeger-Gromov ([ChGr]),
$\{\beta_n(\Gamma_0)\}_{n \geq 0}$, as Gaboriau shows
that $\beta_n(\Gamma_0) = \beta_n(\Cal R_{\Gamma_0})$, for any
countable equivalence relation
$\Cal R_{\Gamma_0}$ implemented by a free, ergodic, measure-preserving
action of the group $\Gamma_0$ on
a standard probability space $(X, \mu)$ ([G2]).

We define in this paper the Betti numbers
$\{\beta^{^{HT}}_n(M)\}_{n\geq 0}$
of a factor $M$ in the class $\Cal H\Cal T$ as the
the $\ell^2$-Betti numbers ([G2]) of the corresponding
equivalence relation $\Cal R^{^{HT}}_M$,
$\{\beta_n({\Cal R^{^{HT}}_M})\}_n$.

Due to the uniqueness of the HT Cartan subalgebra,
the general properties of the Betti numbers for
countable equivalence relations proved in ([G2]) entail
similar properties for the Betti numbers of the factors in the class
$\Cal H\Cal T$. For instance, after proving that
$\Cal H\Cal T$ is closed
under amplifications by arbitrary $t > 0$, we use
the formula $\beta_n(\Cal R^t) = \beta_n(\Cal R)/t$
in ([G2]) to deduce that $\beta^{^{HT}}_n(M^t)
= \beta^{^{HT}}_n(M)/t, \forall n$.
Also, we prove that $\Cal H\Cal T$
is closed under tensor products and that
a K{\"u}nneth type formula holds for
$\beta^{^{HT}}_n(M_1\overline\otimes M_2)$
in terms of the Betti numbers for $M_1, M_2 \in \Cal H\Cal T$,
as a consequence of the similar formula
for groups and equivalence relations
([B], [ChGr], [L\"u], [G2]).

Our main example of a factor in the class
$\Cal H \Cal T$ is the group von Neumann algebra $L(G_0)$
associated with $G_0=\Bbb Z^2 \rtimes SL(2, \Bbb Z)$,
regarded as the group-measure space construction
$L^\infty(\Bbb T^2, \mu)=A_0 \subset A_0
\rtimes_{\sigma_0} SL(2, \Bbb Z)$, where
$\Bbb T^2$ is regarded as the dual of $\Bbb Z^2$
and $\sigma_0$ is the action
implemented by $SL(2, \Bbb Z)$ on it. More
generally, since our HT condition
on the Cartan subalgebra $A$
requires only part of $A$ to be rigid in $M$, we show that
any cross product factor of the form
$A \rtimes_\sigma SL(2, \Bbb Z)$, with
$A=A_0\overline{\otimes} A_1$,
$\sigma=\sigma_0 \otimes \sigma_1$ and
$\sigma_1$ an arbitrary ergodic action of $SL(2, \Bbb Z)$ on
an abelian algebra $A_1$,
is in the class $\Cal H\Cal T$.
By a recent result in ([Hj]), based on
the notion and results on tree-ability in ([G1]),
all these factors are in fact
amplifications of group-measure space factors of the form
$L^\infty(X, \mu) \rtimes \Bbb F_n$, where
$\Bbb F_n$ is the free group on $n$
generators, $n=2, 3, ...$.

To prove that $M$ belongs to the class $\Cal H \Cal T$,
with $A$ its corresponding HT Cartan subalgebra, we use Kazhdan's
rigidity of the inclusion $\Bbb Z^2 \subset \Bbb Z^2 \rtimes SL(2, \Bbb Z)$
and Haagerup's compact approximation property of $SL(2, \Bbb Z)$.
The same arguments are actually used to show  that if $\alpha \in \Bbb C, |\alpha|=1,$ and $L_\alpha(\Bbb Z^2)$ denotes the corresponding
``twisted'' group algebra (or ``quantized'' 2-dimensional thorus),
then $M_\alpha=L_\alpha(\Bbb Z^2) \rtimes SL(2, \Bbb Z)$
is in the class $\Cal H\Cal T$ if and only if $\alpha$ is a root of unity.

Since the orbit equivalence relation $\Cal R^{^{^{HT}}}_M$ implemented by
$SL(2, \Bbb Z)$ on $A$ has exactly one non-zero Betti number, namely
$\beta_1(\Cal R^{^{^{HT}}}_M) = \beta_1(SL(2, \Bbb Z)) = 1/12$ (
[B], [ChGr], [G2]), it follows that the factors $M
=A\rtimes_\sigma SL(2, \Bbb Z)$ satisfy $\beta^{^{HT}}_1(M) =
1/12$ and $\beta^{^{HT}}_n(M) = 0, \forall n \neq 1$.
More generally, if $\alpha$ is a $n$'th primitive root
of 1, then the factors $M_\alpha=
L_\alpha(\Bbb Z^2) \rtimes SL(2, \Bbb Z)$ satisfy $\beta^{^{HT}}_1(M_\alpha)=
n/12, \beta^{^{HT}}_k(M_\alpha)=0, \forall k\neq 1$. We deduce from this
that if $\alpha, \alpha'$ are primitive roots of unity of
order $n$ respectively $n'$ then $M_\alpha \simeq M_{\alpha'}$ if
and only if $n=n'$.

Other examples of factors in the class
$\Cal H\Cal T$ are obtained by taking discrete groups $\Gamma_0$
that can be embedded as
arithmetic lattices in $SU(n,1)$ or $SO(m,1)$, together
with suitable actions $\sigma$ of $\Gamma_0$ on abelian
von Neumann algebras $A\simeq L(\Bbb Z^N)$. Indeed, these groups $\Gamma_0$
have the Haagerup approximation property by ([dCaH], [CowH])
and their action $\sigma$ on $A$ can be taken to be rigid by a
recent result of
Valette ([Va]). In each of these cases, the Betti numbers
have been calculated in
([B]). A yet another
example is offered by the action of $SL(2, \Bbb Q)$ on $\Bbb Q^2$:
Indeed, the rigidity of the action of $SL(2, \Bbb Z)$ (regarded
as a subgroup of $SL(2,
\Bbb Q)$) on $\Bbb Z^2$ (regarded as a subgroup
of $\Bbb Q^2$), as well as the property H of $SL(2, \Bbb Q)$
proved in ([CCJJV]), are enough to insure that $L(\Bbb Q^2 \rtimes
SL(2, \Bbb Q))$ is in the class $\Cal H\Cal T$.

As a consequence of these considerations, we are able
to answer a number
of open questions
in the theory of type II$_1$ factors. Thus, the factors
$M=A\rtimes_\sigma SL(2, \Bbb Z)$ (more generally,
$A\rtimes_\sigma \Gamma_0$ with $\Gamma_0, \sigma$ as above)
provide the first class of type II$_1$ factors
with trivial fundamental group, i.e. $\mycal F(M)
\overset \text{\rm def}
\to = \{ t> 0\mid M^t \simeq M\} = \{1\}$. Indeed,
we mentioned that $\beta^{^{HT}}_n(M^t)
= \beta^{^{HT}}_n(M)/t , \forall n$,
so that if $\beta^{^{HT}}_n(M)\neq 0, \infty$
for some $n$ then
$\mycal F(M)$ is forced to be equal to $\{1\}$.

In particular, the
factors $M$ are not
isomorphic to the algebra of $n$ by $n$ matrices over $M$, for any
$n \geq 2$, thus providing
an answer to Kadison's
Problem 3 in ([K1]) (see also Sakai's Problem
4.4.38 in [S]).
Also, through
appropriate
choice of actions of the form
$\sigma=\sigma_0\otimes \sigma_1$, we obtain factors of the form
$M=A\rtimes_\sigma SL(2, \Bbb Z)$
having the property $\Gamma$ of Murray and von Neumann,
yet trivial fundamental group.

The fundamental group $\mycal F(M)$ of a II$_1$ factor $M$ was
defined by Murray and von Neumann in
the early 40's, in connection with their notion of
continuous dimension. They noticed
that $\mycal F(M) = \Bbb R^*_+$ when
$M$ is isomorphic to the hyperfinite type II$_1$ factor $R$,
and more generally when $M$ ``splits off'' $R$.

The first examples of type II$_1$ factors $M$ with
$\mycal F(M) \neq \Bbb R_+^*$, and the
first occurence of rigidity in the von Neumann
algebra context, were
discovered by Connes in ([C1]). He proved that if
$G_0$ is an infinite conjugacy class discrete group
with the property (T) of Kazhdan then its
group von Neumann algebra $M=L(G_0)$ is a type II$_1$ factor with
countable fundamental group.
It was then proved
in ([Po1]) that this is still the case for
factors $M$ which contain some irreducible copy of such $L(G_0)$.
It was also shown that there exist type II$_1$ factors
$M$ with $\mycal F(M)$ countable and containing any presribed
countable set of numbers ([GoNe], [Po4]).
However, the fundamental group $\mycal F(M)$ could
never be computed exactly, in any of these
examples.

In fact, more than proving that $\mycal F(M)=\{1\}$
for $M=A\rtimes_\sigma SL(2, \Bbb Z)$,
the calculation of the
Betti numbers shows that
$M^{t_1} \overline\otimes M^{t_2} ...
\overline\otimes M^{t_n}$ is isomorphic to
$M^{s_1} \overline\otimes M^{s_2} ... \overline\otimes M^{s_m}$
if and only if $n=m$ and $t_1 t_2 ... t_n = s_1 s_2 ... s_m$.
In particular,
all tensor powers of $M$, $M^{\overline\otimes n},
n=1, 2, 3, ...$, are mutually non-isomorphic and have
trivial trivial fundamental group.
(N.B. The first examples of factors
having non-isomorphic tensor powers were constructed in [C4];
another class of examples was obtained in [CowH]). In fact,
since $\beta^{^{HT}}_k(M^{\overline\otimes n})
\neq 0$ iff $k=n$, the factors
$\{M^{\overline\otimes n}\}_{n\geq 1}$ are not even stably isomorphic.

In particular, since $M^t \simeq L^\infty(X, \mu)\rtimes \Bbb F_n$ for
$t= (12(n-1))^{-1}$ (cf [Hj]), it follows that for each $n\geq 2$
there exist a free ergodic action $\sigma_n$ of $\Bbb F_n$
on the standard probability space
$(X, \mu)$ such that the factors $M_n = L^\infty(X,
\mu)\rtimes_{\sigma_n}\Bbb F_n, n= 2, 3, ...,$
satisfy $M_{k_1} \overline\otimes ... \overline\otimes M_{k_p}
\simeq M_{l_1} \overline\otimes ... \overline\otimes M_{l_r}$ if and
only if $p=r$ and $k_1k_2 ... k_p = l_1 l_2 ... l_r$. Also,
since $\beta^{^{HT}}_1(M_n) \neq 0$,
the K{\"u}nneth formula shows that the factors $M_n$ are prime within the class
of type II$_1$ factors in $\Cal H\Cal T$.

Besides being closed under tensor products and amplifications,
the class $\Cal H\Cal T$ is closed under finite
index extensions/restrictions, i.e.,
if $N \subset M$ are type II$_1$ factors with finite Jones index,
$[M:N] < \infty$,
then $M \in \Cal H\Cal T$ if and only if $N \in \Cal H\Cal T$. In fact,
factors in the class $\Cal H\Cal T$ have a remarkably
rigid ``subfactor picture'':

Thus, if $M \in \Cal H\Cal T$ and $N \subset M$ is an irreducible subfactor
with $[M:N] < \infty$ then $[M:N]$ is an integer.
More than that, the graph of $N\subset M$, $\Gamma=\Gamma_{N,M}$,
has only integer weights $\{v_k\}_k$. Recall that
the weights $v_k$ of the graph of a subfactor $N \subset M$
are given by the ``statistical dimensions''
of the irreducible $M$-bimodules $\Cal H_k$ in the Jones tower
or, equivalently,
as the square roots of the
indices of the corresponding irreducible inclusions of
factors, $M \subset M(\Cal H_k)$. They
give a Perron-Frobenius type eigenvector for $\Gamma$,
satisfying $\Gamma \Gamma^t \vec{v} = [M:N] \vec{v}$.
We prove that if $\beta^{^{HT}}_n(M) \neq 0, \infty$ then
$v_k = \beta_n^{^{HT}}(M(\Cal H_k))/\beta_n^{^{HT}}(M), \forall k$,
i.e., the statistical dimensions are proportional to the Betti numbers.
As an application of this subfactor analysis, we show that the
non-$\Gamma$ factor $L(\Bbb Z^2
\rtimes SL(2, \Bbb Z))$ has two non-conjugate period two automorphims.

We also discuss invariants that can distinguish
between factors in the class $\Cal H\Cal T$
which have the same Betti numbers. Thus, we show that
if $\Gamma_0 = SL(2, \Bbb Z), \Bbb F_n$, or if $\Gamma_0$
is an arithemtic lattice in some $SU(n,1), SO(n,1)$,
for some $n \geq 2$, then there
exist three non-orbit equivalent free ergodic measure preserving
actions $\sigma_i$ of $\Gamma_0$ on $(X, \mu)$, with $M_i = L^{\infty}(X, \mu)
\rtimes_{\sigma_i} \Gamma_0 \in \Cal H\Cal T$ non-isomorphic
for $i=1,2,3.$. Also, we
apply Gaboriau's notion
of approximate dimension to  equivalence relations of the form
$\Cal R^{^{HT}}_M$ to
distinguish between $\Cal H\Cal T$ factors of the form $M_k=\Bbb F_{n_1} \times ...
\times \Bbb F_{n_k} \times S_\infty$, with $S_\infty$ the infinite symmetric
group and $k=1,2,...$, which all have only $0$ Betti numbers.

As for the ``size'' of the class $\Cal H\Cal T$, note that we could
only produce examples of factors
$M=A\rtimes_\sigma \Gamma_0$ in $\Cal H \Cal T$
for certain property H groups $\Gamma_0$, and for certain special
actions $\sigma$
of such groups. We call H$_{_{T}}$ the groups $\Gamma_0$ for which
there exist free ergodic measure preserving actions
$\sigma$ on the standard probability space $(X, \mu)$ such that
$L^\infty(X, \mu) \rtimes_\sigma \Gamma_0 \in \Cal H \Cal T$.
Besides the examples $\Gamma_0
=SL(2, \Bbb Z), SL(2, \Bbb Q), \Bbb F_n$, or $\Gamma_0$
an arithmetic lattice in $SU(n,1), SO(n,1), n \geq 2$, mentioned
above, we show that the class of H$_{_{T}}$
groups is closed under
products by arbitrary property H groups, cross product
by amenable groups and finite index restriction/extension.

On the other hand, we prove that the class
$\Cal H\Cal T$ does not contain factors of the form $M \simeq
M\overline{\otimes} R$, where $R$ is the hyperfinite
II$_1$ factor. In particular, $R \notin \Cal H\Cal T$.
Also, we prove that the factors
$M\in \Cal H\Cal T$ cannot contain
property (T) factors and cannot be
embedded into free group factors (by using arguments
similar to [CJ]). In the same vein, we show that
if $\alpha\in \Bbb T$
is not a root of unity, then
the factors $M_\alpha = L_\alpha(\Bbb Z^2) \rtimes SL(2, \Bbb Z)
=R \rtimes SL(2, \Bbb Z)$ cannot be embedded into any factor
in the class $\Cal H\Cal T$. In fact, such factors
$M_\alpha$ belong to a special
class in its own, that we will study in a forthcoming
paper.

Besides these concrete applications,
our results give a partial answer to a
challanging problem recently raised by Alain Connes,
on defining a notion of Betti numbers $\beta_n(M)$
for type II$_1$ factors $M$, from similar conceptual
grounds as in the case of measure preserving equivalence relations
in ([G2]) (simplicial structure, $\ell^2$ homology/cohomology, etc),
a notion that should satisfy $\beta_n(L(G_0)) = \beta_n(G_0)$ for
group von Neumann factors $L(G_0)$.
In this respect, note that our definition is not the result of a
``conceptual approach'', relying instead on
the uniqueness result
for the HT Cartan subalgebras, which allows reducing the problem
to Gaboriau's work on invariants for
equivalence relations and, through it, to the
results on $\ell^2$-cohomology
for groups in ([ChGr], [B], [L\"u]). Thus, although they are
invariants  for ``global factors'' $M \in
\Cal H\Cal T$, the Betti numbers $\beta_n^{^{HT}}(M)$
are ``relative'' in spirit, a fact that we have indicated by adding
the upper index $^{^{HT}}$. Also, rather than satisfying
$\beta_n(L(G_0)) =  \beta_n(G_0)$, the invariants $\beta_n^{^{HT}}$
satisfy $\beta_n^{^{HT}}(A\rtimes \Gamma_0)
= \beta_n(\Gamma_0)$. In fact, note that
if $A\rtimes \Gamma_0 = L(G_0)$, where $G_0=\Bbb Z^N \rtimes \Gamma_0$, then
$\beta_n(G_0) = 0$,
while $\beta_n^{^{HT}}(L(G_0)) =
\beta_n(\Gamma_0)$ may be different from $0$.

The paper is organized as follows: Section 1 consists
of preliminaries: we first establish some
basic properties of {\it Hilbert bimodules} over
von Neumann algebras and of their
associated {\it completely positive maps}; then
we recall the {\it basic construction} of an
inclusion of finite von Neumann algebras
and study its {\it compact ideal space};
we also recall the definitions of {\it normalizer and quasi-normalizer}
of a subalgebra, as well as the notions of {\it regular, quasi-regular,
dicrete and Cartan subalgebras}, and discuss some
of the results in ([FM]) and ([PoSh]). In Section 2 we consider a relative version of 
Haagerup's compact approximation property
for inclusions of von Neumann algebras, called
{\it relative property H} (cf. also [Bo]), and prove its main
properties. In Section 3 we give examples of
property H inclusions and use ([PoSh]) to
show that if a type II$_1$ factor $M$ has the property H relative
to a maximal abelian subalgebra $A \subset M$ then $A$
is a Cartan subalgebra of $M$. In Section 4 we
define a notion of rigidity (or relative property (T)) 
for inclusions of algebras and investigate 
its basic properties. In Section 5
we give examples of rigid inclusions and relate this property
with the co-rigidity property defined in ([Zi], [A-De], [Po1]). We
also introduce a new notion of property (T) for equivalence relations, called {\it
relative property (T)},
by requiring the associated Cartan subalgebra inclusion to be  rigid.

In Section 6 we define the class $\Cal H\Cal T$ of factors $M$ having
{\it HT Cartan subalgebras} $A\subset M$,
i.e., maximal abelian
$^*$-subalgebras $A\subset M$ such that $M$ has the
property H relative to $A$ and
$A$ contains a subalgebra $A_0 \subset A$
with $A_0'\cap M = A$ and $A_0 \subset M$ rigid. We then prove the main
technical result of the paper, showing that HT Cartan subalgebras
are unique. We show the stability of the class $\Cal H\Cal T$ with respect
to various operations (amplification, tensor product), and prove its
rigidity to perturbations. Section 7
studies the lattice of subfactors of $\Cal H\Cal T$ factors: we
prove the stability of the class $\Cal H\Cal T$ to finite index,
obtain a canonical decomposition for subfactors in $\Cal H\Cal T$ and
prove that the index is always an integer. In Section 8
we define the {\it Betti numbers} $\{\beta^{^{HT}}_n(M)\}_n$
for $M\in \Cal H\Cal T$ and
use the previous sections and
([G2]) to deduce various properties of this invariant.
We also discuss some alternative invariants for factors
$M \in \Cal H \Cal T$, such as the {\it outomorphism group}
Out$_{_{HT}}(M)
\overset\text{\rm def} \to
=$Aut$(\Cal R^{^{HT}}_M)/$Int$(\Cal R^{^{HT}}_M)$,
which we prove is discrete countable, or
$ad_{_{HT}}(M)$, defined to be Gaboriau's
{\it approximate dimension} ([G3]) of $\Cal R^{^{HT}}_M$. We end
with applications, as well as some remarks and open questions.
We have included an Appendix in which we prove
some key technical results on
unitary conjugacy of von Neumann subalgebras in type
II$_1$ factors. The proof uses techniques from ([Chr], [Po2,3,6], [K2]).

{\it Acknowledgement}. I want to thank U.Haagerup, V. Lafforgue
and A. Valette for useful conversations on the
property H and (T) for groups. My special
thanks are due Damien Gaboriau, for keeping
me informed on his beautiful recent
results and for useful comments on the first version of this paper. I am
particularly greatful to Alain Connes
and Dima Shlyakhtenko
for many fruitful conversations and constant
support. I want to express my gratitude
to MSRI and the organizers of the
Operator Algebra year 2000-2001, for their hospitality
and for a most stimulating atmosphere. This article
is an expanded version
of a paper with the same title which appeared as
MSRI preprint 2001/0024.

\heading 1. Preliminaries. \endheading

\noindent
{\bf 1.1. Pointed correspondences}. By using
the GNS construction as a link, a representation
of a group $G_0$ can be viewed in two equivalent ways:
as a group morphism from $G_0$ into the unitary group of a Hilbert space
$\Cal U(\Cal H)$, or as a positive definite function on $G_0$.

The discovery of the appropriate notion of representations for von Neumann
algebras, as so-called {\it correspondences}, is due to
Connes ([C3,7]).
In the vein of group representations,
Connes introduced correspondences in two alternative ways,
both of which use the
idea of ``doubling'' - a genuine conceptual breakthrough.
Thus, correspondences of von Neumann algebras $N$
can be viewed as {\it Hilbert}
$N$-{\it bimodules} $\Cal H$, the
quantized version of group
morphisms into $\Cal U(\Cal H)$;
or as {\it completely positive maps} $\phi:N \rightarrow N$,
the quantized version of
positive definite functions on groups (cf. [C3,7] and [CJ]).
The equivalence of
these two points of view is again realized via a version of
the GNS construction ([CJ], [C7]).

We will in fact need ``pointed''
versions of Connes's correspondences, adapted to the case
of inclusions $B \subset N$, as introduced
in ([Po1]) and ([Po5]). In this Section we detail the two alternative
ways of viewing such pointed correspondences,
in the same spirit as ([C7]): as ``$B$-pointed bimodules''
or as ``$B$-bimodular completely positive maps''. This is a very
important idea, which will be present throughout this paper.
\vskip .1in
\noindent
{\it 1.1.1. Pointed Hilbert bimodules}.
Let $N$ be a finite von Neumann algebra with a
fixed normal faithful tracial state $\tau$ and
$B\subset N$ a von Neumann subalgebra of $N$. A
{\it Hilbert} $(B \subset N)$-{\it bimodule} $(\Cal H, \xi)$
is a Hilbert $N$-bimodule with a fixed unit vector $\xi \in \Cal H$
satisfying $b\xi = \xi b, \forall b\in B$. When $B=\Bbb C$,
we simply call $(\Cal H, \xi)$ a {\it pointed Hilbert} $N$-{\it bimodule}.

If $\Cal H$ is a Hilbert
$N$-{\it bimodule} then $\xi \in \Cal H$ is a {\it cyclic vector}
if $\overline{{\text{\rm sp}}} N\xi N = \Cal H$.

To relate Hilbert $(B \subset N)$-bimodules and $B$-bimodular
completely positive maps on $N$ one uses a generalized version
of the GNS construction, due to Stinespring, which we describe below:
\vskip .1in
\noindent
{\it 1.1.2. From completely positive maps to Hilbert bimodules.}  Let
$\phi$ be a normal completely positive map on $N$, normalized so
that $\tau(\phi(1)) = 1$. We associate to it
the pointed Hilbert $N$-bimodule $(\Cal H_\phi, \xi_\phi)$ in the
following way:

Define on the linear space $\Cal H_0 = N \otimes N$ the sesquilinear form
$\langle x_1 \otimes y_1, x_2 \otimes y_2 \rangle_\phi =
\tau(\phi(x_2^*x_1)y_1y_2^*),
x_{1,2}, y_{1,2} \in N$. The complete positivity of
$\phi$ is easily seen to be
equivalent to the positivity of $\langle \cdot, \cdot \rangle_\phi$.
Let $\Cal H_\phi$
be the completion of $\Cal H_0/\sim$, where $\sim$ is the
equivalence modulo the null
space of $\langle \cdot, \cdot \rangle_\phi$ in $\Cal H_0$. Also, let
$\xi_\phi$ be the class of $1\otimes 1$ in $\Cal H_\phi$. Note that
$\|\xi_\phi\|^2 = \tau(\phi(1)) = 1$.

If $p = \Sigma_i x_i \otimes y_i \in \Cal H_0$, then by using again the complete
positivity of $\phi$ it follows that $N \ni x \rightarrow \Sigma_{i,j}
\tau(\phi(x_j^*xx_i)y_iy_j^*)$ is a positive normal functional
on $N$ of norm $\langle p,p \rangle_\phi$. Similarily,
$N \ni y \rightarrow \Sigma_{i,j} \tau(\phi(x_j^*x_i)y_iyy_j^*)$ is a
positive normal functional on $N$ of norm $\langle p,p \rangle_\phi$. Note that the latter
can alternatively be viewed as a functional on the
opposite algebra $N^{\text{\rm op}}$ (which is the
same as $N$ as a vector space but has
multiplication inverted, $x\cdot y = yx$). Moreover, $N$ acts on $\Cal H_0$
on left and right by $xpy= x (\Sigma_i x_i \otimes y_i)y= \Sigma_ixx_i\otimes y_iy$.
These two actions clearly commute and the complete positivity
of $\phi$ entails:
$$
\langle xp, xp\rangle_\phi = \langle x^*xp,p\rangle_\phi
\leq \|x^*x\| \langle p,p \rangle_\phi = \|x\|^2 \langle p,p\rangle_\phi
$$
Similarly
$$
\langle py, py \rangle_\phi \leq \|y\|^2 \langle p,p \rangle_\phi.
$$
Thus, the above left and right actions of $N$ on $\Cal H_0$ pass to
$\Cal H_0/\sim$ and then extend to commuting left-right actions
on $\Cal H_\phi$. By the normality of the
forms $x\rightarrow \langle xp, p \rangle_\phi$
and $y \rightarrow \langle py,p \rangle_\phi$,
these left-right actions of $N$ on $\Cal H_\phi$ are normal (i.e., weakly
continuous).

This shows that $(\Cal H_\phi, \xi_\phi)$
with the above $N$-bimodule structure is a pointed Hilbert $N$-bimodule,
which in addition is clearly cyclic.
Moreover, if $B \subset N$ is a von Neumann
subalgebra and the completely positive map $\phi$ is $B$-bimodular,
then it is immediate to check that $b\xi_\phi = \xi_\phi b, \forall b \in B$.
Thus, if $\phi$ is $B$-bimodular, then $(\Cal H_\phi, \xi_\phi)$
is a Hilbert $(B \subset N)$-bimodule.

Let us end this paragraph with some useful inequalities which show that
elements that are almost fixed by a $B$-bimodular
completely positive map $\phi$ on $N$
are almost commuting with the associated vector $\xi_\phi\in \Cal H_\phi$:

\proclaim{Lemma} $1^\circ$. $\| \phi (x) \|_2 \leq
\|\phi(1)\|_2, \forall x\in N, \|x\|\leq 1$.
\vskip .1in
$2^\circ$. If $a=1 \vee \phi(1)$ and $\phi'(\cdot) = a^{-1/2}
\phi(\cdot) a^{-1/2}$, then $\phi'$ is completely positive, $B$-bimodular
and satisfies $\phi'(1) \leq 1$,
$\tau \circ \phi' \leq \tau \circ \phi$ and the estimate:
$$
\|\phi'(x)-x\|_2 \leq \|\phi(x) - x\|_2 + 2\|\phi(1) -1 \|_1^{1/2}\|x\|,
\forall x\in N.
$$

$3^\circ$. Assume $\phi(1) \leq 1$
and define $\phi''(x)=\phi(b^{-1/2}xb^{-1/2})$,
where $b=1 \vee ({\text{\rm d}}\tau\circ \phi
/{\text{\rm d}}\tau) \in L^1(N,\tau)_+$.
Then
$\phi''$ is completely positive, $B$-bimodular and
satisfies $\phi''(1)
\leq \phi(1) \leq 1, \tau\circ \phi'' \leq \tau$, as well as the estimate:
$$
\|\phi''(x) - x\|^2_2 \leq
2\|\phi(x)-x\|_2 + 5\|b-1\|_1^{1/2}, \forall x\in N, \|x\|\leq 1.
$$

$4^\circ$. $\|x\xi_\phi - \xi_\phi x \|_2^2
\leq 2\|\phi(x) - x\|_2^2 + 2 \|\phi(1)\|_2 \|\phi(x)-x\|_2, \forall x\in N,
\|x\|\leq 1.$
\endproclaim
\vskip .1in \noindent {\it Proof}. $1^\circ$. Since any $x\in N$
with $\|x\| \leq 1$ is a convex combination of two unitary
elements, it is sufficient to prove the inequality for unitary
elements $u\in N$. By continuity, it is in fact sufficient to
prove it in the case the unitary elements $u$ have finite
spectrum. If $u = \Sigma_i \lambda_i p_i$ for some scalars
$\lambda_i$ with $|\lambda_i| = 1$, $1 \leq i \leq n$, and some
partition of the identity with projections $p_i\in N$, then
$\tau(\phi(p_i)\phi(p_j)) \geq 0, \forall i,j$. Taking this into
account, we get:
$$
\tau(\phi(u)\phi(u^*)) = \Sigma_{i,j} \lambda_i
\overline{\lambda_j} \tau(\phi(p_i)\phi(p_j)) \leq \Sigma_{i,j}
|\lambda_i \overline{\lambda_j}| \tau(\phi(p_i)\phi(p_j))
$$
$$
= \Sigma_{i,j} \tau(\phi(p_i)\phi(p_j)) = \tau(\phi(1)\phi(1)).
$$

$2^\circ$. Since $a \in B'\cap N$, $\phi'$ is
$B$-bimodular. We clearly have $\phi'(1)=a^{-1/2}\phi(1)a^{-1/2} \leq 1$.
Since $a^{-1} \leq 1$, for $x\geq 0$  we get $\tau(\phi'(x))
= \tau(\phi(x) a^{-1}) \leq \tau(\phi(x))$. Also, we have:
$$
\|\phi'(x) - x \|_2\leq \|a^{-1/2}\phi(x)a^{-1/2} - a^{-1/2}xa^{-1/2}\|_2 +
\|a^{-1/2}xa^{-1/2}-x\|_2
$$
$$
\leq \|\phi(x)-x\|_2 + 2\|a^{-1/2} -1\|_2 \|x\|.
$$
But
$$
\|a^{-1/2}-1\|_2 \leq \|a^{-1} -1\|^{1/2}_1=\|a^{-1}-aa^{-1}\|_1
$$
$$
\leq \|a-1\|_1\|a^{-1}\| \leq \|a-1\|_1 \leq \|\phi(1)-1\|_1.
$$
Thus,
$$
\|\phi'(x) - x \|_2\leq \|\phi(x)-x\|_2 + 2\|\phi(1)-1\|_1^{1/2}\|x\|.
$$

$3^\circ$. The first properties are clear by the definitions.
Then note that $\|y\|_2^2 \leq \|y\|\|y\|_1$ and
$\|\phi''(y)\|_1 \leq      \|y\|_1$.
(Indeed, because if ${\phi''}^*$ is as defined
in Lemma 2.1,  then for $z\in N$ with $\|z\| \leq 1$
we have $\|{\phi''}^*(z)\| \leq 1$
so that $\|\phi''(y)\|_1={\text{\rm sup}} \{|\tau(\phi''(y)z)|\mid
z\in N, \|z\|\leq 1\}={\text{\rm sup}} \{|\tau(y{\phi''}^*(z))|\mid
z\in N, \|z\|\leq 1\} \leq {\text{\rm sup}} \{|\tau(yz))|\mid
z\in N, \|z\|\leq 1\}=\|y\|_1$.) Note also that
$\tau(b) \leq 1 + \tau(\phi(1)) \leq 2$. Thus, for
$x\in N, \|x\| \leq 1$, we get:
$$
\|\phi''(x)-x\|^2_2 \leq 2 \|\phi''(x)-x\|_1
$$
$$
\leq  2 \|\phi''(x)-\phi''(b^{1/2}xb^{1/2})\|_1 +  2 \|\phi(x)-x\|_1
$$
$$
\leq 2 \|x-b^{1/2}xb^{1/2}\|_1 + 2 \|\phi(x)-x\|_1.
$$
$$
\leq 2\|x-xb^{1/2}\|_1+2\|xb^{1/2}-b^{1/2}xb^{1/2}\|_1 + 2 \|\phi(x)-x\|_1.
$$

But $\|x\|_2 \leq 1$ and $\|xb^{1/2}\|_2^2 \leq \tau(b)\leq 2$,
so by the Cauchy-Schwartz inequality the above is majorized by:
$$
2\|x\|_2\|1-b^{1/2}\|_2 + 2\|1-b^{1/2}\|_2\|xb^{1/2}\|_2
+2\|\phi(x)-x\|_2
$$
$$
\leq (2+ 2^{3/2}) \|b^{1/2}-1\|_2 + 2 \|\phi(x)-x\|_2 \leq 5\|b-1\|_1^{1/2} +
2\|\phi(x)-x\|_2.
$$
$4^\circ$. Since by the Cauchy-Schwartz inequality we have
$$
\pm {\text{\rm Re}} \tau (\phi(x)(\phi(x)^* - x^*))
\leq \|\phi(x)\|_2 \|\phi(x^*) - x^*\|_2,
$$
it follows that
$$
\|\phi(x) - x\|_2^2 = \tau(\phi(x)\phi(x)^*) + 1 -
2{\text{\rm Re}} \tau(\phi(x)x^*)
$$
$$
= {\text{\rm Re}} \tau(\phi(x)x^*) +
{\text{\rm Re}} \tau(\phi(x)(\phi(x)^* - x^*)) + 1 -
2 {\text{\rm Re}} \tau(\phi(x)x^*)
$$
$$
\geq 1-{\text{\rm Re}} \tau(\phi(x)x^*) - \|\phi(x)-x\|_2\|\phi(x)\|_2
$$
$$
=\|x\xi_\phi - \xi_\phi x \|_2^2/2 - \|\phi(x)-x\|_2\|\phi(x)\|_2,
$$
which by part $1^\circ$ proves the statement.
\hfill Q.E.D.

The inequalities in the previous Lemmas
show in particular that if $\phi$ almost fixes some
$u \in \Cal U(N)$, then $\phi(ux)$ is close
to $u\phi(x)$, uniformly in $x \in N, \|x\| \leq 1$,
whenever we have a control over $\|\phi\|$:

\proclaim{Corollary} For any unitary element $u \in N$ and $x \in N$, we have:
$$
\|\phi(ux) - u \phi(x)\|_2 \leq \|\phi\|^{1/2} \|x\| \|[u, \xi_\phi]\|_2
$$
$$
\leq \|\phi\|^{1/2} \|x\|(2\|\phi(u) - u\|_2^2 +
2\|\phi(1)\|_2\|\phi(u)-u\|_2)^{1/2}.
$$
\endproclaim
\vskip .1in
\noindent
{\it Proof}. By using that
$\|\phi(ux) - u \phi(x)\|_2 =
{\text{\rm sup}} \{ |\tau((\phi(ux)-u\phi(x))y)| \mid y \in N, \|y\|_2\leq 1\}$,
we get:
$$
\|\phi(ux) - u \phi(x)\|_2
$$
$$
={\text{\rm sup}} \{  |\langle ux\xi_\phi y , \xi_\phi \rangle -
\langle x\xi_\phi yu, \xi_\phi \rangle | \mid y \in N, \|y\|_2 \leq 1\}
$$
$$
={\text{\rm sup}} \{  |\langle x\xi_\phi y , [u^*,\xi_\phi] \rangle |
\mid y \in N, \|y\|_2 \leq 1\}
$$
$$
\leq {\text{\rm sup}} \{\|x \xi_\phi y\|_2 \mid y\in N, \|y\|_2 \leq 1\}
\|[u^*, \xi_\phi]\|_2
$$
$$
=\|\phi(x^*x)\|^{1/2} \|[u, \xi_\phi]\|_2 \leq \|\phi\|^{1/2} \|x\| \|[u, \xi_\phi]\|_2.
$$
\hfill Q.E.D.
\vskip .1in
\noindent
{\it 1.1.3. From Hilbert bimodules to completely positive maps.} Conversely,
let $(\Cal H, \xi)$ be a pointed
Hilbert $(B\subset N)$-bimodule,
with $\langle \xi \cdot, \xi \rangle
\leq c \tau$,
for some $c > 0$. Let $T : L^2(N, \tau) \rightarrow \Cal H$ be the unique bounded
operator defined
by $T\hat{y} = \xi y, y \in N$. Thus $\langle \xi y , \xi y \rangle
\leq c \tau (yy^*) = c \|\hat{y}\|_2^2$, so that $\|T\| \leq c^{1/2}$.

It is immediate to check that if we denote
for clarity by $L(x)$ the operator of left multiplication by $x$
on $\Cal H$, then $T$ satisfies:
$$
\langle T^* L(x) T (J_NyJ_N (\hat{y_1})) , \hat{y_2} \rangle_\tau
= \langle L(x) (\xi y_1 y^*), \xi y_2 \rangle_{\Cal H}
$$
$$
=\langle L(x) \xi y_1, \xi y_2 y \rangle_{\Cal H} = \langle
J_Ny J_N (T^* L(x) T)\hat{y_1},
\hat{y_2}\rangle_\tau.
$$
This shows that the operator
$\phi_{(\Cal H, \xi)}(x)\overset\text{\rm def} \to = T^* L(x) T$
commutes with the right multiplication on $L^2(N, \tau)$ by elements $y \in N$.
Thus, $\phi_{(\Cal H, \xi)}(x)$
belongs to $(J_NNJ_N)' \cap \Cal B(L^2(N, \tau)) = N$,  showing
that $\phi_{(\Cal H, \xi)}$ defines a map from $N$ into $N$, which is obviously
completely positive and
$B$-bimodular, by
the definitions. Furthermore, if we denote by $\Cal H'$ the closed linear
span of $N\xi N$ in $\Cal H$, then $U: \Cal H_\phi
\rightarrow \Cal H', U( x \otimes y) = x\xi y$
is easily seen to be an isomorphism of Hilbert $(B\subset N)$-bimodules.

The assumption that $\xi$ is ``bounded from the right'' by $c$ is not really
a restriction for this construction, since if we put
$\Cal H^0 = \{\xi \in \Cal H \mid b\xi = \xi b, \forall b \in B$,
$\xi$ bounded from the left and from the right $\}$, then
it is easy to see
that $\Cal H^0$ is dense in the Hilbert space $\Cal H_0\subset \Cal H$
of all $B$-central vectors in $\Cal H$. This actually implies
that any $(B\subset N)$ Hilbert bimodule $(\Cal H, \xi)$
is a direct sum of some $(B\subset N)$ Hilbert bimodules $(\Cal H_i, \xi_i)$
with $\xi_i$ bounded both from left and right (Hint: just use the above density
and a  maximality argument).

Note that if $(\Cal H, \xi)$ comes itself
from a completely positive $B$-bimodular map $\phi$, i.e.,
$(\Cal H, \xi)=(\Cal H_\phi, \xi_\phi)$
as in 1.1.2, then $\phi_{(\Cal H, \xi)} = \phi.$ Similarly, if $(\Cal H, \xi)$
is a cyclic pointed $(B\subset N)$-Hilbert bimodule and $\phi=\phi_{(\Cal H, \xi)}$,
then $(\Cal H_\phi, \xi_\phi) \simeq (\Cal H, \xi)$.

Let us also note a converse to Lemma 1.1.3,
showing that if $\xi$ almost
commutes with a unitary element $u \in N$ then $u$ is almost
fixed by $\phi=\phi_{(\Cal H, \xi)}$,
provided we have some control over $\|\phi(1)\|_2$:

\proclaim{Lemma} Let
$\xi \in \Cal H$ be a vector bounded from the right and
denote $\phi = \phi_{(\Cal H, \xi)}$.

$1^\circ$. Let $a_0, b_0 \in L^1(N, \tau)_+$ be such that
$\langle \cdot \xi, \xi \rangle
= \tau(\cdot b_0), \langle \xi \cdot, \xi \rangle = \tau(\cdot a_0)$
and put $a = 1\vee a_0, b=1\vee b_0$,
$\xi' = b^{-1/2}\xi a^{-1/2}$. Then $\phi(1)=a_0$ and we have:
$$
\|\xi-\xi'\|^2 \leq 2\|a_0-1\|_1 + 2\|b_0 - 1\|_1.
$$

$2^\circ$. If $u \in \Cal U(N)$,  then we have:
$$
\|\phi(u)-u\|^2_2 \leq \|[u, \xi]\|_2^2 + (\|\phi(1)\|^2_2 - 1).
$$
\endproclaim
\vskip .1in
\noindent
{\it Proof}. $1^\circ$. We have:
$$
\|\xi-\xi'\|^2 \leq  2\|\xi-b^{-1/2}\xi\|^2 + 2\|\xi-\xi a^{-1/2}\|^2
$$
$$
= 2\tau((1-b^{-1/2})^2b_0) + 2\tau((1-a^{-1/2})^2a_0)
$$
$$
\leq 2\|b_0 - 1\|_1 + 2\|a_0-1\|_1.
$$

$2^\circ$. By part $1^\circ$ of Lemma 1.1.2
we have $\tau(\phi(u^*) \phi(u)) \leq \tau(\phi(1)\phi(1))$,
so that:
$$
\|\phi(u) - u \|^2_2=\tau(\phi(u) \phi(u^*)) + 1 - 2{\text{\rm Re}} \tau (\phi(u)u^*)
$$
$$
\leq \tau(\phi(1)\phi(1)) +1 - 2{\text{\rm Re}} \tau (\phi(u)u^*)
$$
$$
= 2 - 2{\text{\rm Re}} \tau (\phi(u)u^*) + (\tau(\phi(1)\phi(1))-1)
$$
$$
= \|[u, \xi]\|_2^2 + (\|\phi(1)\|^2_2 - 1).
$$
\hfill Q.E.D.
\vskip .1in
\noindent
{\it 1.1.4. Correspondences from representations of groups.}
Let $\Gamma_0$ be a discrete
group, $(B, \tau_0)$ a finite von Neumann algebra
with a normal faithful tracial state and $\sigma$ a cocycle
action of $\Gamma_0$ on $(B, \tau_0)$ by $\tau_0$-preserving
automorphisms. Denote by $N = B \rtimes_\sigma \Gamma_0$ the
corresponding cross-product algebra and by
$\{u_g\}_g\subset N$ the canonical
unitaries implementing the action $\sigma$ on $B$.

Let $(\pi_0, \Cal H_0, \xi_0)$ be a pointed, cyclic representation
of the group $\Gamma_0$. We denote by $(\Cal H_{\pi_0}, \xi_{\pi_0})$
the pointed Hilbert space $(\Cal H_0, \xi_0) \overline{\otimes}
(L^2(N, \tau), \hat{1})$. We let $N$ act on the right
on $\Cal H_{\pi_0}$ by $(\xi \otimes \hat{x})y
= \xi \otimes (\hat{xy}), x, y \in N, \xi \in \Cal H_0$ and on the left by
$b (\xi \otimes \hat{x}) = \xi \otimes \hat{bx}$, $u_g(\xi \otimes \hat{x})
= \pi_0(g)(\xi) \otimes \hat{u_gx}, b\in B, x\in N, g \in \Gamma_0, \xi \in \Cal H_0$.

It is easy to check that these are indeed mutually commuting
left-right actions of $N$ on
$\Cal H_{\pi_0}$. Moreover, the vector $\xi_{\pi_0}=\xi_0 \otimes \hat{1}$
implements the trace $\tau$ on $N$, both from left and right. Also,
$\xi_{\pi_0}$ is easily seen to be $B$-central.
Thus, $(\Cal H_{\pi_0}, \xi_{\pi_0})$ is a
Hilbert $(B\subset N)$-bimodule.

Let now $\varphi$ be a positive definite function on $\Gamma_0$
and denote by $(\pi_\varphi, \Cal H_\varphi, \xi_\varphi)$ the
representation obtained from it through the GNS construction. Let
$(\Cal H, \xi)$ denote the $(B \subset B\rtimes \Gamma_0)$-Hilbert
bimodule constructed out of the
representation $\pi_\varphi$ as above and $\phi$ the
completely positive $B$-bimodular map associated with $(\Cal H, \xi)$
as in 1.1.3. An easy calculation shows that $\phi$ acts
on $B \rtimes \Gamma_0$ by $\phi(\Sigma_g b_gu_g) =
\Sigma_g \varphi(g) b_gu_g$.

Conversely, if $(\Cal H, \xi)$ is a $(B\subset N)$ Hilbert bimodule,
then we can associate to it the representation $\pi_0$ on
$\Cal H_0 = \overline{\text{\rm sp}} \{u_g\xi u_g^* \mid g\in \Gamma_0\}$
by $\pi_0(g)\xi' = u_g\xi'u_g^*, \xi'\in \Cal H_0$. Equivalently,
if $\phi$ is the $B$-bimodular completely positive map
associated with $(\Cal H, \xi)$ then $\varphi(g) = \tau(\phi(u_g)u_g^*),
g\in \Gamma_0$, is a positive definite function on $\Gamma_0$.
\vskip .1in
\noindent
{\it 1.1.5. The adjoint of a bimodule}. Let $(\Cal H, \xi_0)$ be a
$(B\subset N)$ Hilbert bimodule. Let $\overline{\Cal H}$ be the
conjugate Hilbert space of $\Cal H$, i.e.,
$\overline{\Cal H}=\Cal H$ as a set, the sum of vectors in $\overline{\Cal H}$
is the same as in $\Cal H$, but the multiplication by scalars is
given by $\lambda \cdot \xi = \overline{\lambda}\xi$ and
$\langle \xi, \eta \rangle_{\overline{\Cal H}} = \langle \eta,\xi \rangle_{\Cal H}$.
Denote by $\overline{\xi}$ the element $\xi$ regarded as a vector in
the Hilbert space  $\overline{\Cal H}$. Define on $\overline{\Cal H}$ the
left and right
multiplication operations
by  $x \cdot \overline{\xi} \cdot y
= \overline{y^*\xi x^*}.$,
for $x, y \in N, \xi \in \Cal H$. It is easy to see that they
define a $N$ Hilbert bimodule structure
on $\overline{\Cal H}$. Moreover, $\overline{\xi_0}$ is
clearly $B$-central. We call $(\overline{\Cal H}, \overline{\xi_0})$
the {\it adjoint of} $(\Cal H, \xi_0)$. Note that we clearly have $(\overline{\overline{\Cal H}}, \overline{\overline{\xi_0}}) = (\Cal H, \xi_0)$.

\proclaim{Lemma} Let $\phi$ be a normal $B$-bimodular completely
positive map on $N$. For each $x\in N$ let $\phi^*(x)\in L^1(N, \tau)$ denote
the Radon-Nykodim derivative of $N \ni y \mapsto \tau(\phi(y)x)$ with
respect to $\tau$.

$1^\circ$. $\phi^*(N) \subset N$ if and only if
$\tau \circ \phi \leq c \tau$ for
some $c > 0$,
i.e., if and only if the Radon-Nykodim derivative
$b_0={\text{\rm d}}\tau\circ \phi/{\text{\rm d}}\tau$
is a bounded operator. Moreover, if the condition is satisfied
then $\phi^*$ defines a normal, $B$-bimodular, completely
positive map of $N$ into $N$ with $\phi^*(1)=b_0$ and
$$
\|\phi^*\| = \|b_0\|
={\text{\rm inf}} \{c > 0 \mid \tau \circ \phi \leq c \tau\}.
$$

$2^\circ$. If $\phi$ satisfies condition $1^\circ$ then $\phi^*$
also satisfies it, and we have $(\phi^*)^* = \phi$. Also,
$(\Cal H_{\phi^*}, \xi_{\phi^*})=
(\overline{\Cal H_{\phi}}, \overline{\xi_{\phi}})$.

$3^\circ$. If $\tau \circ \phi \leq \tau$
then for any unitary element $u\in N$ we have:
$$
\|\phi^*(u)-u\|_2^2 \leq 2 \|\phi(u)-u\|_2.
$$
\endproclaim
\noindent
{\it Proof}. Parts $1^\circ$ and $2^\circ$ are
trivial by the definition of $\phi^*$.

To prove $3^\circ$,
note that by part $1^\circ$,
$\tau \circ \phi \leq \tau$ implies  $\phi^*(1) \leq 1$ and so by
Lemma 1.1.2 we get:
$$
\|\phi^*(u)-u\|_2^2 = \tau(\phi^*(u)\phi^*(u)^*) + 1 -
2{\text{\rm Re}} \tau(\phi^*(u)u^*)
$$
$$
\leq \tau(\phi^*(1)\phi^*(1)) + 1 -
2 {\text{\rm Re}} \tau(\phi(u)u^*) \leq 2 - 2 {\text{\rm Re}} \tau(\phi(u)u^*)
$$
$$
= 2 {\text{\rm Re}} \tau((u-\phi(u))u^*) \leq 2\|\phi(u)-u\|_2.
$$
\hfill Q.E.D.
\vskip .1in
\noindent
{\bf 1.2. Completely positive maps as Hilbert space operators}.
We now show that if a completely positive map $\phi$ on the
finite von Neumann algebra $N$ is sufficiently smooth
with respect to the normal faithful tracial state $\tau$ on $N$,
then it can be extended to the Hilbert space $L^2(N,\tau)$.
In case $\phi$ is $B$-bimodular, for some von Neumann
subalgebra $B\subset N$, these operators belong to the
algebra of the basic construction associated with $B\subset N$, defined
in the next paragraph.

\proclaim{1.2.1. Lemma}
$1^\circ$. If there exists $c > 0$ such that $\|\phi(x)\|_2 \leq
c \|x\|_2, \forall x\in N$, then there exists a bounded operator
$T_\phi$ on $ L^2(N,\tau)$ such that $T_\phi(\hat{x}) = \hat{\phi(x)}$.
The operator
$T_\phi$ commutes with the canonical
conjugation $J_N$. Also, if $B\subset N$
is a von Neumann subalgebra, then
$T_\phi$ commutes with the operators of left
and right multiplication
by elements in $B$ (i.e., $T_\phi \in B'\cap (JBJ)'$)
if and only if
the completely
positive map $\phi$ is $B$-bimodular.

$2^\circ$. If $\tau \circ \phi
\leq c_0 \tau,$ for some constant $c_0 > 0$, then $\phi$ satisfies
condition $1^\circ$ above, and so  there exists a bounded
operator $T_\phi$ on the Hilbert space
$L^2(N, \tau)$ such that $T_\phi(\hat{x}) = \hat{\phi(x)}$, for $x\in N$.
Moreover, if $\phi^*: N \rightarrow N$ is the adjoint of
$\phi$, as defined in $1.1.4$, then
$\|T_\phi\|^2\leq \|\phi(1)\| \|\phi^*(1)\|$. Also,
$\phi^*$ satisfies $\tau \circ \phi^*
\leq \|\phi(1)\| \tau$ and we have $T_{\phi^*} = T_{\phi}^*$.

$3^\circ$. If $\phi$ is $B$-bimodular then $\phi(1) \in B'\cap N$. Thus, if
we assume $B'\cap N = \Cal Z(B)$ then $\phi(1) \in \Cal Z(B)$,
$\tau \circ \phi \leq \|\phi(1)\| \tau$ and
the bounded operator $T_\phi$ exists by $2^\circ$.
If in addition $\phi(1)=1$, then
$\phi$ is trace-preserving as well.
\endproclaim
\vskip .1in
\noindent
{\it Proof}. 1$^\circ$. The existence of $T_\phi$ is trivial. Also,
for $x\in N$ we have
$$
T_\phi(J(\hat{x}))= \hat{\phi(x^*)}=\hat{\phi(x)^*}=J_N(T_\phi(\hat{x})).
$$
If $\phi$ is $B$-bimodular
and $b\in B$ is regarded as an operator of left multiplication
by $b$ on $L^2(N, \tau)$, then we have
$$
bT_\phi(\hat x) =
\hat{b\phi(x)} = \hat{\phi(bx)} = T_\phi(b\hat x).
$$
Thus,
$T_\phi\in B'$.

Similarily,
$$
JbJ (T_\phi(\hat x)) = \phi(x)b = \phi(xb) =
T_\phi(JbJ (\hat x))
$$
showing that $T_\phi\in JBJ'$ as well.
Conversely, if $T_\phi \in B'\cap JBJ'$,
then by exactly the same equalities we get $\phi(bx) =
b\phi(x), \phi(xb) = \phi(x)b, \forall x\in N, b\in B$.

$2^\circ$. By
Kadison's inequality, for $x\in M$ we have
$$
\langle T_\phi(\hat x), T_\phi(\hat x) \rangle =
\tau (\phi(x)^*\phi(x)) \leq \|\phi (1)\| \tau (\phi(x^*x)), \forall x\in N.
$$
Thus, by Lemma 1.1.5 we have $\|T_\phi\|^2 \leq \|\phi(1)\| \|\phi^*(1)\|$.
The last part is now trivial, by 1.1.5 and the
definitions of $T_\phi$, $\phi^*$ and $T_{\phi^*}$.

$3^\circ$. The $B$-bimodularity of $\phi$ implies
$u\phi(1)u^* = \phi(1), \forall u\in \Cal U(B)$, thus $\phi(1)
\in B'\cap N$.

Using again the bimodularity, as well as the
normality of $\phi$, for each fixed $x \in N$ we have
$$
\tau(\phi(x)) = \tau(u\phi(x) u^*)=\tau(\phi(uxu^*)) = \tau (\phi(y))
$$
for all $u \in \Cal U(B)$ and
all $y$ in the weak closure of the convex hull of $\{uxu^* \mid
u\in \Cal U(N)\}$. The latter set contains $E_{B'\cap N}(x)\in
B'\cap N \subset B$
(see e.g. [Po6]), thus
$$
\tau(\phi(x))= \tau(\phi(E_{B'\cap N}(x))) =
\tau(E_{B'\cap N}(x)\phi(1)).
$$
This shows that if $x \geq 0$ then $\tau(\phi(x)) \leq \|\phi(1)\| \tau(x)$.
It also shows that in case $\phi(1)=1$ then
$\tau(\phi(x))=\tau(x), \forall x\in N$.
\hfill Q.E.D.
\vskip .1in
\noindent
{\bf 1.3. The basic construction and its compact ideal space.}
We now recall
from ([Chr], [J1], [Po2,3]) some well known facts about the
{\it basic construction}
for an inclusion of finite von Neumann
algebras $B\subset N$ with a normal
faithful tracial state $\tau$ on it.
Also, we establish some properties of the
ideal generated by finite projections in the semifinite von Neumann
algebra $\langle N, B \rangle$ of the basic construction.
\vskip .1in
\noindent
{\it 1.3.1. Basic construction for} $B\subset N$.
We denote by $\langle N, B \rangle$
the von Neumann algebra generated in $\Cal B (L^2(N, \tau))$
by $N$ (regarded as the algebra of left multiplication operators
by elements in $N$) and by the orthogonal projection $e_B$
of $L^2(M, \tau)$ onto $L^2(B, \tau)$.

Since $e_Bxe_B = E_B(x)e_B, \forall x\in N$,
where $E_B$ is the unique $\tau$-preserving conditional
expectation of $N$ onto $B$, and $\vee \{x(e_B(L^2(N))) \mid
x\in N\}=L^2(N)$, it follows that
sp$Ne_BN$ is a *-algebra with support equal to $1$ in
$\Cal B(L^2(N, \tau))$. Thus,
$\langle N, B \rangle=
{\overline{\text{\rm sp}}}^{\text{\rm w}} \{xe_By \mid
x, y \in N\}$ and
$e_B \langle N, B, \rangle e_B = Be_B$.

One can also readily see that if $J=J_N$ denotes
the canonical conjugation on the
Hilbert space $L^2(N, \tau)$, given
on $\hat{N}$ by $J(\hat{x})=\hat{x^*}$, then $\langle N, B \rangle=
JBJ'\cap \Cal B(L^2(N,\tau)) $.
This shows in particular that $\langle N, B \rangle$
is a semifinite von Neumann algebra. It also shows that
the isomorphism of $N\subset \langle N, B \rangle$ only depends
on $B \subset N$ and not on the trace $\tau$ on $N$
(due to the uniqueness of the standard
representation).

As a consequence, if $\phi$ is a $B$-bimodular completely positive map
on $N$ satisfying $\|\phi(x)\|_2 \leq c\|x\|_2, \forall x\in N$,
for some constant $c > 0$, as in Lemma 1.2.1, then the corresponding
operator $T_\phi$ on $L^2(N, \tau)$ defined by $T_\phi(\hat{x}) =
\hat{\phi(x)}, x\in N$ belongs to $B'\cap \langle N, B \rangle$.

We endow $\langle N, B \rangle$ with the unique normal
semifinite faithful trace $Tr$ satisfying
$Tr(xe_By)=\tau(xy), \forall x, y \in N$.
Note that there exists a unique $N$ bimodule map
$\Phi$ from sp$Ne_BN \subset \langle N, B \rangle $ into
$N$ satisfying $\Phi(xey) = xy, \forall
x, y \in N$,
and $\tau\circ \Phi= Tr$. In particular
this entails $\|\Phi(X)\|_1
\leq \|X\|_{1, Tr}, \forall X \in
{\text{\rm sp}}Ne_BN$. Note that the map $\Phi$ extends uniquely to a
$N$-bimodule map from $L^1(\langle N, B \rangle, Tr)$ onto
$L^1(N, \tau)$, still denoted $\Phi$. This $N$-
bimodule map satisfies the
``pull down'' identity $eX = e\Phi(eX), \forall X \in \langle N, B \rangle$
(see [PiPo], or [Po2]). Note that $\Phi(eX)$ actually belongs
to $L^2(N, \tau)\subset L^1(N, \tau)$,
for $X\in \langle N, B \rangle$.
\vskip .1in
\noindent
{\it 1.3.2. The compact ideal space of a semifinite algebra}.
In order to define the compact ideal space of the semifinite
von Neumann algebra $\langle N, B \rangle$, it will be useful to first
mention some remarks about the compact ideal space
of an arbitrary semifinite von Neumann algebra $\Cal N$.

Thus, we let $\Cal J(\Cal N)$ be the norm-closed
two-sided ideal generated in $\Cal N$ by the finite projections
of $\Cal N$, and call it the {\it compact ideal space} of $\Cal N$
(see e.g., [KafW], [PoRa]).
Note that $T\in \Cal N$ belongs to $\Cal J(\Cal N)$ if and only
if all the spectral projections $e_{[s, \infty)}(|T|), s > 0,$ are finite projections in $\Cal N$.
As a consequence, it follows that the
set $\Cal J^0(\Cal N)$ of all elements supported by
finite projections (i.e., the {\it finite rank} elements in $\Cal J(\Cal N)$) is
a norm dense ideal in $\Cal J(\Cal N)$.

Let further $e\in \Cal N$ be a finite
projection with central support equal to $1$ and denote
by $\Cal J_e(\Cal N)$ the norm closed
two-sided ideal generated by $e$ in $\Cal N$.
It is easy to see that an operator $T \in \Cal N$
belongs to $\Cal J(\Cal N)$ if and only if there exists
a partition of $1$ with projections $\{z_i\}_i$
in $\Cal Z(\Cal N)$  such that $Tz_i \in
\Cal J_e(\Cal N), \forall i$. In particular,
if $p\in \Cal N$ is a finite projection then there exists
a net of projections $z_i \in \Cal Z(\Cal N)$ such that
$z_i \uparrow 1$ and $pz_i \in \Cal J_e(\Cal N), \forall i$
(see e.g., 2.1 in [PoRa]). Also, $T \in \Cal J_e(\Cal N)$
iff $e_{[s, \infty)}(|T|)\in \Cal J_e(\Cal N), \forall s>0$.
In turn, a projection $f\in \Cal N$ lies in $\Cal J_e(\Cal N)$ iff
there exists a constant
$c>0$ such that that $Tr(fz) \leq c Tr (ez)$, for any normal
semifinite trace $Tr$ on $\Cal N$ and any projection $z\in \Cal Z(\Cal N)$.

The next result,
whose proof is very similar to some arguments in ([Po7]),
shows that one can ``push'' elements of $\Cal J(\Cal N)$ into
the commutant of a subalgebra $\Cal B$ of $\Cal N$, while still staying
in the ideal $\Cal J(\Cal N)$, by averaging by unitaries
in $\Cal B$. We include a complete proof, for
convenience.

\proclaim{Proposition} Let $\Cal B \subset \Cal N$ be a von Neumann
subalgebra of $\Cal N$. For $x\in \Cal N$ denote
$K_x=\overline{\text{\rm co}}^w \{uxu^* \mid
u \in \Cal U(\Cal B)\}$. If $x\in \Cal J(\Cal N)$ then
$\Cal B'\cap K_x$ consists of exactly one element,
denoted $\Cal E_{\Cal B'\cap \Cal N}(x)$, and which belongs to $\Cal J(\Cal N)$.
Moreover, the application $x \mapsto \Cal E_{\Cal B'\cap \Cal N}(x)$
is a conditional expectation of $\Cal J(\Cal N)$ onto $\Cal B'\cap \Cal J(\Cal N)$.
Also, if $x\in \Cal J_e(\Cal N)$ for some
finite projection $e\in \Cal N$ of central support $1$, then
$\Cal E_{\Cal B'\cap \Cal N}(x)\in \Cal J_e(\Cal N)$.
\endproclaim
\noindent
{\it Proof}. If $x=f$ is a projection in $\Cal J_e(\Cal N)$
then $Tr(fz) \leq c Tr (ez)$, for any normal
semifinite trace $Tr$ on $\Cal N$ and any projection $z\in \Cal Z(\Cal N)$.
By averaging with unitaries and taking
weak limits, this implies that
for some appropriate constant $c'>0$
we have $Tr(yz) \leq c Tr(ez), \forall y \in K_f$,
so  that $Tr(pz) \leq c' Tr(ez)$, for any spectral projection
$p=e_{[s, \infty)}(y), s>0$. Thus, $K_f \subset \Cal J_e(\Cal N)$. Since
any $x \in \Cal J_e(\Cal N)$ is a norm limit of linear
combinations of projections $f$ in $\Cal J_e(\Cal N)$, this shows that
the very last part of the statement follows from the first part.

To prove the first part, consider first the case when
$\Cal N$ has a normal semifinite faithful trace $Tr$. Assume
first that $x\in \Cal J(\Cal N)$ actually belongs
to $\Cal N \cap L^2(\Cal N, Tr)$ ($\subset \Cal J(\Cal N))$.
Note that in this case all $K_x\subset \Cal N$ is a
subset of the Hilbert space $L^2(\Cal N, Tr)$, where
it is convex and weakly closed.
Let then $x_0\in K_x$
be the unique element of minimal Hilbert norm $\|\quad \|_{2,Tr}$ in $K_x$.
Since $\|ux_0u^*\|_{2,Tr} = \|x_0\|_{2,Tr}, \forall
u\in \Cal U(\Cal B)$, it follows that $ux_0u^* = x_0, \forall
u\in \Cal U(\Cal B)$. Thus, $x_0 \in \Cal B'\cap \Cal N \cap L^2(\Cal N, Tr)$.
In particular, $x_0\in \Cal B'\cap \Cal J(\Cal N)$.

If we now denote by $p$ the orthogonal projection
of $L^2(\Cal N, Tr)$ onto the space of fixed points of the representation
of $\Cal U(\Cal B)$ on it given by $\xi \mapsto u\xi u^*$,
then $x_0$ coincides with $p(x)$. Since $p(uxu^*) = p(x)$,
this shows that $x_0 = p(x)$ is in fact the unique element $y$
in $K_x$
with $uyu^*=y, \forall u\in \Cal U(\Cal B)$.
Thus, if for each $x\in \Cal N \cap L^2(\Cal N, Tr)$ we put
$\Cal E_{\Cal B'\cap \Cal N}(x) \overset \text{\rm def}
\to = p(x)$, then
we have proved the statement for the
subset $\Cal N \cap L^2(\Cal N, Tr)$.

Since $\|y\|\leq \|x\|, \forall y\in K_x$, it follows that if $\{x_n\}_n
\subset \Cal N \cap L^2(\Cal N, Tr)$
is a Cauchy sequence (in the uniform norm), then
so is $\{\Cal E_{\Cal B'\cap \Cal N}(x_n)\}_n$. Thus,
$\Cal E_{\Cal B'\cap \Cal N}$ extends uniquely by continuity to a linear,
norm one projection from $\Cal J(\Cal N)$ onto $\Cal B'\cap \Cal J(\Cal N)$,
which by the above remarks takes the norm dense
subspace $\Cal N \cap L^2(\Cal N, Tr)$ into itself.

Let us now prove that $\Cal B'\cap K_x \neq \emptyset, \forall x\in \Cal J
(\Cal N)$. To this end, let
$x$ be an arbitrary element in $\Cal J(\Cal N)$ and
$\varepsilon > 0$. Let
$x_1 \in \Cal N \cap L^2(\Cal N, Tr)$ with
$\|x-x_1\|\leq \varepsilon$. Write $\Cal E_{\Cal B'\cap \Cal N}(x_1)$
as a weak limit of a net $\{T_{u_\alpha}(x_1)\}_\alpha$,
for some finite tuples $u_\alpha=(u^\alpha_1, ..., u^\alpha_{n_\alpha})
\subset \Cal U(\Cal B)$. By passing to a subnet if necessary, we may assume
$\{T_{u_\alpha}(x)\}_\alpha$ is also weakly convergent, to some
element $x' \in K_{x}$. Since,
$\|T_{u_\alpha}(x)-T_{u_\alpha}(x_1)\|
\leq \|x-x_1\| \leq \varepsilon$, it follows that $\|x'-
\Cal E_{\Cal B'\cap \Cal N}(x_1)\| \leq \varepsilon$. This
shows that the weakly-compact set $K_x$ contains elements which are
arbitrarily close to $\Cal B'\cap \Cal N$. By taking
a weak limit of such elements it follows that
$\Cal B'\cap K_x \neq \emptyset$.

Finally, let $x\in \Cal J(\Cal N)$ and
assume $x^0$ is an element in $\Cal B'\cap K_x$. To prove that
$x^0=\Cal E_{\Cal B'\cap \Cal N}(x)$, let $\varepsilon > 0$ and
$x_1 \in \Cal N \cap L^2(\Cal N, Tr)$ with
$\|x-x_1\|\leq \varepsilon$, as before. Write $x^0$ as a weak limit of a net $\{T_{v_\beta}(x)\}_\beta$,
for some finite tuples $v_\beta=(v_1^\beta, ..., v_{m_\beta}^\beta)
\subset \Cal U(\Cal B)$. By passing to a subnet if necessary, we may assume
$\{T_{v_\beta}(x_1)\}_\beta$ is also weakly convergent, to some
element $x_1^0 \in K_{x_1}$. Since,
$\|T_{v_\beta}(x)-T_{v_\beta}(x_1)\|
\leq \|x-x_1\| \leq \varepsilon$, it follows that $\|x^0-x_1^0\| \leq \varepsilon$.
But $p(x_1^0)=p(x_1)= \Cal E_{\Cal B'\cap \Cal N}(x_1)$,
and $p(x_1^0)$ is obtained as a weak limit of averaging by unitaries
in $\Cal B$, which commute with $x^0$. Thus,
$$
\|x^0- \Cal E_{\Cal B'\cap \Cal N}(x)\|  \leq
\|x^0- \Cal E_{\Cal B'\cap \Cal N}(x_1)\|+
\|\Cal E_{\Cal B'\cap \Cal N}(x_1)-\Cal E_{\Cal B'\cap \Cal N}(x)\|
\leq \varepsilon + \|x_1-x\| \leq 2\varepsilon.
$$
Since $\varepsilon > 0$ was arbitrary, this shows that $x^0 =
\Cal E_{\Cal B'\cap \Cal N}(x)$.

This finishes the proof of the case when $\Cal N$ has a faithful
trace $Tr$. The general case follows now readily, by noticing
that if $\{z_i\}_i$ is an increasing net of projections in
$\Cal Z(\Cal N)$ such that $K_{z_ix} \cap (\Cal Bz_i)'$
consists of exactly one element, which belongs to $\Cal J(\Cal N)z_i=
\Cal J(\Cal Nz_i)$,
$\forall x\in \Cal J(\Cal N)$, then the same holds true for the projection
$\underset i \rightarrow \infty \to \lim z_i$.
\hfill Q.E.D.

\vskip .1in
\noindent
{\it 1.3.3. The compact ideal space of} $\langle N,B \rangle$.
In particular, if $B \subset N$ is an inclusion of finite von Neumann
algebras as
in 1.3.1, then we denote
by $\Cal J(\langle N, B \rangle)$ the compact ideal space of
$\langle N, B \rangle$. Noticing that
$e_B$ has central support $1$ in
$\langle N, B \rangle$, we denote $\Cal J_0(\langle N,B \rangle)$
the norm closed two sided ideal $\Cal J_{e_B}(\langle N,B \rangle)$
generated by $e_B$ in $\langle N, B \rangle$.
Note that if $B = \Bbb C$ then
$\Cal J(\langle N, B \rangle)=\Cal J_0(\langle N, B \rangle)$
is the usual ideal of compact operators $\Cal K(L^2(N))$.

It will be
useful to have the following alternative characterisations
of the compact ideal spaces $\Cal J(\langle N, B \rangle),
\Cal J_0(\langle N, B \rangle)$.

\proclaim{Proposition} Let $N$ be a finite
von Neumann algebra with countably decomposable center and
$B\subset N$ a von Neumann subalgebra. Let $T\in \langle N,B \rangle$.
The following conditions are equivalent:

$1^\circ$. $T\in \Cal J(\langle N,B \rangle)$.

$2^\circ$. For any $\varepsilon > 0$ there exists a finite projection $p\in
\langle N,B \rangle$ such that $\|T(1-p)\| < \varepsilon$.

$3^\circ$. For any $\varepsilon > 0$ there exists $z \in
\Cal P(\Cal Z(J_NBJ_N))$ such that $\tau(1-z) \leq \varepsilon$ and
$Tz \in \Cal J_0(\langle N,B \rangle)$.

$4^\circ$. For any given sequence $\{\eta_n \}_n \in L^2(N)$ with
the properties $E_B(\eta_n^* \eta_n) \leq 1, \forall n \geq 1,$ and
$\underset n \rightarrow \infty \to \lim
\|E_B(\eta_n^*\eta_m)\|_2 = 0, \forall m$, we have
$\underset n \rightarrow \infty \to \lim \|T\eta_n\|_2=0$.

$5^\circ$. For any given sequence $\{x_n \}_n \in N$ with
the properties $E_B(x_n^* x_n) \leq 1, \forall n\geq 1,$ and
$\underset n \rightarrow \infty \to \lim
\|E_B(x_n^* x_m)\|_2 = 0, \forall m$, we have
$\underset n \rightarrow \infty \to \lim \|Tx_n\|_2=0$.

Moreover,
$T \in \Cal J_0(\langle N, B \rangle)$ if and only if condition $2^\circ$
above holds true with projections $p$ in $\Cal J_0(\langle N, B \rangle)$.
\endproclaim
\noindent
{\it Proof}. The equivalence of $1^\circ$ and $2^\circ$
(resp. the equivalence in the last part of the statement) is trivial
by the following fact, noted
in 1.3.2: $T \in \Cal J(\langle N,B \rangle)$
(resp. $T \in \Cal J_0(\langle N,B \rangle)$) iff
$e_{[s, \infty)}(|T|)\in \Cal J(\langle N, B \rangle)$
(resp. $\in \Cal J_0(\langle N,B \rangle)$),
$\forall s > 0$.

$3^\circ \implies 2^\circ$ is trivial by the general remarks in 1.3.2. To
prove  $2^\circ \implies 3^\circ$, for each $n \geq 1$ let
$T_n$ be a linear combination of finite projections
in $\langle N,B \rangle$ such that $\|T-T_n\| \leq 2^{-n}$.
Noting that for any finite
projection $e\in \langle N, B \rangle$
and $\delta > 0$ there exists a projection $z\in
\Cal Z(\langle N,B \rangle) = J_N\Cal Z(B)J_N$ such that $\tau(1-z)
\leq \delta$ and $ez \in \Cal J_0(\langle N,B \rangle)$,
it follows that for each $n$ there exists a projection $z_n \in
J_N\Cal Z(B)J_N$ such that $\tau(1-z_n) \leq 2^{-n} \varepsilon$
and $T_nz_n \in \Cal J_0(\langle N,B \rangle)$.
Let $z= \wedge z_n$. Then we have
$\tau(1-z) \leq \Sigma_n 2^{-n} \varepsilon
\leq \varepsilon$, $T_nz \in
\Cal J_0(\langle N,B \rangle)$ and
$\|(T-T_n)z\| \leq \|T-T_n\| \leq 2^{-n}, \forall n$.
Thus, $Tz \in \Cal J_0(\langle N, B \rangle)$ as well.

$3^\circ \implies 4^\circ$ is just a particular case of (2.5 in [PoRa]).
To prove $4^\circ \implies 1^\circ$, assume by contradiction
that there exists $s > 0$ such that the spectral projection $e=e_s(|T|)$
is properly infinite. It follows that there exist mutually orthogonal,
mutually equivalent projections $p_1, p_2, ... \in \langle N,B \rangle$
such that $\Sigma_n p_n \leq e$ with $p_n$
majorised by $e_B, \forall n$. Thus, for each $n \geq 1$ there exists
$\eta_n \in L^2(N)$ such that $p_n = \eta_n e_B \eta_n^*$. It then
follows that $E_B(\eta_n^*\eta_m) = 0$ for $n \neq m$, with
$E_B(\eta_n^*\eta_n)$ mutually equivalent projections in $B$. In particular,
$\|\eta_n\|_2^2 = \tau(\eta^*_n \eta_n)=c>0$ is constant, $\forall n$.
Thus,
$$s^{-1}\|T\eta_n\|_2\geq
\|e(\eta_n)\|_2 \geq \|p_n(\eta_n)\|_2 = \|\eta_n\|_2 = c^{1/2}, \forall n,
$$
a contradiction.

$4^\circ \implies 5^\circ$ is trivial. To prove
$5^\circ \implies
4^\circ$ assume $5^\circ$
holds true and let $\eta_n$ be  a sequence satisfying
the hypothesis in $4^\circ$.
For each $n$ let $q_n$ be a spectral projection corresponding
to some interval $[0, t_n]$
of $\eta_n\eta_n^*$ (the latter regarded
as a positive, unbounded, summable operator
in $L^1(N)$) such that $\|\eta_n - q_n\eta_n\|_2 < 2^{-n}$.
Thus, $x_n = q_n\eta_n $ lies in $N$. One can
easily check $E_B(x_n^*x_n)
\leq E_B(\eta_n^* \eta_n) \leq 1$ and
$\underset n \rightarrow \infty \to \lim
\|E_B(x_n^* x_m)\|^2_2$ $= \underset n \rightarrow \infty \to \lim
Tr((q_n\eta_ne_B\eta_n^*q_n) (q_m\eta_me_B\eta_m^*q_m))$
$=0$.
Thus
$\underset n \rightarrow \infty \to \lim \|Tx_n\|_2=0$. But
$\|T\eta_n\|_2 \leq \|Tx_n\|_2 + \|T\|\|\eta_n - x_n\|_2
\leq  \|Tx_n\|_2+2^{-n}\|T\|$,
showing that $\underset n \rightarrow \infty \to \lim \|T\eta_n\|_2=0$
as well.
\hfill Q.E.D.

\vskip .1in
\noindent
{\bf 1.4. Discrete embeddings and
bimodule decomposition.} If $B\subset N$ is an inclusion
of finite von Neumann algebras with a faithful normal tracial state $\tau$
as before, then we often consider $N$ as an (algebraic) (bi)module
over $B$ and $L^2(N, \tau)$ as a Hilbert (bi)module over $B$.
In fact any vector subspace $H$ of $N$ which
is invariant under left (resp. right) multiplication
by $B$ is a left (resp. right) module over $B$. Similarly, any
Hilbert subspaces of $L^2(N, \tau)$ which
is invariant under multiplication to the left (resp. right) by elements in $B$
is a left (resp. right) Hilbert module. Also,
the closure in $L^2(N,\tau)$ of a $B$-module $H\subset N$ is a Hilbert $B$-module.
\vskip .1in
\noindent
{\it 1.4.1. Orthonormal basis}. An {\it orthonormal basis} for a
right  (respectively left) Hilbert $B$-module $\Cal H \subset L^2(N, \tau)$
is a subset $\{\eta_i\}_i \subset L^2(N)$ such that
$\Cal H = \overline{\Sigma_k \eta_k B}$
(respectively $\Cal H = \overline{\Sigma_k B \eta_k}$) and
$E_B(\eta_i^*\eta_{i'}) = \delta_{ii'}p_i \in \Cal P(B), \forall i, i',$
(respectively $E_B(\eta_{j'}\eta_j^*) = \delta_{j'j} q_j \in
\Cal P(B), \forall
j, j'$). Note that if this is the case, then we have:
$\xi = \Sigma_i \eta_i
E_B(\eta_i^* \xi),
\forall \xi \in \Cal H$ (resp. $\xi = \Sigma_j E_B(\xi\eta_j^*)\eta_j,
\forall \xi \in \Cal H$).

A set $\{\eta_j\}_j \subset L^2(N, \tau)$ is an
orthonormal basis for $\Cal H_B$ if and only if the
orthogonal projection $f$ of $L^2(N, \tau)$ on $\Cal H$ satisfies
$f = \Sigma_j \eta_j e_B \eta_j^*$ with $\eta_je_B\eta_j^*$ projection
$\forall j$. A simple maximality argument shows that
any left (resp. right) Hilbert $B$-module $\Cal H \subset L^2(N, \tau)$ has an
orthonormal basis (see [Po2] for all this).
The Hilbert module $\Cal H_B$ (resp. $_B\Cal H$) is
{\it finitely generated} if it has a finite orthonormal basis.
\vskip .1in
\noindent
{\it 1.4.2. Quasi-regular subalgebras.}
Recall from ([D]) that if $B\subset N$ is an inclusion
of finite von Neumann algebras
then the {\it normalizer}
of $B$ in $N$ is the set
$\Cal N(B)=\Cal N(B)=\{u\in \Cal U(N)\mid
uBu^* = B \}$. The von Neumann algebra $B$ is called
{\it regular} in $N$ if $\Cal N(B)''= N$.

In the same spirit, the {\it quasi-normalizer}
of $B$ in $N$ is defined to be the set
$q\Cal N(B)\overset\text{def}\to=\{x\in
N\mid\exists\ x_1,x_2,\ldots,x_n\in N$ such that $xB\subset\sum_{i=1}^nBx_i$
and $Bx\subset\sum_{i=1}^nx_iB\}$ (cf. [Po5], [PoSh]).
The condition ``$xB\subset\sum Bx_i$, $Bx\subset\sum x_iB$'' is
equivalent to ``$BxB\subset(\sum_{i=1}^nBx_i)\cap(\sum_{i=1}^nx_iB)$'' and also
to ``sp$BxB$ is finitely generated both as left and as a right $B$-module.'' It
then follows readily that sp$(q\Cal N_N(B))$ is a $^*$-algebra. Thus,
$P\overset\text{def}\to=\overline{\text{sp}}(q\Cal N_N(B))=q\Cal N_N(B)''$ is a
von Neumann subalgebra of $N$ containing $B$.
In case the von Neumann algebra $P=q\Cal N_B(N)''$
is equal to all $N$, then we say that $B$ is {\it quasi-regular} in $N$
([Po5]).

The most interesting case of inclusions $B \subset N$
for which one considers the normalizer $\Cal N(B)$ and
the quasi-normalizer $q\Cal N_N(B)$ of $B$ in $N$ is
when the subalgebra $B$  satisfies the condition $B'\cap N \subset B$,
or equivalently $B'\cap N = \Cal Z(B)$, notably when $B$ and
$N$ are factors
(i.e., when $B'\cap N = \Bbb C$) and
when $B$ is a maximal abelian $^*$-subalgebra
(i.e., when $B'\cap N = B$).

The next lemma lists some useful properties
of $q\Cal N(B)$. In particular,
it shows that if a Hilbert $B$-bimodule $\Cal H \subset L^2(N, \tau)$
is finitely generated
both as a left and as a right Hilbert $B$ module, then it is ``close'' to
a bounded finitely
generated $B$-bimodule $H \subset P$.

\proclaim{Lemma} $(i)$. Let $N$ be a finite von Neumann
algebra with a normal finite faithful trace $\tau$
and $B\subset N$ a von Neumann subalgebra.
Let $p \in B'\cap \langle N, B \rangle$
be a finite projection such that $J_NpJ_N$ is also a finite projection.
Let $\Cal H\subset L^2(N, \tau)$ be the Hibert space on which $p$ projects
(which is thus a Hilbert $B$-bimodule).
Then there exists an increasing
sequence of central projections
$z_n \in \Cal Z(B)$ such that $z_n \uparrow 1$
and such that the Hilbert $B$-bimodules $z_n\Cal Hz_n \subset L^2(N)$
are finitely generated
both as left and as right Hilbert $B$-modules.

$(ii)$. If $B\subset N$ are as in $(i)$ and
$\Cal H^0\subset L^2(N)$
is a Hilbert $B$-bimodule such that $\Cal H_B^0$, $_B\Cal H^0$ are
finitely generated Hilbert modules, with $\{\xi_i \mid 1\leq i \leq n\},
\{\zeta_j \mid 1\leq j \leq m\}$ their corresponding
ortonormal basis, then for any $\varepsilon > 0$ there
exists a projection $q \in B'\cap N$ such that $\tau(1-q) < \varepsilon$ and
$x_i = q\xi_iq\in N , y_j = q\zeta_jq \in N, \forall i,j$, with the
orthogonal projection $p_0$ of $L^2(N)$ onto the closure
of $q\Cal H^0q$ in $L^2(N)$
acting on $N = \hat{N}\subset L^2(N)$ by $p_0(x) =
\Sigma_i x_i E_B(x_i^*x) = \Sigma_j
E_B(xy_j^*)y_j, \forall x \in N$. In particular, $\Sigma_i x_i B
= \Sigma_j By_j = q \Cal H^0 q \cap N$ is dense in $q \Cal H^0 q$
and is finitely generated both as left and right $B$-module.

$(iii)$. If $p$ is a projection as in $(i)$ then
$p \leq e_P$. Also, $B$ is quasiregular in $N$ if and
only if $B$ is discrete in $N$, i.e.,
$B'\cap \langle N,B \rangle$ is generated by projections which are
finite in $\langle N, B \rangle$ $({\text{\rm [ILP]}})$.
\endproclaim
\noindent
{\it Proof}. $(i)$ and
$(ii)$ are trivial consequences of 1.4.1 and of
the definitions.

The first part of $(iii)$ is trivial by $(i), (ii)$. Thus, $e_P$ is
the supremum of all projections $p \in B'\cap \langle N,B \rangle$ such that
both $p$ and $J_NpJ_N$ are finite in $\langle N,B \rangle$. Thus,
if $q\in \langle N,B \rangle$ is a non-zero finite projection
orthogonal to $e_P$ then
any projection $q' \in B'\cap \langle N,B \rangle$
with $q' \leq J_N q J_N$
must be infinite (or else the maximality of $e_P$ would
be contradicted). But if $q$ satsifies this property then
$B'\cap \langle N,B \rangle$ cannot be generated by finite projections.
\hfill Q.E.D.
\vskip .1in
\noindent
{\it 1.4.3. Cartan subalgebras}. Recall from ([D]) that a
maximal abelian $^*$-subalgebra $A$ of a finite von Neumann factor
$M$ is called
{\it semiregular} if $\Cal N(A)$ generates a factor, equivalently, if
$\Cal N(A)'\cap M = \Bbb C$. Also, while
maximal abelian $*$-subalgebras $A$ with $\Cal N(A)''=M$ were called
{\it regular} in ([D]), as mentioned before,
they were later called
{\it Cartan subalgebras} in ([FM]), a terminology that seems to
prevail and which we therefore adopt.

By results of Feldman and Moore
([FM]), in case a type II$_1$ factor $M$
is {\it separable} in the norm $\|\, \|_2$ given by the trace,
to each Cartan subalgebra $A\subset M$
corresponds a countable, measure preserving, ergodic equivalence
relation $\Cal R=\Cal R(A\subset M)$ on the standard probability space
$(X, \mu)$, where $L^\infty(X, \mu) \simeq (A, \tau_{|A})$,
given by orbit equivalence under the action
of $\Cal N(A)$. In fact, $\Cal N(A)$ also gives rise to an
$A$-valued 2-cocycle $v=v(A\subset M)$, reflecting the
associativity mod $A$ of the product
of elements in the normalizing groupoid $\Cal G\Cal N
\overset \text{\rm def} \to =\{pu \mid u \in \Cal N(A), p \in \Cal P(A)\}$.

Conversely, given any pair $(\Cal R, v)$,
consisting of a countable, measure preserving, ergodic equivalence
relation $\Cal R$ on the standard probability space $(X, \mu)$ and
a $L^\infty(X, \mu)$-valued 2-cocycle $v$ for the
corresponding groupoid action (N.B.: $v \equiv 1$ is always
a $2$-cocycle, $\forall \Cal R$), there exists a
type II$_1$ factor with a Cartan subalgebra
$(A\subset M)$ associated with it, via a group-measure
space construction ``\`a la''
Murray-von Neumann.
The association $(A\subset M) \rightarrow (\Cal R, v)
\rightarrow (A\subset M)$ is one to one, modulo
isomorphisms of inclusions $(A\subset M)$ and respectively
measure preserving orbit equivalence of $\Cal R$ with
equivalence of the 2-cocycles $v$ (see [FM] for all this).

Examples of countable, measure preserving, ergodic equivalence
relation $\Cal R$ are obtained by taking free
ergodic measure preserving actions
$\sigma$ of countable groups $\Gamma_0$ on the standard probability
space $(X, \mu)$,
and letting $x \Cal R y$
whenever there exists $g\in \Gamma_0$ such that $y = \sigma_g(x)$.

If $t > 0$ then the {\it amplification of a Cartan subalgebra}
$A \subset M$ by $t$ is the Cartan subalgebra $A^t\subset M^t$
obtained by first choosing
some $n \geq t$ and then
compressing the Cartan subalgebra $A\otimes D
\subset M \otimes M_{n\times n}(\Bbb C)$ by a projection $p\in A \otimes D$
of (normalized) trace equal to $t/n$. (N.B.
This Cartan subalgebra is defined up
to isomorphism.) Also, the {\it amplification of a measurable
equivalence relation} $\Cal R$ by $t$ is the equivalence relation
obtained by reducing the equivalence relation $\Cal R \times \Cal D_n$
to a subset of measure $t/n$, where $\Cal D_n$ is the ergodic
equivalence relation on the $n$ points set. Note that if $A \subset M$
induces the equivalence relation $\Cal R$ then $A^t \subset M^t$ induces
the equivalence realtion $\Cal R^t$. Also, $v_{A\subset M} \equiv 1$
implies $v_{A^t\subset M^t}\equiv 1, \forall t>0$.

By using Lemma 1.4.2, we can reformulate a result from
([PoSh]), based on prior results in ([FM]),
in a form that will be more suitable for us:

\proclaim{Proposition} Let $M$ be a separable type ${\text{\rm II}}_1$ factor.

$(i)$. A maximal abelian
$^*$-subalgebra $A\subset M$ is a Cartan subalgebra
if and only if $A \subset M$ is discrete, i.e., iff
$A'\cap \langle M, A \rangle$
is generated by projections that are finite in $\langle M,A \rangle$.

$(ii)$. Let $A_1, A_2 \subset M$ be two Cartan subalgebras of
$M$. Then $A_1, A_2$ are conjugate by a unitary element of $M$ if and only if
$A_1'\cap \langle M, A_2 \rangle$ is generated by finite projections of
$\langle M, A_2 \rangle$ and $A_2'\cap \langle M, A_1 \rangle$ is
generated by finite projections of
$\langle M, A_1 \rangle$. Equivalently,
$A_1, A_2$ are unitary conjugate if and
only if $_{A_1}L^2(M, \tau)_{A_2}$ is a direct sum $A_1-A_2$
Hilbert bimodules that are finite dimensional both as left $A_1$-Hilbert modules
and as right $A_2$-Hilbert modules.
\endproclaim
\noindent
{\it Proof}. $(i)$. By Lemma 1.4.2, the discreteness
condition on $A$ is equivalent to the
quasi-regularity of $A$ in $N$. By ([PoSh]), the latter is equivalent to
$A$ being Cartan.

$(ii)$. If $A_i'\cap (J_N A_j J_N)'$ is generated by finite projections
of the semifinite von Neumann algebra $(J_N A_j J_N)'$, for $i,j = 1,2$,
and we denote $M = M_2(N)$ the algebra
of $2$ by $2$ matrices over $N$ and $A = A_1 \oplus A_2$
then $A'\cap (J_M A J_M)'$ is also
generated by finite projections of $J_M A J_M$. By part $(i)$, this
implies $A$ is Cartan in $M$. By ([Dy]) this implies there
exists a partial isometry $v \in M$ such that $vv^*=e_{11},
v^*v = e_{22}$, where $\{e_{ij}\}_{i,j= 1,2}$,
is a system of
matrix units for $M_2(\Bbb C)$. Thus, if $u \in N$ is the unitary element
with $ue_{12} = v$ then $uA_1u^* = A_2$.
\hfill Q.E.D.

\heading 2. Relative Property H: Definition and Examples. \endheading

In this Section we consider a ``co-type'' relative version of Haagerup's
compact approximation property for inclusions of von Neumann algebras.
This property can be viewed as a ``weak co-amenability'' property,
as we will comment on in the next Section (see 3.5, 3.6).
It is a property that excludes ``co-rigidity'', as later
explained (see 5.6, 5.7). We first recall the definition
for groups and for single von Neumann algebras, for completeness.
\vskip .1in
\noindent
{\it  2.0.1. Property H for groups.}
In ([H1]) Haagerup proved that the free groups
$\Gamma_0 = \Bbb F_n, 2\leq n \leq \infty,$
satisfy the following condition:
There exist positive definite functions $\varphi_n$ on $\Gamma_0$ such
that
\vskip .1in
\noindent
$(2.0.1')$. $\quad \quad \underset g\rightarrow \infty \to \lim
\varphi_n(g) = 0, \forall n,$ (equivalently, $\varphi_n \in c_0(\Gamma_0)$).
\vskip .1in
\noindent
$(2.0.1'')$. $\quad \quad
\underset n\rightarrow \infty \to \lim
\varphi_n(g) = 1, \forall g \in \Gamma_0 $.
\vskip .1in
Many more groups $\Gamma_0$
were shown to satisfy  conditions (2.0.1) in ([dCaH],
[CowH], [CCJJV]).
This property is often refered to as {\it Haagerup's approximation
property}, or {\it property} H (see e.g., [Cho], [CJ],
[CCJJV]).
By a result of Gromov, a group has property H if and
only if it satisfies a certain embeddability
condition into a Hilbert space, a property he called
a-T{\it -menability} ([Gr]). There has been a lot of
interest in studying these groups lately. We refer the reader to the recent
book ([CCJJV]) for a comprehensive account on this subject. Note that property H
is a hereditary property, so if a group $\Gamma_0$ has it, then any subgroup
$\Gamma_1 \subset \Gamma_0$ has it as well.
\vskip .1in
\noindent
{\it 2.0.2. Property} H {\it for algebras}. A similar {\it property} H,
has been considered for finite von Neumann
algebras $N$ ([C3], [Cho], [CJ]): It
requires the existence of a net of normal completely positive
maps $\phi_\alpha$ on $N$ satisfying the conditions:
\vskip .1in
\noindent
$(2.0.2')$. $\quad  \tau \circ \phi_\alpha
\leq \tau$ and $\phi_\alpha(\{x\in N \mid \|x\|_2 \leq 1\})$
is $\|\quad \|_2$-precompact, $\forall \alpha$;
\vskip .1in
\noindent
$(2.0.2'')$. $\quad  \underset \alpha \rightarrow \infty \to \lim
\|\phi_\alpha(x)-x\|_2 = 0, \forall x\in N$;
\vskip .1in
\noindent
with respect to some fixed normal faithful trace $\tau$ on $N$.
The net can of course be taken to be a sequence in case $N$ is separable in the
$\|\quad \|_2$-topology.
\vskip .1in

It was shown in ([Cho]) that if $N$ is
the group von Neumann algebra $L(\Gamma_0)$ associated to
some group $\Gamma_0$, then $L(\Gamma_0)$ has the property H
(as a von Neumann algebra)
if and only if $\Gamma_0$ has the property H (as a group). It was further
shown in ([Jo1])
that the set of properties $(2.0.2)$  does not depend on the
normal faithful trace $\tau$ on $N$, i.e., if there exists a net
of completely positive maps $\phi_\alpha$ on $N$ satisfying conditions
$(2.0.2'), (2.0.2'')$ with respect to some faithful normal trace $\tau$, then
given any other faithful normal trace $\tau'$ on $N$ there exists
a net of completely positive maps $\phi_\alpha'$
on $N$ satisfying the conditions with
respect to $\tau'$. It was also proved in
([Jo1]) that if $N$ has property H then given any
faithful normal trace $\tau $ on $N$ the completely
positive maps $\phi_\alpha$ on $N$ satisfying
$(2.0.2)$ with respect to $\tau$ can be taken $\tau$-preserving and
unital.

We now extend the definition of the property H
from the above single algebra case to the
relative (``co-type'') case of
inclusions of von Neumann
algebras, by using a similar strategy to the way
the notions of amenabilty and property (T) were extended from
single algebras to inclusions of algebras
in ([Po1,10]; see Remarks 3.5, 3.6, 5.6 hereafter).

\vskip .1in
\noindent
{\it 2.1. Definition}.
Let $N$ be a finite von Neumann algebra with countable decomposable
center and
$B\subset N$ a von Neumann subalgebra.
$N$ has the {\it property} H {\it relative to} $B$ if there exists
a normal faithful tracial
state $\tau$ on $N$ and a
net of normal completely positive $B$-bimodular maps
$\phi_\alpha$ on $N$ satisfying the conditions:
\vskip .1in
\noindent
$(2.1.0.)$ $\quad\quad \tau\circ \phi_\alpha \leq \tau;$
\vskip .1in
\noindent
$(2.1.1)$. $\quad\quad
T_{\phi_\alpha} \in J(\langle N, B \rangle), \forall \alpha;$
\vskip .1in
\noindent
$(2.1.2)$. $\quad\quad \underset \alpha \rightarrow \infty \to \lim
\|\phi_\alpha(x)-x\|_2 = 0, \forall x\in N$,
\vskip .1in
\noindent
where $T_{\phi_\alpha}$ are the operators in the semifinite
von Neumann algebra
$\langle N, B \rangle
\subset \Cal B(L^2(N, \tau))$ defined out of $\phi_\alpha$
and $\tau$, as in 2.1.

Following ([Gr]), one can also use the
terminology: $N$ is a-T-{\it menable relative to} $B$.

Note that the finite von Neumann algebra $N$ has the property H
as a single von Neumann algebra if and only if $N$ has the
property H relative to $B = \Bbb C$.

Note that a similar notion of
``relative Haagerup property'' was considered
by Boca in ([Bo]), to study the behaviour
of the Haagerup property under amalgamated free products.
The definition in ([Bo]) involved a fixed trace and it
required the completely positive maps
to be unital and trace preserving. The next proposition
addresses some of the differences between his definition and 2.1:

\proclaim{2.2. Proposition} Let $N$ be a finite
von Neumann algebra with countably decomposable center and
$B\subset N$ a von Neumann subalgebra.

$1^\circ$. If $N$ has the property ${\text{\rm H}}$ relative to $B$
and $\{\phi_\alpha\}_\alpha$ satisfy $(2.1.0)-(2.1.2)$
with respect to the trace $\tau$ on $N$, then there exists a
net of completely positive maps $\{\phi'_\alpha\}_\alpha$ on $N$,
which still satisfy $(2.1.0)-(2.1.2)$
with respect to the trace $\tau$, but also
$T_{\phi'_\alpha}\in J_0(\langle N, B \rangle)$ and
$\phi'_\alpha(1) \leq 1, \forall \alpha.$

$2^\circ$. Assume $B'\cap N \subset B$.
Then the following conditions are equivalent:

$(i)$. $N$ has the property ${\text{\rm H}}$ relative to $B$.

$(ii)$. Given any faithful normal tracial state
$\tau_0$ on $N$, there exists a net of unital, $\tau_0$-preserving,
$B$-bimodular completely positive maps $\phi_\alpha$ on $N$ such that
$T_{\phi_\alpha}\in \Cal J_0(\langle N,B \rangle), \forall \alpha$,
and such that condition $(2.1.2)$ is satisfied for the
norm $\|\quad\|_2$ given by $\tau_0$.

$(iii)$. There exists a normal faithful tracial state $\tau$
and a net of normal, $B$-bimodular completely
positive maps $\phi_\alpha$ on $N$ such that $\phi_\alpha$ can be extended
to bounded operators $T_{\phi_\alpha}$ on $L^2(N, \tau)$,
such that $T_{\phi_\alpha} \in \Cal J(\langle N, B \rangle)$ and
$(2.1.2)$ is satisfied for the trace $\tau$.

Moreover, in case $N$ is
countably generated as a $B$-module, i.e., there exists a countable set
$S\subset N$ such that $\overline{\text{\rm sp}} SB=N$, the closure
being taken in the norm $\| \quad \|_2$, then the net $\phi_\alpha$
in either $1^\circ$, $2^\circ$ or $3^\circ$ can be taken
to be a sequence.
\endproclaim
\noindent
{\it Proof}. $1^\circ$. By part 3$^\circ$
of Proposition 1.3.3, we can replace if necessary $\phi_\alpha$ by
$\phi_\alpha (z_\alpha \cdot z_\alpha)$,
for some $z_\alpha \in \Cal P(\Cal Z(B))$
with $z_\alpha \uparrow 1$, so that
the corresponding operators on $L^2(N, \tau)$ belong to
$\Cal J_0(\langle N,B \rangle), \forall \alpha$.

By using continuous functional calculus for
$\phi_\alpha(1)$, let $b_\alpha = (1 \vee \phi_\alpha(1))^{-1/2}
\in B'\cap N$.
Then $b_\alpha \leq 1, \|b_\alpha - 1\|_2 \rightarrow 0$ and
$$
\phi'_\alpha (x) = b_\alpha \phi_\alpha(x) b_\alpha, x\in N,
$$
still defines a normal completely positive map on $N$ with
$\|\phi'_\alpha(x) - x\|_2 \rightarrow 0, \forall x\in N.$
Moreover, if $x \geq 0$ then
$$\tau(\phi'_\alpha(x))
=\tau(\phi_\alpha(x) b_\alpha^2) \leq \tau(\phi_\alpha(x)).
$$

Also, since $T_{\phi'_\alpha} = L(b_\alpha) R(b_\alpha) T_{\phi_\alpha}$
and $L(b_\alpha) \in N \subset \langle N,B \rangle$,
$R(b_\alpha) \in J(B'\cap N)J \subset \langle N, B\rangle$ and
$T_{\phi_\alpha} \in \Cal J(\langle N,B \rangle)$,
it follows that $T_{\phi'_\alpha} \in \Cal J(\langle N,B \rangle)$.

$2^\circ$. We clearly have $(ii) \implies (i) \implies (iii)$.

Assume now $(iii)$ holds true for the trace $\tau$ and let $\tau_0$
be an arbitrary normal, faithful tracial state on $N$. Thus,
$\tau_0 = \tau (\cdot a_0)$, for some $a_0 \in \Cal Z(N)_+$
with $\tau(a_0)=1$.
Since $B'\cap N = \Cal Z(B)$, by part $3^\circ$ of
Lemma 1.2.1 we have
$a_\alpha = \phi_\alpha(1) \in \Cal Z(B)$. Also, $(2.1.2)$ implies
$$
\underset \alpha \rightarrow \infty \to \lim
\|a_\alpha-1\|_2 = 0, \tag 2.2.2'
$$
where $\|\quad \|_2$ denotes the norm given by $\tau$.

Let $p_\alpha$ be the spectral
projection of $a_\alpha$ corresponding to $[1/2, \infty)$.
Since $a_\alpha \in
\Cal Z(B)$, $p_\alpha \in \Cal Z(B)$. Also,
condition $(2.2.2')$ implies
$\underset \alpha \rightarrow \infty \to \lim
\|p_\alpha-1\|_2=
\underset \alpha \rightarrow \infty \to \lim
\|a_\alpha^{-1}p_\alpha-p_\alpha\|_2 = 0$. Furthermore,
by condition $3^\circ$ of
Proposition 1.3.3, there
exists $p_\alpha' \in \Cal Z(\Cal B)$.
with $p_\alpha' \leq p_\alpha$, such that
$$
\underset \alpha \rightarrow \infty \to \lim
\|p'_\alpha-1\|_2=0, \underset \alpha \rightarrow \infty \to \lim
\|a_\alpha^{-1}p_\alpha-p_\alpha'\|_2=0. \tag 2.2.2''
$$

Define
$\phi'_\alpha$ on $N$ by
$$
\phi'_\alpha(x) = a_\alpha^{-1/2}p'_\alpha\phi_\alpha (x)
p'_\alpha
a_\alpha^{-1/2} + (1-p'_\alpha)E_B(x)(1-p'_\alpha), x\in N.
$$
Then we clearly have $\phi'_\alpha(1) = 1$, $\phi'_\alpha$
are $B$-bimodular and $T_{\phi'_\alpha}
\in \Cal J_0(\langle N, B \rangle)$.
Since $B'\cap N \subset B$, by part $2^\circ$ in
Lemma 1.2.1, this also implies $\tau \circ \phi'_\alpha
= \tau, \tau_0 \circ \phi'_\alpha
= \tau_0$. Moreover,
since $a_\alpha^{-1}p_\alpha \leq
2$, it follows that for each $x\in N$ we have:
$$
\|\phi'_\alpha(x)-x\|_2
\leq
\|a_\alpha^{-1/2}p'_\alpha \phi_\alpha(x)a_\alpha^{-1/2}p'_\alpha -
p'_\alpha xp'_\alpha\|_2
$$
$$
+\|(1-p'_\alpha)x p'_\alpha\|_2+\|p'_\alpha x(1-p'_\alpha)\|_2
+\|(1-p'_\alpha)(x-E_B(x))(1-p'_\alpha)\|_2
$$
$$
\leq 2\|\phi_\alpha(x) - x\|_2 + 2\|a_\alpha^{-1/2}p'_\alpha x
a_\alpha^{-1/2}p'_\alpha - p'_\alpha x p'_\alpha\|_2
+3\|1-p'_\alpha\|_2 \|x\|,
$$
with the latter tending to $0$
for all $x\in N$, by $(2.2.2'')$. Since this
convergence holds true for
one faithful normal trace, it holds true in the $s$-topology,
thus for the normal trace $\tau_0$ as well.

The last part of $2^\circ$ is trivial.
\hfill Q.E.D.

We now prove some basic properties of the relative
property H, showing that it is well behaved to simple operations
such as tensor products, amplifications,
finite index extensions/restrictions.

\proclaim{2.3. Proposition} $1^\circ$. If $N$ has the property
${\text{\rm H}}$
relative to $B$ and
$B_0 \subset N_0$ is embedded into $B \subset N$
with commuting squares, i.e., $N_0 \subset N, B_0\subset B,
B_0 = N_0 \cap B$ and
$E_{N_0}\circ E_B = E_B\circ E_{N_0} = E_{B_0}$, then
$N_0$ has the property ${\text{\rm H}}$
relative to $B_0$.

$2^\circ$. If $B_1\subset N_1$ and $B_2\subset N_2$ then
$N_1\overline\otimes N_2$ has the property ${\text{\rm H}}$
relative to $B_1\overline\otimes B_2$ if and only if $N_i$ has the
property ${\text{\rm H}}$ relative to $B_i, i=1, 2$.

$3^\circ$. Let $B\subset N_0 \subset N$. If $N$ has the
property ${\text{\rm H}}$ relative to $B$, then
$N_0$ has the
property ${\text{\rm H}}$ relative to $B$.
Conversely, if we assume $N_0\subset N$
has a finite orthonormal basis $\{u_j\}_j$
with $u_j$ unitary elements such that $u_jBu_j^* = B, \forall j$,
and $N_0$ has the property ${\text{\rm H}}$ relative
to $B$, with respect to $\tau_{|N_0}$ for some
normal faithful trace $\tau$ on $N$, then $N$ has the property ${\text{\rm H}}$ relative to $B$, with respect to $\tau$.

$4^\circ$. Assume $B\subset B_0\subset  N$ and $B\subset B_0$ has a
finite orthonormal basis. If $N$ has the property
${\text{\rm H}}$ relative to $B_0$ then
$N$ has property ${\text{\rm H}}$ relative to $B$. If
in addition $B_0'\cap N \subset B_0$ then, conversely, if
$N$ has the property
${\text{\rm H}}$ relative to $B$, then
$N$ has property ${\text{\rm H}}$ relative to $B_0$.
\endproclaim
\noindent
{\it Proof}. $1^\circ$. If $\phi_\alpha : N \rightarrow N$ are $B$-bimodular
completely positive maps approximating the identity on $N$, then
by the commuting squares relation
$E_{N_0}\circ E_B = E_B\circ E_{N_0} = E_{B_0}$, it follows that $\phi'_\alpha
= E_{N_0}\circ\phi_{\alpha \mid N_0}$ approximate the identity on $N_0$ and
are $B_0$-bimodular. Moreover, by commuting squares,
if $T_{\phi_\alpha}$ satisfy condition $5^\circ$ in 1.3.3 then
so do $T_{\phi'_\alpha}$.

$2^\circ$. The implication from left to right follows
by applying $1^\circ$ to $(B\subset N) =
(B_1\overline{\otimes} B_2 \subset N_1 \overline{\otimes} N_2)$
and $(B_0 \subset N_0) = (B_i \otimes \Bbb C \subset N_i \otimes \Bbb C),
i= 1,2$. The implication from right to left follows from the fact that
$T_{\phi^i_\alpha} \in \Cal J(\langle N_i, B_i \rangle), i=1,2,$ implies
$T_{\phi^1_\alpha \otimes \phi^2_\alpha} \in \Cal J(\langle N_1 \overline{\otimes}
N_2, B_1 \overline{\otimes} B_2)$ (by using that the tensor product of finite
projections is a finite projection).

$3^\circ$. For the first implication,
let $\phi_\alpha$ be
completely positive maps on $N$ that satisfy
$(2.1.0)-(2.1.2)$ for $B \subset N$ and for the trace $\tau$ on $N$.
Define $\phi^0_\alpha(x) = E_{N_0}(\phi_\alpha(x)), x\in N_0$. Then
$\phi^0_\alpha$ are completely positive, $B$-bimodular maps
which still satisfy $\tau \circ \phi^0_\alpha \leq \tau$.
Moreover, since $T_{\phi_\alpha}$ satisfy condition
$5^\circ$ in Proposition 1.3.3,
then clearly $\phi^0_\alpha$ do as well.

For the
converse, assuming $\phi^0_\alpha$ are completely
positive maps on $N_0$ that satisfy
$(2.1.0)-(2.1.2)$ for $B \subset N_0$, define $\tilde{\phi}_\alpha$
on $\langle N, e_{N_0} \rangle$ by
$$
\tilde{\phi}_\alpha
(\Sigma_{i,j} u_i x_{ij}e_{N_0}u_j^*)=
\Sigma_{i,j} u_i \phi^0_\alpha(x_{ij})e_{N_0}u_j^*,
$$
where $x_{ij}\in N_0$. It
is then immediate to check that $\tilde{\phi}_\alpha$ are
completely positive, $B$-bimodular and
check (2.1.0)-(2.1.2) with respect to the
canonical trace $\tilde{\tau}$
on $\langle N, e_{N_0} \rangle$ implemented
by the trace $\tau$ on $N$ (which is clearly Markov by hypothesis).
Thus, $\langle N, e_{N_0} \rangle$ has
property H relative to $B$, so by the first part $N$ has
property H relative to $B$ as well (with respect to
$\tilde{\tau}_{|N} = \tau$.

$4^\circ$. For the first implication, note that
the condition that $B_0$ has a finite orthonormal basis
over $B$ implies $\Cal J_0(\langle N,B_0 \rangle)
\subset \Cal J_0(\langle N,B \rangle)$. Indeed, this
follows by first approximating $T\in
\Cal J_0(\langle N,B_0 \rangle)$ by linear
combination of projections in $J_0(\langle N,B_0 \rangle)$ then
noticing that if
dim$(_{B_0}\Cal H)<\infty$
(respectively, dim$(\Cal H_{B_0})<\infty$), then
dim$(_{B}\Cal H)<\infty$
(respectively, dim$(\Cal H_{B})<\infty$).

For the opposite implication, let $\{m'_j\}_j$ be a finite
orthonormal basis of $B_0$ over $B$ and
recall from ([Po2]) that $b=\Sigma_j m'_j {m'_j}^* \in \Cal Z(B_0)$
and $b \geq 1$. Also, since for any $T \in B'\cap \langle N,B\rangle$
we have
$$\Sigma_{i,j} L(m'_j) R({m'_i}^*) \circ T \circ L({m'_j}^*) R(m'_i)
\in B_0'\cap \langle N,B_0 \rangle
$$
(cf. [Po2]), it
follows that if we put $m_j = b^{-1/2}m'_j$ then we have
$$
T^0=\Sigma_{i,j} L(m_j) R(m_i^*) \circ T \circ L(m_j^*) R(m_i)
\in B_0'\cap \langle N,B_0 \rangle.
$$
This shows that if we put $\phi^0_\alpha=
\Sigma_{i,j} m_j \phi_\alpha(m_j^* \cdot m_i)m_i^*$, then
$T^0=T_{\phi^0_\alpha} \in B_0'\cap \langle N,B_0 \rangle$.
Also, if in the above we take $T$ to be a projection
with the property
that $\Cal H=e(L^2(N,\tau))$ is a finitely
generated left-right Hilbert $B$-module, then the support
projection of the corresponding operator $T^0$
is contained in $\Cal H^0=\overline{\Sigma_{i,j} m_i \Cal H m_j^*}$.
To prove that $T^0$ is contained
in $\Cal J_0(\langle N,B_0 \rangle)$ it is sufficient to
show that $\Cal H^0$ is a finitely
generated left-right Hilbert $B_0$-bimodule.

To do this, write first $\Cal H$ as the closure of a finite sum
$\Sigma_k \eta_k B$. Then
$\Cal H^0$ follows the closure of
$$
\Sigma_{i,j} m_i (\Sigma_k \eta_k B) m_j^* = \Sigma_{i,k} (m_i \eta_k
(\Sigma_j B m_j^*) = \Sigma_{i,k} m_i \eta_k B_0.
$$
This shows that dim$_{B_0}\Cal H^0 < \infty$. Similarly,
dim$\Cal H^0_{B_0} < \infty$.

Taking linear combinations and
norm limits, we get that $T \in \Cal J_0(\langle N,B \rangle)$
implies $T^0 \in \Cal J_0(\langle N, B_0 \rangle)$.

Finally, since $\Sigma_j m_j m_j^*=1$, by Corollary 1.1.2 the
convergence to $id_N$ of $\phi_\alpha$ implies convergence to $id_N$
of $\phi^0_\alpha$.
By condition $(iii)$ in 2.3.2$^\circ$,
this implies $N$ has the property H relative to $B_0$.
\hfill Q.E.D.

\proclaim{2.4. Proposition}
$1^\circ$. If $N$ has the property ${\text{\rm H}}$
relative to $B$ and $p\in \Cal P(B)$ or $p\in \Cal P(B'\cap N)$,
then $pNp$ has the property ${\text{\rm H}}$ relative to $pBp$.

$2^\circ$. If $\{p_n\}_n \subset \Cal P(B)$ or
$\{p_n\}_n \subset \Cal P(B'\cap N)$ are such that
$p_n \uparrow 1$ and $p_nNp_n$ has the
property ${\text{\rm H}}$ relative to $p_nBp_n$, $\forall n$, then
$N$ has the property ${\text{\rm H}}$ relative to $B$.

$3^\circ$. Assume there exist partial isometries
$\{v_n\}_{n \geq 0}\subset N$ such that $v_n^*v_n \in pBp,$
$v_nv_n^* \in B,$ $v_nBv_n^* =
v_nv_n^* B v_nv_n^*, \forall n\geq 0, \Sigma_n v_nv_n^* = 1$ and
$B \subset (\{v_n\}_n \cup pBp)''$.  If $pNp$ has
property ${\text{\rm H}}$ relative to $pBp$ then
$N$ has property ${\text{\rm H}}$ relative to $B$.

$4^\circ$. If $B \subset N_0 \subset N_1 \subset ...$, then
$N = \overline{\cup_k N_k}$ has the property ${\text{\rm H}}$
relative to $B$ (with respect to a trace $\tau$ on $N$)
iff $N_k$ has the
property ${\text{\rm H}}$ relative to $B$
(with respect to $\tau_{|N_k}$), $\forall k$.
\endproclaim
\noindent
{\it Proof}. $1^\circ$. In both cases, if $\phi$ is
$B$-bimodular completely positive on $N$ then
$p\phi(p \cdot p)p$ is a $pBp$-bimodular
completely positive map
on $pNp$. Also, $\tau \circ \phi \leq \tau$ implies $\tau_p \circ
(p\phi(p \cdot p)p) \leq \tau_p$, where
$\tau_p(x) = \tau(x)/\tau(p), x \in pNp$.
Finally, if $T_{\phi}$ satisfies condition $5^\circ$ in
1.3.3 as an element in $\langle N, B \rangle$ then clearly
$T_{p\phi(p\cdot p)p}$ satisfies the condition as an
element in $\langle pNp, pBp\rangle$.

The case
$\{p_n\}_n \subset \Cal P(B'\cap N)$ of $2^\circ$ follows by noticing that
if $p\in \Cal P(B'\cap N)$
and $\phi_p$ is $Bp$-bimodular completely positive map
on $pNp$, with $\tau_p \circ \phi_p \leq \tau_p$, $\tau(1-p)\leq \delta$,
$\|\phi_p(x)-x\|_2
\leq \delta$,
$\forall x\in pFp$, for some finite set $F\subset N$, and
$T_{\phi_p} \in \Cal J_0(\langle pNp, Bp\rangle)$, then $\phi(y)
\overset\text{\rm def} \to = \phi_p(pyp) + E_B((1-p)y(1-p))$,
$\forall y\in N$ is $B$-bimodular and satisfies
$\tau \circ \phi \leq \tau$, $\|\phi(x)-x\|_2 \leq \varepsilon(\delta),
\forall x\in F$ and $T_{\phi} \in \Cal J_0(\langle N, B\rangle)$,
where $\underset \delta \rightarrow 0 \to \lim
\varepsilon(\delta) = 0$.

To prove $3^\circ$, let $\phi^p_\alpha$ be
$pBp$-bimodular, completely positive maps on
$pNp$ with $\tau_p \circ \phi^p_\alpha \leq \tau_p$, $T_{\phi^p_\alpha}
\in \Cal J_0(\langle pNp, pBp \rangle)$ and
$\phi^p_\alpha \rightarrow id_{pNp}$.
Define $\phi_\alpha$ on $N$ by
$$
\phi_\alpha(x) =
\Sigma_{i,j} v_i\phi^p_\alpha(v_i^*xv_j)v_j^*, x\in N.
$$

It is immediate to check that $\tau \circ \phi_\alpha \leq \tau$ and that
$\phi_\alpha \rightarrow id_N$. Also,
if $b \in pBp$ or $ b = v_i v_j^*$ then $ b\phi_\alpha(x) =
\phi_\alpha(bx), \phi_\alpha(x)b = \phi_\alpha (xb), \forall x\in N$.
Thus, if we denote by $B_1$ the von Neumann algebra
generated by $pBp$ and $\{v_n \}_n$ then $\phi_\alpha$
is $B_1$-bimodular.

Also, the same argument as in the last
part of the proof of 2.3.4$^\circ$ shows that $T_{\phi^p_\alpha} \in
\Cal J_0(\langle pNp, pBp \rangle)$ implies $T_{\phi_\alpha(p_n \cdot p_n)} \in
\Cal J_0(\langle p_nNp_n,p_nB_1p_n \rangle)$,
where $p_n = \Sigma_{0 \leq k \leq n} v_k^*v_k$. Thus, $p_nNp_n$
has property H relative to $p_nB_1p_n$. Since
$p_nBp_n \subset p_nB_1p_n$ and $p_nB_1p_n$ has finite orthonormal basis
over $p_nBp_n$, by 2.4.1$^\circ$ above and the
first implication in 2.3.4$^\circ$ it follows that
$p_nNp_n$ has property H relative to $p_nBp_n, \forall n$.

For each n let $\{z^n_k\}_{k}$ be a partition of the identity with
projections in $\Cal Z(B)$ such that $z_k^n$ has a finite partition
into projections in $B$ that are majorized by $p_nz^n_k$. Thus, there
exist finitely many partial isometries $v^n_0=p_nz^n_k, v^n_1, v^n_2, ...$
in $B$ such that ${v_i^n}^*v_i^n \geq {v_{i+1}^n}^*v_{i+1}^n, \forall i \geq 0$
and such that $\Sigma_i v^n_i{v^n_i}^* = z^n_k$. By the first part of the proof,
$z^n_k N z^n_k$ has property H
relative to $Bz^n_k$. By the case of $2^\circ$ that we have already
proved, it follows that $N$ has the property H relative to $B$.

The case $\{p_n\}_n \subset B$ in $2^\circ$ now follows by using 3$^\circ$,
to reduce the problem to the case $p_n$ are central in $B$ (as in the
proof of the last part of $3^\circ$).

$4^\circ$. The implication $\implies$ follows by
condition 2.3.3$^\circ$.
The reverse implication follows immediately once we note that if
$\phi$ is a completely positive map on $N_k$ such that
$\tau \circ \phi \leq \tau$ and $T_\phi \in \Cal J(\langle N_k, B \rangle)$,
then the completely positive map $\phi^k = \phi \circ E_{N_k}$ on $N$
satisfies $\tau \circ \phi^k \leq \tau$ and
$T_{\phi^k} \in \Cal J(\langle N,B \rangle)$
(for instance, by $5^\circ$ in 1.3.3).
\hfill Q.E.D.
\proclaim{2.5. Corollary} Let $A \subset M$ be a Cartan subalgebra
of the type ${\text{\rm II}}_1$ factor $M$. If $t > 0$ then
$M^t$ has the property ${\text{\rm H}}$ relative to $A^t$
if and only if $M$ has the property ${\text{\rm H}}$
relative to $A$ (see $1.4$ for the definition
of the amplification by $t$ of a Cartan subalgebra).
\endproclaim
\noindent {\it Proof}. Since the amplification by $1/t$ of $A^t
\subset M^t$ is $A \subset M$, it is sufficient to prove one of
the implications. Assume $M$ has the property H relative to $A$
and let $n \geq t$. By 2.3.2$^\circ$ it follows that $M \otimes
M_{n\times n}(\Bbb C)$ has property H relative to $A \otimes D_n$,
where $D_n$ is the diagonal algebra in $M_{n\times n}(\Bbb C)$. If
$p \in A \otimes D_n$ is a projection with $\tau(p) =t/n$ then, by
$2.4.1^\circ$, $M^t=p(M \otimes M_{n\times n}(\Bbb C))p$ has the
property H relative to $A^t=(A\otimes D_n)p$. \hfill Q.E.D. \vskip
.1in \noindent {\bf 2.6. Remark}. We do not know whether the
``smoothness'' condition $(2.1.0)$ on the $B$-bimodular,
completely positive, compact maps $\phi_n$ approximating the
identity on $N$ in Definition 2.1 can be removed. This is not
known even in the case $B=\Bbb C1_N$. In this respect, we mention
that in fact, for all later applications, the following weaker
``property H''-type condition will be sufficient: \vskip .1in
\noindent $(2.6.1)$. There exists a net of completely positive
$B$-bimodular maps $\phi_\alpha$ on $N$ which satisfy condition
$(2.2.2)$ and so that for all $\{u_n\}_n \subset \Cal U(N)$ with
$\underset n \rightarrow \infty \to \lim \|E_B(u_n^*u_m)\|_2 = 0,
\forall m$, we have $\underset n \rightarrow \infty \to \lim
\|\phi_\alpha(u_n)\|_2=0$. \vskip .1in

We do not know whether $(2.6.1)$
implies conditions $(2.1.0)-(2.1.2)$, not
even in the case $N$ is a factor and $B=\Bbb C1_N$.

We mention however that for type II$_1$
factors $N$ without the property $\Gamma$
of Murray and von Neumann
([MvN]), the smothness condition $(2.1.0)$ is automatically
satisfied, in case the completely positive map $\phi$ is sufficiently
close to the identity, thus making condition
(2.1.0) redundant. Indeed, we have the following observation, essentially
due to Connes and Jones ([CJ]):

\proclaim{2.7. Lemma} If $N$ is a non-$\Gamma$ type ${\text{\rm II}}_1$
factor then for any $\varepsilon > 0$
there exist $\delta > 0$ and a finite subset $F \subset \Cal U(N)$
such that the following conditions hold true:

$1^\circ$. If $\phi$ is a completely positive map
satisfying $\|\phi(u)-u\|_2 \leq \delta, \forall u\in F$, then
there exists a normal completely positive map $\phi''$ on $N$
such that $\phi''(1) \leq 1,
\tau \circ \phi'' \leq \tau$, $\|\tau \circ \phi'' - \tau\|
\leq \varepsilon$ and $\|\phi''(x)-x\|_2
\leq \|\phi(x) - x\|_2 + \varepsilon, \forall x\in N, \|x\| \leq 1$.
Moreover, if $\phi$ is $B$-bimodular for some $B \subset N$, then
$\phi''$ can be taken $B$-bimodular.

$2^\circ$. If $(\Cal H, \xi)$ is a $(B\subset N)$ Hilbert
bimodule with $\|u \xi - \xi u\| \leq \delta, \forall u\in F$
then $\|\langle \cdot \xi, \xi \rangle - \tau\| \leq \varepsilon$ ,
$\|\langle \xi \cdot, \xi \rangle - \tau\| \leq \varepsilon$.
\endproclaim
\noindent
{\it Proof}. $1^\circ$. Since $N$ is non-$\Gamma$,
by ([C2]) there exist unitary elements $u_1, u_2, ..., u_n$
in $N$ such  that if a state
$\varphi \in N^*$ satisfies $\|\varphi - \varphi(u\cdot u^*)\| \leq \delta$
then $\|\varphi - \tau\| \leq \varepsilon^2/9$.

Let $F = \{1\} \cup \{u_i\}_i$. Assume $\phi$ is a
completely positive map on $N$ such that $\|\phi(u)-u\|_2
\leq \delta^4/200, \forall u \in F$.
Let $a = 1 \vee \phi(1)$ and first define $\phi'$ on $N$
as in part $2^\circ$ of Lemma 1.1.2, i.e.,
$\phi'(x) = a^{-1/2} \phi (x) a^{-1/2}, x\in N$.
By 1.1.2, $\phi'(1) \leq 1$ and
$$
\|\phi'(x)-x\|_2
\leq \|\phi(x)-x\|_2 + 2 \|\phi(1) -1 \|_2^{1/2}\|x\|.
$$
Thus, by Corollary 1.1.2 we have for all $x\in N$ with $\|x\|\leq 1$
the estimates:
$$
\|\phi'(uxu^*) - u \phi'(x)u^* \|_2 \leq 2 (2\|\phi'(u)-u\|_2^2 +
2\|\phi'(u)-u\|_2)^{1/2} \leq \delta.
$$
Thus, if $\varphi = \tau \circ \phi'$ then
$\|\varphi - \varphi (u_i\cdot u_i^*)\| \leq \delta, \forall i$,
implying that $\|\varphi - \tau\| \leq \varepsilon^2/9$.

Thus, if we now take $\phi_1$ to be the normal part
of $\phi'$ then we still have $\phi_1(1) \leq 1,$
$\|\tau \circ \phi_1 - \tau \|
\leq \varepsilon^2/9$ and
$$\|\phi_1(x)-x\|_2 \leq
\|\phi(x)-x\|_2 + 2 \|\phi(1) -1 \|^{1/2}_2 \leq \|\phi(x) - x
\|_2 + \delta^2/6,
$$
for all $x\in N, \|x\|\leq 1$.  Finally, let $b_1 \in L^1(N, \tau)$
be the Radon-Nykodim derivative of $\tau \circ \phi_1$ with respect to
$\tau$ and define $b= 1 \vee b_1$,
$\phi'' = \phi_1( b^{-1/2}\cdot b^{-1/2})$,
as in Lemma 1.1.2. Thus, by part $3^\circ$ of that Lemma, all
the required conditions are satisfied.

$2^\circ$. This part is now trivial, by part $1^\circ$ above and 1.1.3.
\hfill Q.E.D.

\heading 3. More on property H. \endheading

In this Section we provide examples of inclusions of
finite von Neumann algebras with the
property H. We also prove that if a type II$_1$ factor $N$ has the
property H relative to a maximal abelian $*$-subalgebra $B$ then
$B$ is necessarily a Cartan subalgebra of $N$. Finally, we relate
relative property H with notions
of relative amenability considered in ([Po1,5]).

The examples we construct
arise from cross product constructions, being a consequence of the following
relation between groups and inclusions of algebras with property H:

\proclaim{3.1. Proposition} Let $\Gamma_0$ be a discrete
group and $(B, \tau_0)$ a finite von Neumann algebra
with a normal faithful tracial state. Let $\sigma$ be a cocycle
action of $\Gamma_0$ on $(B, \tau_0)$ by $\tau_0$-preserving
automorphisms. Then $N = B \rtimes_\sigma \Gamma_0$
has the property ${\text{\rm H}}$ relative to $B$
if and only if the group $\Gamma_0$ has the property ${\text{\rm H}}$.
\endproclaim
\noindent
{\it Proof}. First assume that $\Gamma_0$ has property H and let
$\varphi_\alpha: \Gamma_0 \rightarrow \Bbb C$ be unital
positive definite functions such that $\varphi_\alpha \in c_0(\Gamma_0)$
and $\varphi_\alpha(g) \rightarrow 1, \forall g\in \Gamma_0$. Also,
without loss of generality, we may assume $\varphi_\alpha(e)=1, \forall \alpha$.
For each $\alpha$, let $\phi_\alpha$ be the associated completely
positive map on $N=B \rtimes \Gamma_0$ defined as in Section 1.4, by
$\phi( \Sigma_g b_g u_g) = \Sigma_g \varphi(g)b_gu_g$. Note that
$\phi_\alpha$ are unital, trace preserving and $B$-bimodular (cf. 1.4).

Also, since $T_{\phi_\alpha}= \Sigma \varphi(g)
u_ge_Bu_g^*$, it follows that
$T_{\phi_\alpha} \in \Cal J(\langle N,B \rangle)$ if and only
if $\varphi_\alpha \in c_0(\Gamma_0)$. Finally, since
$|1 - \varphi_\alpha(g)| = \|\phi(u_g) - u_g\|_2,$ it follows that
$\underset \alpha \rightarrow \infty \to \lim
\varphi_\alpha(g) = 1, \forall g\in \Gamma_0,$ if and only
if $\underset \alpha \rightarrow \infty \to \lim
\|\phi_\alpha(x) - x \|_2 = 0, \forall x \in N$.

In particular, this shows that $N$ has property H
relative to $B$.

Conversely, assume $N$ has property H relative to
$B$ and let $\phi_\alpha : N \rightarrow N$ be a net of
completely positive maps satisfying $(2.3.0)-(2.3.2)$.
Let $\varphi_\alpha: \Gamma_0
\rightarrow \Bbb C$ be defined out
of $\phi_\alpha$, as in Section 1.4, i.e.,
by $\varphi_\alpha(g) =
\tau(\phi_\alpha(u_g)u_g^*), \forall g\in \Gamma_0$. By 2.4.5$^\circ$ we have
$$
\underset g
\rightarrow \infty \to \lim\|\phi_\alpha(u_g)\|_2 = 0, \forall \alpha.
$$

Thus, by using the Cauchy-Schwartz inequality  it follows that
$$
\underset g
\rightarrow \infty \to \lim \varphi_\alpha(g) = 0, \forall \alpha.
$$

Similarly, $\underset \alpha  \to \lim\|\phi_\alpha(u_g) - u_g\|_2=0$
implies $\underset \alpha  \to \lim\varphi_\alpha(g)=1$, thus
showing that $\Gamma_0$ has the property H.
\hfill Q.E.D.
\vskip .1in
\noindent
{\bf 3.2. Examples of groups with property H.}
The following groups $\Gamma_0$
(and thus, by hereditarity, any of their subgroups as well)
are known to have the property H,
thus giving rise to property H inclusions
$B \subset B \rtimes \Gamma_0$ whenever
acting (possibly with a cocycle)
on a finite von Neumann algebra $(B, \tau_0)$, by
trace preserving automorphisms, as in 3.1:
\vskip .1in
\noindent
{\it 3.2.0}. Any amenable group $\Gamma_0$ (cf. [BCV]; see also
3.5 below).
\vskip .1in
\noindent
{\it 3.2.1.} $G=\Bbb F_n,$
for some $2\leq n\leq \infty$, more generally $\Bbb F_S$, for $S$ an
arbitrary set of generators (cf. [H]).
\vskip .1in
\noindent
{\it 3.2.2.} $\Gamma_0$ a discrete subgroup of $SO(n, 1)$, for some $n\geq 2$
(cf. [dCaH]).
\vskip .1in
\noindent
{\it 3.2.3.} $\Gamma_0$ a discrete subgroup of $SU(n, 1)$, for some $n\geq 2$
(cf. [CowH]).
\vskip .1in
\noindent
{\it 3.2.4.} $SL(2, \Bbb Q)$, more generally $SL(2, \Bbb K)$ for any
field $\Bbb K\subset \Bbb R$ which is a finite extension over $\Bbb Q$
(by a result
of Jolissaint, Julg and Valette, cf. [CCJJV]).
\vskip .1in
\noindent
{\it 3.2.5.} $\Gamma_0=G_1 *_H G_2$, where $G_1, G_2$ have the property ${\text{\rm
H}}$ and $H\subset G_1, H\subset G_2$ is a common
finite subgroup (cf. [CCJJV]). In particular $\Gamma_0=SL(2, \Bbb Z)$.
\vskip .1in
\noindent
{\it 3.2.6.} $\Gamma= \Gamma_0 \times \Gamma_1$,
with $\Gamma_0, \Gamma_1$ property H groups. Also,
$\Gamma=\Gamma_0\rtimes_\gamma \Gamma_1$, with $\Gamma_0$ a
property H group
and $\Gamma_1$ an amenable group acting on it
(cf. [CCJJV]).
\vskip .1in
\noindent
We refer the reader to the book ([CCJJV]) for a more comprehensive list
of groups with the property H.
As pointed out in ([CCJJV]), the only
known examples of groups which do not have Haagerup property are the groups
$G_0$ containing infinite
subgroups $G \subset G_0$ such that $(G_0, G)$ has the relative property (T)
in the sense of ([Ma, dHVa]; see also next Section).
\vskip .1in
\noindent
{\bf 3.3. Examples of actions.} We are interested in constructing
examples of cocycle actions
$\sigma$ of (property H) groups $\Gamma_0$ on
finite von Neumann algebras $(B,\tau)$
(see e.g. [CJ] for the def. of cocycle actions)
that are {\it ergodic}
(i.e., $\sigma_g(b)=b, \forall g\in \Gamma_0$
implies $b\in \Bbb C1$) and {\it properly outer} (i.e., $\sigma_g(b)b_0 =
b_0b, \forall b\in B,$ implies $g=e$ or $b_0=0$).
Also, we consider the condition of {\it weak mixing},
which requires that $\forall F \subset B$ finite and
$\forall \varepsilon > 0$, $\exists g\in \Gamma_0$
such that $|\tau(\sigma_g(x)y)
- \tau(x) \tau(y)| \leq \varepsilon, \forall x, y \in F$. 
Weakly mixing actions are clearly ergodic.

Recall that the proper outernes of $\sigma$
is equivalent to
the condition $B'\cap B\rtimes_\sigma \Gamma_0
= \Cal Z(B)$. Also,
if $\sigma$ is a properly outer action,
then $\sigma$ acts ergodically on the center of $B$ if and only
if $B\rtimes_\sigma \Gamma_0$ is a factor. Finally, weak-mixing
is equivalent to the fact that $L^2(B, \tau)$ has no
$\sigma$-invariant finite dimensional subspaces other than $\Bbb C1$.

A yet another property of actions that we consider is
the following: the action $\sigma$
of $\Gamma_0$ on $(B,\tau)$ is {\it strongly ergodic} if
$B$ has no non-trivial
approximately $\sigma$-invariant sequences,
i.e., if $(b_n)_n \in \ell^\infty(\Bbb N, B)$ satisfies
$\underset n
\rightarrow \infty \to \lim \|\sigma_g(b_n)-b_n\|_2=0, \forall g\in \Gamma_0$
then $\underset n
\rightarrow \infty \to \lim \|b_n-\tau(b_n)1\|_2=0$. Note that if
we denote $N = B\rtimes_\sigma \Gamma_0$ and take $\omega$ to be a
free ultrafilter on $\Bbb N$, then this condition
is equivalent to $N'\cap B^\omega = \Bbb C$.

\vskip .1in
\noindent
{\it 3.3.1. Bernoulli shifts}. Given any
countable discrete group $\Gamma_0$ and any
finite von Neumann algebra $(B_0, \tau_0)$, $\Gamma_0$
acts on $(B, \tau) = (B,\tau)={\overline{\underset
g \in \Gamma_0 \to \otimes}} (B_0, \tau_0)_g$ by Bernoulli
shifts $\sigma_g$, namely $\sigma_g ( \otimes_h x_h )
= x'_h$, where $x'_h = x_{g^{-1}h}$.

If $B_0$ has no atoms or if $\Gamma_0$ is an infinite group,
then $\sigma$ is known to be
properly outer. Also, if $\Gamma_0$ is
infinite, then $\sigma$ is ergodic, in fact
even mixing. A Bernoulli shift action
is strongly ergodic if and only if $\Gamma_0$ is
non-amenable (cf. [J2]).
\vskip .1in
\noindent
{\it 3.3.2. Actions induced by automorphisms of groups}.
Let $\gamma$ be an action of
an infinite group $\Gamma_0$ on a
group $G$, by automorphisms. Let also
$\nu$ be a (normalized) scalar 2-cocycle on $\Gamma_0$ such that
$\nu_{\gamma(g),\gamma(h)}=\nu_{(g,h)}, \forall g,h\in \Gamma_0$.
Then $\gamma$ implements an action on the
``twisted'' group von Neumann algebra
$L_\nu(G)$, that we denote by $\sigma_\gamma$,
by taking $\sigma_\gamma(\lambda(g))=\lambda(\gamma(g)),
\forall g\in G$.
Note that $\sigma_\gamma$ preserves the canonical trace
$\tau$ of $L_\nu(G)$. Then we have:

\proclaim{Lemma} $(i)$. The following conditions are equivalent:
$(a)$. $\sigma_\gamma$ is ergodic; $(b)$. $\sigma_\gamma$ is weakly mixing;
$(c)$. $\gamma$ has no finite invariant subsets $\neq \{e\}$;
$(d)$. For any finite subset $S \subset G$ there exists $h \in \Gamma_0$
such that $\gamma_h(S) \cap S = \emptyset$.

$(ii)$. If $G_1 \subset G$ is
so that $\{ g_1^{-1} g_0 \gamma_h (g_1) \mid g_1 \in G_1\}$ is
infinite, $\forall h \in \Gamma_0 \setminus \{e\}$,
$\forall g_0 \in G$ then $L_\nu(G_1)' \cap L_\nu(G)\rtimes_\sigma \Gamma_0
\subset L_\nu(G)$. In particular, if
this holds true for $G_1=G$ then $\sigma$ is properly outer. If $\nu =1$
then the converse holds true as well.

$(iii)$. Let $\Gamma_1 \subset \Gamma_0, G_1 \subset G$
be subgroups of finite index such that $G_1$ is invariant to the restriction
of $\gamma$ to $\Gamma_1$. If $\gamma, \Gamma_0, G$
satisfy either of the conditions $(c), (d)$ in $(i)$, or $(ii)$
then $\gamma_{|\Gamma_1},
\Gamma_1, G_1$ satisfy that condition as well.
\endproclaim
\noindent
{\it Proof}. $(i)$. $(b) \implies (a)$ is trivial.

$(a) \implies (c)$.
If $\gamma_h(S) = S, \forall h \in
\Gamma_0$ for some finite set $S \subset G$ with $e \not\in S$, then
$x = \Sigma_{g \in S} \lambda(g)\notin
\Bbb C1$ satisfies $\gamma_h(x) = x, \forall h\in
\Gamma_0$, implying that $\sigma$ is not ergodic.

$(c) \implies (d)$. If $\gamma_h(S) \cap S
\neq \emptyset, \forall h\in \Gamma_0$,
for some finite set $S \subset G\backslash \{e\}$, then
denote by $f$ the characteristic function of $S$
regarded as an element of $\ell^2(G)$. If we denote
by $\tilde{\gamma}$ the action (=representation) of $\Gamma_0$ on
$\ell^2(G)$ implemented by $\gamma$, then we have
$\langle \tilde{\gamma}_h(f), f\rangle \geq 1/|S|, \forall h\in \Gamma_0$.
Thus, the element $a$ of minimal norm $\|\quad\|_2$ in the weak closure of
co$\{\tilde{\gamma}_h(f) \mid h \in \Gamma_0\} \subset \ell^2(G)$
is non zero. But then any ``level set'' of $a \geq 0$ is
invariant to $\gamma$, showing that $(c)$ doesn't hold true.

$(d) \implies (b)$. Let $E_0$ be a finite set in the unit
ball of $L_\nu(G)$, $\varepsilon > 0$ and
$F_0 \subset \Gamma_0 \backslash \{e\}$ a finite set as well.
Let $S_0 \subset G\backslash \{e\}$ be such that $\|(x -
\tau(x)1)-x_{S_0}\|_2
\leq \varepsilon/2, \forall x\in E_0$. By applying the
hypothesis to $S = \cup \{\gamma_g(S_0) \mid g\in F_0\}$, it follows that
there exists $g\in \Gamma_0$ such that
$\gamma_g(S) \cap S = \emptyset$. But then $g\notin F_0$ and
$\gamma_g(S_0) \cap S_0 = \emptyset$. Also, by Cauchy-Schwartz,
for each $x,y \in E_0$ we have:
$$
|\tau(\sigma_{\gamma(g)}(x)y) - \tau(x) \tau(y)|
$$
$$
\leq \|(x-\tau(x)1)
-x_{S_0}\|_2 \|y\|_2 +
\|(y-\tau(y)1)-y_{S_0}\|_2 \|x\|_2 + |\tau(\sigma_{\gamma(g)}(x_{S_0})y_{S_0})|
$$
$$
=\|(x-\tau(x)1)
-x_{S_0}\|_2 \|y\|_2 +
\|(y-\tau(y)1)-y_{S_0}\|_2 \|x\|_2 \leq \varepsilon.
$$

$(ii)$. If $y_0 \in L_\nu(G) \rtimes_\sigma \Gamma_0$ satisfies
$y_0 x = y_0 x, \forall x\in L_\nu(G_1)$ and $y_0 \notin L_\nu(G)$
then there exists $h \in \Gamma_0, h\neq e,$ such that
$\sigma_{\gamma(h)}(x) a = a x, \forall x \in L_\nu(G),$ for some
$a \in L_\nu(G), a \neq 0$. This implies
$\lambda(\gamma_h(g_1)) a \lambda(g_1^{-1}) = a, \forall
g_1 \in G_1$. But if this holds true then $\{\gamma_h(g_1)g'g_1^{-1} \mid
g_1 \in G_1\}$ must be finite, for any $g' \in G$ in the support of $a$.
When $G_1=G$ and $\nu=1$, reversing the implications proves the converse.

$(iii)$. Note first that if
$S \subset G_1$ is a finite
subset such that $\gamma_h(S) = S,
\forall h\in \Gamma_1$, the set $\cup_{h \in \Gamma_0} \gamma_h(S)$
follows finite as well. Thus,
if $\gamma, \Gamma_0, G$ checks $(c)$ in $(i)$ so does $\gamma_{|\Gamma_1},
\Gamma_1, G_1$.

Then note that if $\gamma, \Gamma_0, G$
verifies $(ii)$ and for some $g_1 \in G_1$ the set
$\{\gamma_h(g) g_1 g^{-1} \mid g \in G_1\}$ is finite, then the set
$\{\gamma_h(g) g_1 g^{-1} \mid g \in G\}$ follows finite,
contradiction.
\hfill Q.E.D.

\proclaim{Corollary} Let $\tilde{\gamma}$ be the action
of the group $SL(2, \Bbb R)$ on $\Bbb R^2$. For each $\alpha =e^{2\pi i t}
\in \Bbb T,$ let $\tilde{\nu}=\tilde{\nu}(\alpha)$ be the
unique normalized scalar $2$-cocycle
on $\Bbb R^2$ satisfying the relation $u_xv_y = exp(2\pi i txy)v_yu_x$,
where $u_x = (x,0), v_y = (0,y)$ for $x, y \in \Bbb R$.
Then $\tilde{\nu}$ is $\tilde{\gamma}$-invariant.
Moreover, the following restrictions $(\gamma, \Gamma_0, G, \nu)$
of $(\tilde{\gamma}, SL(2, \Bbb R), \Bbb R^2, \tilde{\nu})$
are stronly ergodic and satisfy conditions
$(i), (ii)$ in the previous ${\text{\rm Lemma}}$ (so
the corresponding actions $\sigma_\gamma$ of $\Gamma_0$ are free and
weakly mixing on $L_\nu(G)$):

$(a)$. $\Gamma_0=SL(2, \Bbb Z),
G=\Bbb Z^2$, or any other subgroup $G$ of $\Bbb R^2$ which
is $SL(2, \Bbb Z)$-invariant, with $\gamma$ the
appropriate restriction of $\tilde{\gamma}$
(and of $\tilde{\nu}$).

$(b)$. $\Gamma_0=SL(2, \Bbb Q)$, $G=\Bbb Q^2$
(or any other $SL(2, \Bbb Q)$-invariant
subgroup of $\Bbb R^2$), with $\gamma$ the appropriate
restriction of $\tilde{\gamma}$.

$(c)$. $\Gamma_0 \simeq \Bbb F_n$,
regarded as a subgroup of finite index in $SL(2, \Bbb Z)$
(see e.g., ${\text{\rm [dHVa]}}$), and $G=L((k \Bbb Z)^2)$,
for some $k \geq 1$.
\endproclaim
\noindent
{\it Proof}. Both conditions
$(i)$ and $(ii)$ of the
Lemma are trivial to check in case $(a)$
and $(b)$. Then $(c)$ is just a simple consequence of part $(iii)$
of the Lemma. The strong ergodicity
of these actions was proved in ([S1], [Va]).
\hfill Q.E.D.
\vskip .1in
\noindent
{\it 3.3.3. Tensor products of actions}. We'll  often need
to take tensor products of actions $\sigma_i$ of the same
group $\Gamma_0$ on $(B_i, \tau_i)$, $i=1, 2, ...$,
thus getting
an action $\sigma= \sigma_1\otimes \sigma_2\otimes...$ of $\Gamma_0$ on
$(B, \tau) = (B_1, \tau_1) \overline{\otimes} (B_2, \tau_2)
\overline{\otimes}...$.

It is easy to see that the tensor product
of a properly outer action
$\sigma$ of a group $\Gamma_0$ with any other
action $\sigma_0$ of $\Gamma_0$ gives a properly outer action.
In fact, if $\sigma$ is an action
of $\Gamma_0$ on $(B, \tau)$ and $A_0 \subset B$ is
so that $A_0'\cap B\rtimes_\sigma \Gamma_0 \subset B$
then given any action $\sigma_0$ of
$\Gamma_0$ on some $(B_0, \tau_0)$,
we have $(A_0 \otimes 1)' \cap (B\overline{\otimes} B_0
\rtimes_{\sigma\otimes\sigma_0} \Gamma_0) = (A_0'\cap B)
\overline{\otimes} B_0$.

While ergodicity
does not always behave well with respect to tensor products,
weak-mixing does: If $\sigma$ is weakly mixing and $\sigma_0$
is ergodic then $\sigma \otimes \sigma_0$
is ergodic. If $\sigma_i, i\geq 1,$ are weakly
mixing then $\otimes_i \sigma_i$
is weakly mixing.

If $\sigma_0$ is not strongly ergodic, then $\sigma\otimes \sigma_0$
is not strongly ergodic $\forall \sigma$. Note that by ([CW]),
if $\Gamma_0$ is an infinite  property H group then there always exist
free ergodic measure preserving actions $\sigma_0$ of $\Gamma_0$
on $L^\infty(X, \mu)$ which are not strongly ergodic.
Thus, given
any $\sigma$, $\sigma\otimes \sigma_0$ is not strongly ergodic either.

The following combination of
Bernoulli shifts and
tensor products of actions will be of interest to us: Let $\sigma_0$
be an action of $\Gamma_0$ on $(B_0, \tau_0)$. Let also
$\Gamma_1$ be another discrete group and $\gamma$ an action
of $\Gamma_1$ on $\Gamma_0$ by group automorphisms. (N.B.:
The action $\gamma$ may be trivial.) Let
$\sigma_1$ be the Bernoulli shift action
of $\Gamma_1$ on $(B,\tau)={\overline{\underset
g_1 \in \Gamma_1 \to \otimes}} (B_0, \tau_0)_{g_1}$. Let also
$\sigma_0^\gamma$ be the action of $\Gamma_0$ on $(B, \tau)$
given by $\sigma_0^\gamma =
\otimes_{g_1} \sigma_0 \circ \gamma(g_1)$.

\proclaim{Lemma} $1^\circ$. We have $\sigma_1(g_1) \sigma^\gamma_0(g_0)
\sigma_1(g_1^{-1})= \sigma^\gamma_0(\gamma(g_1)(g_0))$,
for any $g_0 \in \Gamma_0$ and $g_1\in \Gamma_1$. Thus,
$(g_0, g_1) \mapsto \sigma^\gamma_0(g_0) \sigma_1(g_1)$
implements an action
$\sigma=\sigma_0\rtimes_\gamma \sigma_1$ of
$\Gamma_0 \rtimes_\gamma \Gamma_1$
on $(B,\tau)$.

$2^\circ$. If the group $\Gamma_0$ is infinite and the action
$\sigma_0$ is properly outer then
the action $\sigma$ defined in $1^\circ$ is properly outer.
Moreover, if $B_1\subset B_0$ satisfies
$B_1'\cap (B_0 \rtimes_{\sigma_0} \Gamma_0) \subset B_0$,
and we identify $B_1$
with $...\otimes \Bbb C \otimes B_1 \otimes \Bbb C...
\subset B$, then $B_1'\cap (B\rtimes_\sigma (\Gamma_0 \rtimes \Gamma_1))
=B_1'\cap B$.

$3^\circ$. If the action $\sigma_0$ is weakly mixing, or if the group
$\Gamma_1$ is infinite,  then $\sigma$ is weakly mixing
(thus ergodic).

$4^\circ$. If the group $\Gamma_1$ is non-amenable, then $\sigma$ is
strongly ergodic.
\endproclaim
\noindent
{\it Proof}. $1^\circ$ is straightforward direct calculation.

$2^\circ$ follows once we notice that if $\Gamma_0$ is infinite and
$\sigma_0$ is properly outer, it automatically follows that
$B_0$ has no atomic part. This in turn implies the Bernoulli shift
of $\Gamma_1$ on
$(B_0, \tau_0)^{\otimes\Gamma_1}$ is a properly outer action,
even when $\Gamma_1$
is a finite group.

$3^\circ$. This follows by the observations at the beginning of 3.3.3
and 3.3.1.

$4^\circ$. This follows from the properties of the
Bernoulli shift listed in  3.3.1 (cf. [J2]).
\hfill Q.E.D.

\proclaim{3.4. Proposition} If the finite
von Neumann algebra $N$ has the property
${\text{\rm H}}$ relative to its
von Neumann subalgebra $B \subset N$,
then $B$ is quasiregular in $N$. If in addition $N$ is
a type ${\text{\rm II}}_1$ factor $M$ and $B=A$ is maximal abelian
in $M$, then $A$ is a Cartan subalgebra of $M$.
\endproclaim
\noindent
{\it Proof}. By Proposition 2.3, given any $x_1, x_2, ..., x_n \in N$,
with $\|x_i\|_2 \leq 1$, and
any $\varepsilon > 0$, there exists an operator $T \in B'\cap
J(\langle N,B \rangle)$ such that $\|T\| \leq 1$ and $\|T(\hat{x_i}) -
\hat{x_i}\|_2 < \varepsilon^2/32, \forall i$. Since $\|T\| \leq 1$, this implies
$$
\|T^*(\hat{x_i})-\hat{x_i}\|^2_2 = \|T^*(\hat{x_i})\|^*_2 - 2{\text{\rm Re}}
\langle T^*(\hat{x_i}),\hat{x_i} \rangle + \|x_i\|_2^2
$$
$$
\leq 2\|x_i\|^2_2 - 2 {\text{\rm Re}}\langle T^*(\hat{x_i}),\hat{x_i} \rangle
=2 {\text{\rm Re}}\langle \hat{x_i}, (\hat{x_i}-T(\hat{x_i}))\rangle
$$
$$
\leq 2 \|x_i\|_2 \|\hat{x_i} - T(\hat{x_i})\|_2 < \varepsilon^2/16.
$$
As a consequence, we get:
$$
\|T^*T(\hat{x_i}) - \hat{x_i}\|_2 \leq \|T^*\| \|T(\hat{x_i}) - \hat{x_i}\|_2
+ \|T^*(\hat{x_i}) - \hat{x_i}\|_2 < \varepsilon/2.
$$
Thus, if we let $e$ be the spectral projection of $T^*T$ corresponding to
$[1-\delta, 1]$ then $\|T^*T - T^*Te\| \leq \delta$, yielding
$$
\|e(\hat{x_i}) -
\hat{x_i}\|_2 \leq  \|T^*T(\hat{x_i}) - \hat{x_i}\|_2 +
\|e(T^*T(\hat{x_i}) - (\hat{x_i}))\|_2 + \|T^*T - T^*Te\|
$$
$$
\leq 2\|T^*T(\hat{x_i}) - \hat{x_i}\|_2 + \delta.
$$
But for $\delta$ sufficiently small the latter
follows less than $\varepsilon, \forall i$. Since
the projection $e$ lies in $B'\cap J(\langle N,B \rangle)$,
this proves that $\vee \{f\mid f \in \Cal P(B'\cap \langle N,B \rangle), f$
finite projection in $\langle N,B\rangle \}=1$. By part
$(iii)$ of Lemma 1.4.2,
this implies $B$ is
quasiregular in $N$. If in
addition $B$ is a maximal abelian subalgebra then
$B$ follows Cartan by ([PoSh]; see also part $(i)$ in Proposition
1.4.3).
\hfill Q.E.D.
\vskip .1in
\noindent
{\bf 3.5. Remarks}.
$0^\circ$. It is interesting to note that
in most known examples
of groups $\Gamma_0$ with the property H, the positive definite
functions $\varphi_n \in c_0(\Gamma_0)$ approximating the identity
can be chosen in $\ell^p(\Gamma_0)$, for some $p=p(n)$. This is the case,
for instance, with the free groups
$\Bbb F_m$ (cf. [H]), the
arithmetic
lattices in $SO(m,1), SU(m,1)$, etc. It is a known fact that
if all $\varphi_n$ can be taken in the same $\ell^p(\Gamma_0)$,
(which is easily seen to imply they can be taken in 
$\ell^2(\Gamma_0), \forall n$) then $\Gamma_0$
follows amenable. This fact,
along with many other similar observations, justifies regarding Haagerup's
approximating property as a ``weak amenability'' property.

$1^\circ$. The same proof
as in ([Cho]) shows that
if $G \subset G_0$ is an inclusion of discrete groups
with the property that there exists a net of positive definite
functions $\varphi_\alpha$ on $G_0$ which are constant on double
cosets $Gg_0G, \forall g_0\in G_0$
(thus factoring out to bounded functions on $G\backslash G_0/G$) and satisfy
\vskip .1in
\noindent
(3.5.1'). $\quad$ $G$ is quasi-normal in  $G_0$ and
$\varphi_\alpha \in c_0(G\backslash G_0/G), \forall \alpha$;
\vskip .1in
\noindent
(3.5.1''). $\quad \underset \alpha \rightarrow \infty \to \lim
\varphi_\alpha(g_0) = 1, \forall g_0 \in \Gamma_0$.
\vskip .1in
\noindent
then $L_\nu(G_0)$ has the property H relative to $L_\nu(G)$
for any scalar 2-cocycle $\nu$ for $G_0$.

When $G\subset G_0$ satisfies the set of conditions $(3.5.1)$ we say that
$G_0$ has the {\it property} H {\it relative to} $G$. Note that in the case
$G$ is normal in $G_0$ this is equivalent to $G_0/G$ having the property H
as a group. (See 3.18-3.20 in [Bo] for similar such considerations).

$2^\circ$. The relative property H for
inclusions of finite von Neumann algebras is related to the following
notion of relative amenability considered
in ([Po1,5]): If $B \subset N$ is an
inclusion of finite
von Neumann algebras then $N$ is {\it amenable relative to}
$B$ if there exists a norm one projection of $\langle N,B \rangle =
(J_NBJ_N)' \cap \Cal B(L^2(N))$
onto $N$, where $L^2(N)$ is the standard representation of $N$ and $J_N$
is the corresponding canonical conjugation.

It is easy
to see that if $B \subset N$ is a cross-product
inclusion $B \subset B\rtimes_\sigma \Gamma_0$ for some cocycle action
$\sigma$ of a discrete group $\Gamma_0$ on $(B, \tau_0)$,
with $\tau_0$ a faithful normal trace on $B$, then $N$
is amenable relative to $B$ in the above sense
if and only if $\Gamma_0$ is amenable,
a fact that justifies the terminology.
Thus, in this case $N$ amenable
relative to $B$ implies $N$ has the property H relative to $B$.

If $N$ is an arbitrary finite von Neumann algebra
with a normal faithful tracial state
$\tau$ and $B \subset N$ is a von Neumann subalgebra, then
the amenability of $N$ relative to $B$ is equivalent to the
existence of a $N$-{\it hypertrace}
on $\langle N,B \rangle$, i.e., a state $\varphi$ on $\langle N,B \rangle$
with $N$ in its centralizer: $\varphi (xT) = \varphi (Tx), \forall x\in N,
T\in \langle N,B \rangle$ (cf. [Po1]). It is also
easily seen to be equivalent
(by using the standard Day-Namioka-Connes trick)
to the following
F\o lner type condition: 
$\forall F\subset \Cal U(N)$ finite and
$\varepsilon > 0$, $\exists$ 
$e \in \Cal P(\langle N,B \rangle)$ with Tr$e < \infty$ such that
\vskip .1in
\noindent
$(3.5.2)$. $\| u_0 e - e u_0\|_{2,{\text{\rm Tr}}}
< \varepsilon \|e\|_{2,{\text{\rm Tr}}},
\forall u_0\in F.$
\vskip .1in

Note that in case $(B\subset N)
= (L_\nu(G) \subset L_\nu(G_0))$ for some inclusion
of discrete groups $G\subset G_0$ and a scalar $2$-cocycle $\nu$ on $G_0$,
condition $(3.5.2)$ amounts to the following: $\forall F\subset G_0$
finite and $\varepsilon > 0$,  $\exists E \subset G_0/G$
finite such that
\vskip .1in
\noindent
$(3.5.2$'$)$. $|g_0E - E| < \varepsilon |E|, \forall g_0 \in F$.
\vskip .1in

This condition for inclusions of groups, for which the terminology
used is ``$G$ {\it co-F\o lner in} $G_0$'', was first considered
in ([Ey]). It has been used in ([CCJJV]) to prove that if $G
\subset G_0$ is an inclusion of groups, $G_0$ is amenable relative
to $G$ and $G$ has Haagerup property, then $G_0$ has Haagerup's
property. It would be interesting to know whether a similar result
holds true in the case of inclusions of finite von Neumann
algebras $B\subset N$.

$3^\circ$. A stronger version of relative amenability
for inclusions of finite von Neumann algebras
$B\subset N$ was considered in ([Po5]), as follows: $N$ is
{\it s-amenable relative to} $B$ if
given any finite set
of unitaries $F \subset \Cal U(N)$ and any $\varepsilon > 0$ there
exists a projection $e \in B' \cap \langle N,B \rangle$, with Tr$e < \infty$,
such that $e$ satisfies the F\o lner condition $(3.5.2)$ and
$\|Tr(\cdot e)/Tr(e) - \tau\| \leq \varepsilon$.
(No specific
terminology is in fact used in [Po5] 
to nominate this amenability property.)
Note that in case $B'\cap N = \Bbb C$,
we actually have $Tr(\cdot e)/Tr(e) = \tau$ for any finite projection
$e$ in $B'\cap \langle N, B\rangle$, so the second condition is redundant. The 
$s$-amenability of $N$ relative to $B$ 
is easily seen to be 
equivalent to: There exists a net of $B$-bimodular completely positive maps
$\phi_\alpha$ on $N$ such that $\tau \circ \phi_\alpha \leq \tau$,
$T_{\phi_\alpha}$
belong to the (algebraic) ideal generated in $\langle N, B \rangle$ by $e_B$ and
$\underset \alpha \rightarrow \infty \to \lim
\|\phi_\alpha(x)-x\|_2 = 0, \forall x\in N.$ Thus, 
$N$ $s$-amenable relative to $B$ 
implies $N$ has property H relative to $B$. Also, 
one can check that 
if $N = B \rtimes_\sigma \Gamma_0$ for some
cocycle action $\sigma$ of a discrete
group $\Gamma_0$ on $(B, \tau)$, then 
$N$ is $s$-amenable relative to
$B$ iff $N$ is amenable relative to $B$ and iff 
$\Gamma_0$
is an amenable group.

$4^\circ$. Let $N \subset M$ be an
extremal inclusion of
type II$_1$ factors with finite Jones index and let
$T=M\vee M^{{\text{\rm op}}} \subset
M\bt_{e_N}M^{{\text{\rm op}}}=S$
be its associated {\it symmetric enveloping inclusion}, as
defined in ([Po5]). It was shown in ([Po5]) that
$T$ is quasiregular in $S$.
It was also shown that
$S$ is amenable relative to
$T$ iff
$S$ is s-amenable relative to
$T$ and iff
$N\subset M$ has amenable graph $\Gamma_{N,M}$ (or, equivalently,
$N \subset M$ has amenable standard invariant $\Cal G_{N,M}$).

By (Sec. 3 in [Po5]), if 
$N\subset M$ is the subfactor associated to a properly
outer cocycle action $\sigma$
of a finitely generated group $\Gamma_0$ on a factor $\simeq M$,
then the corresponding
symmetric enveloping inclusion $T=M\vee M^{{\text{\rm op}}} \subset
M\bt_{e_N}M^{{\text{\rm op}}}=
S$ is isomorphic to
$M\overline{\otimes} M^{op} \subset M\overline{\otimes} M^{op}
\rtimes_{\sigma\otimes\sigma^{op}} \Gamma_0$,
so $T$ is regular in $S$.
But if $N\subset M$ has index
$\lambda^{-1}
\geq 4$ and Temperley-Lieb-Jones (TLJ) standard invariant
$\Cal G_{N,M}=\Cal G^\lambda$,
then the corresponding symmetric enveloping inclusion
$T\subset S$ is
quasi-regular but not regular. In particular, if $\lambda^{-1}=4$
then $[S:T] = \infty$ and
$S$ has property H 
relative to $T$ (because $\Cal G_{N,M}$ is amenable by [Po3]), 
while $T$ is quasi-regular but
not regular in $S$.

$5^\circ$. By using exactly the same
arguments as in the case of the property (T)
for standard lattices considered in ([Po5]),
it can be shown that for an extremal standard lattice $\Cal G$
the following conditions are equivalent: $(i)$. There exists
an irreducible subfactor $N \subset M$ with $\Cal G_{N,M} = \Cal G$
such that $M\bt_{e_N}M^{{\text{\rm op}}}$ has property H relative
to $M\vee M^{{\text{\rm op}}}$; $(ii)$. Given any subfactor
$N \subset M$ with $\Cal G_{N,M} = \Cal G$,
$M\bt_{e_N}M^{{\text{\rm op}}}$ has property H relative
to $M\vee M^{{\text{\rm op}}}$. If $\Cal G$ satisfies either of these
conditions, we say that the standard lattice $\Cal G$ {\it has the
property} H. By 3.6, any amenable $\Cal G$ has the property H.
We will prove in a forthcoming paper that
the TLJ standard lattices $\Cal G^\lambda$ have the
property H, $\forall \lambda^{-1} \geq 4$, while they are known to be
amenable iff $\lambda^{-1}=4$ ([Po2], [Po5]).

$6^\circ$. When applied to the case of Cartan subalgebras
$A \subset M$ coming from
standard equivalence relations $\Cal R$
(i.e., countable, free, ergodic, measure preserving) and having
trivial 2-cocycle $v\equiv 1$, definition 2.2 gives the
following: A standard equivalence relation $\Cal R$
has
{\it the property} H (or is of {\it Haagerup-type})
if $M$ has the property H relative to $A$. Note that in case
$\Cal R$ comes from an action $\sigma$ of a group $\Gamma_0$
then the property H of the corresponding $\Cal R$ depends entirely
on the group $\Gamma_0$, and not on the action (cf. 3.1).
Since in addition $A \rtimes \Gamma_0$ has the
property H relative to $A$ if and only if
$p(A \rtimes \Gamma_0)p$ has the property H relative to $Ap$,
for $p\in \Cal P(A)$ (cf. 2.5), it follows that property H for groups is
invariant to stable orbit equivalence (this fact was
independently noticed by Jolissaint;
see [Fu] for a reformulation of stable orbit equivalence as
Gromov's ``measure equivalence'', abreviated ME).

\heading 4. Rigid embeddings: definitions and properties. \endheading

In this section we consider a notion of
rigid embeddings for finite von Neumann
algebras, inspired by
Kazhdan's example of the rigid embedding of groups $\Bbb Z^2 \subset
\Bbb Z^2 \rtimes SL(2, \Bbb Z)$. Our
definition will be the operator algebraic version of the
notion of property (T) for pairs of groups in
([Ma], [dHVa]),
in the same spirit Connes and Jones
defined the property (T) for single von Neumann algebras starting from the
property (T) of groups, in
([CJ]). Thus, like in ([CJ]), to formulate
the definition we use Connes's idea ([C3]) of
regarding Hilbert bimodules as an operator algebra substitute
for unitary representations of groups, and
completely positive maps as an operator algebra
substitute for positive definite functions on groups
(see Section 1.1).
For convenience (and comparison),
we first recall the definition of property (T) for
inclusions of groups and for single II$_1$ factors:
\vskip .1in
\noindent
{\it 4.0.1. Relative property} (T) {\it for pairs of groups}. The key part
in Kazhdan's proof that the groups $SL(n, \Bbb R)$
(resp. $SL(n, \Bbb Z)$), $n\geq 3$,
have the property (T) consists in showing that representations
of $\Bbb R^2 \rtimes SL(2, \Bbb R)$ that are
close to the trivial representation contain copies of the trivial 
representation of $\Bbb R^2$. This
type of ``relative
rigidity'' property was later emphasized as a notion in its own right by
Margulis ([Ma]; see also [dHVa]), as follows:

Let $G \subset G_0$ be an inclusion of discrete groups. {\it The
pair} $(G_0, G)$ {\it has the relative property} (T) if the following
condition holds true: 
\vskip .1in 
\noindent (4.0.1). There exist
finitely many elements $g_1, g_2, ..., g_n \in G_0$ and
$\varepsilon > 0$, such that if $\pi : G_0 \rightarrow \Cal U(\Cal
H)$ is a unitary representation of the group $G_0$ on the Hilbert
space $\Cal H$ with a unit vector $\xi \in \Cal H$ satisfying
$\|\pi(g_i)\xi - \xi \| < \varepsilon, \forall i$, then there
exists a non-zero vector $\xi_0 \in \Cal H$ such that $\pi(h)
\xi_0 = \xi_0, \forall h\in G$. 
\vskip .1in 
Due to a recent result
of Jolissaint ([Jo2]), the above condition is equivalent to:
\vskip .1in 
\noindent (4.0.1'). For any $\varepsilon > 0$, there
exist a finite subset $E' \subset G_0$ and $\delta' > 0$ such that
if $(\pi, \Cal H)$ is a unitary representation of $G_0$ on the
Hilbert space $\Cal H$ and $\xi \in \Cal H$ is a unit vector
satisfying $\|\pi(h)\xi - \xi \| \leq \delta', \forall h\in E'$,
then $\|\pi(g)\xi - \xi \| \leq \varepsilon, \forall g \in G.$
\vskip .1in 
Note that the equivalence of (4.0.1) and (4.0.1') is
easy to establish in case $G$ is a normal subgroup of $G_0$
(exactly the same argument as in [DeKi] will do), but it is less
simple in general (cf. [Jo2]). On the other hand, condition
(4.0.1') is easily seen to be equivalent to: 
\vskip .1in 
\noindent
(4.0.1''). For any $\varepsilon > 0$, there exist a finite subset
$E \subset G_0$ and $\delta > 0$ such that if $\varphi$ is a
positive definite function on $G_0$ with $|\varphi(h)-1| \leq
\delta, \forall h\in E$ then $|\varphi(g) - 1|\leq \varepsilon,
\forall g\in G$. 
\vskip .1in Note that in the case $G=G_0$,
condition (4.0.1) amounts to the usual {\it property} T {\it of
Kazhdan} for the group $G_0$ ([Kaz]; see also [DeKi], [Zi]). We
will in fact also use the following alternative terminologies to
designate property (T) pairs: $G \subset G_0$ {\it is a property} (T)
(or {\it rigid}) {\it embedding}, or $G$ is a {\it relatively
rigid subgroup} of $G_0$. 
\vskip .1in 
\noindent 
{\it 4.0.2. Property} (T) {\it for factors}. The abstract 
definition of property
(T) for a single von Neumann factors is due to Connes and Jones
([CJ]): A type II$_1$ factor $N$ {\it has the property} (T) if the
following condition holds true \vskip .1in \noindent (4.0.2).
There exist finitely many elements $x_1, x_2, ..., x_n \in N$ and
$\varepsilon_0 > 0$ such that if $\Cal H$ is a $N$ Hilbert
bimodule with a unit vector $\xi \in \Cal H$ such that $\|x_i \xi
- \xi x_i\| \leq \varepsilon_0, \forall i$, then $\Cal H$ contains
a non-zero vector $\xi_0$ such that $x\xi_0 = \xi_0 x, \forall
x\in N$. \vskip .1in Connes and Jones have also proved that the
fixed vector $\xi_0$ can be taken close to the initial $\xi$, if
the ``critical set'' in $N$ is taken sufficiently large and the
``commutation constant'' sufficiently small ([CJ]), by showing
that (4.0.2) is equivalent to the following: 
\vskip .1in 
\noindent
(4.0.2'). For any $\varepsilon > 0$, there exist a finite subset
$F' \subset N$ and $\delta' > 0$ such that if $\Cal H$ is a
Hilbert $N$-bimodule and $\xi \in \Cal H$ is a unit vector
satisfying $\|y \xi - \xi y\| \leq \delta', \forall y\in F'$, then
there exists $\xi_0 \in \Cal H$ such that $x\xi_0=\xi_0 x, \forall
x \in N$ and $\|\xi - \xi_0\| \leq \varepsilon.$ 
\vskip .1in

For inclusions of finite von Neumann algebras, we first establish the
equivalence of several conditions:

\proclaim{4.1. Proposition} Let
$N$ be a finite von Neumann algebra
with countable decomposable center (i.e., which has normal faithful
tracial states). Let $B \subset N$ be a von
Neumann subalgebra. The following conditions are equivalent:

$1^\circ$. There exists a normal faithful tracial state $\tau$ on
$N$ such that: $\forall \varepsilon > 0$, $\exists F'\subset N$
finite and $\delta' > 0$ such that if $\Cal H$ is a Hilbert
$N$-bimodule with a vector $\xi \in \Cal H$ satisfying the
conditions $\|\langle \cdot\xi, \xi\rangle - \tau\| \leq \delta',
\|\langle \xi \cdot, \xi\rangle - \tau\| \leq \delta'$ and $\|y
\xi - \xi y\| \leq \delta', \forall y\in F'$ then $\exists
\xi_0\in \Cal H$ such that $\|\xi_0 - \xi\| \leq \varepsilon$ and
$b \xi_0 = \xi_0 b, \forall b\in B$.

$2^\circ$. There exists a normal faithful tracial state $\tau$ on
$N$ such that: $\forall \varepsilon > 0$, $\exists F\subset  N$
finite and $\delta > 0$ such that if $\phi: N \rightarrow N$ is a
normal, completely positive map with $\tau\circ \phi \leq \tau,
\phi(1)\leq 1$ and $\|\phi(x) - x \|_2 \leq \delta, \forall x\in
F$, then $\|\phi(b)-b\|_2 \leq \varepsilon, \forall b\in B$,
$\|b\| \leq 1$.

$3^\circ$. Condition $1^\circ$ above is satisfied for any normal faithful
tracial state $\tau$ on $N$.

$4^\circ$. Condition $2^\circ$ above is satisfied for any normal
faithful tracial state $\tau$ on $N$.
\endproclaim
\noindent
{\it Proof}. We first prove that
condition $1^\circ$ holds true for a specific
normal faithful tracial state $\tau$ if and only if
condition $2^\circ$ holds true for that same trace.
Then we prove $1^\circ \Leftrightarrow 3^\circ$,
which due to the equivalence of $1^\circ$ and $2^\circ$
ends the proof of the Proposition.

$2^\circ \implies 1^\circ$. By part
$1^\circ$ of Lemma 1.1.3, we may assume the vectors $\xi \in \Cal H$
in condition 4.1.1$^\circ$
also satisfy $\langle \cdot \xi, \xi\rangle \leq \tau$ and
$\langle \xi\cdot, \xi \rangle \leq \tau$, in addition to the
given properties. We take
$x_1, x_2, ..., x_n$ to be an enumeration of the finite
set $F$ and for any given $\varepsilon' > 0$ let $\delta'$ be
the $\delta$ given by condition $2^\circ$ for $\varepsilon =
{\varepsilon'}^2/4.$
By part $2^\circ$ of Lemma 1.1.3, such a vector
$\xi$ gives rise to a completely positive map $\phi=\phi_{(\Cal H, \xi)}$
on $N$ which satisfies condition 4.1.2$^\circ$.
Thus, $\|\phi(b)-b\|_2 \leq \varepsilon, \forall b\in B,
\|b\|\leq 1$. By Lemma 1.1.2, this implies
$\xi$ (which is equal to $\xi_\phi$) satisfies
$\|u\xi u^* - \xi \| \leq 2 \varepsilon^{1/2} \leq \varepsilon',
\forall u\in \Cal U(B)$. By averaging over the
unitaries $u\in \Cal U(B)$, it follows that there
exists $\xi_0 \in \Cal H$ such that $\|\xi_0 - \xi\|
\leq \varepsilon'$ and $\xi_0$ commutes with $B$.

$1^\circ \implies 2^\circ$. Let $\varepsilon > 0$. Define
$F(\varepsilon) = F'(\varepsilon^2/8), \delta(\varepsilon)
=\delta'(\varepsilon^2/8)^2/4$. Let then $\phi:N\rightarrow N$ be a
completely positive map satisfying the conditions $2^\circ$ for this
$F(\varepsilon)$ and $\delta(\varepsilon)$. Let $(\Cal H_\phi, \xi_\phi)$
be constructed as in 1.1.2. By part $4^\circ$
of Lemma 1.1.2, we have for $x\in F(\varepsilon)$ the inequality
$$
\|x \xi_\phi - \xi_\phi x\| \leq 2 \|\phi(x)-x\|_2^{1/2} \leq
\delta'(\varepsilon^2/8).
$$

Thus, there exists $\xi_0 \in \Cal H_\phi$ such that $\|\xi_\phi - \xi_0\|
\leq \varepsilon^2/8$ and $b\xi_0 = \xi_0 b, \forall b\in B$. But then,
if $u \in \Cal U(B)$ we get
$$
\|\phi(u)-u\|_2^2 \leq 2 - 2 {\text{\rm Re}}
\langle u\xi_\phi u^*, \xi_\phi \rangle
$$
$$
\leq 2 - 2 \|\xi_0\|^2 + 4 \|\xi_0 - \xi_\phi\| \leq 2 -
2(1-\varepsilon^2/8)^2 + 4 \varepsilon^2/8 < \varepsilon^2.
$$
Since any $b \in B$, $\|b\|\leq 1$, is a convex combination of
unitary elements, we are done.

$3^\circ \implies 1^\circ$ is trivial. To prove $1^\circ \implies
3^\circ$, let $\tau_0$ be a normal faithful tracial state on $N$.
We have to show that $\forall \varepsilon > 0$, $\exists F_0
\subset N$ finite and $\delta_0 > 0$ such that if $\Cal H$ is a
Hilbert $N$-bimodule with $\eta \in \Cal H$ satisfying $\|\langle
\cdot\eta, \eta\rangle - \tau_0\| \leq \delta_0, \|\langle \eta
\cdot, \eta\rangle - \tau_0\| \leq \delta_0$ and $\|y \eta - \eta
y\| \leq \delta_0, \forall y\in F_0$ then $\exists \eta_0\in \Cal
H$ such that $\|\eta_0 - \eta\| \leq \varepsilon$ and $b \eta_0 =
\eta_0 b, \forall b\in B$.

By Sakai's Radon-Nykodim theorem,
$\tau_0$ is of the form $\tau_0 = \tau(\cdot a_0)$ for
some $a_0 \in L^1(\Cal Z(N), \tau)_+$ with $\tau(a_0)=1$.
It is clearly sufficient to prove the statement in the case
$a_0$ is
bounded and with finite spectrum (thus bounded away from 0 as well). Also,
by taking the spectral projections of $a_0$ to be in $F_0$ and slightly
perturbing $\eta$,
we may assume $\eta$ commutes with $a_0$.
We take $F_0=F'(\varepsilon/\|a_0\|)$ and
$\delta_0=\delta'(\varepsilon/\|a_0\|)/\|a_0^{-1}\|$, as
given by condition $1^\circ$ for $\tau$.

Let $\xi= a_0^{-1/2}\eta=\eta a_0^{-1/2}$. Then we have
$$
\|\langle \cdot \xi, \xi \rangle -\tau\| =
\|\langle \cdot a_0^{-1} \eta, \eta \rangle - \tau_0(\cdot a_0^{-1})\|
\leq \|a_0^{-1}\|(\delta'/\|a_0^{-1}\|) = \delta'.
$$
Similarly, $\|\langle \xi \cdot, \xi \rangle -\tau\|\leq \delta'$.
Also, for $y \in F_0$ we have:
$$
\|[y, \xi]\| = \|[y, a_0^{-1/2}\eta]\|
\leq \|a_0^{-1/2}\|(\delta'/\|a_0^{-1}\|)\leq \delta'.
$$
Thus, by $1^\circ$,
there exists $\xi_0 \in \Cal H$ such that $b\xi_0 = \xi_0 b, \forall b\in B$
and $\|\xi_0 - \xi\| \leq \varepsilon/\|a_0\|$. In addition, since $\xi$
commutes with $a_0$, we may assume $\xi_0$ also does. Let
$\eta_0 = a_0^{1/2} \xi_0$. Then $\eta_0$ still commutes with $B$ and
we have the estimates:
$$
\|\eta_0 - \eta\| = \|a_0^{1/2} \xi_0 - a_0^{1/2} \xi\| \leq
\|a_0^{1/2}\| \|\xi_0 - \xi\| \leq \|a_0^{1/2}\|
(\varepsilon/\|a_0\|) \leq \varepsilon.
$$
\hfill Q.E.D.

\vskip .1in \noindent {\it 4.2. Definitions}. Let $N$ be a
countable decomposable finite von Neumann algebra and $B\subset N$
a von Neumann subalgebra. \vskip .05in \noindent {\it 4.2.1}.
$B\subset N$ is {\it a rigid} (or {\it property} (T)) {\it
embedding} (or, $B$ is a {\it relatively rigid subalgebra} of $N$,
or {\it the pair} $(N,B)$ {\it has the relative property} (T)) if $B\subset
N$ satisfies the equivalent conditions $4.1$. \vskip .05in
\noindent {\it 4.2.2.} If $N$ is a finite factor and $
\varepsilon_0 > 0$ then $B \subset N$ is $\varepsilon_0$-{\it
rigid} if $\exists F \subset N$ finite and $\delta > 0$ such that
if $\phi$ is a completely positive map on $N$ with $\phi(1) \leq
1$, $\tau \circ \phi \leq \tau$ and $\|\phi(x)-x\|_2 \leq \delta,
\forall x \in F$ then $\|\phi(b)-b\|_2 \leq \varepsilon_0, \forall
b \in B, \|b\|\leq 1$.
\vskip .05in
Note that if $N$ is a finite 
factor then an 
embedding $B \subset N$ is rigid if and only if it is
$\varepsilon_0$-rigid $\forall \varepsilon_0 > 0$. We'll see that
if some additional conditions are satisfied (e.g.: $B$ regular
in $N$, in $4.3.2^\circ$; $B, N$
group algebras coming from a group-subgroup situation, in
$5.1$;) then $B \subset N$ $\varepsilon_0$-rigid, for $\varepsilon_0
=1/3$, is in fact sufficient to insure that $B \subset N$ is rigid. 

\proclaim{4.3. Theorem} Let $N$ be a separable type ${\text{\rm
II}}_1$ factor and $B \subset N$ a von Neumann subalgebra.

$1^\circ$. Assume $B \subset N$ is either rigid or
$\varepsilon_0$-rigid, for some $\varepsilon_0 < 1$, with $B$
semi-regular. Then $N'\cap N^\omega = N'\cap (B'\cap N)^\omega$,
for any free ultrafilter $\omega$ on $\Bbb N$. If in addition to
either of the above conditions $B$ also satisfies $B'\cap N = \Cal
Z(B)$ (resp. $B'\cap N = \Bbb C$) then $N$ is non-McDuff (resp.
non-$\Gamma$).

$2^\circ$. Assume that either $B$ is regular in $N$ or that $\Cal
N_N(B)'\cap N^\omega = \Bbb C$. Then $B \subset N$ is rigid if and
only if it is $\varepsilon_0$-rigid for some $\varepsilon_0 \leq
1/3$.
\endproclaim
\noindent {\it Proof}. $1^\circ$. Assume first that $B \subset N$
is rigid. By applying 4.1.2$^\circ$ to the completely positive
maps $\phi = {\text{\rm Ad}}u$ for $u \in \Cal U(N)$, it follows
that for any $\varepsilon > 0$ there exist $\delta > 0$ and $x_1,
x_2, ..., x_n \in N$ such that if $u\in \Cal U(N)$ satisfies
$$
\|ux_i - x_i u\|_2 \leq \delta, \forall i,
$$
then
$$
\|ub - bu\|_2 \leq \varepsilon, \forall b\in B, \|b\|\leq 1.
$$
In particular, $\|vuv^* - u\|_2 \leq \varepsilon, \forall v \in
\Cal U(B)$. Thus, by taking averages over the unitaries $v \in B$,
it follows that $\|E_{B'\cap N}(u)-u\|_2 \leq \varepsilon$. Thus,
if $(u_n) \subset \Cal U(N)$ is a central sequence of unitary
elements in $N$, i.e.,
$$
\underset n \rightarrow \infty \to \lim\|[x, u_n]\|_2 = 0, \forall
x\in N,$$ then
$$
\underset n \rightarrow \infty \to \lim\|u_n - E_{B'\cap
N}(u_n)\|_2 = 0.
$$

Assume now that $B \subset N$ is $\varepsilon_0$-rigid, with
$\varepsilon_0 < 1$, and that $\Cal N(B)'\cap N = \Bbb C$. We
proceed  by contradiction, assuming there exists $u=(u_n)_n \in
\Cal U(N'\cap N^\omega)$ such that $u \not\in (B'\cap N)^\omega$.
By taking a suitable subsequence of $(u_n)$, it follows that there
exists $(v_n)_n \subset \Cal U(N)$ such that $\underset n
\rightarrow \infty \to \lim \|[v_n, x]\|_2=0$, $\forall x\in N$,
and $\|E_{B'\cap N}(v_n)\|_2 \leq c, \forall n$, for some $c < 1$.
It further follows that given any separable von Neumann subalgebra
$P \subset N^\omega$ there exist $k_1 \ll k_2 \ll  ...$ such that
$\underset n \rightarrow \infty \to \lim \|[v_{k_n}, y_n]\|_2=0$,
$\forall y=(y_n)_n \in P$.

Moreover, if $P \subset \Cal N_{N^\omega}(B^\omega)''$, then the
subsequence $v'=(v_{k_n})_n$ can be taken so that to also have
$[E_{{B^\omega}'\cap N^\omega}(v'),y]=0$, $\forall y\in P$. To see
this, let $S \subset \Cal N(B^\omega)$ be a countable set
such that the von Neumann algebra $P_0$ generated by $S$ contains $P$. 
Choose $k_n \uparrow \infty$ so
that $\underset n \rightarrow \infty \to \lim \|[v_{k_n},
w_n]\|_2=0$, $\forall w=(w_n)_n \in S$. We then have
$wE_{{B^\omega}'\cap N^\omega}(v')w^*=wE_{{B^\omega}'\cap
N^\omega}(w^*v'w)w^*=E_{{B^\omega}'\cap N^\omega}(v')$, $\forall
w\in S$. Thus $[E_{{B^\omega}'\cap N^\omega}(v'), S]=0$ implying
$[E_{{B^\omega}'\cap N^\omega}(v'), P_0]=0$ as well.

Now notice that $(B'\cap N)^\omega = {B^\omega}'\cap N^\omega$
(see e.g. [Po2]). As a consequence, since $E_{{B^\omega}'\cap
N^\omega} (x)$ is the element of minimal norm $\| \quad\|_2$ in
$\overline{\text{\rm co}}^w\{w x w^* \mid w\in \Cal U(B^\omega)\}$, 
which in turn can be
realized as a $\|\quad \|_2$-limit of convex combinations of the
form $wxw^*$ with $w$ in a suitable countable subset of $\Cal
U(B^\omega)$, it follows that for any $x \in N^\omega$ there
exists a separable von Neumann subalgebra $P \in B^\omega$ such
that $E_{P'\cap N^\omega}(x)=E_{{B^\omega}'\cap N^\omega}(x)$.
Also, since $\Cal N_{N^\omega}(B^\omega) \supset  \underset n \rightarrow
\omega \to \Pi \Cal N_N(B)$, $\Cal N(B^\omega)''$ follows a factor
and for any $x'\in N^\omega$ there exists a separable von Neumann
subalgebra $P_0$ generated by a countable subset in $\Cal
N(B^\omega)$ such that $P_0 \supset P$ and $E_{P_0'\cap
N^\omega}(x') = \tau(x')1$.

Using all the above, we'll prove the following statement: \vskip
.05in \noindent $(4.3.1')$. If $x\in N^\omega$ then there exists a
subsequence $(v_{k_n})_n$ of $(v_n)_n$ such that $v'=(v_{k_n})_n
\in N^\omega$ satisfies $\|E_{{B^\omega}'\cap N^\omega}(xv')\|_2
=\|E_{{B^\omega}'\cap N^\omega}(x)\|_2 \|E_{{B^\omega}'\cap
N^\omega}(v')\|_2$. \vskip .05in To see this, take first a
separable von Neumann subalgebra $P \subset B^\omega$ such that
$E_{{B^\omega}'\cap N^\omega}(x)=E_{P'\cap N^\omega}(x)$. Then
take $P_0$ a von Neumann algebra generated by a countable subset
in $\Cal N(B^\omega)$ such that $P_0\supset P$ and $E_{P_0'\cap
N^\omega}(x') =\tau(x')1$ where $x'=E_{{B^\omega}'\cap
N^\omega}(x)^*E_{{B^\omega}'\cap N^\omega}(x)$. Since
${B^\omega}'\cap N^\omega \subset P'\cap N^\omega$, if the
subsequence $(v_{k_n})_n$ is chosen such that $[v', P_0]=0$ then
$[v',P]=0$ and we have
$$
E_{{B^\omega}'\cap N^\omega}(xv')=E_{{B^\omega}'\cap N^\omega}
(E_{P'\cap N^\omega}(xv')) =E_{{B^\omega}'\cap N^\omega}
(E_{P'\cap N^\omega}(x)v')
$$
$$
=E_{{B^\omega}'\cap N^\omega} (E_{{B^\omega}'\cap N^\omega}(x)v')
=E_{{B^\omega}'\cap N^\omega}(x)E_{{B^\omega}'\cap N^\omega}(v').
$$
Also, since $y'=E_{{B^\omega}'\cap N^\omega}(v')E_{{B^\omega}'\cap
N^\omega}(v')^*$ satisfies $[y',P_0]=0$, we get
$$
\|E_{{B^\omega}'\cap N^\omega}(xv')\|_2^2=\|E_{{B^\omega}'\cap
N^\omega}(x)E_{{B^\omega}'\cap N^\omega}(v')\|_2^2
=\tau(x'y')=\tau(E_{P_0'\cap N^\omega}(x'y'))
$$
$$
=\tau(E_{P_0'\cap N^\omega}(x')y')=\tau(x')\tau(y')
=\|E_{{B^\omega}'\cap N^\omega}(x)\|_2^2 \|E_{{B^\omega}'\cap
N^\omega}(v')\|_2^2.
$$

Now, by applying recursively $(4.3.1')$, it follows that we can
choose a subsequence $v^1$ of $v=(v_n)_n$, then  $v^2$ of $v^1$,
etc, such that
$$
\|E_{{B^\omega}'\cap N^\omega}(\Pi_j v^j)\|_2 = \Pi_j
\|E_{{B^\omega}'\cap N^\omega}(v^j)\|_2 =\|E_{{B^\omega}'\cap
N^\omega}(v)\|_2^m \leq c^m.
$$

Take $m$ so that $c^m < 1-\varepsilon_0$ and put $w=v^1v^2...
v^m$, $w=(w_n)_n$, with $w_n \in \Cal U(N)$, and
$\phi_n={\text{\rm Ad}}(w_n)$. It follows that
$$
\underset n \rightarrow \omega \to \lim \|E_{B'\cap N}(w_n)\|_2 <
1-\varepsilon_0, \underset n \rightarrow \infty \to \lim
\|\phi_n(x)-x\|_2 =0, \forall x\in N \tag 4.3.1''
$$
By the $\varepsilon_0$-rigidity of $B \subset N$ the second
condition in $(4.3.1'')$ implies that for large enough $n$ we have
$$
\|uw_nu^*-w_n\|_2=\|w_nuw_n^*-u\|_2=\|\phi_n(u)-u\|_2 \leq
\varepsilon_0, \forall u\in \Cal U(B).
$$
Taking convex combinations over $u$, this yields $\|E_{B'\cap
N}(w_n)-w_n\|_2 \leq \varepsilon_0$. Thus $\|E_{B'\cap N}(w_n)\|_2
\geq 1-\varepsilon_0$ for all large enough $n$, contradicting the
first condition in $(4.3.1'')$.

$2^\circ$. We need to show that if $(\psi_n)_n$ are completely
positive maps on $N$ satisfying
$$
\tau \circ \psi_n \leq \tau, \psi_n(1) \leq 1, \forall n,
\underset n \rightarrow \infty \to \lim \|\psi_n(x) - x\|_2=0,
\forall x\in N, \tag a
$$
then $\underset n \rightarrow \infty \to \lim {\text{\rm sup}}
(\{\|\psi_n(b)-b\|_2 \mid b \in B, \|b\|\leq 1\}) =0$. Assume by
contradiction that there exist $(\psi_n)_n$ satisfying $(a)$ but
$$
{\text{\rm inf}}_n ( {\text{\rm sup}} \{\|\psi_n(b)-b\|_2 \mid b
\in B, \|b\|\leq 1\})
>0. \tag b
$$
Note that by the $\varepsilon_0$-rigidity of $B \subset N$, $(a)$
implies
$$
\underset n \rightarrow \infty \to \limsup ({\text{\rm sup}}
\{\|\psi_n(b)-b\|_2 \mid b \in B, \|b\|\leq 1\}) \leq
\varepsilon_0. \tag c
$$

If $(\psi_n)_n$ satisfies $\tau \circ \psi_n \leq \tau, \psi_n(1)
\leq 1, \forall n$ in $(a)$ then
$$
\Psi((x_n)_n) \overset\text{\rm def} \to = (\phi_n(x_n))_n,
(x_n)_n \in N^\omega, \tag d
$$
gives a well defined completely positive map $\Psi$ on $N^\omega$
with $\tau\circ \Psi \leq \tau$, $\Psi(1) \leq 1$. Thus, the fixed
point set $(N^\omega)^\Psi \overset\text{\rm def} \to =\{x \in
N^\omega \mid \Psi(x) = x \}$ is a von Neumann algebra. If
$(\psi_n)_n$ also satisfies the last condition in $(a)$, then $N
\subset (N^\omega)^\Psi$. In particular $\Psi(1)=1$ which together
with $\|T_\Psi\| \leq 1$ implies ${T_\Psi}^*(\hat{1})=\hat{1}$,
equivalently $\Psi^*(1)=1$, i.e., $\tau \circ \Psi = \tau$.

If in addition to $(a)$ the sequence $(\psi_n)_n$ satisfies $(b)$,
then $B^\omega \not\subset (N^\omega)^\Psi$. Let us prove that the
$\varepsilon_0$-rigidity of $B \subset N$ entails
$$
B^\omega \subset_{\varepsilon_0} (N^\omega)^\Psi. \tag e
$$
For $\psi$ a map on an algebra denote by $\psi^m$ the $m$-time
composition $\psi \circ \psi ... \circ \psi$. Then note that for
each $m \geq 1$ the sequence $(\psi^m_n)_n$ still satisfies $(a)$,
and thus, by $\varepsilon_0$-rigidity, $(c)$ as well. Thus
$$
\|\Psi^k(b)-b\|_2 \leq \varepsilon_0, \forall b \in B^\omega,
\|b\|\leq 1.
$$
But by von Neumann's ergodic theorem applied to $\Psi$ and $x \in
N^\omega$, we have
$$
\underset n \rightarrow \infty \to \lim \|m^{-1}\Sigma_{k=1}^m
\Psi^k(x) - E_{(N^\omega)^\Psi}(x)\|_2 = 0, \tag f
$$
which together with the previous estimate shows that for $x = b
\in B^\omega, \|b\|\leq 1,$ we have
$\|E_{(N^\omega)^\Psi}(b)-b\|_2 \leq \varepsilon_0$, i.e., $(e)$.

The assumption $\Cal N(B)'\cap N^\omega=\Bbb C$ implies in particular
that $N'\cap N^\omega = \Bbb C \subset (N^\omega)^\Psi$. We next
prove that $B$ regular in $N$ implies $N'\cap N^\omega \subset
(N^\omega)^\Psi$ as well, for any $\Psi$ on $N^\omega$ associated
as in $(d)$ to a sequence $(\psi_n)_n$ satisfying $(a)$. Denote
$P=(N^\omega)^\Psi$ and assume by contradiction that $N'\cap
N^\omega \nsubseteq P$. Since $N'\cap N^\omega$ and $P$ make a
commuting square, this implies there exists $x \in N'\cap
N^\omega$, $x\neq 0$, such that $E_P(x)=0$. Moreover, we may
assume $x=(x_n)_n$ satisfies $x_n=x_n^*, \|x_n\|_2 =1, \forall n$.

By using $(f)$, we can choose ``rapidly'' increasing $k_1 \ll k_2
\ll ...$ and ``slowly'' non-decreasing $m_1 \leq m_2 \leq...$ such
that the sequence of completely positive maps $\psi'_n=(m_n)^{-1}
\Sigma_{j=1}^{m_n} \psi_{k_n}^j$ satisfies $(a)$ and $\underset n
\rightarrow \infty \to \lim \|\psi'_n(x'_n)\|_2=0$, with
$\underset n \rightarrow \infty \to \lim \|[x'_n, y]\|_2=0,
\forall y\in N$, $\underset n \rightarrow \infty \to \lim
\tau((x'_n)^k) = \tau (x^k), \forall k$, where $x'_n = x_{k_n}$.

Denote by $\Psi_1$ the completely positive map on $N^\omega$
associated with $(\psi'_n)_n$, as in $(d)$, and put $X=X_1 =
(x'_n)_n \in N^\omega$. Since each separable von Neumann
subalgebra of $N^\omega$ is contained in a separable factor and
since for each separable $Q \subset N^\omega$ there exists $j_1 \ll 
j_2 \ll ...$ such that $X'=(x'_{j_n})_n \in Q'\cap N^\omega$, it
follows that there exist separable factors $Q_0=N \subset Q_1
\subset ... \subset Q_{m-1}$ in $N^\omega$ and consecutive
subsequences of indices $(j,1) < (j,2) < ...$, for $j=1, 2, ...,
m$, with $(1,n)=n$, such that $X_j=(x'_{j,n})_n \in N^\omega$
satisfy $X_1, X_2, ..., X_j \in Q_j$, $[Q_j, X_{j+1}]=0$, for $0
\leq j \leq m-1$. Denote by $\Psi_j$ the completely positive map
on $N^\omega$ associated with $(\psi'_{j,n})_n$, noticing that
each one of these sequences checks $(a)$. Thus for each $j = 1, 2,
..., m$ we have $\Psi_j(x)=x, \forall x\in N$ and $\Psi_j(X_j)=0$.
Moreover, the von Neumann algebra generated by $X_1, X_2, ...,
X_m$ in $N^\omega$ is isomorphic to the tensor power $(A(X),
\tau)^{\otimes m}$, where $A(X)$ is the von Neumann algebra
generated by $X\in N^\omega$.

Let $\tilde{X} = m^{-1/2} \Sigma_{j=1}^m X_j$ and $\tilde{\Psi} =
m^{-1} \Sigma_{j=1}^m \Psi_j$. Let $P_j=(N^\omega)^{\Psi_j}, 1\leq
j \leq m,$ and $\tilde{P}=(N^\omega)^{\tilde{P}}$. By $(a)-(e)$,
$\tilde{P}, P_j$ are von Neumann algebras containing $N$ and
$B^\omega \subset_{\varepsilon_0} P_j, \tilde{P}$. Moreover, since
by convexity we have $\tilde{\Psi}(Y)=Y$ iff $\Psi_j(Y)=Y$,
$\forall j$, it follows that $\tilde{P}=\cap_j P_j$. Thus, since
$\Psi_j(X_j)=0$ implies $E_{P_j}(X_j)=0$, it follows that
$E_{\tilde{P}}(\tilde{X}) = 0$.

But by the central limit theorem, as $m \rightarrow \infty$,
$\tilde{X}$ gets closer and closer (in distribution) to an element
$Y=Y^*$ with Gaussian spectral distribution, independently of $X$.
Let  $Y' = Ye_{[-2,2]}(Y)$ and 
$\|Y\|_2^2 =t$. By using Mathematica, one finds $t> 0.731$). 
Thus, for large enough $m$,
$\tilde{X}'=\tilde{X}e_{[-2,2]}(\tilde{X})$ satisfies
$\|\tilde{X}'\|_2^2 = t_-$ with $t_-$ close to $t$. Let $\tilde{X}'' =
\tilde{X}-\tilde{X}'$ and note that $\tilde{X}' \tilde{X}''=0$,
so $\|\tilde{X}'\|_2^2+ \|\tilde{X}''\|_2^2 = \|\tilde{X}\|_2^2
=1$. Also,
$$
E_{\tilde{P}}(\tilde{X}')=E_{\tilde{P}}(\tilde{X}-\tilde{X}'')
=-E_{\tilde{P}}(\tilde{X}''),
$$
implying that $\|E_{\tilde{P}}(\tilde{X}')\|^2_2 \leq
\|\tilde{X}''\|_2^2 = 1-t_-$. Altogether
$$
\|\tilde{X}' - E_{\tilde{P}}(\tilde{X}')\|_2^2 =
\|\tilde{X}'\|_2^2 - \|E_{\tilde{P}}(\tilde{X}')\|_2^2  
\geq  2t_- - 1 .
$$
Since $\tilde{X}_1 \in N'\cap N^\omega \subset B^\omega$ and
$\|\tilde{X}_1\| = 2$, if we take $\tilde{X}_0 = \tilde{X}'/2$ then 
$\|\tilde{X}_0\| = 1$ and $\|\tilde{X}_0 - E_{\tilde{P}}(\tilde{X}_0\|_2^2 
= (2t_- - 1)/4 > (1/3)^2$, this contradicts $B^\omega
\subset_{1/3} \tilde{P}$.

This finishes the proof of the fact that $N'\cap N^\omega \subset
(N^\omega)^\Psi$, independently of $\Psi$, for arbitrary
$(\psi_n)_n$ checking $(a)$. Thus $P=\cap_i (N^\omega)^{\Psi_i}$,
where $\Psi_i, i\in \Cal I,$ is the family of all completely
positive maps on $N^\omega$ coming from sequences $(\psi_{i,n})_n$
satisfying $(a)$, still satisfies $N, N'\cap N^\omega \subset P$.
Let us show that this newly designated $P$ still satisfies
$B^\omega \subset_{\varepsilon_0} P$. To see this, take a finite
subset $I \subset \Cal I$ and consider the sequence $\psi_{I,n} =
|I|^{-1} \Sigma_i \psi_{i,n}$, which clearly satisfies $(a)$.
Thus, the associated completely positive map $\Psi_I$ on
$N^\omega$ satisfies
$$
\|E_{P_I}(b)-b\|_2 \leq \varepsilon_0, \forall b\in B^\omega,
\|b\|\leq 1 .
$$
where $P_I\overset\text{\rm def} \to 
=(N^\omega)^{\Psi_I}$. Since $|I|^{-1} \Sigma_i \Psi_i (x) = x$
iff $\Psi_i(x)=x$, $\forall i \in I$, we have $P_I = \underset i
\in I \to \cap (N^\omega)^{\Psi_i}$. But $P_I \downarrow P$ as $I
\uparrow \Cal I$, implying that $\|E_P(b)-b\|_2 \leq
\varepsilon_0, \forall b$, as well.

Denote $\Cal U_0=\Cal N(B)\cup \Cal U(\Cal N(B)'\cap
(B^\omega)'\cap N^\omega)$, $N_0=\Cal U_0''$ and notice that
$v(B^\omega)v^* = B^\omega, \forall v \in \Cal U_0$. Also, if we
let $M=N^\omega$, $Q=B^\omega$, then by $1^\circ$ both the
assumption $\Cal N(B)'\cap N^\omega = \Bbb C$ and $\Cal
N_N(B)''=N$ imply that $\Cal U_0 \subset P$ and $N_0'\cap M =
Q'\cap \Cal Z(N_0)$ are satisfied. Thus, A.3 applies to get a
non-zero projection $p \in Q'\cap \Cal Z(N_0)$ such that $Qp
\subset P$. In the case $\Cal N(B)'\cap N^\omega=\Bbb C$, this
implies $p=1$ and we get $B^\omega \subset P$, a contradiction
which finishes the proof under this assumption.

If $B$ is regular in $N$, then the group $\Cal N(B)=\Cal N(B \vee
B'\cap N)$ generates the factor $N$, a fact that is easily seen to
imply $\Cal N_{N^\omega}(B^\omega)' \cap N^\omega=\Bbb C$. This
implies there exists a countable subgroup $\Cal U_1 \subset \Cal
N(B^\omega)$ such that $\tau(p)1$ is a limit in the norm-$\|\quad
\|_2$ of convex combinations of elements of the form $u_1pu_1^*$,
$u_1 \in \Cal U_1$. Let then $(\psi_n)_n$ be the sequence of
completely positive maps satisfying $(a)-(b)$ at the beginning of
the proof, with $b_n \in B, \|b_n\|\leq 1$, $\|\psi_n(b_n)-b_n\|_2
\geq c >0$, $\forall n$. If we choose a sufficiently rapidly
increasing $k_1\ll k_2 \ll ...$, then the completely positive map
$\Psi'$ associated with $(\psi_{k_n})_n$ as in $(d)$ has both $N$
and $\Cal U_1$ in the fixed point algebra $(N^\omega)^{\Psi'}$.
But since $P \subset (N^\omega)^{\Psi'}$, it follows that
$(N^\omega)^{\Psi'}$ contains $B^\omega p$, and thus $u_1(B^\omega
p)u_1^*=B^\omega(u_1pu_1^*), \forall u_1 \in \Cal U_1$ as well. This
implies $B^\omega \subset (N^\omega)^{\Psi'}$, contradicting
$\|\Psi'(b')-b'\|_2 \geq c >0$, where $b'=(b_{k_n})_n \in
B^\omega$. \hfill Q.E.D.

\proclaim{4.4. Theorem} Let $N$ be a type ${\text{\rm II}}_1$
factor and $B\subset N$ a von Neumann subalgebra such that $B'\cap
N = \Cal Z(B)$ and such that the normalizer of $B$ in $N$, $\Cal
N(B)$, acts ergodically on the center of $B$. Let $\Cal G_B
\subset {\text{\rm Aut}}N$ be the group generated by ${\text{\rm
Int}}N$ and by the automorphisms of $N$ that leave all elements of
$B$ fixed. If $B\subset N$ is $\varepsilon_0$-rigid for some
$\varepsilon_0 < 1$ then $\Cal G_B$ is open and closed in
${\text{\rm Aut}}N$. Thus, ${\text{\rm Aut}}N/\Cal G_B$ is
countable.
\endproclaim
\noindent {\it Proof}. By applying condition 4.2.2$^\circ$ to the
completely positive maps $\theta \in {\text{\rm Aut}}N$, it
follows that there exist $\delta > 0$ and $x_1, x_2, ..., x_n \in
N$ such that if $\|\theta(x_i) - x_i\|_2 \leq \delta$ then
$$
\|\theta(u) - u\|_2 \leq \varepsilon_0, \forall u\in \Cal U(B).
$$

Thus, if $k$ denotes the unique element of minimal norm
$\|\quad\|_2$ in $K=\overline{\text{\rm co}}^w \{ \theta(u)u^*\mid
u\in \Cal U(B)\}$ then $\|k-1\|_2\leq \varepsilon_0$ and thus
$k\neq 0$. Also, since $\theta(u)Ku^* \subset K$ and
$\|\theta(u)ku^*\|_2 = \|k\|_2, \forall u \in \Cal U(B)$, by the
uniqueness of $k$ it follows that $\theta(u)ku^*=u$, or
equivalently $\theta(u)k = ku$, for all $u\in \Cal U(B)$. By a
standard trick, if $v\in N$ is the (non-zero) partial isometry in
the polar decomposition of $k$, then $\theta(u)v=vu, \forall u\in
\Cal U(B)$, $v^*v \in B'\cap N = \Cal Z(B), vv^* \in
\theta(B)'\cap N = \theta(\Cal Z(B))$. Since $\Cal N(B)$ acts
ergodically on $\Cal Z(B)$ (equivalently, $\Cal N(B)'\cap N = \Bbb
C$), there exist finitely many partial isometries $v_0 = v^*v,
v_1, v_2, ..., v_n \in N$ such that $v_i^*v_i = v^*v, 0\leq i \leq
n-1$, $v_n^*v_n \in \Cal Z(B)v^*v$ and $v_iv_i^* \in \Cal Z(B),
v_iBv_i^* = Bv_iv_i^*, \forall i$.

If we then define $w = \Sigma_i \theta(v_i) v v_i^*$, an easy
calculation shows that $w$ is a unitary element and $wbw^* =
\theta(b), \forall b\in B$. \hfill Q.E.D.

\proclaim{4.5. Proposition} Let $N$ be a type $\text{\rm II}_1$
factor and $B\subset N$ a rigid embedding.

$1^\circ$. For any $\varepsilon_0 >0$ there exist $F_0\subset N$
and $\delta_0 > 0$ such that if $N_0\subset N$ is a subfactor with
$B\subset N_0$ and $F_0 \subset_{\delta_0} N_0$, then $B\subset
N_0$ is $\varepsilon_0$-rigid. In particular, if $N_k\subset N,
k\geq 1$ is an increasing sequence of subfactors such that
$B\subset N_k, \forall k$, and $\overline{\cup_k N_k} = N$, then
for any $\varepsilon_0 > 0$ there exists $k_0$ such that $B\subset
N_k$ is $\varepsilon_0$-rigid $\forall k \geq k_0$.

$2^\circ$. Assume in addition that $B$ is regular in $N$ and
$B'\cap N = \Cal Z(B)$. For any $\varepsilon > 0$ there exist a
finite subset $F \subset N$ and $\delta > 0$ such that if $N_0
\subset N$ is a subfactor with $N_0'\cap N = \Bbb C$ and $F
\subset_\delta N_0$ then there exists $u\in \Cal U(N)$ such that
$\|u-1\|_2 \leq \varepsilon$ and $uBu^* \subset N_0$, with
$uBu^*\subset N_0$ rigid embedding. If in addition $N_0
\supset B$ then one can take $u=1$. In particular, if $N_k \subset N$ is an
increasing sequence of subfactors with $N_k'\cap N = \Bbb C$ and
$N_k \uparrow N$ then there exist $k_0$ such that $u_kBu_k^*
\subset N_k$ rigid, $\forall k \geq k_0$, for some $u_k\in \Cal
U(N)$, $\|u_k-1\|_2 \rightarrow 0$, and such that if 
$N_k \supset B, \forall k$, then $B \subset N_k$ rigid $\forall k
\geq k_0$.
\endproclaim
\noindent {\it Proof}. $1^\circ$. With the notations of
$4.1.2^\circ$, for the critical sets $F(\varepsilon')$ and
constants $\delta(\varepsilon')$ for $B \subset N$, let
$F_0=F(\varepsilon_0)$ and $\delta_0 = \delta(\varepsilon_0)/2$.
Let $N_0\subset N$ be a von Neumann algebra with $B\subset N_0$,
$\|E_{N_0}(y)-y\|_2 \leq \delta_0, \forall y \in F_0$. We want to
prove that $B \subset N_0$ is $\varepsilon_0$-rigid by showing
that if $\phi_0$ is a completely positive map on $N_0$ with
$\phi_0(1) \leq 1, \tau\circ \phi_0 \leq \tau$ and
$$
\|\phi_0(y_0)-y_0\|_2 \leq \delta(\varepsilon_0)/2, \forall y_0
\in E_{N_0}(F_0),
$$
then $\|\phi_0(b)-b\|_2 \leq \varepsilon_0, \forall b \in B,
\|b\|\leq 1$. To this end let $\phi = \phi_0 \circ E_{N_0}$, which
we regard as a completely positive map from $N$ into $N$ ($\supset
N_0$). Clearly $\phi(1) \leq 1, \tau \circ \phi \leq \tau$. Also,
for $y \in F(\varepsilon_0)$ we have
$$
\|\phi(y)-y\|_2 \leq \|\phi_0(E_{N_0}(y)) - E_{N_0}(y)\|_2 +
\|E_{N_0}(y)-y\|_2 \leq \delta(\varepsilon_0).
$$
Thus, $\|\phi(b)-b\|_2 \leq \varepsilon_0, \forall b\in B,
\|b\|\leq 1$. Since for $b \in B$ we have $\phi(b)=\phi_0(b)$, we
are done.

$2^\circ$. By applying condition $4.1.2^\circ$ to the completely
positive maps $E_{N_0}$, it follows that if we denote
$\varepsilon(N_0) = {\text{\rm sup}} \{\|E_{N_0}(b)-b\|_2\mid b\in
B, \|b\|\leq 1\}$, then $\varepsilon(N_0) \rightarrow 0$ as
$E_{N_0} \rightarrow  id_N$. Thus, by Theorem A.2 it follows that
there exist unitary elements $u=u(N_0) \in N$ such that $uBu^*
\subset N_0$ and $\|u-1\|_2 \rightarrow 0$. Moreover, by $1^\circ$
above and 4.3.2$^\circ$, it follows that $uBu^* \subset N_0$
(equivalently, $B \subset uN_0u^*$) is a rigid embedding when
$N_0$ is close enough to $N$ on an appropriate finite set of
elements. The fact that $B$ is still regular in $N_0$ is a
consequence of ([JPo]). The last part is now trivial. \hfill
Q.E.D.

\proclaim{4.6. Proposition} $1^\circ$. $(B_i \subset N_i)$ are
rigid embeddings for $i=1, 2$ if and only if $(B_1
\overline{\otimes} B_2 \subset N_1 \overline{\otimes} N_2)$ is a
rigid embedding.

$2^\circ$. Let
$B \subset N_0 \subset N$. If $B \subset N_0$
is a rigid embedding
then $B \subset N$
is a rigid embedding.
Conversely, if we assume
$N_0\subset N$ is a $\lambda$-Markov inclusion (${\text{\rm [Po2]}}$),
i.e., $N$ has an orthonormal basis $\{m_j\}_j$ with $\Sigma m_jm_j^* =
\lambda^{-1}$ for some constant $\lambda > 0$
(e.g., if $N, N_0$ are factors and $[N:N_0] < \infty$) then
$B \subset N$ rigid embedding,
implies $B\subset N_0$
is a rigid embedding.

$3^\circ$. Let $B \subset B_0 \subset N$.
If $B_0 \subset N$ is a rigid embedding, then
$B\subset N$ is a rigid
embedding. Conversely, if $B_0$ has a finite
orthonormal basis with respect to $B$ and $B\subset N$ is a rigid
embedding,
then $B_0 \subset N$ is a rigid embedding.
\endproclaim
\noindent {\it Proof}. $1^\circ$. Assume first that $(B_i \subset
N_i)$ are rigid embeddings$/\tau_i$, for $i=1,2$. Let $\varepsilon
> 0$ and $F'_i(\varepsilon/2), \delta_i'(\varepsilon/2)$ be the
critical sets and constants for $B_i \subset N_i$, as given by
4.1.1$^\circ$, for $\varepsilon/2$. Define $F' = F'_1 \otimes 1
\cup 1 \otimes F'_2, \delta' = {\text{\rm min}} \{\delta_1',
\delta_2'\}$.

Put $N=N_1\overline{\otimes} N_2, B=B_1\overline{\otimes} B_2$.
Let $\Cal H$ be a Hilbert $N$-bimodule with a vector $\xi \in \Cal H$
which satisfies conditions 4.1.1$^\circ$
with respect to the trace $\tau_1\otimes \tau_2$,
for $F', \delta'$. In particular, $\Cal H$ is a Hilbert
$N_i$ bimodule, for $i=1,2$. Thus, if we denote by $p_i$ the
orthogonal projection of $\Cal H$ onto the Hilbert
subspace of all vectors in $\Cal H$ that commute with $B_i$,
then $\|\xi - p_i(\xi)\|_2 \leq \varepsilon/2,
i=1,2$, for any vector $\xi \in \Cal H$ that satisfies
4.1.1$^\circ$ for the above $F', \delta'$. But $p_1$ and $p_2$
are commuting projections and $p_1p_2$ projects onto the
Hilbert subspace of vectors commuting
with both $B_1$ and $B_2$, i.e., onto the Hilbert
subspace of vectors commuting with $B$. Since
$$
\|\xi - p_1p_2(\xi)\|
\leq \|\xi - p_1(\xi)\| + \|p_1(\xi) - p_1(p_2(\xi))\|
$$
$$
\leq  \|\xi-p_1(\xi)\| + \|\xi - p_2(\xi)\| \leq \varepsilon,
$$
it follows that $B\subset N$ satisfies $4.1.1^\circ$.

Assume now that $B\subset N$ satisfies 4.1.2$^\circ$
for some trace $\tau$. Since
$N_1 \otimes N_2$ is a dense
$*$-subalgebra in $N$, by using Kaplanski's density theorem
and the fact that in 4.1.2$^\circ$ we only have to deal
with completely positive maps $\phi$
satisfying $\tau\circ \phi \leq \tau, \phi(1) \leq 1$,
it follows that we may assume the critical set $F'(\varepsilon)$ is contained
in $N_1\otimes N_2$ (by diminishing if necessary
the corresponding $\delta'(\varepsilon)$).

Let $F'_i \subset N_i$ be finite subsets such that
$F'\subset {\text{\rm sp}} F'_1\otimes F'_2$.
There clearly exist $\delta_i' > 0$ such that
if $\phi_i$ are completely positive maps
on $N_i$ with $\tau \circ \phi_i \leq \tau, \phi_i(1) \leq 1$ and
$\|\phi_i(x_i)-x_i\|_2 \leq \delta'_i, \forall x_i \in F'_i, i=1,2$,
then $\phi=\phi_1 \otimes \phi_2$ satisfies $\|\phi(x) - x\|_2
\leq \delta', \forall x\in F'$. Thus, $\|\phi(b) - b\|_2 \leq \varepsilon,
\forall b\in B, \|b\|\leq 1$. Taking $b\in B_i$, we get
$\|\phi_i(b) - b\|_2 \leq \varepsilon,
\forall b\in B_i, \|b\|\leq 1, i=1,2.$

$2^\circ$. The implication $\implies$ follows by noticing
that if $\phi$ is a completely positive map on $N$ such that
$\phi(1) \leq 1$ and $\tau\circ \phi \leq \tau$ then for
$x \in N_0$ we have $\|E_{N_0}(\phi(x)) - x\|_2 \leq \|\phi(x)-x\|_2$
while for $b\in B, \|b\|\leq 1$, we have
$$
\|\phi(b) - b \|_2^2 \leq  \|E_{N_0}(\phi(b)) - b \|_2^2
+ 2 \|E_{N_0}(\phi(b)) - b \|_2.
$$
Thus, if 4.1.2$^\circ$ is satisfied for $B\subset N_0$
with critical set $F_0(\varepsilon)$ and constant $\delta_0(\varepsilon)$,
then 4.1.2$^\circ$ holds true for $B \subset N$ for the same
set $F_0$ but constant $\delta(\varepsilon) =
\delta_0(\varepsilon)^2/3$.

To prove the opposite implication,
let $e=e_{N_0}$ be the Jones
projection corresponding to $N_0 \subset N$ and
$N_1 = \langle N, e \rangle$ the basic construction. Since $N_0 \subset N$
is $\lambda$-Markov, there exists a unique trace $\tau$ on
$N_1$ extending the trace $\tau$ of $N$ and such that
$E^\tau_N(e)=\lambda 1$.

We may assume  $1$ belongs to the orthonormal basis
$\{m_j\}_j$ of $N$ over $N_0$.
Note that $x=\Sigma_j m_jE_N({m_j}^*x),
\forall x\in N$. Any element $X \in N_1$ can be
uniquely written in the form
$X=\Sigma_{i,j} m_i x_{ij} em_j^*$ for some $x_{ij} \in p_iN_0p_j$,
where $p_i=E_{N_0}(m_i^*m_i) \in \Cal P(N_0)$. Also,
if $x\in N$ then
$$
x = (\Sigma_i m_iem_i^*) x (\Sigma_j m_jem_j^*)= \Sigma_{i,j}
m_iE_{N_0}(m_i^*x m_j)e{m_j}^* \tag 4.6.2'
$$

For each completely positive map $\phi$ on $N_0$
define $\tilde{\phi}$ on $N_1$ by
$$
\tilde{\phi}(\Sigma_{i,j} m_i x_{ij}em_j^*) = \Sigma_{i,j} m_i
\phi(x_{ij})e m_j^* \tag 4.6.2''
$$

Note that if $X = \Sigma_{i,j} m_i x_{ij}e m_j^*\geq 0$
and $\tau \circ \phi \leq \tau$ then
$$
\tau(\tilde{\phi}(X)= \tau(\tilde{\phi}(\Sigma_{i,j} m_i x_{ij}e m_j^*))
= \lambda \Sigma_{i,j} \tau (m_i \phi(x_{ij})m_j^*)
$$
$$
= \lambda \Sigma_{i,j} \tau (m_i \phi(x_{ij})m_j^*)
= \lambda \Sigma_{i} \tau ( \phi(x_{ii})p_i)
$$
$$
\leq  \lambda \Sigma_{i} \tau ( \phi(x_{ii})) \leq \lambda \Sigma_i
\tau(x_{ii}) = \tau(X).
$$
Similarly, if $\phi(1) \leq 1$ then $\tilde{\phi}(1) \leq 1$.

Let now $\varepsilon > 0$ be given. Let $F=F(\lambda
\varepsilon^2/3), \delta= \delta(\lambda \varepsilon^2/3)$ be the
critical set and constant for $B\subset N$, corresponding to
$\lambda \varepsilon^2/3$. Let $F_0 = \{E_{N_0}(m_i^*x m_j) \mid
\forall i,j, \forall x\in F \}$. Formulas $(4.6.2'), (4.6.2'')$
above show that there exists $\delta_0 > 0$ such that if
$\|\phi(x)-x\|_2 \leq \delta_0, \forall x\in F_0$ then
$\|\tilde{\phi}(x)-x\|_2 \leq \delta, \forall x\in F$.

We claim that $F_0, \delta_0$ give the critical set and constant
for $B\subset N_0$, corresponding to $\varepsilon$. To see this,
note first that by the proof of $\implies$ above we get
$\|\tilde{\phi}(b)-b\|_2 \leq \lambda^{1/2}\varepsilon, \forall
b\in B, \|b\|\leq 1$. By $(4.6.2'')$ this gives
$$
\lambda^{1/2} \|\phi(b)-b\|_2 = \|(\phi(b)-b)e\|_2
$$
$$
\leq \|\tilde{\phi}(b)-b\|_2 \leq \lambda^{1/2} \varepsilon.
$$

$3^\circ$. The first implication is trivial. The opposite
implication is equally evident, if
we take the critical set $F_0(\varepsilon)$
and constant $\delta_0(\varepsilon)$ for $B_0
\subset N$ to be
defined as follows:
We first choose $\delta_1>0$ with the
property that if $\phi$ is a completely positive map
on $N$ with $\tau \circ \phi \leq \tau, \phi(1) \leq 1$ and
$\|\phi(b)-b\|_2 \leq \delta_1, \forall b\in B, \|b\|\leq 1$
and $\|\phi(b^0_j)-b_j^0\|_2 \leq \delta_1$, then
$\|\phi(b_0)-b_0\|_2 \leq \varepsilon, \forall b_0\in B_0, \|b_0\|\leq 1$
($\{b^0_j\}_j$ denotes here the orthonormal basis of $B_0$ over $B$).
We then define
$F_0(\varepsilon)=F(\delta_1) \cup \{b^0_j\}_j$ and put
$\delta_0(\varepsilon)=\delta_1$.
\hfill Q.E.D.

\proclaim{4.7. Proposition} $1^\circ$. If $B\subset N$ and
$\{p_n\}_n$ is an increasing sequence of projections in $N$, with
$p_n \uparrow 1$, which lye either in $B$ or in $B'\cap N$, and
with the property that $p_nBp_n \subset p_nNp_n$ are rigid
embeddings, $\forall n$, then $B \subset N$ is a rigid embedding.
In particular, if $B$ is atomic then $B \subset N$ is rigid.

$2^\circ$. If
$B\subset N$ is a rigid embedding and
$p\in \Cal P(B)$ or $p\in \Cal P(B'\cap N)$ then $pBp
\subset pNp$ is a rigid embedding.

$3^\circ$. Let $B\subset N$ and $p\in \Cal P(B)$. Assume there
exist partial isometries $\{v_n\}_{n \geq 0}\subset N$ such that
$v_n^*v_n \in pBp,$ $v_nv_n^* \in B,$ $v_nBv_n^* = v_nv_n^* B
v_nv_n^*, \forall n\geq 0, \Sigma_n v_nv_n^* = 1$ and $B \subset
(\{v_n\}_n \cup pBp)''$. If $pBp \subset pNp$ is a rigid embedding
then $B\subset N$ is a rigid embedding.
\endproclaim
\noindent
{\it Proof}.
$1^\circ$. Notice first
that if $\phi$ is completely positive on $N$ and $\tau\circ
\phi \leq \tau, \phi(1) \leq 1$ then
$\tau(p_n\phi(p_n x p_n)p_n) \leq \tau(\phi(p_nxp_n))
\leq \tau(p_nxp_n), \forall x \geq 0,$ and $p_n\phi(p_n)p_n \leq p_n$.
Then we simply take the critical set and constant
for $B \subset N$ to be the critical set and constant for
$p_nBp_n \subset p_n N p_n$, with $n$ sufficiently large,
and apply the above to deduce that for $\phi$
satisfying the conditions for this set and constant,
$p_n \phi(p_n \cdot p_n)p_n$ follows uniformly close to the
identity on the unit ball of $p_nBp_n$.

The case when $B$ is atomic is now trivial, by first applying
4.6.3$^\circ$ and then the first part of the proof.

$2^\circ$. The statement is clearly true in
case $p\in \Cal Z(\Cal N)$. Assume next that $p\in \Cal P(B)$.
By part $1^\circ$ above, we may suppose $pBp$ has some non-atomic
part.

By noticing that there exist projections $z_n \in \Cal Z(N)$
with $z_n \uparrow 1$ such that each $z_n$ is a sum of finitely many
projections in $Bz_n$ which are majorized by $pz_n$ in $B$,
by $1^\circ$ above
it is sufficient to prove the case when there exist partial
isometries $v_0 = p, v_1, v_2, ..., v_n \in B$ such that
$v_i^*v_i \leq p, \forall i$,
$\Sigma_i v_iv_i^* = 1$.

Let then $\varepsilon > 0$. Let $F=F(\varepsilon \tau(p))$ and
$\delta=\delta(\varepsilon \tau(p))$ be given by 4.1.2$^\circ$ for
the inclusion $B \subset N$. Let also $F_0 = \{v_i^*xv_j \mid
1\leq i,j \leq n, x\in F\}$. We show that $F_0$ and
$\delta_0=\delta$ are good for $pBp \subset pNp$. Thus, let $\phi$
be a completely positive map on $pNp$ such that $\phi(p) \leq p$,
$\tau_p \circ \phi \leq \tau_p$ and $\|\phi(y) - y \|_2 \leq
\delta_0, \forall y\in F_0$. Define $\tilde{\phi}(x) =
\Sigma_{i,j} v_i\phi(v_i^*xv_j)v_j^*$. Like in the proof of
4.6.1$^\circ$, we get $\tau\circ \tilde{\phi}(x) \leq \tau(x),
\forall x\in N$ and $\tilde{\phi}(1) \leq 1$.

An easy calculation shows that $\|\tilde{\phi}(x) - x\|_2
\leq \delta$ for $x\in F$. Thus, $\|\tilde{\phi}(b)-b\|_2
\leq \varepsilon \tau(p), \forall b\in B, \|b\| \leq 1$.
But this implies $\|\phi(pbp) - pbp\|_2 \leq \varepsilon \|p\|_2, \forall b\in
B, \|b\|\leq 1$ as well.

If the projection $p$ lies in $B'\cap N$ then by the last part of
4.6.3$^\circ$ the subalgebra $B_0\subset N$ generated by $B$ and
$\{1, p\}$ is rigid in $N$. But then we apply the first part to
get $pBp = pB_0p $ is rigid in $pNp$.

$3^\circ$. By 1$^\circ$ above, it is sufficient to prove the case
when the set $\{v_i\}_i$ is finite. Let $\varepsilon > 0$ and
$F_p=F(\varepsilon'), \delta_p=\delta(\varepsilon')$ be given by
condition 4.1.2$^\circ$, for $pBp \subset pNp$ and $\varepsilon' =
\varepsilon ({\text{\rm min}}_i\tau(v_iv_i^*)/2)^2$. Then define
$F_0 = F_p \cup \{v_i\}_{0\leq i\leq n}$. If $\phi$ is a
completely positive map on $N$ such that $\|\phi(x)-x\|_2 \leq
\delta_0$ with $\delta_0 \leq \delta_p \tau(p)^{1/2}, \forall x\in
F_0,$ then in particular we have $\|\phi(x)-x\|_{2,p} \leq
\delta_p, \forall x\in F_p$. Thus, $\|p\phi(b)p-b \|_2 \leq
\varepsilon ({\text{\rm min}}_i \|v_iv_i^*\|_2/2)^2, \forall b\in
pBp, \|pbp\|\leq 1$. This easily gives $\|\phi(b)-b\|_2 \leq
\varepsilon$ for all $b$ in the von Neumann algebra $B_0 =
\Sigma_{i,j} v_iBv_j^*,$ generated by $pBp$ and $\{v_i\}_{0\leq
i\leq n}$, with $\|b\|\leq 1$ (in fact, even for all $b \in B_0$
that satisfy $\|v_i^*bv_j\|\leq 1, \forall i,j$). Thus, $B_0
\subset N$ is rigid, so by $4.6.3^\circ$, $B\subset N$ is rigid as
well. \hfill Q.E.D.

\heading 5. More on rigid embeddings. \endheading

In this Section we produce examples of rigid
inclusions of algebras, by using
results of Kazhdan ([Kaz]) and Valette
([Va]), which provide examples
of property (T) inclusions of groups, and the result below,
which establishes the link between the property (T)
for an inclusion of groups and the property (T) (rigidity) for
the inclusion of the corresponding group von Neumann
algebras (as defined in $(4.2)$).

\proclaim{5.1. Proposition} Let $G \subset G_0$ be an inclusion of
discrete groups and $\nu$ a scalar $2$-cocycle for $G_0$. Denote
$(B\subset N)=(L_\nu(G) \subset L_\nu(G_0))$. Conditions $(a)-(d)$
are equivalent. If in addition $L_\nu(G_0)$ is a factor then
$(a)-(e)$ are equivalent.

$(a)$. $(G_0, G)$ is a property ${\text{\rm (T)}}$ pair, i.e.,
$G \subset G_0$ checks the equivalent conditions $(4.0.1)$, $(4.0.1')$,
$(4.0.1'')$.

$(b)$. $B\subset N$ is a rigid embedding of algebras.

$(c)$. For any $\varepsilon > 0$ there exist a finite set
$F'\subset N$ and $\delta' > 0$ such that if
$\Cal H$
is a Hilbert $N$-bimodule with a unit vector $\xi \in \Cal H$
satisfying $\|x_i \xi - \xi x_i\| \leq \delta', \forall i$ then
there exists a vector $\xi_0\in \Cal H$ such that
$\|\xi_0 - \xi\| \leq \varepsilon$ and $b \xi_0 = \xi_0 b, \forall b\in B$.

$(d)$. For any $\varepsilon > 0$ there exist a finite set
$F\subset  N$ and $\delta > 0$ such that if
$\phi: N \rightarrow N$
is a normal completely positive map with
$\|\phi(x) - x \|_2 \leq \delta, \forall x\in F$,
then $\|\phi(b)-b\|_2 \leq \varepsilon, \forall b\in B$, $\|b\| \leq 1$.

$(e)$. $L_\nu(G) \subset L_\nu(G_0)$ is $\varepsilon_0$-rigid for
some $\varepsilon_0 < 1$.
\endproclaim
\noindent
{\it Proof}. To prove $(a) \implies (c)$, we prove $(4.0.1') \implies  (c)$.
Let $\varepsilon > 0$ and let $E\subset G_0$,
$\delta'> 0$ be given by $(4.0.1')$, for this $\varepsilon$.
Let $\Cal H$ be a Hilbert $N$ bimodule with $\xi \in \Cal H$,
$\|\xi\|=1$, $\|u_h \xi - \xi u_h\| \leq \delta',
\forall h\in E'$. Taking $\pi(g)\eta = u_g\eta u_g^*, \eta \in \Cal H,
g\in G_0$, gives a representation of $G_0$ on $\Cal H$,
with $\|\pi(h)\xi - \xi\|=\|u_h\xi - \xi u_h\| \leq \delta'$.
Thus, there exists $\xi_0 \in \Cal H$ fixed by $\pi(G)$
(equivalently, $u_g\xi_0 = \xi_0 u_g, \forall g\in G$)
and such that $\|\xi_0 - \xi\| \leq \varepsilon$.

$(b) \implies (a)$. We prove that 4.1.1$^\circ$ implies $(4.0.1')$.
Let $\varepsilon > 0$.
By part $1^\circ$ in Lemma 1.1.3 and by
Kaplanski's density theorem (which implies that the unit ball of the group
algebra $\Bbb C_\nu G_0$ is dense in the unit ball of $L_\nu(G_0)$
in the norm $\|\quad\|_2$), it follows that given
any $\varepsilon$ there exist a finite set $E_0\subset
G_0$ and $\delta_0 > 0, \delta_0\leq
\varepsilon,$ such that if $\Cal H$ is a
$L_\nu(G_0)$ Hilbert bimodule with $\xi \in \Cal H$ a unit
vector which is left and right $\delta_0$-tracial and satisfies
$\|u_h \xi - \xi u_h\| \leq \delta_0, \forall h\in E_0$,
then there exists $\xi_1 \in \Cal H$ such that
$\| \xi_1 - \xi\| \leq \varepsilon/2$ and $b \xi_1 = \xi_1 b, \forall
b\in L_\nu(G), \|b\|\leq 1$.

Let then $(\pi_0, \Cal H_0, \xi_0)$ be a cyclic representation
of $G_0$ such that $\|\pi_0(h)\xi_0 -\xi_0\| \leq \delta_0, \forall h \in E_0$.
Let $(\Cal H_{\pi_0}, \xi_{\pi_0})$ be the pointed Hilbert
$L_\nu(G_0)$ bimodule, as defined in 1.4. We clearly have
$\|u_h \xi_{\pi_0} - \xi_{\pi_0} u_h\|= \|\pi_0(h)\xi_0 - \xi_0\|
\leq \delta_0, \forall h\in E_0$, by the definitions. Thus, there exists
$\xi_1\in \Cal H_{\pi_0}$ such that $\|\xi_1-\xi_{\pi_0}\|\leq
\varepsilon/2$ and $\xi_1$ commutes with $L_\nu(G)$.
But this implies that for all $g\in G$ we have
$$
\|\pi_0(g)\xi_0-\xi_0\| = \|u_g \xi_{\pi_0}
- \xi_{\pi_0} u_g\|
$$
$$
\leq  \|[u_g, (\xi_{\pi_0}-\xi_1)]\| + \|[u_g, \xi_1]\|
\leq 2\varepsilon/2 = \varepsilon.
$$

Taking the element of minimal norm $\xi_2$ in the weak closure
of co$\{\pi_0(g)\xi_1 \mid g\in G\}$, it follows that
$\xi_2$ is fixed by $\pi_0$ and $\|\xi_2-\xi_0\|
\leq \varepsilon$.

The implications $(c) \implies (b)$, $(d) \implies (b)$, $(b)
\implies (e)$ (the latter for factorial $L_\nu(G_0)$) are trivial.

To prove $(a) \implies (d)$, we prove (4.0.1') $\implies  (d)$.
Let $\varepsilon > 0$ and let $E' \subset G_0$,
$\delta' > 0$ be given by (4.0.1'), for $\varepsilon/2$. Also, we take
$E'$ to contain the unit $e$ of the group $G_0$.

Let $\phi$
be a completely positive map on $L_\nu(G_0)$ such that $\|\phi(u_h)-u_h\|_2
\leq \delta', \forall h\in E'$, where the norm $\|\quad \|_2$
is given by some trace $\tau$ on $L_\nu(G_0)$.
Let $F=\{u_h \mid h\in E'\}$.

Let $(\Cal H_{\phi}, \xi_{\phi})$ be the pointed Hilbert
$N$-bimodule defined out of $\phi$ as in 1.1.2. Let $\pi$ be the
associated representation of $G_0$ on $\Cal H_{\phi}$, as in the
last part of  1.1.4. It follows that there exists $\xi_0\in \Cal
H_{\phi}$ such that $b\xi_0=\xi_0b, \forall b\in  L_\nu(G)$ and
$\|\xi_{\phi}-\xi_0\|\leq \varepsilon/2$. Since $1 \in F$, part
$2^\circ$ of Lemma 1.1.2 shows that we may assume $\phi(1)\leq 1$.
By part 1$^\circ$ of Lemma 1.1.2 it then follows that for any
$u\in \Cal U(B)$ we have
$$
\|\phi(u)-u\|_2^2
\leq 2 - 2{\text{\rm Re}}\tau(\phi(u)u^*)
= \|u\xi_{\phi}-\xi_{\phi} u\|^2
$$
$$
= \|u(\xi_{\phi}-\xi_0)-(\xi_{\phi}-\xi_0) u\|^2
\leq 4\|\xi_{\phi}-\xi_0\|^2 \leq \varepsilon^2.
$$

$(e)\implies (a)$. As in the proof of $(b) \implies (a)$, by
Kaplanski's density theorem it follows that there exists $\delta >
0$ and $E \subset G_0$ such that if $\phi$ is completely positive
on $N= L_\nu(G_0)$, with $\phi(1) \leq 1, \tau \circ \phi \leq
\tau$ and $\|\phi(u_h) - u_h \|_2 \leq \delta, \forall h \in E$,
then $\|\phi(b)-b\|_2 \leq \varepsilon_0,$ for all $b$ in the unit
ball of $B=L_\nu(G)$.

Let $(\pi_0, \Cal H_0, \xi_0)$ be a cyclic representation of $G_0$
such that $\|\pi_0(h)\xi_0-\xi_0\| \leq \delta, \forall h \in E$.
Define $\phi_0$ on $N$ by $\phi_0(\Sigma_g \alpha_g u_g) =
\Sigma_g \langle \pi_0(g)\xi_0, \xi_0 \rangle \alpha_g u_g$. We
clearly have $\phi_0(1), \tau\circ \phi_0 = \tau$,
$\|\phi_0(u_h)-u_h\|\leq \delta, \forall h \in E$. Thus,
$\|\phi_0(u_g)-u_g\|_2 \leq \varepsilon_0, \forall g \in G$,
yielding $|\langle \pi_0(g)\xi_0, \xi_0 \rangle - 1| \leq
\varepsilon_0 < 1, \forall g\in G$. Taking the vector $\xi$ of
minimal norm in $\overline{\text{\rm co}}\{\pi_0(g) \mid g\in G\}
\subset \Cal H_0$, it follows that $\xi\neq 0$ and $\pi_0(g)(\xi)
= \xi, \forall g\in G$. This shows that the pair $(G_0, G)$
satisfies $(4.0.1)$, i.e., it has the relative property (T).
\hfill
Q.E.D.
\vskip .1in
For the first part of the next Corollary recall
that any (normalized, unitary, multiplicative) scalar $2$-cocycle
$\nu$ on $\Bbb Z^2$ is given by a bicharacter, and it is uniquely
determined by a relation of the form $uv=\alpha vu$ between the
generators $u=(1,0), v=(0,1)$ of $\Bbb Z^2$, where $\alpha$ is
some scalar with $|\alpha|=1$. We already considered such
2-cocycles in Corollary 3.3.2, where we pointed out that they are
$SL(2, \Bbb Z)$-invariant. Thus, if we denote by $L_\alpha(\Bbb
Z^2)$ the twisted group algebra $L_\nu(\Bbb Z^2)$, then the action
$\sigma$ of $SL(2, \Bbb Z)$ on $\Bbb Z^2$ 
induces an action still denoted $\sigma$ of
$SL(2, \Bbb Z)$ on $L_\alpha(\Bbb Z^2)$, preserving the canonical
trace (cf. 3.3.2). We have:

\proclaim{5.2. Corollary} $1^\circ$. The inclusion $\Bbb Z^2
\subset Z^2 \rtimes SL(2, \Bbb Z)$ is rigid. Thus, given any
$\alpha \in \Bbb T$, $L_\alpha(\Bbb Z^2) \subset L_\alpha(\Bbb
Z^2) \rtimes SL(2, \Bbb Z)$ is a rigid embedding of algebras.
Moreover, if $\alpha$ is not a root of unity, then the
``$2$-dimensional non-commutative torus'' $L_\alpha(\Bbb Z^2)$ is
isomorphic to the hyperfinite ${\text{\rm II}}_1$ factor $R$, thus
giving rigid embeddings $R \subset R\rtimes_\sigma SL(2, \Bbb Z)$.
If $\alpha$ is a primitive root of unity of order $n$, then
$$
(L_\alpha(\Bbb Z^2) \subset L_\alpha(\Bbb Z^2) \rtimes
SL(2, \Bbb Z)) =(L((n\Bbb Z)^2) \subset L((n\Bbb Z)^2) \rtimes
SL(2, \Bbb Z)) \otimes M_{n\times n}(\Bbb C)
$$
$$
\simeq (L(\Bbb Z^2) \subset L(\Bbb Z^2) \rtimes SL(2, \Bbb Z))
\otimes M_{n\times n}(\Bbb C) =(L^\infty(\Bbb T^2, \lambda)
\subset L^\infty(\Bbb T^2, \lambda) \rtimes SL(2, \Bbb Z))\otimes
M_{n\times n}(\Bbb C).
$$

$2^\circ$. If $n\geq 2$ and $\Bbb F_n \subset SL(2, \Bbb Z)$ has
finite index, then the restriction to $\Bbb F_n$ of the canonical
action of $SL(2, \Bbb Z)$ on $\Bbb T^2=\hat{\Bbb Z^2}$ (resp. on
$L_\alpha(\Bbb Z^2)\simeq R$, for $\alpha$ not a root of unity) is
free, weakly mixing, measure preserving, with $L^\infty(\Bbb T^2,
\mu) \subset L^\infty(\Bbb T^2, \lambda) \rtimes \Bbb F_n$
rigid (resp. $R \subset R \rtimes \Bbb F_n$ rigid).

$3^\circ$. For each $n\geq 2$ and each arithmetic lattice
$\Gamma_0$ in $SO(n, 1)$ (resp. in $SU(n, 1)$) there exist
free weakly mixing measure preserving actions of
$\Gamma_0$ on $A\simeq L^\infty(X, \mu)$ such that the 
corresponding cross-product inclusions $A \subset
A\rtimes \Gamma_0$ are rigid.

$4^\circ$. Let $\sigma_0$ be a properly outer, weakly mixing
action of some group $\Gamma_0$ on
$(B_0, \tau_0)$ such that
$B_0 \subset B_0 \rtimes_{\sigma_0} \Gamma_0$ be rigid
(e.g., like in $1^\circ$, $2^\circ$ or $3^\circ$).
Let $\sigma_1$ be any action of $\Gamma_0$ on some finite
von Neumann algebra $(B_1, \tau_1)$, which acts ergodically
on the center of $B_1$. If we denote
$B=B_0\overline{\otimes} B_1$ and
$M=(B_0\overline{\otimes} B_1)
\rtimes_{\sigma_0\otimes \sigma_1} \Gamma_0$, then
$M$ is a factor, $B_0'\cap M \subset B$, and
$B_0 \subset M$ is a rigid embedding.
\endproclaim
\noindent {\it Proof}. $1^\circ$. The rigidity of $\Bbb Z^2
\subset Z^2 \rtimes SL(2, \Bbb Z)$ is a well known result in
([Kaz]; see also [Bu], [Sha] for more elegant proofs). The fact
that $L_\alpha(\Bbb Z^2) \simeq R$ if  $\alpha$ is not a root of
unity and that $L_\alpha(\Bbb Z^2) \simeq A \otimes M_{n\times
n}(\Bbb C)$, with $A = \Cal Z(L_\alpha(\Bbb Z^2)) \simeq L((n\Bbb
Z)^2)$, if $\alpha$ is a primitive root of order $n$, are folklore
type results (see [Ri] and [HkS]).

In the latter case, if $p\in 1\otimes M_{n\times n}(\Bbb C)
\subset L_\alpha(\Bbb Z^2)$ is a projection of central trace $1/n$
then $\sigma_g(p)$ has central trace $1/n$ as well, so there
exists $v_g \in \Cal U(L_\alpha(\Bbb Z^2))$ such that $v_g
\sigma_g(p)v_g^* = p$. Thus, since $v_g$ commute with the center
$A$, if we denote by $\sigma'_g$ the action implemented by the
restriction of Ad$v_g\circ \sigma_g$ to $p(L_\alpha(\Bbb Z^2))p
=Ap \simeq A \simeq L((n\Bbb Z)^2)$, then $\sigma_g'$ coincides
with the restriction of $\sigma_g$ to $A\simeq L((n\Bbb Z)^2)$.

Moreover, if $u_g \in L_\alpha(\Bbb Z^2) \rtimes SL(2, \Bbb Z)$
are the canonical unitaries implementing $\sigma_g$ on $L_\nu(\Bbb
Z^2)$, then $u'_g=v_gu_gp$ implement the action
$\sigma'_g=\sigma_{|A}$ on $A$, but with an $A$-valued 2-cocycle
$v'$, i.e., $p(L_\alpha(\Bbb Z^2) \subset L_\alpha(\Bbb Z^2)
\rtimes_\sigma SL(2, \Bbb Z))p \simeq (A \subset A
\rtimes_{\sigma',v'} SL(2, \Bbb Z))$. But by ([Hj]), $A\subset
(A\rtimes_{\sigma',v'} SL(2, \Bbb Z)$ is the amplification by 12
of an inclusion of the form $A_0 \subset A_0 \rtimes \Bbb F_2$,
for some free ergodic action of $\Bbb F_2$ on $A_0$. Since any
action by the free group has trivial cocycle, $A_0 \subset  A_0
\rtimes \Bbb F_2$ is associated to the bare equivalence relation
it induces on the probability space, with trivial cocycle. Thus,
so does its $1/12$ reduction (see 1.4), i.e., $(A \subset A
\rtimes_{\sigma'} SL(2, \Bbb Z))=(L^\infty(\Bbb T^2, \lambda)
\subset L^\infty(\Bbb T^2, \lambda) \rtimes_\sigma SL(2, \Bbb Z))$.

The rest of the statement follows from part $(a)$ of Corollary
3.3.2$^\circ$.

$2^\circ$ follows from part $1^\circ$ above, Proposition
4.6.2$^\circ$ and part $(c)$ of Corollary 3.3.2$^\circ$.

$3^\circ$ follows by a recent result in ([Va]),
showing that there exist actions $\gamma$ of such $\Gamma_0$ on some
appropriate $\Bbb Z^N$ which give rise to rigid embeddings
$\Bbb Z^N \subset \Bbb Z^N\rtimes \Gamma_0$.
It is easy to
see that the actions $\gamma$ in ([Va]) can be taken to satisfy conditions
$(i), (ii)$ in Lemma 3.3.2.

$4^\circ$. By 3.3.3, since $\sigma_0$ is
properly outer, it follows that
$\sigma_0\otimes \sigma_1$ is properly outer and
$B_0'\cap M =\Cal Z(B_0) \otimes B_1$.
Also, since $\sigma_0$ is weakly mixing and $\sigma_1$ is
ergodic, it follows that $\sigma_0\otimes \sigma_1$
is ergodic and $M$ is a factor.
\hfill Q.E.D.

\proclaim{5.3. Corollary} $1^\circ$. Let $\Gamma_0$ be an
arbitrary discrete, countable group. Denote by $\sigma_1$ the Bernoulli
shift action of $\Gamma_0$ on $(A_1, \tau_1) =
\overline{\otimes}_{g\in \Gamma_0}
(L^\infty(\Bbb T, \lambda))_g$ and let
$\sigma_0$ be an ergodic action  of $\Gamma_0$ on
an abelian von Neumann algebra
$(A_0, \tau_0)$.
If we denote $A =A_0\overline{\otimes} A_1, \sigma=\sigma_0\otimes \sigma_1$
then $\sigma$ is free ergodic and the inclusion $A \subset A
\rtimes_\sigma \Gamma_0$ is not rigid.

$2^\circ$. $L(\Bbb Q^2)=A \subset M=L(\Bbb Q^2) \rtimes SL(2, \Bbb Q)$
is not a rigid inclusion but $A_0=L(\Bbb Z^2)\subset A$
satisfies $A_0 \subset M$ rigid and $A_0'\cap M=A$.

$3^\circ$. If $\Gamma_0$ is equal to
$SL(2, \Bbb Z)$, or to $\Bbb F_n$, for some $n \geq 2$,
or to an
arithemtic lattice in some $SO(n,1)$, $SU(n,1)$, $n\geq 2,$ then
there exist three non orbit equivalent
free ergodic measure preserving actions
$\sigma_i, 1\leq i\leq 3,$ of $\Gamma_0$ on the probability
space $(X,\mu)$. Moreover, each $\sigma_i$
can be taken such that $A=L^\infty(X, \mu)$
contains a subalgebra $A_i$ with $A_i \subset A\rtimes_{\sigma_i}
\Gamma_0$ rigid and $A_i'\cap A\rtimes_{\sigma_i}
\Gamma_0=A$.
\endproclaim
\noindent {\it Proof}. $1^\circ$. Write $L^\infty(\Bbb T, \lambda)
= \overline{\cup_n A^n}$, with $A^n$ an increasing sequence of
finite dimensional subalgebra and denote $A^n_1 =
\overline{\otimes}_g (A^n)_g \subset A_1$. Then $A^n_1 \uparrow
A_1$ and $\sigma_g(A^n_1)=A^n_1, \forall g \in \Gamma_0, \forall
n$. Thus, if we denote by $N_n = (A_0\overline{\otimes} A^n_1 \cup
\{u_g\}_g)''$ then $N_n \uparrow N = A\rtimes_\sigma \Gamma_0$. So
if we assume $A \subset N$ is rigid, then by 4.5 there exists $n$
such that $\|E_{N_n}(a)-a\|_2 \leq 1/2, \forall a\in A, \|a\|\leq
1$. But if $a\in 1\otimes A_1$ then $E_{N_n}(a)=E_{A^n_1}(a)$. Or,
since $A^n$ is finite dimensional and $L^\infty(\Bbb T, \lambda)$
is diffuse, there exists a unitary element $u_0\in L^\infty(\Bbb
T, \lambda)$ such that $E_{A_1^n}(u_0)=0$. Taking  $u= ... \otimes
1\otimes u_0 \otimes 1... \in A$, it follows that $E_{A_n}(u)=0$,
so that $1=\|E_{A^n_1}(u)-u\|_2 = \|E_{N_n}(u)-u\|_2 \leq 1/2$, a
contradiction.

$2^\circ$. For each $n$ let $\Bbb Q_n$ be the ring of rationals
with the denominator having prime decomposition with only the
first $n$ prime numbers appearing. Then $A\supset A_n = L(\Bbb
Q_n) \subset L(\Bbb Q_n) \rtimes SL(2, \Bbb Q_n)=M_n\subset M$ and
we have $E_{M_n}\circ E_A = E_{A_n}, \forall n$. If $A \subset M$
would be rigid, then by 4.5 there exists $n$ such that
$\|E_{M_n}(a)-a\|_2 \leq 1/2, \forall a\in A, \|a\|\leq 1$. But
any unitary element $u\in A=L(\Bbb Q^2)$ corresponding to a group
element in $\Bbb Q\setminus \Bbb Q_n$ satisfies $E_{A_n}(u)=0$, a
contradiction.

$3^\circ$. We take
$\sigma_1$ to be the action of $\Gamma_0$ on $A=L^\infty(X, \mu)$
given by 5.2.1$^\circ$-5.2.3$^\circ$.

We then take $\sigma_2$ to be the tensor product of $\sigma_1$
with the Bernoulli shift action of $\Gamma_0$ on
$\overline{\otimes}_{g\in \Gamma_0} (L^\infty(\Bbb T, \lambda))_g$.

Finally, we take $\sigma_3$ to be the tensor product of $\sigma_1$
with a free ergodic measure preserving action of $\Gamma_0$
which is not strongly ergodic, as
provided by the Connes-Weiss Theorem ([CW]; this is possible
because $\Gamma_0$ has the property H, so it does not have the
property (T)).

By part 1$^\circ$ we have $(A \subset A\rtimes_{\sigma_1} \Gamma_0)
\not\simeq (A \subset A\rtimes_{\sigma_2} \Gamma_0)$.
By results of Klaus Schmidt ([Sc]; see also [J2]) $\sigma_1, \sigma_2$
are strongly ergodic, while $\sigma_3$
is not. Thus, $(A \subset A\rtimes_{\sigma_3} \Gamma_0)
\not\simeq (A \subset A\rtimes_{\sigma_i} \Gamma_0), i=1,2.$

Since all these Cartan subalgebras have trivial 2-cocycle by construction,
their non-isomorphism implies the non-equivalence of the corresponding
orbit equivalence relations.

The existence of ``large'' subalgebras $A_i\subset A$
with $A_i \subset A\rtimes_{\sigma_i} \Gamma_0$ rigid
follows by construction and by 3.3.3.
\hfill Q.E.D.

\proclaim{5.4. Theorem} $1^\circ$. If $N$ is a type
${\text{\rm II}}_1$ factor with the
property ${\text{\rm H}}$ (as defined in $2.0.2$), then $N$
contains no diffuse relatively rigid
subalgebras $B \subset N$.

$2^\circ$. If $N$ has the property ${\text{\rm H}}$ relative to
a type ${\text{\rm I}}$ von Neumann algebra $B_0 \subset N$ then $N$ contains
no relatively rigid type ${\text{\rm II}}_1$ von Neumann
subalgebras $B \subset N$.
\endproclaim
\noindent
{\it Proof}. $1^\circ$. Let $\phi_n$ be completely positive maps
on $N$ such that $\phi_n \rightarrow id_N$, $\tau\circ\phi_n \leq \tau$ and
$T_{\phi_n} \in \Cal K(L^2(N,\tau))$. If $B \subset N$
is a rigid inclusion then by 4.1.2$^\circ$ it follows
that there exists $n$ such that $\phi=\phi_n$ satisfies
$\|\phi(u)-u\|_2 \leq 1/2,
\forall u\in \Cal U(B)$. If in addition $B$ has no atoms,
then any maximal abelian
subalgebra $A$ of $B$ is diffuse. Thus, such $A$ contains unitary
elements $v$ with $\tau(v^m)=0, \forall m\neq 0$. Since
the sequence $\{\hat{v^m}\}_m\subset L^2(N,\tau)$
is weakly convergent to $0$ and $T_{\phi}$ is compact,
it follows that $\|\phi(v^m)\|_2=\|T_{\phi}(\hat{v^m})\|_2 \rightarrow 0$.
Thus,
$$
\underset m \rightarrow \infty \to \lim \|\phi(v^m)-v^m\|_2 =
\underset m \rightarrow \infty \to \lim \|v^m\|_2 =1,
$$
contradicting $\|\phi(v^m)-v^m\|_2 \leq 1/2, \forall m$.

$2^\circ$. Assume $N$ does contain a
relatively rigid type ${\text{\rm II}}_1$ von Neumann
subalgebra $B \subset N$. Let $\phi_n$ be completely positive
$B_0$ bimodular maps
on $N$ such that $\phi_n \rightarrow id_N$, $\tau\circ\phi_n \leq \tau$ and
$T_{\phi_n} \in \Cal J_0(\langle N, B_0 \rangle)$.
By the rigidity of $B\subset N$ it follows that $\varepsilon_n=$
sup$\{\|\phi_n(u) - u\|_2 \mid u\in \Cal U(B)\} \rightarrow 0$. Since
$$
\|u^*T_{\phi_n}u(\hat{1})-\hat{1}\|_2=\|u^*\phi_n(u) -1\|_2
=\|\phi_n(u) - u\|_2,
$$
by taking convex combinations and weak limits of elements
of the form $uT_{\phi_n}u^*$, by Proposition 1.3.2 it follows that
there exists $T_n\in K_{T_{\phi_n}} \cap
(B'\cap \Cal J(\langle N,B \rangle))$ such that
$\|T_n(\hat{1})-\hat{1}\|_2\rightarrow 0$.
Thus, $T_n\neq 0$ for $n$ large enough,
so $B'\cap \langle N,B_0 \rangle$ contains
non-zero projections of finite trace. By ([Chr]), this implies
there exist non-zero projections $p\in B, q\in B_0$ and a
unital isomorphism
$\theta$ of $pBp$ into $qB_0q$. But $qB_0q$ is type I and $pBp$ is not,
a contradiction.
\hfill Q.E.D.

\proclaim{5.5. Corollary} $1^\circ$. If $N$ has a diffuse relatively
rigid subalgebra $B \subset N$ then $N$ cannot be embedded into
a free group factor $L(\Bbb F_n)$. In particular, the factors constructed
in ${\text{\rm Corollary}}$ $5.2$ cannot be embedded into $L(\Bbb F_n)$.

$2^\circ$. The factors $L_\alpha (\Bbb Z^2) \rtimes SL(2, \Bbb Z)$,
constructed in $5.2.1^\circ$
for $\alpha$
irrational, cannot be embedded
into $L_{\alpha'}(\Bbb Z^2) \rtimes SL(2, \Bbb Z)$
for $\alpha'$ rational.
\endproclaim
\noindent
{\it Proof}. Part $1^\circ$ is a consequence
of  5.4.1, while part $2^\circ$ follows trivially from 5.4.2.
\hfill Q.E.D.
\vskip .1in
\noindent
{\bf 5.6. Remarks}. $1^\circ$. In the case when $N$ is a finite factor,
a different notion of ``relative property
T'' for inclusions $B\subset N$,
was considered in ([A-De], [Po1]), as follows:
\vskip .1in
\noindent
(5.6.1). $N$ {\it has property} T {\it relative to} $B$
(or $B$ {\it is co-rigid in} $N$) if
there exists a finite set $F_1 \subset N$ and $\varepsilon_1 > 0$ such that
if $(\Cal H, \xi)$ is a $(B\subset N)$ Hilbert bimodule
(recall that by definition this requires $[B, \xi]=0$) such that $\|x \xi
- \xi x\| \leq \varepsilon, \forall x\in F$, then there exists
$\xi_0 \in \Cal H, \xi_0\neq 0,$ with $x\xi_0 = \xi_0 x, \forall x\in N$.
\vskip .1in

In the case $B$ is a Cartan subalgebra $A$ of a type II$_1$ factor
$N=M$, this definition is easily seen to be equivalent to Zimmer's
property (T) ([Zi2]) for the countable, measurable,
measure-preserving equivalence relation $\Cal R_{A \subset M}$,
which it thus generalizes to the case of arbitrary inclusions of
von Neumann algebras (cf. Section 4.8 in [Po1]). Thus, in this
re-formulation, a standard equivalence relation $\Cal R$ satisfies
Zimmer's relative property (T) iff the Cartan subalgebra $A \subset
M$, constructed as in ([FM]) out of $\Cal R$ and the trivial
2-cocycle $v\equiv 1$, is co-rigid in the sense of ([Po1],
[A-De]). We will in fact call such equivalence relations $\Cal R$
{\it co-rigid}.

$2^\circ$. It is easy to see that in case
$(B\subset N)= (B \subset B \rtimes_\sigma \Gamma_0)$,
for some cocycle action $\sigma$ of a group $\Gamma_0$ on $(B, \tau)$
then $N$ has the property (T) relative
to $B$ (i.e., $B$ is co-rigid in
$N$) if and only if $\Gamma_0$ has the property (T) of Kazhdan
(cf. [A-De], [Po1]; also
[Zi] for the
Cartan subalgebra case). In particular,
if $H\subset G_0$ is a normal subgroup of $G_0$ then $L(G_0)$ has
the property (T) relative to $L(H)$ if and only if the
quotient group $G_0/H$ has the
property (T). In fact, it is easy to see that if $H\subset G_0$
is an inclusion of discrete groups
then $L(G_0)$ has property (T) relative
to $L(H)$ iff the following holds true:
\vskip .1in
\noindent
(5.6.2). There exist a finite set $E \subset G_0$
and $\varepsilon > 0$ such that
if $\pi$ is a unitary representation
of $G_0$ on a Hilbert space $\Cal H$ with a
unit vector $\xi \in \Cal H$ such that $\pi(h)\xi = \xi, \forall h\in H$
and $\|\pi(g)\xi - \xi \| \leq \varepsilon, \forall g\in E$, then
$\Cal H$ contains a non-zero vector $\xi_0$ such that $\pi(g)\xi_0 = \xi_0,
\forall g\in G_0$.
\vskip .1in

A sufficient condition for an
inclusion of groups $H \subset G_0$ to satisfy 5.6.2$^\circ$ is when
$G_0$ has {\it finite length over} $H$, i.e.,
when the following holds true:
\vskip .1in
\noindent
$(5.6.2')$. There exists $n \geq 1$ and
a finite set $E \subset G_0$ such that any
element $g\in G_0$ can be written as
$g=h_1f_1h_2f_2... h_nf_n$,
for some $f_i \in E, h_j\in H$.
\vskip .1in
Indeed, because
then $\pi(h)\xi = \xi, \forall h\in H$
and  $\xi$ almost fixed by $\pi(f), f\in E$, implies that
$\xi$ is almost fixed by $\pi(g)$, uniformly
for all $g\in G_0$. This, of course, shows that $\Cal H$ has
a non-zero vector fixed by $\pi(G_0)$. (N.B.
Finite length was exploited in relation to rigidity in [Sha]).

An example of inclusion of groups $H \subset G_0$ satisfying
$(5.6.2')$ is obtained by
taking $G_0$ to be the group of all affine transformations
of $\Bbb Q$ and $H$ to be the subgroup of all homotheties of $\Bbb Q$.
Indeed, because if we take $E$ to be the single element set
consisting of the
translation by $1$ on $\Bbb Q$, then we clearly have $G_0 = HEH$.
Thus, $L(G_0)$, which is isomorphic to the hyperfinite type
II$_1$ factor
$R$, has the property (T) relative to $L(H)$, which is a singular maximal
abelian subalgebra in $L(G_0)$ (cf. [D]).

\proclaim{5.7. Proposition} Let $N$ be a type $II_1$ factor
and $B\subset N$ a von Neumann subalgebra.

$1^\circ$. If $\langle N, B \rangle$ is finite then
$N$ has both the property (T) relative to $B$
(in the sense of $(5.6.1)$) and the property H relative to $B$.

$2^\circ$. If $N$ has both property (T) and H relative to $B$ then
there exists a non-zero $q\in \Cal P(B'\cap N)$ such that $qNq$ is a
finitely generated $Bq$-module. Thus, if
in addition $B$ is a subfactor with $B'\cap N=\Bbb C$
then $[N:B] < \infty$ and if $B$ is a maximal
abelian von Neumann subalgebra in $N$ then
${\text{\rm dim}} N < \infty$.
\endproclaim
\noindent {\it Proof}. $1^\circ$. If $\langle N,B \rangle$ is
finite, then there exists a sequence of projections $p_n \in \Cal
Z(B), p_n \uparrow 1$, such that $p_nNp_n$ has finite orthonormal
basis over $Bp_n$. By 2.3.4$^\circ$, this implies $p_nNp_n$ has
the property H relative to $Bp_n$ and by 4.6.3$^\circ$, $Bp_n
\subset p_nNp_n$ follows rigid. By 2.4.2$^\circ$ this implies $N$
has the property H relative to $B$ and by 4.7.1$^\circ$, $B
\subset N$ is rigid.

$2^\circ$. Note first that if there exist no
$q\in \Cal P(B'\cap N)$ such that $qNq$ is a
finitely generated $Bq$-module, then
$N'\cap \langle N,B \rangle$ contains no finite
projections of $\langle N,B \rangle$.

On the other hand, if $N$ has the property H relative to
$B$ then by 2.2.1$^\circ$ there exist unital
completely positive, $B$-bimodular maps $\phi_n$ on $N$ such that
$\tau\circ \phi_n \leq \tau$, $\phi_n(1) \leq 1$,
$\phi_n \rightarrow id_N$ and
$T_{\phi_n} \in \Cal J_0(\langle N, B\rangle)$. If in addition
$N$ has the property (T) relative to $B$, then $\exists n$ such that
$\|\phi_n(u)-u\|_2 \leq 1/4,
\forall u\in \Cal U(N)$. By 1.3.3,
$\exists$  a spectral projection $p\in B'\cap \Cal J_0(\langle N, B\rangle)$
of $T^*_{\phi_n} T_{\phi_n}$ such that $\|T_{\phi_n}(1-p)\| < 1/4$.
If we now assume $N'\cap \langle N, B\rangle$ has no finite projections,
then there exists a unitary element $u\in \Cal U(N)$ such that
$Tr(pue_Bu^*) < 1/4$.
But $Tr(pue_Bu^*) =\|p(\hat{u})\|_2^2$ (see the proof of 6.2
in the next Section).
Altogether, since $\|p(\hat{u})\|_2 \geq \|T_{\phi_n}(\hat{u})\|_2 -
\|T_{\phi_n}((1-p)(\hat{u}))\|_2 \geq 1/2$, it
follows that $1/4 > Tr(pue_Bu^*) \geq 1/4$, a contradiction.
The last part of 2$^\circ$ follows trivially from ([PiPo]).
\hfill Q.E.D.
\vskip .1in
\noindent
{\bf 5.8. Remarks}. $1^\circ$. Both
the notion 4.2 considered here and the notion considered in ([A-De],
[Po1]) are in some sense ``relative property (T)'' notions for an
inclusion $B\subset N$, but while the notion in
([A-De], [Po1]) means ``$N$ has the property (T) relative to $B$'',
thus being a ``co''-type property (T,)
the notion considered in this paper is a ``property (T) of $B$ relative
to its embedding into $N$''.  The two notions are complementary
one to the other, and together they imply (and are implied by)
the property (T) of the global factor (see Proposition 5.9 below).

$2^\circ$. An
interesting relation between these two complementary notions
of  property (T) is the following:  If a group $\Gamma_0$ acts on
$(B, \tau)$ such that
$B \subset N = B \rtimes \Gamma_0$ is a rigid embedding, then
$N$ has the property (T) relative to its
group von Neumann subalgebra $L(\Gamma_0)$ (i.e.,
$L(\Gamma_0)$ is co-rigid
in $N$). Indeed, because
if $(\Cal H, \xi)$ is a $(L(\Gamma_0)\subset N)$-Hilbert bimodule
with $\xi$ almost commuting with all $u \in \Cal U(B)$, uniformly,
then $\xi$ almost commutes with the group of elements $\Cal G=\{u u_g\mid
u\in \Cal U(B), g\in \Gamma_0\}$, thus $\xi$ is close to a vector commuting
with all $v\in \Cal G$, thus with all $x\in N$. For instance,
the factor $L(\Bbb Z^2 \rtimes SL(2, \Bbb Z))$ has the property
T relative to its subalgebra $L(SL(2, \Bbb Z))$
(in the sense of definition (5.6.1)).

\proclaim{5.9. Proposition} Let $N$ be a type ${\text{\rm II}}_1$ factor
and $B \subset N$ a von Neumann subalgebra. The following
conditions are equivalent:

$1^\circ$. $N$ has
the property ${\text{\rm T}}$ in the sense of Connes and Jones
(i.e., of the equivalent conditions $(4.0.2), (4.0.2')$).

$2^\circ$. The identity embedding $N\subset N$ is rigid,
i.e., for any $\varepsilon> 0$ there exists a finite
subset $x_1, x_2, ..., x_n \in N$ and
$\delta > 0$ such that if $\Cal H$ is a
Hilbert $N$-bimodule with a unit vector $\xi \in \Cal H$
satisfying $\|\langle \cdot \xi, \xi\rangle -\tau\| \leq \delta$,
$\|\langle  \xi, \xi \cdot \rangle -\tau\| \leq \delta$ and
$\| x_i \xi - \xi x_i \|
\leq \delta, \forall i$, then there exists a vector $\xi_0
\in \Cal H$ such that $\|\xi-\xi_0\| \leq
\varepsilon$ and $x\xi_0 = \xi_0 x, \forall x\in N$.

$3^\circ$. $B \subset N$ is a rigid embedding (in the sense
of definition $4.2$) and $N$ has the property ${\text{\rm T}}$ relative to $B$
(in the sense of $(5.6.1)$).
\endproclaim
\noindent {\it Proof}. $1^\circ \implies 3^\circ$ and $1^\circ
\implies 2^\circ$ are trivial, by the characterization $(4.0.2')$
of the property (T) for $N$.

To prove $3^\circ \implies 1^\circ$ let $F_1\subset N$ and
$\varepsilon_1$ give the critical set and constant for the
property (T) of $N$ relative to $B$ and
$F' \subset N, \delta' >0$ be the critical set and constant
for the rigidity of $B\subset N$, corresponding to $\varepsilon_1/4$.
Let $F= F'\cup F_1$ and let $\Cal H$ be a Hilbert $N$ bimodule
with a unit vector $\xi$ which is left and right $\delta'$-tracial
and satisfies $\|y \xi - \xi y\| \leq \delta', \forall y\in F$.
By the rigidity of $B\subset N$ it follows that there
exists $\xi_0 \in \Cal H$ such that $b\xi_0 = \xi_0 b, \forall b\in B$
and $\|\xi_0-\xi\| \leq \varepsilon_1/4.$ Thus,
if we assume $\varepsilon_1 \leq 1/4$
from the beginning and denote $\xi_1=\xi_0/\|\xi_0\|$, then
$\|\xi_1\|=1$, $b\xi_1 = \xi_1 b, \forall b\in B$,
and $\|y \xi_1 - \xi_1 y\| \leq \varepsilon_1,
\forall y\in F$, in particular for all $y\in F_1$. Thus,
by the property (T) of $N$ relative to $B$, $\Cal H$ has a non-zero $N$-central
vector.

$2^\circ \implies 1^\circ$. By part $1^\circ$ of Theorem 4.3, $N$
follows non-$\Gamma$. Thus, by Lemma 2.9 it is sufficient to check
that any Hilbert $N$ bimodule with a vector that's  almost
left-right tracial and almost central  has a non-zero central
vector for $N$. But this does hold true by the fact that $N$
checks condition 2$^\circ$. \hfill Q.E.D. \vskip .1in \noindent
{\bf 5.10. Remark.} When applied to the case of Cartan subalgebras
coming from standard equivalence relations with trivial 2-cocycle,
the definition of rigid embeddings 4.2 gives the following new
property for equivalence relations: \vskip .1in \noindent {\it
5.10.1. Definition.} A countable, ergodic, measure preserving
equivalence relation $\Cal R$ has the {\it relative property} (T) if
its associated Cartan subalgebra $A\subset M$, constructed out of
$\Cal R$ and the trivial 2-cocycle $v\equiv 1$ as in ([FM]), is a
rigid embedding (definition 4.2). 
\vskip .1in 
Since the rigidity
for Cartan subalgebras is an invariant for the isomorphism class
of $A \subset M$, this relative property (T) is an orbit equivalence
invariant for equivalence relations $\Cal R$. Also, when applied
to the particular case of Cartan subalgebras with trivial
2-cocycle, all the results on rigid embeddings of algebras in
Sections 4 and 5 translate into corresponding results about
standard equivalence relations $\Cal R$. For instance, by 4.6,
4.7, if $\Cal R$ has the relative property (T) then $\Cal R^t$ has
the relative property (T), $\forall t>0$, and if $\Cal R_1, \Cal
R_2$ have the relative property (T) then so does $\Cal R_1 \times
\Cal R_2$. Also, if $\Cal R$ has relative property (T) then
Out$(\Cal R)\overset\text{\rm def} \to ={\text{\rm Aut}}(\Cal
R)/{\text{\rm Int}}(\Cal R)$ is discrete (cf. 4.4) and if we
further have $\Cal R= \cup_n \Cal R_n$ for some increasing
sequence of ergodic sub-equivalence relations, then $\Cal R_n$
have the relative property (T) for all large enough $n$.

We have proved that equivalence relations implemented by
Bernoulli shift actions of a group $\Gamma_0$  cannot have the
relative property (T), no matter the group $\Gamma_0$ (cf. 5.3).
Thus, equivalence relations coming from actions
of the same group $\Gamma_0$ may
or may not have the relative property (T), depending on the action. 
While by ([Zi]; see also [A-De], [Po1]), $A
\rtimes_\sigma \Gamma_0$ has the property (T) relative to $A$, in
the sense of definition (5.6.1) if and only if $\Gamma_0$ has
Kazhdan's property (T), thus being a property entirely depending on
the group. Even more: since by ([Po1]) if $A \subset M$ is a
Cartan subalgebra in a II$_1$ factor and $p\in \Cal P(A)$ then
$pMp$ has property (T) relative to $Ap$ if and only if $M$ has
property (T) relative to $A$, it follows that the property (T) for
groups is invariant to stable orbit equivalence, or equivalently,
it is a ME invariant (see [Fu] for an ``ergodic theory'' proof of
this fact).

Proposition 5.9 shows that when the relative property (T)  (5.10.1)
for $\Cal R$ is combined with the co-rigidity property (5.6.1) for
$\Cal R$ they imply, and are implied by, the ``full'' {\it
property} T {\it of} $\Cal R$, which by definition requires that
the finite factor $M=M(\Cal R)$ has the property (T) in the sense
(4.0.2), of Connes-Jones. It is thus of great interest to answer
the following: \vskip .1in \noindent {\bf 5.10.2. Problem}.
Characterize the countable discrete groups $\Gamma_0$ that can act
rigidly on the probability space, i.e., for which there exist free
ergodic measure preserving actions $\sigma$ on $(X, \mu)$ such
that $L^\infty(X, \mu) \subset L^\infty(X, \mu) \rtimes_\sigma
\Gamma_0$ is a rigid embedding. Do all property (T) groups
$\Gamma_0$ admit such rigid actions (i.e., in view of the above,
actions $\sigma$ with the property that the II$_1$ factor
$L^\infty(X, \mu) \rtimes_\sigma \Gamma_0$ has the property (T) in
the sense of (4.0.2)) ?

\heading 6. HT subalgebras and the class $\Cal H \Cal T$. \endheading

\noindent
{\it 6.1. Definition}. Let $N$ be a finite
von Neumann algebra with a
normal faithful tracial state and $B \subset N$
a von Neumann subalgebra. $B$ is a HT {\it subalgebra} of
$N$ (or $B\subset N$ is a HT {\it inclusion})
if the following two conditions are met:
\vskip .1in
\noindent
(6.1.1). $N$ has the property H relative to $B$ (as defined in Section 2).
\vskip .1in
\noindent
(6.1.2). There exists a von Neumann subalgebra $B_0 \subset B$ such that
$B_0'\cap N \subset B$ and $B_0 \subset N$ is a rigid
(or property (T)) embedding.
\vskip .1in
Also, $B$ is a HT$_{_{s}}$ {\it subalgebra} of $N$ if conditions (6.1.1)
and (6.1.2) hold true with $B_0=B$, i.e., if $N$ has the
property H relative to $B$ and $B \subset N$ is itself a rigid embedding.

If $A\subset M$ is a Cartan subalgebra of
a finite factor $M$ and $A \subset M$ satisfies the conditions
(6.1.1) and (6.1.2), then we call it a HT {\it Cartan subalgebra}.
Similarly, if a Cartan subalgebra $A \subset M$
satisfies (6.1.1) and is a rigid embedding then it is called a HT$_{_{s}}$
{\it Cartan subalgebra}.

Note that condition (6.1.2) implies that $B'\cap N \subset B$ and
(6.1.1) implies $B$ is quasi-regular in $N$ (cf.3.4). In particular, by
Proposition 3.4, for $A\subset M$
a maximal abelian $*$-subalgebra of type II$_1$ factor
$M$, the condition that $A$
is an HT (resp. HT$_{_{s}}$) subalgebra of $M$ is sufficient to
insure that $A$ is an HT (resp. HT$_{_{s}}$)
Cartan subalgebra of $M$.

\proclaim{6.2. Theorem} Let $M$ be a type ${\text{\rm II}}_1$ factor
with two abelian von Neumann subalgebras $A, A_0$ such that $A$,
$A_0'\cap M$ are maximal abelian in $M$.
Assume that
$M$ has property ${\text{\rm H}}$ relative to $A$
and that
$A_0\subset M$ is a rigid inclusion. Then both $A$ and
$A_0'\cap M$ are ${\text{\rm HT}}$
Cartan subalgebras of $M$ and there exists a
unitary element $u$ in $M$ such that $uA_0u^* \subset A$,
and thus $u(A_0'\cap M)u^* = A$. In particular, if $A_1, A_2$ are
${\text{\rm HT}}$ Cartan subalgebras of a type
${\text{\rm II}}_1$ factor
$M$ then there
exists a unitary element $u \in \Cal U(M)$ such that $uA_1u^* = A_2$.
\endproclaim
\noindent {\it Proof}. We first prove that there exists a non-zero
partial isometry $v\in M$ such that $v^*v \in A_0 '\cap M$,
$vv^*\in A$ and $vA_0v^* \subset Avv^*$. If we  assume by
contradiction that this is not the case, then Theorem A.1 implies
$0\in K_{\Cal U(A_0)}(e_A) \subset \langle M, A \rangle $. This in
turn implies that given any finite projection $f \in \langle M, A
\rangle$, with $Tr(f) < \infty$, and any $\varepsilon > 0$, there
exists a unitary element $u\in \Cal U(A_0)$ such that
$Tr(fue_Au^*) < \varepsilon$. Indeed, because if for some $c_0 >
0$ we would have $Tr(fue_Au^*) \geq c_0, \forall u\in \Cal
U(A_0),$ then by taking appropriate convex combinations and weak
limits, we would get that $0 = Tr(f 0) \geq c_0>0$, a
contradiction.

By the property H of $M$ relative to $A$,
there exist completely positive, unital,
$A$-bimodular maps $\phi_n : M \rightarrow M$ which
tend strongly to the identity
and satisfy $\phi_n(1) \leq 1, \tau \circ \phi_n \leq \tau$,
$T_{\phi_n} \in \Cal J_0(\langle M, A \rangle)$.

Let $0< \varepsilon_0 < 1$. By the rigidity of the embedding $A_0
\subset M$, there exists $n$ large enough such that $\phi=\phi_n$
satisfies
$$
\|\phi(v)-v\|_2 \leq \varepsilon_0, \forall v\in \Cal U(A_0) \tag
6.2.1
$$
On the other hand,
since $T_\phi\in \Cal J_0(\langle M, A \rangle)$,
it follows that there exists a finite projection $f\in
\Cal J_0(\langle M, A \rangle)$ such that Tr$(f) < \infty$ and
$$
\|T_\phi(1-f)\| \leq (1-\varepsilon_0)/2   \tag 6.2.2
$$

Let then $u\in \Cal U(A_0)$ satisfy the condition
$$
Tr(fue_Au^*) < (1-\varepsilon_0)^2/4 \tag 6.2.3
$$

Let $\{m_j\}_j \subset L^2(M, \tau)$ be
such that $\Sigma_j m_je_Am_j^* = f$. Equivalently,
$\oplus_j L^2(m_jA) = f L^2(M, \tau)$. Thus, if $x\in N = \hat{N}
\subset L^2(M, \tau)$ then $f(\hat{x}) = \Sigma_j m_j E_A(m_j^*x)$ and
$\|f(\hat{x})\|_2^2 = \Sigma_j \|m_j E_A(m_j^*x)\|_2^2$.

It follows that we have:
$$
Tr(fue_Au^*)= Tr(fue_Au^*f)
$$
$$
= Tr ((\Sigma_j m_je_Am_j^*) ue_Au^*(\Sigma_i m_ie_A m_i^*))
$$
$$
=\Sigma_{i,j} \tau(m_jE_A(m_j^*u)E_A(u^*m_i)m_i^*)
= \|f(\hat u)\|_2^2
$$
By (6.2.3) this implies
$$
\|f(\hat{u})\|_2 < (1-\varepsilon_0)/2 \tag 6.2.4
$$
Thus, by taking into account that $\|T_\phi\|\leq 1$, (6.2.2) and (6.2.4)
entail:
$$
\|T_\phi(\hat u)\|_2 \leq \|T_\phi((1-f)(\hat u))\|_2 + \|f(\hat u)\|_2
$$
$$
\leq (1-\varepsilon_0)/2 + \|f(\hat{u})\|_2 < 1-\varepsilon_0.
$$

But by (6.2.1), this implies:
$$
\|u\|_2 \leq \|T_\phi(\hat u)\|_2 + \|\phi(u)-u\|_2
$$
$$
< 1-\varepsilon_0 + \varepsilon_0 = 1.
$$
Thus $1=\tau(uu^*) < 1,
$
a contradiction.

Let now $(\Cal V, \leq)$ denote the set of partial isometries $v \in M$
with $v^*v \in
A_0 '\cap M$, $vv^*\in A$ and $vA_0v^* \subset Avv^*$, endowed
with the order $\leq $ given by restriction, i.e., $v\leq v'$
if $v=vv^*v'$. $(\Cal V, \leq)$ is clearly inductively ordered.
Let $v_0 \in \Cal V$ be a maximal element. Assume $v_0$ is not
a unitary element.

By 2.4.1$^\circ$, $(1-v_0v_0^*)M(1-v_0v_0^*)$ has the property H
relative to $A(1-v_0v_0^*)$ and by 4.7.2$^\circ$ the inclusion
$A_0(1-v_0^*v_0)\subset (1-v_0^*v_0)M(1-v_0^*v_0)$ is rigid. Let
$u_0\in M$ be a unitary element extending $v_0$ and denote $M^0 =
(1-v_0v_0^*)M(1-v_0v_0^*)$, $A_0^0=u_0(A_0(1-v_0^*v_0))u_0^*$,
$A^0=A(1-v_0v_0^*)$. Thus, $M^0$ has the property H relative to
$A^0$ and $A_0^0\subset M_0$ is rigid. By the first part it
follows that there exists a non-zero partial isometry $v\in M^0$
such that $v^*v\in (A_0^0) '\cap M$, $vv^*\in A^0$ and $vA^0_0v^*
\subset A^0vv^*$. But then $v'=v_0 + vu_0^* \in \Cal V$, $v'\geq
v_0$ and $v'\neq v_0$, contradicting the maximality of $v_0$.

We conclude that $v_0$ is a unitary element, so that $A, A_0'\cap M$
are conjugate in $M$. The last part follows now by Proposition 3.4.
\hfill Q.E.D.

\vskip .1in \noindent {\bf 6.3. Remarks}. $1^\circ$. If in the
last part of Theorem 6.2 we restrict ourselves to the case $A_1,
A_2$ are HT$_{_{s}}$ Cartan subalgebras of the type II$_1$ factor
$M$, then we can give the following alternative proof of the
statement, by using part $(ii)$ of Proposition 1.4.3 ``in lieu''
of Theorem A.1 and an argument similar to the proof of
5.4.2$^\circ$:  By the property H of $M$ relative to $A_1$ there
exists completely positive $A_1$ bimodular trace preserving unital
maps  $\phi_n$ on $M$ such that $\phi_n \rightarrow id_M$ and
$T_{\phi_n} \in \Cal J_0(\langle N, A_0 \rangle)$. By the rigidity
of $A_2\subset M$ it follows that $\varepsilon_n=$
sup$\{\|\phi_n(u) - u\|_2 \mid u\in \Cal U(A_2)\} \rightarrow 0$.
Fix $x\in M$ and note that by Corollary 1.1.2 we have
$$
\|u^*T_{\phi_n}u(\hat{x})-\hat{x}\|_2=\|\phi_n(ux) -ux\|_2
$$
$$
\leq \|\phi_n(ux) -u\phi_n(x)\|_2 + \|\phi_n(x)-x\|_2 \leq
2\varepsilon^{1/2} + \|\phi_n(x)-x\|_2.
$$
Thus, by taking weak limits of appropriate convex combinations
of elements
of the form $u^*T_{\phi_n}u$ with $u\in \Cal A_2,$ and
using Proposition 1.3.2 it follows that
$T_n=\Cal E_{A_2'\cap \langle M, A_1\rangle} (T_{\phi_n})
\in K_{T_{\phi_n}} \cap
(A'\cap J_0(\langle N,A \rangle))$ satisfy
$\underset n \rightarrow \infty \to \lim
\|T_n(\hat{x})-\hat{x}\|_2 = 0$. But $x\in M$ was arbitrary.
This shows that the right supports of $T_n$ span all the identity
of $\langle M, A_1 \rangle$. Since $T_n$ are compact,
this shows that $A_2'\cap \langle M, A_1\rangle$,
is generated by finite projections of $\langle M, A_1\rangle$.
Thus, $A_2$ is discrete over $A_1$. Similarly, $A_1$ is discrete over $A_2$
and $A_1$ follows conjugate to $A_2$ by a result in ([PoSh];
see part $(ii)$ of Proposition 1.4.3).

$2^\circ$. The above argument uses the fact that two Cartan
subalgebras $A_1, A_2$ in $M$ are unitary equivalent in $M$ if and
only if the $A_1-A_2$ Hilbert bimodule $L^2(M, \tau)$ is a direct
sum of Hilbert bimodules that are finite dimensional both as left
$A_1$ modules and as right $A_2$ modules. The proof of Theorem 6.2
uses instead Theorem A.1, which shows that in order for an abelian
von Neumann algebra $A_0\subset M$ to be unitary conjugate to a
subalgebra of a semi-regular maximal abelian $*$-subalgebra $A$ of
$M$ it is sufficient that $A_0'\cap M$ be  semi-regular abelian
and that $_{A_0}L^2(M,\tau)_{A}$ contains a non-zero $A_0-A$
Hilbert bimodule which is finite dimensional as a right $A$-module
(a much weaker requirement).

$3^\circ$. Note that by 3.4 and 4.3.2$^\circ$, $A \subset M$ is
HT$_{_{s}}$ Cartan iff $A \subset M$ is maximal abelian, $M$ has
property H relative to $A$ and $A \subset M$ is
$\varepsilon_0$-rigid for some $\varepsilon_0 \leq 1/3$.

$4^\circ$. Note that the proof of Theorem 6.2 shows in fact that
if $A, A_0$ are abelian von Neumann subalgebras of a finite factor
$M$ such that $A$ is maximal abelian, $M$ has property H relative
to $A$, $A_0'\cap M$ is semi-regular abelian and $A_0 \subset M$
is $\varepsilon_0$-rigid, for some $\varepsilon_0<1$, then there
exists $u\in \Cal U(M)$ such that $u(A_0'\cap M)u^*=A$. In particular, if one calls HT$_w$ the 
Cartan subalgebras $A\subset M$ with the properties 
that $M$ has property H
relative to $A$ and there exists $A_0\subset A$ with
$A_0'\cap M=A$, $A_0 \subset M$ $\varepsilon_0$-rigid, for some
$\varepsilon_0<1$, then any two HT$_w$ Cartan 
subalgebras of a II$_1$ factor are unitary conjugate.
 \vskip .1in \noindent {\it 6.4. Notation}. We denote by $\Cal
H\Cal T$ the class of finite separable (in norm $\|\quad \|_2$)
factors with HT Cartan subalgebras and by $\Cal H\Cal T_{_{s}}$
the class of finite separable factors with HT$_{_{s}}$ Cartan
subalgebras. Note that $\Cal  H\Cal T_{_{s}}\subset \Cal H\Cal T$
and that Theorem 6.2 shows the uniqueness up to unitary conjugacy
of HT Cartan subalgebras in factors $M \in \Cal H\Cal T$.

\proclaim{6.5. Corollary} If $A_i\subset M_i,
i=1,2,$ are ${\text{\rm HT}}$
Cartan subalgebras and $\theta$ is an isomorphism from $M_1$ onto
$M_2$ then there exists a unitary element $u\in M_2$ such
that $u\theta(A_1)u^* = A_2$. Thus, there exists a unique
(up to isomorphism) standard
equivalence relation $\Cal R_M^{^{HT}}$
on the standard probability space, implemented by the normalizer of the
${\text{\rm HT}}$ Cartan subalgebra of $M$.
\endproclaim

The next result shows that $\Cal H\Cal T$ is closed to
amplifications and tensor products and that it has good
``continuity'' properties. The proof of part 3$^\circ$ below, like
the proof of 4.5.2$^\circ$, uses A.2 and is inspired by the proofs
of (4.5.1, 4.5.6 in [Po1]).

\proclaim{6.6. Theorem} $1^\circ$. If $M\in \Cal H\Cal T$
(resp. $M\in \Cal  H\Cal T_{_{s}}$)
and $t > 0$ then $M^t\in \Cal H\Cal T$ (resp. $M^t\in \Cal  H\Cal T_{_{s}}$).

$2^\circ$. If $M_1, M_2 \in \Cal H\Cal T$
(resp. $M_1, M_2 \in \Cal  H\Cal T_{_{s}}$) then
$M_1 \overline{\otimes} M_2 \in \Cal H\Cal T$
(resp. $M_1 \overline{\otimes} M_2 \in \Cal  H\Cal T_{_{s}}$).

$3^\circ$. If $M\in \Cal  H\Cal T_{_{s}}$ then there exist a finite set
$F \subset M$ and $\delta > 0$ such that
if $N\subset M$ is a subfactor with $F \subset_\delta N$ then
$N\in \Cal  H\Cal T_{_{s}}$. In particular, if $N_k \subset M$ are
subfactors with $N_k \uparrow M$, then there exists $k_0$
such that $N_k \in \Cal  H\Cal T_{_{s}}, \forall k\geq k_0$. If in
addition we have $N_k'\cap M = \Bbb C$, then
all the $N_k, k\geq k_0,$
contain the same ${\text{\rm HT}}_{_{s}}$ Cartan subalgebra of $M$.
\endproclaim
\noindent
{\it Proof}. $1^\circ$. Let $A\subset M$ be an HT Cartan
subalgebra and $A_0 \subset A$ be so
that $A_0 \subset M$ is a rigid embedding and
$A_0'\cap M = A$. Choose some
integer $n \geq t$. By 2.3.2$^\circ$  it follows that
if $D$ denotes the
diagonal of $M_0=M_{n\times n}(\Bbb C)$
then $A\otimes D \subset M_n(M)$ has the property H. Also,
$(A_0\otimes D)'\cap M \otimes M_{n\times n}(\Bbb C)
= A\otimes D$ and by 4.4.1$^\circ$, $A_0\otimes D
\subset M\otimes M_{n\times n}(\Bbb C)$ is a rigid embedding.

If we now take $p \in A_0\otimes D$ to be a projection of trace
$\tau(p) = t/n$,
then by 2.4.1$^\circ$ and 4.5.2$^\circ$, it follows that
$A_0^t = (A_0 \otimes D)p \subset M^t =pM_{n\times n}(\Bbb C)p$
is rigid and $M^t$ has the property H relative to $A^t$.
Thus, $M^t\in \Cal H\Cal T$. In case $A_0 = A$, then we get $A_0^t = A^t$,
so that $M^t$ follows in $\Cal H\Cal T_{_{s}}$.

$2^\circ$. This follows trivially
by applying 2.3.2$^\circ$ and 4.4.1$^\circ$, once we notice that if
$A_i \subset M_i$ are maximal abelian $*$-subalgebras and
$A_0^i \subset A_i$ satisfy $(A_o^i)'\cap M_i = A_i$,
then $(A_0^1 \overline{\otimes} A^2_0)'\cap M_1 \overline{\otimes} M_2
=A_1 \overline{\otimes} A_2$.

$3^\circ$. Let $A \subset M$ be a fixed HT$_{_s}$ Cartan
subalgebra of $M$. By 4.5.2$\circ$, it follows that there exist a
finite subset $F$ in the unit ball of $M$ and $\varepsilon > 0$
such that if a subfactor $N_0\subset M$ satisfies $F
\subset_\varepsilon N_0$ and $N_0'\cap M = \Bbb C$ then $N_0$
contains a unitary conjugate $A_0=uAu^*$ of $A$ with  $A_0 \subset
N_0$ rigid and Cartan. Moreover, $N_0$ has property H relative to
$A_0$ by 2.3.3$^\circ$ (since $M$ has property H relative to
$A_0$). Thus, $A_0 \subset N_0$ is HT$_{_s}$ Cartan, proving the
statement in the case of subfactors with trivial relative
commutant.

To prove the general case, note first
that by Step 1 in the proof of A.2, for the above
given $\varepsilon>0$ there exists $\delta_0> 0$,
with $\delta_0 \leq \varepsilon/4$,  such that if
$N \subset M$ is a subfactor with $A \subset_{\delta_0} N$ then
there exist projections $p\in A$, $q\in N$, a unital isomorphism
$\theta : Ap \rightarrow qNq$ and a partial
isometry $v\in M$ such that $\tau(p) \geq 1-\varepsilon/4$,
$v^*v=p$, $vv^* = qq'$,
for some projection $q' \in \theta(Ap)'\cap qMq$, and
$va = \theta(a)v, \forall a\in Ap$.

Since $Ap$ is maximal abelian in $pMp$, by spatiality it follows
that $\theta(Ap)q'$ is maximal abelian in $q'qMq'q$. Thus,
if $x\in \theta(Ap)'\cap qMq$ then $q'xq' \in \theta(Ap)q'
\simeq \theta(Ap)$. Thus, there exists a unique normal conditional
expectation $E$ of $\theta(Ap)'\cap qMq$ onto $\theta(Ap)$
satisfying $q'xq'=E(x)q', \forall x\in \theta(Ap)'\cap qMq$.

Let $q'_0\in N'\cap M$ be the support projection of $E_{N'\cap M}(q')$.
Thus, $q_0'\geq q'$ and if $b \in q_0'(N'\cap M)q_0'$ is so that
$q'b=0$ then $b=0$. Since $E$ is
implemented by $q'$, $E$ follows
faithful on $q_0'(N'\cap M)q_0'q$, implying that if
$b\in q_0'(N'\cap M)q_0'q$ and $a\in \theta(Ap)$ are positive elements
with $E(b)a=0$ then $ba=0$. But if $ba=0$ then
$0=E_N(ba)=E_N(b)a = (\tau(b)/\tau(q)) a$ (because $b$ commutes with
the factor $qNq$). This shows that $E(b) \in \theta(Ap)$
must have support equal to $q$ for any $b \in q_0'(N'\cap M)q_0'q$,
with $b\geq 0,
b \neq 0$.
Thus, if $f$ is a non-zero projection in $q_0'(N'\cap M)q_0'q$ then
$q'fq' = E(f)q'$ has supposrt $q'$. This implies that
any projection $f\neq 0$ in $q_0'(N'\cap M)q_0'q$ must have
trace $\tau(f)
\geq \tau(q') \geq 1-\varepsilon/4$, showing that $N'\cap M$ has an
atom $q_1'$ of trace $\tau(q_1') \geq 1-\varepsilon/4$.

An easy calculation shows that if we denote by $\tilde{N} \subset M$
a unital subfactor with $q_1'\in \tilde{N}$ and
$q_1'\tilde{N}q_1'=Nq_1'$
(N.B.: $\tilde{N}$ is obtained by amplifying $Nq_1'$ by $1/\tau(q_1')$), then
$F\subset_\varepsilon \tilde{N}$. Also, $\tilde{N}'\cap M = \Bbb C$
by construction. Thus, by the first part
of the proof, $\tilde{N} \in \Cal H\Cal T_s$.
Since $N$ is isomorphic to a reduction of $\tilde{N}$
by a projection, by part $1^\circ$ it
follows that $N\in \Cal H\Cal T_s$ as well.
\hfill Q.E.D.

\proclaim{6.7. Corollary} $1^\circ$. If $A \subset M$ is a ${\text{\rm HT}}$
Cartan subalgebra then any automorphism of $M$ can
be perturbed by an inner automorphism
to an automorphism that leaves $A$ invariant, i.e.,
$$
 {\text{\rm Aut}}M/{\text{\rm Int}}M=
{\text{\rm Aut}}(M, A)/{\text{\rm Int}}(M, A).
$$

$2^\circ$. Let $M \in \Cal H\Cal T_{_{s}}$
with $A \subset M$ its ${\text{\rm HT}}_{_{s}}$ Cartan subalgebra.
Denote by $\Cal G_{_{HT}}(M)$
the subgroup of ${\text{\rm Aut}}(M)$ generated
by the inner automorphisms and by the automorphisms
leaving all elements of $A$ fixed. Then $\Cal G_{_{HT}}(M)$
is an open-closed normal subgroup of ${\text{\rm Aut}}(M)$,
the quotient
group ${\text{\rm Out}}_{_{HT}}(M)\overset\text{\rm def} \to =
{\text{\rm Aut}}(M)/\Cal G_{_{HT}}(M)$ is countable and is
naturally isomorphic to the group of outer automorphisms
of $\Cal R^{^{HT}}_M$,
${\text{\rm Out}}(\Cal R^{^{HT}}_M)
\overset\text{\rm def} \to =
{\text{\rm Aut}}(\Cal R^{^{HT}}_M)/{\text{\rm Int}}(\Cal R^{^{HT}}_M)$.
\endproclaim
\noindent
{\it Proof}. $1^\circ$. If $\theta \in$Aut$(M)$ then $\theta(A)$
is HT Cartan, so by Theorem 6.2 there exists a
unitary element $u\in M$ such that $u\theta(A)u^* = A$.

$2^\circ$. This is trivial by 4.4.
\hfill Q.E.D.

\proclaim{6.8. Corollary} If $M \in \Cal H\Cal T$ then any central sequence
of $M$ is contained in the ${\text{\rm HT}}$ Cartan subalgebra of $M$.
Thus, $M'\cap M^\omega$ is always abelian and
$M$ is non$-\Gamma$ if and only if the equivalence relation
$\Cal R_M^{^{HT}}$ is strongly
ergodic. In particular, $M\not\simeq M\overline{\otimes} R$.
\endproclaim
\noindent{\it Proof}. If $A \subset M$ is the HT Cartan subalgebra
of $M$ and $A_0\subset A$ is so that $A_0 \subset M$ is rigid and
$A_0'\cap M=A$ then by  4.3.1$^\circ$ we have $M'\cap M^\omega =
M'\cap (A_0'\cap M)^\omega = M'\cap A^\omega$. \hfill Q.E.D.

\vskip .1in 
\noindent 
{\bf 6.9. Examples}. We now give a list of
examples of HT inclusions of the form $B \subset B
\rtimes_{\sigma} \Gamma_0$ and of factors in the class $\Cal H\Cal
T$ of the form $L^{\infty}(X, \mu) \rtimes \Gamma_0$, based on the
examples in 5.2, 5.3.2$^\circ$, 5.3.3$^\circ$. Note that if
$B\subset B\rtimes_\sigma \Gamma_0$ is an HT inclusion then
$\Gamma_0$ must have the property H (cf. 3.1), but that in Section
5 we were able to provide examples of inclusions $B \subset
B\rtimes_\sigma \Gamma_0$ satisfying the rigidity condition
(6.1.2) only for certain property H groups $\Gamma_0$ and for
certain actions of such groups (see Problem 6.12 below). Note also
that by Theorem 6.2 if $M=L^{\infty}(X, \mu) \rtimes_\sigma
\Gamma_0$ belongs to the class $\Cal H\Cal T$ and $\Gamma_0$ is a
property H group then $A=L^{\infty}(X, \mu)$ is automatically the
(unique) HT Cartan subalgebra of $M$, i.e., $A\subset M$ must
satisfy the rigidity condition (6.1.2) as well. 

\vskip .05in
\noindent 
{\it 6.9.1}. Let $\Gamma_0=SL(2, \Bbb Z)$,
$B_0=L_\alpha(\Bbb Z^2)$, for some $\alpha \in \Bbb T\subset \Bbb
C$, and $\sigma_0$ be the action of the group $SL(2, \Bbb Z)$ on
$B_0$ induced by its action on $\Bbb Z^2$. Then $B_0 \subset
M_\alpha \overset\text{\rm def} \to = B_0 \rtimes_{\sigma_0} SL(2,
\Bbb Z)$ is a HT$_{_{s}}$ inclusion with $M_\alpha$ a type II$_1$
factor. In case $\alpha$ is not a root of 1, this gives
HT$_{_{s}}$ inclusions $R=B_0 \subset M_\alpha$ and when $\alpha$
is a $n$'th primitive root of $1$, this gives HT$_{_{s}}$
inclusions $B_0\subset M_\alpha$, with $B_0$ homogeneous of type
I$_n$ and diffuse center. Indeed, in all these examples the
property (6.1.1) is satisfied by 3.2, and property (6.1.2) is
satisfied by 5.1. Moreover, by the isomorphism in 5.2.1$^\circ$,
if $\alpha$ is a root of $1$ then $M_\alpha \in \Cal H\Cal
T_{_{s}}$ and any maximal abelian subalgebra of $B_0=L_\alpha(\Bbb
Z^2)$ is Cartan in $M_\alpha$. 

\vskip .05in 
\noindent 
{\it 6.9.1}'. If we take the inclusion $A=L(\Bbb Z^2) \subset L(\Bbb
Z^2) \rtimes SL(2, \Bbb Z)=M$ from the previous example, which we
regard as the group measure space construction $L^\infty(\Bbb T^2,
\lambda)\subset L^\infty(\Bbb T^2, \lambda)\rtimes SL(2, \Bbb Z)$,
through the usual identification of $\Bbb T^2$ with the dual of
$\Bbb Z^2$ and of $L^\infty(\Bbb T^2, \lambda)$ with $L(\Bbb
Z^2)$, and we ``cut it in half'' with a projection $p \in A$ of
trace $1/2$, then we obtain the inclusion $(Ap \subset pMp) \simeq
(L^\infty(\Bbb S^2, \lambda) \subset L^\infty(\Bbb S^2, \lambda)
\rtimes PSL(2, \Bbb Z))$, where $\Bbb S^2$ is the 2-sphere. Thus,
by 6.9.1 and Theorem 6.6, it follows that $L^\infty(\Bbb S^2,
\lambda) \rtimes PSL(2, \Bbb Z)\in \Cal H\Cal T_{_{s}}$. 

\vskip.05 in 
\noindent 
{\it 6.9.2}. If $\Bbb F_n \subset SL(2, \Bbb Z)$
is an embedding with finite index and $\sigma_0$ is the
restriction to $\Bbb F_n$ of the action $\sigma_0$ on
$B_0=L_\alpha(\Bbb Z^2)$ considered in $1^\circ$, then $B_0
\subset B_0 \rtimes_{\sigma_0} \Bbb F_n$ is a HT$_{_{s}}$
inclusion, which in case $\alpha=1$ is a  HT$_{_{s}}$ Cartan
subalgebra. Also, if $p\in L(\Bbb Z^2)$ has trace $(12(n-1))^{-1}$
then the inclusion $(L(\Bbb Z^2)p \subset p(L(\Bbb Z^2 \rtimes
SL(2, \Bbb Z))p)$ is a HT$_{_{s}}$ Cartan subalgebra of the form
$(A \subset A\rtimes \Bbb F_n)$. In all these cases, again,
property (6.1.1) is satisfied by 3.2, and property (6.1.2) is
satisfied by 5.2.2$^\circ$. 

\vskip .05in 
\noindent 
{\it 6.9.3}. If
$\Gamma_0$ is an arithmetic lattice in $SU(n,1), SO(n,1), n \geq
2,$ then there exist free weakly mixing trace preseving actions
$\sigma_0$ of $\Gamma_0$ on $A=L^\infty(X, \mu)$ such that $A
\subset M=A\rtimes_{\sigma_0} \Gamma_0$ is HT$_{_{s}}$ Cartan (all
this cf. 3.2 and 5.2.3$^\circ$). 

\vskip .05in 
\noindent 
{\it 6.9.4}. If $\Gamma_0=SL(2, \Bbb Q)$, $A =L(\Bbb Q^2)$ and
$M=L(\Bbb Q^2 \rtimes SL(2, \Bbb Q)) = A\rtimes SL(2, \Bbb Q)$,
then $A\subset M$ is HT Cartan but not HT$_{_{s}}$ Cartan (cf. 3.2
and 5.3.2$^\circ$). 

\vskip .05in 
\noindent 
{\it 6.9.5}. Let
$\Gamma_0, \sigma_0, (B_0, \tau)$ be as in 6.9.1, 6.9.2 or 6.9.3.
Let $n\geq 1$ and $B=B_0^{\otimes n}, \sigma = \sigma_0^{\otimes
n}$. Then $B \subset B \rtimes_{\sigma} \Gamma_0$ is an
HT$_{_{s}}$ inclusion (cf. 3.2, 3.3.3 and 5.2). Moreover, if
$B_0=A_0$ is abelian, then $A_0^{\otimes n}= A \subset A
\rtimes_{\sigma} \Gamma_0$ is HT$_{_{s}}$ Cartan. 

\vskip .05in
\noindent 
{\it 6.9.6}. Let $\Gamma_0, \sigma_0, (B_0, \tau)$ be
any of the actions considered above.  Let $\sigma_1$ be an ergodic
action of $\Gamma_0$ on a von Neumann algebra $B_1\simeq
L^\infty(X, \mu)$. If $B = B_0\overline{\otimes} B_1$ and $M=B
\rtimes_{\sigma_0\otimes \sigma_1} \Gamma_0$, then $B \subset M$
is a HT inclusion (cf. 3.2 and 5.2.4$^\circ$). In particular, if
$B_0=A_0, B_1=A_1$ are abelian and $A = A_0\overline{\otimes}
A_1$, then $A \subset M$ is a HT Cartan subalgebra. If $\sigma_1$
is taken to be a Bernoulli shift, then $A \subset M$ is not
HT$_{_{s}}$ Cartan. For any such group $\Gamma_0$ the action
$\sigma_1$ can be taken non-strongly ergodic by ([CW]). In this
case, the resulting factor $M$ has the property $\Gamma$ of Murray
and von Neumann, with $M'\cap M^\omega = M'\cap A^\omega$ abelian.
Note that for each of the groups $\Gamma_0$ this gives three
distinct HT Cartan subalgebras of the form $A \subset A \rtimes
\Gamma_0$ (cf. 5.3.3$^\circ$). 

\vskip .05in 
\noindent 
{\it 6.9.7}.
Let $\Gamma_0, \sigma_0, (B_0, \tau)$ be any of the actions
considered above (so that $B_0 \subset B_0 \rtimes_{\sigma_0}
\Gamma_0$ is an HT inclusion). Let also $\Gamma_1$ be a property H
group and $\gamma$ an action of $\Gamma_1$ on $\Gamma_0$ such that
$\Gamma=\Gamma_0 \rtimes_\gamma \Gamma_1$ has the property H (for
instance, if $\Gamma_1$ is amenable or if $\gamma$ is the trivial
action, giving $\Gamma= \Gamma_0 \times \Gamma_1$). Let $\sigma$
denote the action $\sigma_0\rtimes \sigma_1$ on
$B=\overline{\otimes}_{g\in \Gamma_1} (B_0, \tau_0)_g$ constructed
in 3.3.3. Then $B \subset B\rtimes_\sigma \Gamma$ is an HT
inclusion, which follows an HT Cartan subalgebra in case $B_0$ is
abelian (cf. 3.1, 3.3.3, and the definitions).

\proclaim{6.10. Corollary} $1^\circ$. If $M$ is a
McDuff factor, i.e., $M \simeq M\overline{\otimes} R$,
then $M \notin \Cal H\Cal T$. In particular,
$R \notin \Cal H\Cal T$.

$2^\circ$. If $M$ contains a relatively rigid type ${\text{\rm II}}_1$ von
Neumann subalgebra then $M \notin \Cal H\Cal T$. In particular,
if $M$ contains $L(G)$ for some infinite property ${\text{\rm T}}$
group $G$, or if $M$ contains a property ${\text{\rm T}}$ factor, then
$M\notin \Cal H\Cal T$.

$3^\circ$. If $M$ contains a copy
of some $L_\alpha(\Bbb Z^2) \rtimes_\sigma \Gamma_1$,
with $\Gamma_1$ a subgroup of finite index in $SL(2, \Bbb Z)$ and
$\alpha$ an irrational rotation, then
$M\notin \Cal H\Cal T$.

$4^\circ$. If $M$ has the property ${\text{\rm H}}$
(e.g., $M \simeq L(\Bbb F_n)$ for some $2\leq n \leq \infty$)
then $M \notin \Cal H\Cal T$. In fact such factors do not
even contain subfactors in the class $\Cal H\Cal T$.
\endproclaim
\noindent
{\it Proof}. $1^\circ$ is trivial by 6.8,
$2^\circ$ and $3^\circ$ are clear by 5.4.2$^\circ$ and $4^\circ$
follows from 5.4.1$^\circ$.
\hfill Q.E.D.

\vskip .1in
\noindent
{\it 6.11. Definition}. A countable discrete group $\Gamma_0$
is a H$_{_{T}}$ (resp. H$_{{_{T}}_{_{s}}}$)
{\it group} if there exists a free ergodic measure
preserving action $\sigma_0$ of $\Gamma_0$ on the
the standard probability space $(X, \mu)$
such that $L^\infty(X, \mu) \subset L^{\infty} (X, \mu) \rtimes_{\sigma_0}
\Gamma_0$ is an HT (resp. HT$_{_{s}}$) Cartan subalgebra. Note
that a H$_{_{T}}$ group has the property H but is not amenable.
\vskip .1in
\noindent
{\bf 6.12. Problems}. 1$^\circ$. Characterise the class of all
H$_{_{T}}$ (resp. H$_{{_{T}}_{_{s}}}$) groups.

$2^\circ$. Construct examples of free ergodic
measure preserving actions $\sigma$ of $\Gamma_0=\Bbb F_n$
(or of any other non-inner amenable property H group $\Gamma_0$)
on $A=L^\infty(X, \mu)$ such that $A \subset M= A \rtimes_\sigma
\Gamma_0$ is not HT Cartan. Is this the case if $\sigma$ is a Bernoulli
shift ?

\proclaim{6.13. Corollary} $1^\circ$. $SL(2, \Bbb Z), \Bbb F_n, n\geq 2$
as well as any
arithmetic lattice in $SU(n, 1)$ or $SO(n, 1), n\geq 2$, are
${\text{\rm HT}}_{_{s}}$ groups.

$2^\circ$. Let $\Gamma \subset \Gamma_0$ be an inclusion of groups
with $[\Gamma_0: \Gamma] < \infty$. Then $\Gamma_0$ is an ${\text{\rm HT}}$
(resp. ${\text{\rm HT}}_{_{s}}$) group if and only
if $\Gamma$ is an ${\text{\rm HT}}$ (resp. ${\text{\rm HT}}_{_{s}}$) group.

$3^\circ$. If $\Gamma_0$ is an ${\text{\rm HT}}$ group and $\Gamma_1$
has the property ${\text{\rm H}}$
(for instance, if $\Gamma_1$ is amenable) then $\Gamma_0 \times \Gamma_1$
is an ${\text{\rm HT}}$ group.

$4^\circ$. If $\Gamma_0$ is an ${\text{\rm HT}}$ group and
$\Gamma_1$ is amenable and acts on $\Gamma_0$ then $\Gamma_0\rtimes \Gamma_1$
is an ${\text{\rm HT}}$ group.
\endproclaim
\noindent
{\it Proof}. Part 1$^\circ$ follows from 6.9.1$^\circ-3^\circ$,
while parts 3$^\circ$ and $4^\circ$ follow from 6.9.7.

To prove $2^\circ$ note first that by 3.1 and
2.3.3$^\circ$, $\Gamma_0$ has the
property H iff $\Gamma$ has the property H (this
result can be easily proved directly, see e.g.
[CCJJV]).

If $\Gamma_0$ is an H$_{_{T}}$ group
and $A \subset A \rtimes_\sigma
\Gamma_0$ is HT Cartan and
$A_0\subset A$ is so that
$A_0 \subset M$ is rigid and $A_0'\cap M=A$
then $A_0 \subset A \rtimes_\sigma \Gamma$ is also rigid,
by 4.4.2$^\circ$. Moreover,
the fixed point algebra $A^{\Gamma}$
is atomic (because $[\Gamma_0:\Gamma] < \infty$),
so if $p$ is any minimal projection in $A^{\Gamma}$ then
$p(A \rtimes_\sigma \Gamma)p$ is a factor and
$Ap \subset  p(A \rtimes_\sigma \Gamma)p$ is an HT Cartan subalgebra.
Thus, $\Gamma$ is an H$_{_{T}}$ group.

Conversely, if $\Gamma$ is an H$_{_{T}}$ group,
then let $\Gamma_1 \subset \Gamma$ be
a subgroup of finite index so that $\Gamma_1 \subset
\Gamma_0$ is normal. By the first part, $\Gamma_1$ is an H$_{_{T}}$  group.
By part 4$^\circ$, it follows that $\Gamma_0$ is an H$_{_{T}}$ group.
\hfill Q.E.D.

\heading 7. Subfactors of an
$\Cal H\Cal T$ factor.\endheading

In this section we prove that the class $\Cal H\Cal T$ is closed
under extension/restrictions of finite Jones index.
More than that, we show that the lattice of
subfactors of finite index of a factor in the class $\Cal H\Cal T$
is extremely rigid.

\proclaim{7.1. Lemma} Let $N \subset M$ be an irreducible inclusion of factors
with $[M:N] < \infty$ and $A\subset N$
a  Cartan subalgebra of $N$. Denote by
$\Cal N=\Cal N_N(A)$ the normalizer of $A$ in $N$.
Then we have:

$1^\circ$. $A'\cap M$ is a homogeneous type ${\text{\rm I}}_m$
algebra, for some $1\leq m < \infty$, and if we denote $A_1 = \Cal Z(A'\cap M)$
then there exists a partition of the identity $q_1, q_2, ..., q_n\in
\Cal P(A_1)$ such that $A_1
= \Sigma_i Aq_i$ and $E_N(q_i)= E_A(q_i) = 1/n, \forall i$.

$2^\circ$. $\Cal N$
normalizes $A_1$ and
$Q\overset \text{\rm def} \to = {\text{\rm sp}} A_1N =
\overline{\text{\rm sp}} A_1\Cal N$ is a
type ${\text{\rm II}}_1$ factor containing $N$,
with $[Q:N]=n$. Moreover,  $A_1 \subset Q$
is a Cartan subalgebra and we have the non-degenerate
commuting square:
$$
\matrix
N &\subset& Q&\\
\cup && \cup && \\
A & \subset& A_1&
\endmatrix
$$

$3^\circ$. $\Cal N$ normalizes $A'\cap M= A_1'\cap M
\simeq M_{m\times m}(A_1)$ and $P\overset \text{\rm def} \to =
{\text{\rm sp}} (A_1'\cap M)N =
\overline{\text{\rm sp}} (A_1'\cap M)\Cal N$ is a
type ${\text{\rm II}}_1$ factor containing $Q$,
with $[P:Q] = m^2$. Moreover,
we have the non-degenerate commuting square
$$
\matrix
Q &\subset& P&\\
\cup && \cup && \\
A_1 & \subset& A_1'\cap M&
\endmatrix
$$

$4^\circ$. Any maximal abelian $*$-subalgebra $A_2$ of
$A'\cap M=A_1'\cap M$ is a Cartan
subalgebra in $P$, with
$A_2p \subset pPp$ implementing the same equivalence relation
as $A_1 \subset Q$, $\forall p\in \Cal P(A_2)$,
$\tau(p)=1/m$, i.e., $\Cal R_{A_2p \subset pPp} \simeq
\Cal R_{A_1\subset Q}$ (equivalently, $\Cal R_{A_2 \subset P}
\simeq (\Cal R_{A_1\subset Q})^m$), but with the
two Cartan subalgebras possibly differing by their $2$-cocycles.
\endproclaim
\noindent
{\it Proof}. Since $\Cal N$ normalizes $A$,
it also normalizes $A'\cap M$, and thus
$\Cal Z(A'\cap M) = A_1$ as well. In particular, $A_1\Cal N = \Cal N A_1$
and $(A'\cap M)\Cal N=\Cal N (A'\cap M)$,
showing that sp$A_1\Cal N$ and sp$(A'\cap M)\Cal N$ are $*$-algebras.
Since $\Cal N'\cap M=N'\cap M = \Bbb C$, this implies that $Q, P$ are factors.
In particular, this shows that the squares
of inclusions in 2$^\circ$ and $3^\circ$ are commuting
and non-degenerate. Also, by definitions, $A_1$
is Cartan in $Q$.

Since $N \subset Q$ is a
$\lambda$-Markov inclusion, for $\lambda^{-1}=[Q:N]$ (see e.g.,
[Po2] for the definition), it follows that $A\subset A_1$,
with the trace $\tau$ inherited from $M$, is
$\lambda$-Markov. Thus, $e=e^Q_N$ implements
the conditional espectation $E^{A_1}_A$ and
$A_1\subset
B=\langle A_1, A \rangle = \langle A_1, e \rangle$ gives the basic
construction for $A \subset A_1$. Moreover, since $A,A_1$ are abelian,
it follows that $\Cal Z(B)=A=J_{A_1} A J_{A_1}$ and that
$$
A_1'\cap B =
J_{A_1} A_1 J_{A_1} \cap (J_{A_1}AJ_{A_1})'=
J_{A_1} (A_1\cap A') J_{A_1} = J_{A_1}A_1J_{A_1}=A_1.
$$
Thus, $A_1$ is maximal abelian in $B$, implying that the Markov
expectation of $B$ onto $A_1$ given by $E(xey)=\lambda xy,$ for $x, y\in A_1$,
is the unique expectation of $B$ onto $A_1$.

Also, for each $u \in \Cal N$, Ad$u$ acts
on $A\subset A_1$ $\tau$-preservingly.
Thus, Ad$u$ extends uniquely
to an automorphism $\theta_u$ on $B =
\langle A_1, e^{A_1}_A\rangle =
\langle A_1, e_N^Q \rangle$ by $\theta_u(e^{A_1}_A)=e^{A_1}_A$.
This automorphism leaves invariant the Markov trace on $B$.
Also, since $\theta_u, u\in \Cal N,$ act ergodically on
$A=\Cal Z(B)$, it follows that $B$ is homogeneous of type I$_n$, for some $n$.
By ([K2]), it follows that there exists a matrix units system
$\{e_{ij}\}_{1\leq i,j \leq n}$ in $B$ such that $B = A\vee {\text{\rm sp}}
\{e_{ij}\}_{i,j}$ with $A_1 =\Sigma_i Ae_{ii}$.

By the uniqueness of the conditional expectation $E$ of $B$ onto $A_1$,
it follows that if we put $q_i=e_{ii}$ then $E(X)=\Sigma_i q_i X q_i,
\forall X\in B$. In particular, the index of
$A_1 \subset B$ is given by $\lambda^{-1} =n = \tau(e)^{-1}$ and
by Markovianity we have
$1/n = E(e)=\Sigma_i q_i e q_i$. Thus, $q_ieq_i=n^{-1} q_i$,
and so $eq_ie = n^{-1} e = E(q_i)e$ as well, since $\tau(e)=\tau(q_i)$.
This ends the proof of $1^\circ$ and $2^\circ$.

Now, since $A_1$ is the center of $B_1=A'\cap M = A_1'\cap M$
and Ad$u, u\in \Cal N,$ act ergodically on $A_1$, it also follows that
$B_1$ is homogeneous of type I$_m$, for some $m \geq 1$. This
clearly implies $3^\circ$.

To prove $4^\circ$, let $\{f_{ij}\}_{1\leq i,j\leq m}\subset B_1$
be a matrix units system in
$B_1$ such that $A_2 = \Sigma_j A_1 f_{jj}$ and $B_1 = \Sigma_{i,j} A_1 f_{ij}$
(cf. [K2]). $A_2$ is Cartan in $P$ because by
construction $f_{ij}$ are in the normalizing groupoid of $A_2$ in $P$.

For each
$u\in \Cal N$ let $v(u)$ be a unitary element in $B_1$
such that $v(u)(uf_{jj}u^*)v(u)^* = f_{jj}, \forall j$ (this is
possible because $uf_{jj}u^*$ and $f_{jj}$ have the same central trace
$1/m$ in $B_1$). Since $v(u)$ commute with $A_1=\Cal Z(B_1)$,
$\forall u \in \Cal N$, it
follows that $A_1f_{11}$ with the
action implemented on it by $\{v(u) u \mid u\in \Cal N\}$
is isomorphic to $A_1$ with the action implemented on it by $\Cal N$.
Thus, the equivalence relation $\Cal R_{A_1f_{11} \subset f_{11}Mf_{11}}$
is the same as the
equivalence relation $\Cal R_{A_1 \subset Q}$, but with the
2-cocycle coming from the multiplication between the unitaries $v(u)u,
u\in \Cal N$ (for $A_1f_{11} \subset f_{11}Mf_{11}$)
possibly differing from
the 2-cocycle given by the multiplication of the corresponding $u\in \Cal N$
(for $A_1 \subset Q$).  \hfill Q.E.D.

\proclaim{7.2. Lemma} $1^\circ$.
Let $A^1\subset M_1$ be a maximal abelian $*$-subalgebra in the type
${\text{\rm II}}_1$ factor $M_1$.
If there exists a von Neumann
subalgebra $A^0 \subset A^1$ such that $A^0 \subset M_1$ is rigid
and $A^1\subset {A^0}' \cap M_1$ has finite index (in the sense of
${\text{\rm [PiPo]}}$), then $A^1$ contains a von Neumann
subalgebra $A^1_0$ such that $A^1_0 \subset M_1$ is
rigid and ${A^1_0}'\cap M_1 = A^1$.

$2^\circ$. Let $M_0 \subset M_1$ be a subfactor of
finite index with an ${\text{\rm HT}}$
(resp. ${\text{\rm HT}}_{_{s}}$) Cartan subalgebra $A\subset M_0$. If
$A^1 \subset M_1$ is a maximal abelian $*$-subalgebra of $M_1$
such that $A^1 \supset A$ and $M_1$ has property ${\text{\rm H}}$
relative to $A^1$ then $A^1 \subset M_1$ is an ${\text{\rm HT}}$
(resp. ${\text{\rm HT}}_{_{s}}$) Cartan subalgebra.
\endproclaim
\noindent
{\it Proof}. $1^\circ$. Since $A^1\subset {A^0}'\cap M_1$ has finite index,
it follows that ${A^0}'\cap M_1$ is a type I$_{fin}$ von Neumann
algebra and $A^1$ is maximal abelian in it
(see e.g., [Po7]). It follows that there exists a
finite partition of the identity with projections
$\{f_k\}_k$ in $A^1$ such that $\{f_k\}_k' \cap {A^0}'\cap M_1 \subset A^1$.
Thus, if we let $A^1_0 \overset \text{\rm def} \to = \Sigma_k A^0f_k$,
then ${A^1_0}'\cap M_1 \subset A^1$. By
4.4.3$^\circ$ it follows that $A_0 \subset M_1$ is a rigid embedding.

$2^\circ$. This is an immediate application of $1^\circ$, once we notice that
if $A^0 \subset A$ is so that $A^0 \subset M_0$ is rigid
and ${A^0}'\cap M_0=A$ then $A \subset {A^0}'\cap M_1$ has index
majorized by $[M_1:M_0]$, implying that
$A^1 \subset {A^0}'\cap M_1$ has finite index as well.
\hfill Q.E.D.

\proclaim{7.3. Theorem} Let $N\subset M$ be an inclusion of type
${\text{\rm II}}_1$ factors with $[M:N] < \infty$. Then we have:

$1^\circ$. $N \in \Cal H\Cal T$ (resp.
$N \in \Cal H\Cal T_{_{s}}$) if and only if
$M \in \Cal H\Cal T$ (resp.
$M \in \Cal H\Cal T_{_{s}}$).

$2^\circ$. Assume $N'\cap M = \Bbb C$ and $N, M \in \Cal H\Cal T$.
If  $Q , P\subset M$
are the intermediate subfactors constructed out of
an ${\text{\rm HT}}$ Cartan subalgebra of $N$, as in $7.1$, then
$Q, P \in \Cal H\Cal T$ and
the triple inclusion $N \subset Q \subset P \subset M$ is canonical.
Moreover, the ${\text{\rm HT}}$ Cartan subalgebra of $P$ is an
${\text{\rm HT}}$ Cartan subalgebra in $M$.

$3^\circ$. If $M\in \Cal H\Cal T$ and $N \subset M$ is an
irreducible subfactor then $[M:N]$ is an integer. Moreover,
the canonical weights
of the graph $\Gamma_{N,M}$ of $N \subset M$ are integer numbers.
\endproclaim
\noindent
{\it Proof}. 1$^\circ$. Since the
algebra $\langle M, N \rangle$ in the
basic construction $N\subset M \subset \langle M, N \rangle$
is an amplification of $N$, by Theorem 6.6 it
follows that it is sufficient to
prove that if $N\in \Cal H\Cal T$ (resp.
$N \in \Cal H\Cal T_{_{s}}$) then $M\in \Cal H\Cal T$ (resp.
$M \in \Cal H\Cal T_{_{s}}$). By 6.6.1$^\circ$, it is in fact sufficient
to prove this implication in the case $N'\cap M = \Bbb C$.

Let $A \subset N$ be an HT Cartan subalgebra and
$A_1 = \Cal Z(A'\cap M) \subset Q$ be constructed out of $A\subset N$ as in
Lemma 7.1. We begin by showing that
$A_1 \subset Q$ is an HT Cartan subalgebra.
Let $q_1, q_2, ..., q_n \in A_1\subset Q$ be so that $A_1 = \Sigma_i Aq_i$,
$E_N(q_i)=E_A(q_i)= 1/n$,
as in Lemma 7.1. By the last part of 2.3.3$^\circ$, it
follows that $Q$ has property H relative to $A$. But by
the last part of 2.3.4$^\circ$ this implies
$Q$ has property H relative to $A_1$. Also,
$A_1 \subset Q$ satisfies the conditions in part 2$^\circ$
of Lemma 7.2,
implying that it is HT Cartan.

Next we prove that if $A_2$ is constructed as in part $3^\circ$ of Lemma 7.1,
then $A_2 \subset P$ is an HT Cartan subalgebra. Let
$\{e_{ij}\}_{1\leq i,j\leq m}\subset A_1'\cap M$ be a matrix units system
which together with $A_1$ generates $A_1'\cap M$ and such that
$A_2 = \Sigma_j A_1 e_{jj}$. Since $P$ has an
orthonormal basis made up of
unitary elements commuting with $A_1$, by the last part of
2.3.3$^\circ$ it follows that $P$ has property H
relative to $A_1$. By applying the last part of 2.3.4$^\circ$,
it then follows that
$P$ has property H relative to $A_2$. Then 7.2.2$^\circ$
applies to deduce that $A_2 \subset P$ is an HT Cartan subalgebra,
which is even HT$_{_{s}}$ when $A\subset N$
(and thus $A_1 \subset Q$) is HT$_{_{s}}$.

Having proved that $A_2 \subset P$ is
an HT Cartan subalgebra,
we now prove that $A_2$ is HT Cartan in $M$ as well.
By noticing that $A_2$ is maximal abelian in $M$,
7.2.2$^\circ$ shows that it
is sufficient to prove that $M$ has property H relative to $A_2$.
To do this, we prove that if $A_3$ is any maximal abelian
subalgebra in $A_2'\cap M_1$, where $M_1=
\langle M, P\rangle$, then $A_3 \subset M_1$ is HT Cartan.
This would finish the proof, because by the first part of 2.3.4$^\circ$
$M_1$ would have the property H relative to $A_2$, and then by the first part
of 2.3.3$^\circ$ this would imply $M$ has the property H relative to $A_2$.

Since $M_1$
is an amplification of $P \in \Cal H\Cal T$, by Theorem 6.6 it
follows that $M_1$, as well as any reduced of $M_1$ by
projections in $M_1$, belong to $\Cal H\Cal T$.
Let $\Cal N_1$ be the normalizer of $A_2$ in $P$. Since $A_2$ is regular in $P$,
$\Cal N_1'' = P$ and $\Cal N_1'\cap M_1 = P'\cap M_1$. Let
$\{p'_t\}_t$ be a partition of the identity with
minimal projections in $P'\cap M_1$.
For each $t$, the inclusion $A_2 p'_t \subset Pp'_t \subset p'_t M_1 p'_t$
satisfies the hypothesis of Lemma 7.1. Thus, if $A_2^t$
is a maximal abelian $*$-subalgebra of $(A_2p'_t)'\cap p_t'M_1p_t'$,
then $A_2p'_t$ is included in $A_2^t$ and by 7.1.4$^\circ$,
$A_2^t$ is semi-regular
in $p'_tM_1p'_t$. In addition, by 7.2.1$^\circ$ it follows that
$A_2^t$ contains a von Neumann subalgebra $A_0^t$ with ${A^t_0}'\cap
p'_tM_1p_t' = A_2^t$ and $A_0^t \subset p'_tM_1 p'_t$
rigid. Since $p'_tM_1p'_t\in\Cal H\Cal T$, by Theorem 6.2
it follwos that $A^t_2 \subset p'_tM_1p'_t$ is HT Cartan.
Moreover, $M_1 \in \Cal H\Cal T$ implies $A_3=\Sigma_t A_2^t$
is HT Cartan in $M_1$, while clearly $A_2 \subset A_3$, by construction.

$2^\circ$. The triple inclusion $(N\subset Q \subset P \subset M)$
depends on the choice of the Cartan subalgebra $A\subset N$. But such $A$ is unique
up to conjugacy by unitaries in $N$, which leave fixed $Q$ and $P$.
The fact that the HT Cartan subalgebra of $P$ is HT  Cartan in $M$
has been proved in part $1^\circ$.

$3^\circ$. With the notations in $1^\circ$, we have $[M:N]= nm^2[M:P]$,
with $[M:P]$ being itself an integer, since $P$ contains a Cartan subalgebra
of $M$ (see e.g., [Po8]).

The weights $v_k$ of $\Gamma=\Gamma_{N,M}$ 
are square roots of indices of irreducible subfactors
appearing in the Jones tower for $N \subset M$. Thus, $v_k$ are square roots of integers.
Since $v_*=1$, $[M:N]\in \Bbb N$
and $\Gamma$ is irreducible and
has non-negative integer
entries, by the relations comming from $\Gamma\Gamma^t
\vec{v} = [M:N] \vec{v}$,
it follows recursively that all $v_k$ must be integers.
\hfill Q.E.D.
\vskip .1in
\noindent
{\it 7.4. Definitions}. Let $N\subset M$ be an irreducible inclusion of
factors in the class $\Cal H\Cal T$ with $[M:N] < \infty$ and let
$N \subset Q \subset P \subset M$ be the canonical triple inclusion
defined in part $2^\circ$ of Theorem 7.3.
\vskip .05in
\noindent
{\it 7.4.1}. $N \subset Q \subset P \subset M$
is called the {\it canonical decomposition} of $N\subset M$.
\vskip .05in
\noindent
{\it 7.4.2}. If $M=Q$, i.e., if the HT Cartan
subalgebra $A$ of $N$ is so that $A'\cap M$ is abelian (thus HT Cartan in $M$)
and  $M={\text{\rm sp}}AN=M$, then
$N \subset M$ is a {\it type} $C_-$ inclusion (or subfactor). If
$N=P$, i.e., if $A'\cap M=A$ (so that $A$ is Cartan in both $N$ and $M$)
then $N\subset M$ is of {\it type} $C_+$. If
$P=Q$, i.e., if $A'\cap M$ is abelian, then $N\subset M$
is of {\it type} $C_{\pm}$.
\vskip .05in
\noindent
{\it 7.4.3}. If $N=Q, P=M$ then $N\subset M$ is
of {\it type} $C_0$. More generally, an
extremal inclusion $N\subset M$
of factors in the class $\Cal H\Cal T$ is of {\it type} $C_0$ if
the HT Cartan subalgebra $A$ of $N$ satisfies
$A'\cap M= A \vee P_0$, with
$P_0\simeq M_{m\times m}(\Bbb C), m = [M:N]^{1/2},$ and
$M={\text{\rm sp}} (A'\cap M) N ={\text{\rm sp}}P_0 N$.

\proclaim{7.5. Theorem} $1^\circ$. Let
$N\subset M$ be an irreducible inclusion
of factors in the class $\Cal H\Cal T$, with $[M:N] < \infty$.
$N\subset M$ is of type $C_-$
(resp. $C_+, C_{\pm}, C_0$) if and only if its dual inclusion
$M \subset \langle M, N \rangle$ is of type $C_+$ (resp. $C_-,
C_{\pm}, C_0$).

$2^\circ$. If $N\subset M$ and $M \subset L$ are
irreducible inclusions
of factors in the class $\Cal H\Cal T$ with finite index
and both of type $C_-$ (resp. $C_+$), then $N \subset L$ is an
irreducible inclusion of type $C_-$ (resp. $C_+$).

$3^\circ$. If $N\subset M$ and $M \subset L$ are
extremal inclusions of factors in the class $\Cal H\Cal T$, both of
type $C_0$, then $N \subset L$ is of type $C_0$ and
so are all subfactors of the form $Np \subset pLp$, with
$p\in \Cal P(N'\cap L)$.

$4^\circ$. Let $N\subset M$ and $M \subset L$ be
irreducible inclusions
of factors in the class $\Cal H\Cal T$ with finite index
and such that $N\subset M$ is of type $C_+$ and $M\subset L$
is of type $C_-$. If $A \subset N$
is a $\text{\rm HT}$ Cartan subalgebra then $A'\cap L$ is abelian
and each irreducible inclusion
$Np \subset pLp$ for $p$ minimal projection in $N'\cap L$
is of type $C_{\pm}$. In particular this is the case
if $(M\subset L)=(M\subset \langle M, N \rangle)$.

$5^\circ$. Let $N\subset M$ be an inclusion
of factors in the class $\Cal H\Cal T$ with
$[M:N] < \infty$. If $N\subset M$ is either irreducible
of type $C_-$ or
extremal of type $C_0$ then $N \subset \langle M, N \rangle$
is a type $C_0$ inclusion, and so are all subfactors
of the form $Np \subset p\langle M, N \rangle p$,
for $p$ projection in $N'\cap \langle M,N \rangle$.
\endproclaim
\noindent
{\it Proof}. $1^\circ$. Let $A\subset N$ be an HT Cartan subalgebra of $N$.
If $N\subset M$ is of type $C_-$ then let $A'\cap M = \Sigma_i Aq_i$,
where $\{q_i\}_{1\leq i\leq n}
\subset  A'\cap M$ is a partition of the identity
with projections satisfying $E_N(q_i)=1/n, \forall i$.
Let $\alpha=e^{2\pi i/n}$ and denote
$u = n \Sigma_i q_i e_N q_{i+1}$. We clearly have $[u,A]=0$,
$uq_iu^*=q_{i+1}$ and $E_N(u^j)=0, \forall j\leq n-1$.
Thus, the HT Cartan subalgebra $A_1
= A'\cap M$ of $M$ is maximal abelian in $\langle M, N \rangle $ and
is normalized by
$u^j$, with $\langle M, N \rangle=\Sigma_j u^jM$, i.e., $A_1$
is the HT Cartan subalgebra in $\langle M, N\rangle$ as well,
showing that $M \subset \langle M, N \rangle$ is of type $C_+$.

If $N \subset M$ is of type $C_+$, $A \subset N \subset M$ is HT Cartan
in both factors and $u_1, u_2, ..., u_n \in \Cal N_M(A)$
are unitary elements such that $M= \Sigma_i u_i N$
and $E_N(u_i^* u_j)=\delta_{ij}$ then
$q_j = u_je_Nu_j^*$ is a partition of the identity
with projections in  $ \langle M,N \rangle$ and we have
$A'\cap  \langle M,N \rangle = \Sigma_j q_j A,
\langle M,N \rangle = \Sigma_j q_j M$. Thus, $M \subset
\langle M,N \rangle$ is of type $C_-$.

If $N\subset P \subset M$ is so that $N\subset P$ is $C_-$, $P\subset M$
is $C_+$ then we have the irreducible
inclusions $M \subset \langle M, P \rangle$, which is $C_-$, and
$\langle M, P \rangle\subset
\langle M, N \rangle$, which is an amplification
of $P \subset \langle P, N \rangle$, thus of type $C_+$. This shows that
$M \subset \langle M, N \rangle$ is $C_{\pm}$.

If $N \subset M$ is of type $C_0$ and $A\subset N$ is an
HT Cartan subalgebra with $A'\cap M = \Sigma_{i,j} e_{ij} A$
for some
matrix units system $\{e_{ij}\}_{1\leq i,j \leq m} \subset A'\cap M$,
then let $e'_{ij}= m \Sigma_k e_{ki}e_Ne_{jk}, 1\leq i,j \leq m$.
It is immediate to check that $\{e'_{ij}\}_{i,j}$ is a
matrix units system which commutes with $A$ and with
$\{e_{kl}\}_{k,l}$, that
$\{e'_{ij}\}_{i,j}$
is an orthonormal basis of $\langle M,N \rangle$ over $M$
and that $\{e'_{ij}e_{kl}\}_{i,j,k,l}$ is an
orthonormal basis
of $\langle M, N \rangle$ over $N$. It follows that
$A'\cap \langle M,N \rangle = {\text{\rm sp}}\{e'_{ij}e_{kl}\}_{i,j,k,l}A$.
Thus, if $A_2 \subset A'\cap M$ is a maximal abelian subalgebra, then
$A_2'\cap \langle M, N \rangle = \Sigma_{i,j} e'_{ij} A_2$.
This shows that $M \subset \langle M, N \rangle$ is of type $C_0$.

$2^\circ$. By duality in the Jones tower ([PiPo]) and
part 1$^\circ$, it
is sufficient to prove that if $N\subset M, M\subset L$ are
of type $C_+$ then so is $N \subset L$. But this is trivial,
since if $A\subset N$ is HT Cartan in $N$ then
it first follows Cartan in $M$,
then in $L$.

$3^\circ$. Let $\{e_{ij}\}_{1\leq i,j \leq m} \subset A'\cap M$
be a matrix units system such that $A'\cap M = \Sigma_{i,j} e_{ij} A$,
as in the proof of the last part of $1^\circ$
(thus, $[M:N]=m^2$). Let $A_2 = \Sigma_j e_{jj} A$,
which follows HT Cartan in $M$. Let $\{f'_{kl}\}_{1\leq k,l \leq m'}
\subset A_2'\cap L$ be a matrix units system such that $A_2'\cap L
= \Sigma_{k,l} f'_{kl}A_2$, with ${m'}^2 = [L:M]$.
Then $\{f_{ts}\}_{t,s} = \{ e_{i1}f'_{kl}e_{1j} \mid
1\leq i, j\leq m, 1\leq k,l \leq m'\}$ is a matrix units system in
$A'\cap L$ and if we denote $P_0\simeq
M_{mm' \times mm'}(\Bbb C)$ the algebra it generates, then
clearly $E_N(f_{st}) = \delta_{st}/mm'$. Since
$[L:N] = (mm')^2$, and since we have the commuting square
$$
\matrix
N &\subset& L&\\
\cup && \cup && \\
A & \subset& A'\cap L&
\endmatrix
$$
as well as
$$
\matrix
N &\subset& L&\\
\cup && \cup && \\
A & \subset& A \vee P_0&
\endmatrix
$$
with $A\vee P_0 \subset A'\cap L$ and with $P_0$
containing an orthonormal system of $L$ over $N$ made up of
$mm'$ elements, it follows that $A\vee P_0 = A'\cap L$,
thus showing that $N \subset L$ is of type $C_0$.

Finally, if $p \in \Cal P(N'\cap L)$ then in particular
$p \in A \vee P_0$. By the above
commuting squares, we have $E_A(p)=E_N(p)=\tau (p)1$.
But $A=\Cal Z(A \vee P_0)$, implying that $p$ has
scalar central trace in $A \vee P_0$. Thus,
$(Ap)'\cap pLp=p(A \vee P_0)p$ is homogeneous of type I.
Since we also have $pLp=p({\text{\rm sp}} P_0 N)p
= p({\text{\rm sp}} P_0)p Np$, this shows that
$Np \subset pLp$ is of type C$_0$.

4$^\circ$. Let $A\subset N$ be the HT Cartan subalgera of $N$,
which is thus
HT Cartan in $M$ as well.
Thus $A_1=A'\cap L$ is abelian
with $L = {\text{\rm sp}} A_1 M$. Since
any irreducible projection $p\in N'\cap L$
lies in $A_1$, by cutting these
relations with $p$ we obtain that $(Ap)'\cap pLp$
is abelian, which by Lemma 7.1 means that $Np
\subset pLp$ has only type $C_-$ and $C_+$ components
in its canonical decomposition.

$5^\circ$. This is immediate from
the proofs in $1^\circ$ and the last part of 3$^\circ$.
\hfill Q.E.D.
\vskip .1in
\noindent
{\bf 7.6. Examples}. $1^\circ$. Let $\Gamma_0$ be a property H group
and $\sigma$ a free weakly mixing
measure preserving action of $\Gamma_0$ on the probability space
$(X, \mu)$ such that the
Cartan subalgebra $L^\infty(X, \mu) = A \subset N = L^\infty(X, \mu)
\rtimes_\sigma \Gamma_0$ contains a von Neumann subalgebra
$A_1 \subset A$ with $A_1'\cap N =A$ and $A_1 \subset N$ rigid.
Let $\Gamma_1 \subset \Gamma_0$ be a subgroup of finite index and
$\sigma_0$ the left action of $\Gamma_0$ on $\Gamma_0/\Gamma_1$.
Let $A_0 =
\ell^\infty(\Gamma_0/\Gamma_1)$ and
$M = A\otimes A_0 \rtimes_{\sigma\otimes\sigma_0} \Gamma_0$.

Then $N, M \in \Cal H\Cal T$ and if we identify
$N$ with the subfactor of $M$ generated by $A=A\otimes \Bbb C$
and by the canonical unitaries $\{u_g\}_g \subset M$ implementing
the action $\sigma\otimes \sigma_0$ on $A\otimes A_0$,
then $N \subset M$ is an irreducible type $C_-$ inclusion.
Moreover, if we denote $N_1 = A \vee
\{u_g\}_{g\in \Gamma_1} \simeq A \rtimes_\sigma \Gamma_1 \subset N$
then $N_1 \subset N$ is a type $C_+$ inclusion and $N_1 \subset N \subset M$
is a basic construction.

We have $[M:N]=[N:N_1]=[\Gamma_0:\Gamma_1]$, the standard
invariant of $N_1 \subset N$ coincides with the standard invariant
$\Cal G_{\Gamma_1 \subset \Gamma_0}$
of $R \rtimes \Gamma_1 \subset R \rtimes \Gamma_0$
studied in ([KoYa]) and the
standard invariant of $N \subset M$ is the dual
of $\Cal G_{\Gamma_1 \subset \Gamma_0}$. In particular,
$N_1\subset N \subset M$ are finite depth inclusions.

$2^\circ$. Let $\Gamma_0, \sigma, A$ be as in example $1^\circ$ above and
let $\pi_0$ be a finite dimensional irreducible projective representation
of $\Gamma_0$ on the Hilbert space $\Cal H_0$,
with scalar 2-cocycle $v$. Let $B_0 = \Cal B(\Cal H_0)$
and $\sigma_0(g)= {\text{\rm Ad}} \pi_0(g)$ be the action
of $\Gamma_0$ on $B_0$ implemented by $\pi_0$. Denote
$M=M_{\pi_0} = A\otimes B_0 \rtimes_{\sigma\otimes \sigma_0} \Gamma_0$
and let $N$  be the subfactor of $M$ generated by
$A \otimes 1 = \Cal Z(A \otimes B_0)$ and
by the canonical unitaries $\{u_g\}_{g\in \Gamma_0}\subset M$ implementing
the action $\sigma\otimes \sigma_0$. Thus,
$N \simeq A \rtimes_\sigma \Gamma_0$, $M \simeq M_{n\times n}
(A \rtimes_{\sigma, v} \Gamma_0)$ and both belong to the class $\Cal H\Cal T$.

Moreover, $N \subset M$ is an irreducible type $C_0$ inclusion
and its standard invariant coincides with the standard invariant
of the generalized
Wassermann-type subfactor corresponding to the
projective representation
$\pi_0$, i.e.:
$$
\matrix
\Bbb C &\subset&  End(\Cal H_0)^{\sigma_0} &
\subset& End(\Cal H_0\otimes
\overline{\Cal H_0})^{\sigma_0\otimes \overline{\sigma_0}}
& \subset& \dots \\ && \cup && \cup \\
&&\Bbb C & \subset&  \Bbb C \otimes
End(\Cal H_0)^{\overline{\sigma_0}} &\subset& \dots
\endmatrix
$$

$3^\circ$. Let $\sigma$ be the action of
$SL(2, \Bbb Z)$ on $L_\alpha(\Bbb Z^2)$
implemented by the action of $SL(2, \Bbb Z)$ on $\Bbb Z^2$,
as in 5.2.1$^\circ$ and 6.9.1$^\circ$, for $\alpha$ a
primitive root of 1 of order $n$. Let $M_\alpha = L_\alpha(\Bbb Z^2)
\rtimes_\sigma
SL(2, \Bbb Z)$,
$A = \Cal Z(L_\nu(\Bbb Z^2))$
and $N = A \vee \{u_g\}_g$ the von Neumann algebra generated by
$A$ and the canonical unitaries in $M_\alpha$ implementing the action $\sigma$.
Then $N, M_\alpha \in \Cal H\Cal T_{_{s}}$ and $N \subset M_\alpha$
is an irreducible
inclusion of type $C_0$ with $[M_\alpha:N]=n^2$. Indeed,
we have already noticed in 6.9.1$^\circ$ that $N \in \Cal H\Cal T_{_{s}}$,
so by 7.3 we have $M_\alpha \in \Cal H\Cal T_{_{s}}$. Also, by
construction we have $A'\cap M_\alpha = L_{\alpha}(\Bbb Z^2)
= A \otimes B_0$, with $B_0 \simeq M_{n\times n}(\Bbb C)$, and
$M_\alpha={\text{\rm sp}} L_\alpha(\Bbb Z^2) N$.

One can show that $N \subset M_\alpha$
is isomorphic to a type $C_0$ inclusion $N\subset M_{\pi_0}$ 
as in example $2^\circ$, when taking $\Gamma_0=SL(2, \Bbb Z)$, with 
$\sigma, \sigma_0$ the actions of $SL(2, \Bbb Z)$ on $A = \Cal Z(L_\alpha(\Bbb Z^2))
\simeq L((n\Bbb Z)^2)$, $B_0=L_\alpha((\Bbb Z/n\Bbb Z)^2)
\simeq M_{n\times n}(\Bbb C)$. Note that the standard invariant  
([Po3]) of $N \subset M_\alpha$ depends only on the order $n$ of $\alpha$,   
because if $\pi_0, \pi_0'$ are representations corresponding to  
primitive roots $\alpha, \alpha'$ of order $n$ then there exists 
an automorphism $\gamma$ of the group $(\Bbb Z/n\Bbb Z)^2$ such that 
$\pi'=\pi \circ \gamma$. But we do not know whether the 
isomorphism class of $N\subset M_\alpha$ depends only on $n$. 

\vskip .05in
We now reformulate the results in Theorem 7.5 in terms
of correspondences. For the definition of Connes' general $N-M$
{\it correspondences} (or $N-M$ {\it Hilbert bimodules}) $\Cal H$ $=_N\Cal H_M$, of
the {\it adjoint}
$\overline{\Cal H}=_M\overline{\Cal H}_N$ of
$\Cal H$, as well as for the definition of the {\it composition}
$\Cal H \circ \Cal K$
(also called {\it tensor product}, or {\it fusion}) of
correspondences $\Cal H=_N\Cal H_M, \Cal K=_M\Cal K_P$
see ([C7], [Po1], [Sa]).
\vskip .1in
\noindent
{\it 7.7. Definition}. Let $N,M \in \Cal H\Cal T$ and $\Cal K$ be a $N-M$
correspondence, viewed as a Hilbert $N-M$ bimodule. Assume that
dim$_N\Cal K_M \overset\text{\rm def} \to = {\text{\rm dim}}_N\Cal K
\cdot {\text{\rm dim}}\Cal K_M < \infty$ and that $\Cal K$ is
{\it irreducible}, i.e., $N \vee (M^{op})' = \Cal B(\Cal K)$.
$\Cal K$ is of {\it type} $C_-$ (resp. $C_+, C_{\pm}, C_0$) if
the inclusion $N \subset (M^{op})'$ is of type $C_-$ (resp. $C_+, C_{\pm}, C_0$),
in the sense of Definitions 7.4.

Finite index
correspondences (resp. bimodules)
between factors in the class $\Cal H\Cal T$
will also be called HT {\it correspondences} (resp. HT {\it bimodules}).

\proclaim{7.8. Corollary} Let $_N\Cal H_M, _M\Cal K_L$ be
irreducible ${\text{\rm HT}}$ bimodules.

$1^\circ$. $\Cal H$ is of type $C_-$ (resp. $C_+, C_{\pm}, C_0$) iff
$\overline{\Cal H}$ is of type $C_+$ (resp. $C_-, C_{\pm}, C_0$).

$2^\circ$. If both $\Cal H, \Cal K$ are of type $C_-$ (resp. $C_+$, resp.
$C_0$) then  $\Cal H\circ \Cal K$ is irreducible of type $C_-$ (resp.
irred. $C_+$, resp. a sum of irreducible $C_0$). In particular,
the class of $\text{\rm HT}$ bimodules (or correspondences)
of type $C_0$ over a
$\text{\rm HT}$ factor forms a selfadjoint tensor category.

$3^\circ$. If $\Cal H$ is of type $C_+$ and $\Cal K$ is of type $C_-$
then $\Cal H \circ \Cal K$ is a direct sum of irreducible type $C_{\pm}$
bimodules. Also, $\Cal K \circ \overline{\Cal K}$ is a direct sum
of irreducible $C_0$ bimodules.
\endproclaim
\noindent
{\it Proof}. Part $1^\circ$ is a reformulation of $7.5.1^\circ$, while
$2^\circ$ and $3^\circ$ are direct consequences of $7.5.2^\circ-5^\circ$.
\hfill Q.E.D.
\vskip .1in
\noindent
{\it 7.9. Definition}.
Let $M \in \Cal H\Cal T$ and $\theta \in $Aut$M$ be a periodic
automorphism of $M$, with $\theta^n = id$ and
$\theta^k$ outer $\forall 0 < k < n$.
Then $\theta$ {\it is of type} $C_-$ (resp. $C_+$) if
the inclusion $M \subset M\rtimes_\theta \Bbb Z/n\Bbb Z$
is of type $C_-$
(resp. $C_+$). By the uniqueness
of the HT Cartan subalgebra, this property is clearly
a conjugacy invariant for $\theta$.

\proclaim{7.10. Corollary} The factor $N = L(\Bbb Z^2 \rtimes SL(2, \Bbb Z))$
has two non-conjugate period two automorphisms, one of type
$C_-$ and one of type $C_+$.
\endproclaim
\noindent
{\it Proof}. In example 7.6.1$^\circ$, take $\Gamma_1 \subset
\Gamma_0 = SL(2, \Bbb Z)$
a subgroup of index $2$ and
$(X, \mu)=(\Bbb T^2, \mu)$ with
$SL(2, \Bbb Z)$ acting on it in  the usual
way. Then $N=L(\Bbb Z^2 \rtimes SL(2,\Bbb Z))$
and the resulting type $C_-$ inclusion
$N \subset M$ given by the construction
7.6.1$^\circ$ is of index 2. Thus, by Goldman's Theorem, it is
given by a period 2 automorphism $\theta_-$, which is thus of type $C_-$.
Alternatively,
we can take $\theta_-$ to be the
automorphism given by the non-trivial
character $\gamma$ of
$\Bbb Z^2 \rtimes SL(2, \Bbb Z)$ with $\gamma^2=1$, defined by
$\gamma(a)=-a$, $\gamma(b)=b$, on the generators
$a, b$ of period $4$, resp. $6$ of $SL(2, \Bbb Z)$,
and $\gamma(\Bbb Z^2)=1$.

Now take $\theta_+$ to be the automorphism of $N$
implemented by $\left( \matrix 1 & 0 \\
0 & -1 \endmatrix \right) \in GL(2, \Bbb Z)$. Thus, $N \subset M=
N\rtimes_{\theta_+} \Bbb Z/2\Bbb Z$ coincides with
$L(\Bbb Z^2 \rtimes SL(2, \Bbb Z)) \subset
L(\Bbb Z^2 \rtimes GL(2, \Bbb Z))$, and since $GL(2, \Bbb Z)$
acts freely on $\Bbb Z^2$, it follows that $L(\Bbb Z^2)'\cap M =
L(\Bbb Z^2)$, so that $N \subset M$ is of type $C_+$.
\hfill Q.E.D.
\vskip .1in
\noindent
{\bf 7.11. Question}. Let $N \simeq L(\Bbb Z^2 \rtimes SL(2, \Bbb Z))$.
Is then any irreducible type $C_-, C_+$ or $C_0$
inclusion
of factors $N \subset M$
isomorphic to a ``model''
inclusion 7.6.$1^\circ-7.6.2^\circ$ ?

\heading 8. Betti numbers for
$\Cal H\Cal T$ factors.
\endheading

\noindent {\it 8.1. Definition}. Let $M\in \Cal H\Cal T$  and
$\Cal R^{^{HT}}_M$ be the standard equivalence relation
implemented by the normalizer of the HT Cartan subalgebra of $M$,
as in Corollary 6.5. Let $\{\beta_n(\Cal R^{^{HT}}_M)\}_{n\geq 0}$
be the $\ell^2$-Betti numbers of $\Cal R_M^{^{HT}}$, as defined by
Gaboriau in ([G2]). For each $n = 0, 1, 2, ...$, we denote
$\beta^{^{HT}}_n(M) \overset\text{\rm def} \to = \beta_n(\Cal
R^{^{HT}}_M)$ and call it the $n$'th $\ell^2_{_{HT}}$-{\it Betti
number} (or simply the $n$'th Betti number) of $M$. By 6.5,
$\beta^{^{HT}}_n(M)$ are isomorphism invariants for $M$.

From the results in Section 6 and the properties proved by Gaboriau for
$\ell^2$-Betti numbers of standard equivalence relations, one
immediately gets:

\proclaim{8.2. Corollary} $0^\circ$. If $M$ is of type II$_1$ then
$\beta^{^{HT}}_0(M)=0$ and if $M=M_{n\times n}(\Bbb C)$ then
$\beta^{^{HT}}_0(M)=1/n$.

$1^\circ$. If $A\subset M=A\rtimes_\sigma \Gamma_0$
is a ${\text{\rm HT}}$ Cartan
subalgebra, for some countable discrete
group $\Gamma_0$ acting freely and ergodically on $A\simeq L^\infty(X, \mu)$,
then $\beta^{^{HT}}_n(M)$ is equal to the $n$'th $\ell^2$-Betti number
of $\Gamma_0$, $\beta_n(\Gamma_0)$,
as defined in ${\text{\rm [ChG])}}$.

$2^\circ$. If $M \in \Cal H\Cal T$
and $t > 0$ then
$\beta_n^{^{HT}}(M^t) = \beta(M)/t, \forall n$.

$3^\circ$. If $M_1, M_2 \in \Cal H\Cal T$ then for each $n\geq 0$
we have the K{\"u}nneth-type formula:
$$
\beta_n^{^{HT}} (M_1 \overline{\otimes} M_2) = \underset i+j=n \to \Sigma
\beta_i^{^{HT}}(M_1) \beta_j^{^{HT}}(M_2),
$$
where $0 \cdot \infty = 0$ and $b \cdot \infty = \infty$ if $b \neq 0$.

$4^\circ$. Let $M \in \Cal H\Cal T_{_{s}}$ and let
$N_k \subset M, k\geq 1,$
be an increasing sequence of subfactors with  $N_k \uparrow M$
(so that $N_k \in \Cal H\Cal T_{_{s}}$, for $k$ large enough,
by $6.8.3^\circ$). Then
$\underset k \rightarrow \infty \to \liminf \beta^{^{HT}}_n(N_k)
\geq \beta^{^{HT}}_n(M)$.
\endproclaim
\noindent
{\it Proof}. $0^\circ$. This is trivial by the definitions and ([G2]).

1$^\circ$. By 8.1, we have $\beta^{^{HT}}_n(M) =
\beta_n(\Cal R^{^{HT}}_M)$.
But $\Cal R^{^{HT}}_M= \Cal R_{\Gamma_0}$, and by Gaboriau's theorem the latter
has Betti numbers $\beta_n(\Cal R_{\Gamma_0})$ equal to the Cheeger-Gromov
$\ell^2$-Betti numbers $\beta_n(\Gamma_0)$ of the group $\Gamma_0$.

2$^\circ$. By Section 6 we know that the class $\Cal H\Cal T$ is
closed under amplifications and tensor products. Moreover,
by 1.4.3 the ``amplification''
by $t$ of a Cartan
subalgebra $A\subset M$ has a normalizer that gives
rise to the standard
equivalence relation $(\Cal R^{^{HT}}_M)^t$.
Then formula $2^\circ$ is a consequence of Gaboriau's
similar result for standard equivalence relations.

Part $3^\circ$ follows similarily, by taking into account that
if $A_1 \subset M_1, A_2 \subset M_2$ are Cartan subalgebras
then $\Cal N(A_1 \otimes A_2)'' = (\Cal N(A_1) \otimes
\Cal N(A_2))''$.

$4^\circ$. By 6.8.3$^\circ$, there exists $k_0$ and an HT$_{_{s}}$
Cartan subalgebra $A$ of $M$ such that $A\subset N_k, \forall k\geq k_0$.
Then the statement follows from Theorem 5.13 in ([G2]).
\hfill Q.E.D.

\proclaim{8.3. Corollary} $1^\circ$.
If $M\in \Cal H\Cal T$ has at least one non-zero, finite
Betti number then $\mycal F(M)=\{1\}$ and in fact
$M^{t_1} \overline{\otimes}
...\otimes M^{t_n}$ is isomorphic to
$M^{s_1} \overline{\otimes}
...\otimes M^{s_m}$ if and
only if $n=m$ and $t_1... t_n =
s_1 ... s_m$. Equivalently, $\{M^{\overline{\otimes}m}\}_{m \geq 1}$ are stably
non-isomorphic and
all the tensor powers $M^{\overline{\otimes}m}$ have trivial fundamental group,
$\mycal F(M^{\overline{\otimes}m})=\{1\}, \forall m \geq 1$.

$2^\circ$. If $M \in \Cal H\Cal T$ and $\beta_1^{^{HT}}(M) \neq 0, \infty$,
then $M$ is not the tensor product of two factors $M_1, M_2$
in the class $\Cal H\Cal T$. More generally if
$\beta_k^{^{HT}}(M)$ is the first non-zero finite Betti number
for $M$, then $M^{\overline{\otimes}m}$
cannot be expressed as the tensor
product of $km+1$ or more factors in the class $\Cal H\Cal T$.
\endproclaim
\noindent
{\it Proof}. $1^\circ$. First note that if $M$ has $\beta_k^{^{HT}}(M)$
as first non-zero
Betti number, then the formula
$\beta_k^{^{HT}}(M^t)=\beta_k^{^{HT}}(M)/t$ implies
that $M \not\simeq M^t$ if $t \neq 1$. Thus, $\mycal F(M)=\{1\}$.

Also, by the K{\"u}nneth formula 8.2.2$^\circ$,
if $\beta^{^{HT}}_{n_i}(M_i)$ is the first non-zero finite Betti number
for $M_i \in \Cal H\Cal T, i=1,2,$, and we put $n = n_1 + n_2$, then
$\beta_n^{^{HT}}(M_1 \overline{\otimes} M_2)
=\beta_n^{^{HT}}(M_1) \beta_n^{^{HT}}(M_2)$,
is the first non-zero finite Betti number
for $M_1 \overline{\otimes} M_2$.

Thus,
$\beta_{km}^{^{HT}}(M^{\overline{\otimes}m})$ is the first non-zero finite
Betti number for $M^{\overline{\otimes}m}, m\geq 1$, showing that
$\{M^{\overline{\otimes}m}\}_{m \geq 1}$ are stably
non-isomorphic.

$2^\circ$. This is trivial by the
first part of the proof and the K{\"u}nneth formula 8.2.2$^\circ$.
\hfill Q.E.D.

\proclaim{8.4. Corollary} $1^\circ$. Let $N \subset M$
be an irreducible inclusion of factors in the class $\Cal H\Cal T$ with
$[M:N] < \infty$. If $N\subset M$ is of type $C_-$ then
$\beta_n^{^{HT}}(M)=\beta_n^{^{HT}}(N), \forall n$.
If $N\subset M$ is of type $C_+$ then $\beta_n^{^{HT}}(M)
=[M:N] \beta_n^{^{HT}}(N)$.

$2^\circ$. Let $N \subset M$
be an extremal inclusion of factors in the class $\Cal H\Cal T$.
If $N\subset M$ is of type $C_0$ then
$\beta_n^{^{HT}}(M)=[M:N]^{1/2}\beta_n^{^{HT}}(N), \forall n$.

$3^\circ$. If $N \subset Q \subset P \subset M$ is the canonical
decomposition of an irreducible
inclusion of factors $N\subset M$ in the class $\Cal H\Cal T$, then
$\beta_n^{^{HT}}(Q)=\beta_n^{^{HT}}(N)$, $\beta_n^{^{HT}}(P)
=[P:Q]^{1/2} \beta_n^{^{HT}}(N)$ and $\beta_n^{^{HT}}(M)
=[M:P] \beta_n^{^{HT}}(P)$.

$4^\circ$. Let $M\in \Cal H\Cal T$,
$N\subset M$ a
subfactor of finite index, $(\Gamma_{N,M}, (v_k)_k)$
the graph of $N\subset M$, with its standard weights.
Let also
$\{\Cal H_k\}_k$ be the list of irreducible Hilbert $M$-bimodules
appearing in some
$L^2(M_n, \tau), n=0, 1, 2, ...$, with
$\{M\subset M(\Cal H_k)\}_k$ the corresponding irreducible
inclusions of factors.
If $\beta^{^{HT}}_n(M) \neq 0, \infty$ for some $n \geq 1$ then
$v_k = \beta_n^{^{HT}}(M(\Cal H_k))/\beta_n^{^{HT}}(M), \forall k$.
Thus,
$$
\Gamma_{N,M} \Gamma^t_{N,M} (\beta_n^{^{HT}}(M(\Cal H_k)))_k =
[M:N] (\beta_n^{^{HT}}(M(\Cal H_k)))_k.
$$
\endproclaim
\noindent
{\it Proof}. 1$^\circ$. If $N\subset M$ is of type
$C_+$ then $\Cal R^{^{HT}}_N$
is a subequivalence relation of index $[M:N]$ in $\Cal R^{^{HT}}_M$,
so by ([G2]) we have
$$
\beta_n^{^{HT}}(M) = \beta_n(\Cal R^{^{HT}}_M) =
[M:N]\beta_n(\Cal R^{^{HT}}_N) = [M:N]\beta_n^{^{HT}}(N).
$$

If $N\subset M$ is of type $C_-$ then by part $1^\circ$ of Theorem 7.5
$M \subset \langle M,N \rangle$ is of type $C_+$. Since
$\langle M,N \rangle$ is the $[M:N]$-amplification
of $N$, by the first part and by formula 8.2.2, we get:
$$
\beta_n^{^{HT}}(N) = [M:N]^{-1} \beta_n^{^{HT}}(\langle M,N \rangle)
=[M:N]^{-1}[M:N]\beta_n^{^{HT}}(M)
$$

$2^\circ$. If $N\subset M$ is of type $C_0$ then by 7.1
the equivalence relation
$\Cal R_M^{^{HT}}$ is a $[M:N]^{1/2}$-amplification of $\Cal R_N^{^{HT}}$.
Thus, $\beta_n^{^{HT}}(M)=[M:N]^{1/2}\beta_n^{^{HT}}(N).
$

$3^\circ$. This is just a combination of $1^\circ$ and $2^\circ$.

$4^\circ$. Note that all subfactors $M\subset M(\Cal H_k)$ appear
as irreducibe inclusions of factors in some $M \subset M_{2n}$.
By Jones' formula for the local indices ([J1]), if $p$
is a minimal projection in  $M'\cap M_{2n}$ with $(Mp \subset pM_{2n}p)
\simeq (M\subset M(\Cal H_k))$ then $[M(\Cal H_k):M]/\tau(p)^2
=[M_{2n}:M]$. On the other
hand, since $M_{2n}$ is the $[M:N]^n$-amplification of $M$ and
since $M(\Cal H_k)\simeq pM_{2n}p$,
it follows that $M(\Cal H_k)$ is the $\tau(p)[M:N]^n $-
amplification of $M$. By 8.2.2$^\circ$, this yields
$\beta_n^{^{HT}}(M(\Cal H_k))=[M(\Cal H_k):M]^{1/2}\beta_n^{^{HT}}(M)=v_k\beta_n^{^{HT}}(M)$.
\hfill Q.E.D.
\vskip .1in
Using the inventory of examples 6.9 of factors
in the class $\Cal H\Cal T$, and the calculations of
$\ell^2$-Betti numbers for groups in ([ChGr], [B]),
from 8.2.1$^\circ$ above we get the following list
of Betti numbers for factors:

\proclaim{8.5. Corollary} $1^\circ$. If
$\alpha \in \Bbb T$ is a primitive root of unity of order $n$,
then $M_\alpha=
L_\alpha(\Bbb Z^2) \rtimes SL(2, \Bbb Z)
\in \Cal H\Cal T_{_{s}}$ (cf. $6.9.1$) and
$\beta_1^{^{HT}}(M_\alpha) = (12n)^{-1}$, while
$\beta_k^{^{HT}}(M_\alpha)= 0, \forall k\neq 1$.

$2^\circ$. If $\alpha,
\alpha'$ are primitive roots of unity of order $n$
respectively $n'$ then  $M_\alpha \simeq M_{\alpha'}$
if and only if $n=n'$.
\endproclaim
\noindent
{\it Proof}. $1^\circ$. By 5.2.1$^\circ$,
8.2.1$^\circ$ and 8.2.2$^\circ$, $\beta_k^{^{HT}}(M_\alpha)
=\beta_k(SL(2, \Bbb Z))/n$. But by ([B]) we have 
$\beta_1(SL(2, \Bbb Z))=1/12$,
$\beta_k(SL(2, \Bbb Z))=0$ if $k\neq 1$.

$2^\circ$. By 5.2.1$^\circ$, if $n=n'$ then $M_\alpha \simeq
M_{\alpha'}$, while if $n\neq n'$ then
$\beta_1^{^{HT}}(M_\alpha) \neq \beta_1^{^{HT}}(M_{\alpha'})$,
so $M_\alpha \not\simeq M_{\alpha'}$.
\hfill Q.E.D.

\proclaim{8.6. Corollary} $1^\circ$. If $M=L^\infty(\Bbb S^2, \lambda)
\rtimes PSL(2, \Bbb Z²)$ as in $6.9.1'$ then
$\beta_1^{^{HT}}(M)=1/6$ and $\beta_n^{^{HT}}(M)=0, \forall n \neq
1$.

$2^\circ$. Let $\sigma$ be any of the actions $6.9.2$ or $6.9.6$
of the free group $\Bbb F_n$
on the diffuse abelian von Neumann
algebra $(A, \tau)$, and $M = A \rtimes_\sigma \Bbb F_n$
the corresponding factor in the class $\Cal H\Cal T$.
Then $\beta_1^{^{HT}}(M) = (n-1)$,
$\beta_k^{^{HT}}(M)= 0, \forall k\neq 1$.

$3^\circ$. Let $\Gamma_0$ be an
arithmetic lattice in $SU(n,1), n \geq 2,$ or
in $SO(2n,1)$,
$n \geq 1,$ and $\sigma$ a free ergodic
trace preseving action  of $\Gamma_0$ on the diffuse abelian
von Neumann algebra
$A$ as in $6.9.3$
or $6.9.6$. Let $M=A\rtimes_\sigma \Gamma_0\in \Cal H\Cal T$
be the corresponding $\Cal H\Cal T$ factor.
Then $\beta_n^{^{HT}}(M) \neq 0$ and
$\beta_k^{^{HT}}(M)= 0, \forall k\neq n$. Also, if
$\Gamma_0$ is an arithmetic lattice in some $SO(2n+1,1), n\geq 1$, then
the corresponding $\Cal H\Cal T$ factors constructed in $6.9.3$
satisfy $\beta_k^{^{HT}}(M)= 0, \forall k\geq 0$.

$4^\circ$. Let $\Gamma_0$ be an ${\text{\rm HT}}$ group (in the
sense of definition $6.11$; e.g., any
of the groups listed in $6.13$)
and $\Gamma_1$ an infinite amenable group. Let
$M \in \Cal H\Cal T$ be of the form $M = L^\infty(X, \mu) \rtimes
(\Gamma_0 \times \Gamma_1)$ (cf. $6.13.3^\circ$). Then
$\beta_k^{^{HT}}(M)= 0, \forall k\geq 0$.
\endproclaim
\noindent
{\it Proof}. For each
of the groups in 1$^\circ, 2^\circ$
the $\ell^2$-Betti numbers for certain
specific co-compact actions were calculated in ([B]). Then
the statements follow by ([G2], [ChGr]) and 8.2.1$^\circ$.
Similarly for $3^\circ$.
\hfill Q.E.D.

\proclaim{8.7. Corollary} If $\Gamma_0 = SL(2, \Bbb Z),
\Bbb F_n$ or if $\Gamma_0$ is an arithmetic lattice in $SU(n,1),
SO(n,1)$, for some $n \geq 2$, then there exist three
non-isomorphic factors
$M_i=L^{\infty}(X, \mu) \rtimes_{\sigma_i} \Gamma_0, 1\leq i\leq 3,$
in the class $\Cal H\Cal T$,
with $M_1 \in \Cal H\Cal T_{_{s}}$, $M_{2,3} \notin \Cal H\Cal T_{_{s}}$,
$M_{1,2}$ non-$\Gamma$ and $M_3$ with the property $\Gamma$.
\endproclaim
\noindent
{\it Proof}. All the groups mentioned have the property H (see 3.2).
The statement then follows from the last
part of 5.3.3$^\circ$.
\hfill Q.E.D.

\proclaim{8.8. Corollary} There exist both property $\Gamma$ and non-$\Gamma$
type ${\text{\rm II}}_1$ factors $M$ with trivial
fundamental group, $\mycal F(M)=\{1\}$. Moreover, such
factors $M$ can be taken to have non stably-isomorphic
tensor powers, all with trivial fundamental group.
\endproclaim

\vskip .1in
\noindent
{\it 8.9. Definition}. Let $M \in \Cal H\Cal T$. The
HT-{\it approximate dimension}
of $M$, denoted $ad_{_{HT}}(M)$, is by definition
Gaboriau's approximate dimension ([G2], [G3]) of the equivalence relation
$\Cal R^{^{HT}}_M$ associated with the HT Cartan subalgebra
of $M$. Note that $ad_{_{HT}}(M^t)=ad_{_{HT}}(M), \forall t>0$.

\proclaim{8.10. Corollary} Let $M \in \Cal H\Cal T$ be
of the form $M_k = L^{\infty}(X, \mu)\rtimes \Gamma_k$,
where $\Gamma_k=\Gamma_0 \times \Bbb F_{n_1}\times ... \times \Bbb F_{n_k}$,
for some $2 \leq n_i < \infty, \forall 1\leq i\leq k$,
with $\Gamma_0$ an increasing  union of finite groups.
Then $ad_{_{HT}}(M_k)=k$, so the factors $M_k, k\geq 1$, are
non stably-isomorphic.
\endproclaim
\noindent
{\it Proof}. By ([G3]), the approximate dimension
of the group $\Gamma_k$, and thus of $\Cal R^{^{HT}}_{M_k}$, is equal to $k$.
\hfill Q.E.D.
\vskip .1in
\noindent
{\it 8.11. Definition}. Let $M \in \Cal H\Cal T_{_{s}}$
and Out$_{_{HT}}(M)$ be the
countable discrete group
defined in Corollary 6.7.2$^\circ$. We call it
the HT-{\it outomorphism group} of $M$. As noted in 6.7,
Out$_{_{HT}}(M)$ can be identified with
the outer automorphism group of the equivalence relation
$\Cal R^{^{HT}}_M$, Out$(\Cal R^{^{HT}}_M)=$
Aut$(\Cal R^{^{HT}}_M)/$Int$(\Cal R^{^{HT}}_M)$.
Note that Out$_{_{HT}}(M^t)=$Out$_{_{HT}}(M), \forall t>0$.
The outer automorphism group of an equivalence relation $\Cal R$ has been
first considered by I.M. Singer in ([Si]), and
was also studied in ([FM]). By 6.7
this group is discrete (with the quotient topology) and countable.
Thus, it seems likely that Out$_{_{HT}}(M)$
can be computed in certain
specific examples. In this respect we mention the following:
\vskip .1in
\noindent
{\bf 8.12. Problem}. Calculate Out$_{_{HT}}(M)$ for
$M = L(\Bbb Z^2 \rtimes SL(2, \Bbb Z))$,
more generally for $M_n =
L((\Bbb Z^2)^n \rtimes SL(2, \Bbb Z))$,
with $SL(2, \Bbb Z)$ acting
diagonally on $(\Bbb Z^2)^n=\Bbb Z^2 \oplus ... \oplus \Bbb Z^2$. 
Let $\Cal G_n$
be the normalizer of $SL(2, \Bbb Z)$ in $GL(2n, \Bbb Z)$,
where $SL(2, \Bbb Z)$ is embedded in $GL(2n, \Bbb Z)$
block-diagonally. Is Out$_{_{HT}}(M_n)$ equal
to the quotient group $\Cal G_n/SL(2, \Bbb Z)$, in
particular is Out$_{_{HT}}(M_1)$ equal to $\{\theta_+, id\}$,
for $\theta_+$ the $C_+$ period
2 automorphism in Corollary 7.7 ?
\vskip .1in
\noindent
{\bf 8.13. Remarks}. $1^\circ$. Note that the above Corollary
8.8 (and also 8.5-8.7)
solves Problem 3 from Kadison's Baton Rouge list, providing
lots of examples of factors $M$ with the property
that the algebra of $n$ by $n$ matrices over $M$ is
not isomorphic to $M$, for any $n \geq 2$.

$2^\circ$.
We could extend the definition of
$\beta_n^{^{HT}}(M)$ to arbitrary II$_1$ factors
$M$, by simply letting $\beta_n^{^{HT}}(M)=0, \forall n$,
whenever $M$ does not belong to the class $\Cal H\Cal T$. This definition
would still be consistent with the property $\beta_n^{^{HT}}(M^t) =
\beta_n^{^{HT}}(M)/t, \forall t > 0$. However,
in order for this definition to also
satisfy the K{\"u}nneth formula (an
imperative!), one needs to solve the
following:
\vskip .05in
\noindent
{\bf 8.13.2. Problem}. Does
$M_1 \overline{\otimes} M_2 \in \Cal H\Cal T$ imply
$M_1, M_2 \in \Cal H\Cal T$ ?

\vskip .05in
Note that if this problem would have an
affirmative answer,
our factors $A \rtimes \Bbb F_n \in \Cal H\Cal T$
would follow prime,
i.e., $A \rtimes \Bbb F_n$
would not be expressable as a tensor product of type II$_1$
factors $M_1 \overline{\otimes} M_2$.
Indeed, this is because $\beta_1^{^{HT}}(M_1 \overline{\otimes}
M_2)= 0$ for $M_1, M_2 \in \Cal H \Cal T$,
by the K{\"u}nneth formula, while $\beta_1^{^{HT}}(A \rtimes \Bbb F_n)
= n-1 \neq 0$.

$3^\circ$. It would be interesting to extend the class of factors
in the ``good class'' for which a certain uniqueness result can
be proved for some special type of Cartan subalgebras,
beyond the HT factors we considered here. Such generalizations can
go two ways: by either extending the class of groups $\Gamma_0$
for which $A \subset A\rtimes_\sigma \Gamma_0$ works,
for certain $\sigma$, or by
showing that for the groups $\Gamma_0$ we already considered here
(e.g., the free groups) any action $\sigma$ works (see Problems
6.12.1$^\circ$ and respectively 6.12.2$^\circ$, in
this respect).

$4^\circ$. During a Conference at MSRI
in May 2001 ([C6]), Alain Connes posed
the problem of constructing $\ell^2$-type Betti number invariants
$\beta_k(M)$ for type
II$_1$ factors $M$, building on similar
conceptual grounds as in ([A], [C4], [ChGr], [G2,3]),
through appropriate definitions of simplicial complexes,
$\ell^2$-homology/cohomology for $M$, and which
should satisfy $\beta_k(L(G_0))=\beta_k(G_0)$
for von Neumann factors $M=L(G_0)$ associated
to discrete groups $G_0$. Thus,
since $\beta_k(\Bbb Z^2 \rtimes
SL(2, \Bbb Z))=0, \forall k$ (cf. [ChGr]), such Betti numbers would give
$\beta_k(L(\Bbb Z^2 \rtimes
SL(2, \Bbb Z)))=0, \forall k$.

Instead,
our approach to defining $\ell^2$-Betti numbers invariants was to
restrict our attention
to a class of factors $M$ having
a special type of Cartan subalgebras $A$,
the HT ones, for which we could prove
a uniqueness result, thus being able to use
the notion of Betti numbers for equivalence relations in ([G2]).
Thus, our Betti numbers are defined
``relative'' to HT Cartan subalgebras,
a fact we emphasized by using
the terminology ``$\ell^2_{_{HT}}$-Betti numbers''
and the notation ``$\beta^{^{HT}}_n(M)$''.
When $M=A \rtimes G_0$ these $\ell^2_{_{HT}}$-Betti numbers satisfy
$\beta_k^{^{HT}}(M)=\beta_k(G_0)$. In particular, if
$M=L(\Bbb Z^2\rtimes SL(2, \Bbb Z))$ then
$\beta_1^{^{HT}}(M)= \beta_1(SL(2, \Bbb Z))
\neq 0$. Thus $\beta_1^{^{HT}}(M)\neq \beta_1(M)$, if $\beta_k(M)$
could be
defined as asked in ([C6]).

Moreover, if such $\beta_k(M)$ are possible,
then according to Voiculescu's formula ([V1])
for the number of generators of the amplifications/compressions $M^t$
of the free group factors $M=L(\Bbb F_n)$
(cf. also [Ra], [Dy], [Sh]), the first Betti number $\beta_1 (M^t)$
(= (number of generators of $M^t$) $-1$)
should satisfy a formula of the type $\beta_1(M^t) = \beta_1(M)/t^2$,
rather than $\beta^{^{HT}}_1(M^t) = \beta^{^{HT}}_1(M)/t$, as we have
in this paper!

\heading Appendix: Some conjugacy results. \endheading

We prove here several conjugacy results for subalgebras in type
II$_1$ factors. The first one, Theorem A.1, plays a key role in
the proof of 6.2. The starting point in its proof is the following
simple observation: If $B_0, B$ are finite von Neumann algebras
for which there exists a $B_0-B$ Hilbert bimodule $\Cal H$ with
dim$\Cal H_B < \infty$ then a suitable reduced algebra of $B_0$ is
isomorphic to a subalgebra of some reduced of $B$. In the context
of C$^*$-algebras, this is reminiscent of the fact that
imprimitivity bimodules entail Morrita equivalence.  In the von
Neumann context, if both $B_0, B$ are subalgebras in some finite
factor $M$ then existence of  Hilbert $B_0-B$ bimodules $\Cal
H\subset L^2(M, \tau)$ with dim$\Cal H_B < \infty$ amounts to
existence of finite projections in $B_0'\cap \langle M, B\rangle$ 
($\langle M, B \rangle$ being the basic construction algebra) and the corresponding isomorphism of
$B_0$ into $B$ follows implemented by an element in $M$.

The basic construction  was first used in conjugacy problems by
Christensen ([Chr]), to
study ``small perturbations'' of subalgebras of type II$_1$ factors.
Although in A.1 we deal with conjugacy of subalgebras for
which no ``small distance'' assumption is made,
we still use the basic construction
as a set-up for the proof.
This framework allows us to use a trick
inspired from ([Chr]), and then
to utilise techniques
from ``subfactor theory'', notably  the pull down
identity ([PiPo], [Po2,3]).
We also use von Neumann algebra analysis of projections,
with repeated use of results from ([K2]). For notations and
elementary properties of the basic construction, see Section 1.3
and ([J1], [PiPo], [Po2,3]).

To state A.1, let $M$ be a finite factor, $B \subset M$
a von Neumann subalgebra and $\Cal U_0 \subset M$ be a subgroup of
unitary elements. Let $B_0 = \Cal U_0''$ be the
von Neumann algebra it generates in $M$. For each $b \in 
\langle M, B \rangle$, $Tr(b^*b) < \infty$, 
we denote by
$K_{\Cal U_0}(b)$ the weak
closure of the convex hull
of $\{u_0 b u_0^* \mid u_0 \in \Cal U_0 \}$,
i.e., $K_{\Cal U_0}(b) = {\overline{\text{\rm co}}}^{\text{\rm w}}
\{u_0 b u_0^* \mid u_0 \in \Cal U_0 \}$.
Note that
$K_{\Cal U_0}(b)$ is also contained in the Hilbert space
$L^2(\langle M, B \rangle, Tr)$, where it is still
weakly closed.

Let $h=h_{\Cal U_0}(b)
\in K_{\Cal U_0}(b)$ be the unique element of minimal
norm $\| \, \|_{2,Tr}$ in $K_{\Cal U_0}(b)$.
Since $u K_{\Cal U_0}(b) u^* = K_{\Cal U_0}(b)$ and
$\|uhu^*\|_{2, Tr} = \|h\|_{2, Tr}, \forall u\in \Cal U_0$,
by the uniqueness of $h$
it follows that $uhu^*=h, \forall u\in \Cal U_0$.
Thus $h\in \Cal U_0' \cap \langle M, B \rangle =B_0'\cap
\langle M, B \rangle$. Moreover, by the definitions, we see that
if $0\leq b \leq 1$ then $0\leq k \leq 1$ and
Tr$(k) \leq {\text{\rm Tr}}(b)$, for all
$k\in K_{\Cal U_0}(b)$.

\proclaim{A.1. Theorem} Let $M, B, B_0, \Cal U_0$
be as above. Assume the von Neumann subalgebra
$B\subset M$ is maximal abelian
in $M$ and $B_0$ is abelian with
$B_{01}\overset \text{\rm def}
\to = B_0'\cap M$ still abelian (thus maximal abelian
in $M$). Then
the following conditions are equivalent:

$1^\circ$. The element $h_{\Cal U_0}(e_B)$ is non-zero.

$2^\circ$. There exists a non-zero projection $e_0\in B_0'\cap
\langle M, B \rangle$ with $Tr(e_0) < \infty$.

$3^\circ$. There exist non-zero projections $q_0\in
B_0 '\cap M$, $q\in B$ and a
partial isometry $v\in M$ such that $v^*v =q_0,
vv^* = q$ and $vB_0v^* \subset Bq$.
\endproclaim
\vskip .1in
\noindent
{\it Proof}. $3^\circ \implies 1^\circ$. If $v$ satisfies
condition 3$^\circ$ then $B_0q_0$
is contained in $v^*Bv$. Since $e_B$ commutes
with $B$, it follows that $e_0 = v^*e_Bv$ commutes
with $B_0$, i.e.,
$e_0\in B_0'\cap \langle M, B \rangle$. Also,
Tr$e_0={\text{\rm Tr}}(v^*e_Bv)
\leq {\text{\rm Tr}}(e_B)=1$.

$1^\circ \implies 2^\circ$. Let $e_0$ be the
spectral projection of $h = h_{\Cal U_0}(e_B)$
corresponding to the interval $(\|h\|/2, \infty)$. Then $e_0 \neq 0$
and $h \geq {1\over 2} e_0$. Thus,
$$
Tr(e_0) \leq 2 Tr(h) \leq 2 Tr(e_B)
<\infty.
$$

Thus, $e_0$ is a finite projection in $\langle M, B \rangle$ and
$e_0$ commutes with $B_0$ (since $h$ does).

2$^\circ \implies 3^\circ$. Denote $M_1 = \langle M, B\rangle$.
Since $B_0e_0$ is abelian, it is contained in a maximal abelian
subalgebra $B_1$ of $e_0M_1e_0$. Since $M_1 = (JBJ)'\cap \Cal B(L^2M)$,
it is a type I von Neumann algebra. Thus,
by a result of Kadison ([K2]),
$B_1$ contains a non-zero abelian projection $e_1$ of
$M_1$ (i.e., $e_1M_1e_1$ is abelian). Since $e_B$ is
a maximal abelian projection in $M_1$ and
has central support 1 in $M_1$, it
follows that $e_B$ majorizes $e_1$. Thus, $e_1$
satisfies $e_1(L^2(M,\tau)) = \overline{\xi B}$
for some $\xi \in L^2(M, \tau)$.

Let $V\in M_1$ be a partial isometry such that
$V^*V = e_1 \leq e_0$ and
$VV^* \leq e_B$. It follows that $VB_1e_1V^*$ is
a subalgebra of $e_BM_1e_B =
Be_B$. Since $e_1$ commutes with $B_0$, it follows that if we denote
by $f'$ the maximal projection in $B_0$ such that $f'e_1=0$
and let $f_0=1-f'$, then there exists a unique isomorphism
$\alpha$ from $B_0f_0$ into $B$ such that $\alpha(b)e_B
= VbV^*, \forall b\in B_0f_0$. Let $f = \alpha(f_0)
\in B$.

It follows that
$\alpha(b)e_BV
= e_BVb, \forall b\in B_0f_0$. By applying $\Phi$ to both sides and
denoting $a$ the square
integrable operator $a=\Phi(e_BV)\in L^2(M, \tau)$, it follows that
$\alpha(b)a=ab, \forall b\in B_0$. Since $e_Ba = e_BV = V$,
it follows that $a\neq 0$.

By the usual trick, if we denote by $v_0\in M$ the unique
partial isometry in the polar decomposition of $a$ such
that the right supports of $a$ and $v_0$ coincide, then
$p_0=v_0^*v_0$ belongs
to the algebra $B_0'\cap M = B_{01}$, which is abelian
by hypothesis, $p
=v_0v_0^*$ belongs to $(\alpha(B_0)f)'\cap
fMf$ and
$\alpha(b)v_0=v_0b, \forall b\in B_0f_0$.

But $B_{01}=B_0'\cap M$ maximal abelian in $M$ implies
$B_{01}f_0$
maximal abelian in $f_0Mf_0$. Moreover, since
$v_0B_0v_0^* = \alpha(B_0)p$, if we denote
$B_{11}=v_0B_{01}v_0^*$, then by spatiality we have
$$
B_{11}=v_0B_{01}v_0^*=v_0(B_0'\cap M)v_0^*
$$
$$
={v_0B_0v_0^*}'\cap pMp = (\alpha(B_0)p)'\cap pMp =
p((\alpha(B_0)f)'\cap
fMf)p.
$$
Tis implies that $p$ is an abelian projection
in $(\alpha(B_0)f)'\cap fMf$. Thus, if $z$ is
the central projection of $p$ in
$(\alpha(B_0)f)'\cap fMf$ then $((\alpha(B_0)f)'\cap fMf)z=
((\alpha(B_0)z)'\cap zMz$
is finite of type I.

Since $Bf$ is maximal abelian in $fMf$ it follows that $z \in Bf$ and
$Bz$ is maximal abelian in the type
I$_{fin}$ algebra $((\alpha(B_0)z)'\cap zMz$.
By ([K2]) it follows that there
exists a projection
$f_{11} \in Bz$ such that $f_{11}$ is equivalent to $p$
in $((\alpha(B_0)z)'\cap zMz$. Let $v_1\in
((\alpha(B_0)z)'\cap zMz$ be such that
$v_1v_1^*=f_{11}, v_1^*v_1 = p$ and denote $v=v_1v_0\in M$.
Then we have
$v^*v=p_0\in B_0', vv^* = f_{11} \in B$
and $vB_0v^* = \alpha(B_0)f_{11} \subset Bf_{11}$.
\hfill Q.E.D.

Our second conjugacy result, A.2, is a ``small perturbation''-type
result, needed in the proofs of 4.5 and 6.6.3$^\circ$. The
starting point in its proof is a trick from ([Chr]). Then, like in
A.1, we use techniques from ([Po2,3,7], [PiPo]). Note that the
proof of Step 1 below is a refinement of the proof of (4.4.2 in
[Po1]), while the proof of Step 2 is a refinement of an argument
used in proving (4.5.1, 4.5.6 and 4.7.3 in [Po1]).

\proclaim{A.2. Theorem} For any
$\varepsilon_0 > 0$ there exists $\delta > 0$
such that if $M$ is a type ${\text{\rm
II}}_1$ factor, $B\subset M$ is a subfactor with
$B'\cap M=\Bbb C$, $B_0\subset M$
is a von Neumann subalgebra with
$B_0'\cap M = \Cal Z(B_0)$, $\Cal N_M(B_0)''=M$
and $B_0 \subset_\delta B$ then there exists a
unitary element
$u\in M$ such that
$\|u-1\|_2 \leq \varepsilon_0$ and $uB_0u^* \subset B$.
\endproclaim
\vskip .1in
\noindent
{\it Proof}. Step 1. Let $\varepsilon = \varepsilon_0^2/4$.
We first prove that
$\exists \delta > 0$ such that if $B_0, B\subset M$
satisfy $B_0'\cap M=\Cal Z(B_0)$ and
$B_0 \subset_\delta B$ then
$\exists p_0\in \Cal P(B_0), p\in \Cal P(B)$,
a unital isomorphism $\theta$ of
$p_0B_0p_0$ into $pBp$, a
projection $q\in \theta(p_0B_0p_0)'\cap pMp$ and a partial isometry
$v\in M$ such that $v^*v = p_0, vv^* = q\leq p$, $\|v-1\|_2
\leq \varepsilon$,
$\tau(q) \geq 1-\varepsilon$
and $vb_0 = \theta(b_0)v, \forall b_0\in p_0B_0p_0$.

To do this
note first that if $u_0 \in \Cal U(B_0)$ then
$\|u_0e_Bu_0^* - e_B\|^2_{2,Tr}/2 = 1-Tr(e_Bu_0e_Be_0^*)$
$=\|u_0 - E_B(u_0)\|^2_2$ (see e.g., line 17 on page 322 in [Po9]).
So if $\|u_0 - E_B(u_0)\|_2 \leq \delta$,
$\forall u_0 \in \Cal U_0 = \Cal U(B_0)$, then
with the notations in A.1 we get
$h=h_{\Cal U_0}(e_B) \in B_0'\cap \langle M, B \rangle$, with
$h \leq 1$, $Tr(h) \leq 1$
and $\|h - e_B\|_{2, Tr} \leq 2^{1/2}\delta$.
Thus, by (1.1 in [C2])
there exists $s > 0$ such that the spectral projection $e$ of
$h$ corresponding to the interval $[s, \infty)$
satisfies $\|e-e_B\|_{2,Tr} \leq (2\delta)^{1/2}$. Note that
$e \in B_0'\cap \langle M, B \rangle$ as well. We next want to show that
by slightly shrinking $e$ we may assume in addition $(B_0e)'\cap
e\langle M, B \rangle e=\Cal Z(B_0)e$.

So let $u \in \Cal U(C)$, where $C=(B_0e)'\cap e\langle M, B \rangle e$.
Since $e_B\langle M, B \rangle e_B = Be_B$ and $e$ is $(2\delta)^{1/2}$-close
to $e_B$ in the norm $\|\quad \|_{2, Tr}$, if we
denote by $b$ the unique element
in $B$ with $be_B=e_Bue_B,$ then $u$ is close to $ebe$ in the
norm $\|\quad \|_{2,tr}$ implemented by
the normalized trace $tr = Tr(e)^{-1} Tr$ on $e\langle M, B \rangle e$.
This implies that $\|[ebe, v]\|_{2,tr}  \leq
\varepsilon(\delta)$, $\forall v\in \Cal U(B_0e)$,
in which $\varepsilon(\delta)$ denotes
from now on a constant depending on
$\delta$,
with $\underset \delta \rightarrow 0 \to \lim
\varepsilon(\delta) = 0$ (but
$\varepsilon(\delta)$ possibly changing in each
of the subsequent estimates).
Since $B_0'\cap M = \Cal Z(B_0)$, if we average $ebe$ by
unitaries in $B_0e$, we see that
$u$ is $\varepsilon(\delta)$-close to an element in $\Cal Z(B_0)e$. Thus
$C \subset_{\varepsilon(\delta)} A_0$, where
$A_0 = \Cal Z(B_0)e$. Noticing
that $A_0 \subset \Cal Z(C)$, we infer that
this implies $\exists e'\in \Cal Z(C)$,
with $tr(e') \geq 1-\varepsilon(\delta)$ and $Ce' = A_0e'$,
i.e., $(B_0e)'\cap
e\langle M, B \rangle e=\Cal Z(B_0)e$. Indeed, for
if $q'\in \Cal Z(C)$ is the maximal projection with $Cq'$ abelian
and $A\subset C$ is a maximal abelian $^*$-subalgebra with $A_0
\subset A$ then $q'\in A$ and there exists $u\in \Cal U(B(1-q'))$
with $E_A(u)=0$. Since $q'+u\in \Cal U(C)$ we have:
$$
tr(1-q')=\|u\|_{2,tr}^2=\|(q'+u)-E_A(q'+u)\|^2_{2,tr}
$$
$$
\leq \|(q'+u)-E_{A_0}(q'+u)\|^2_{2,tr} \leq \varepsilon(\delta)^2.
$$
This reduces the problem to the case $C$ is abelian, which is an easy
exercise (e.g., use the argument on page 745 in [Po7]).

Taking $e'$ for $e$ in the above, this shows that
if $B_0 \subset_\delta B$ then $\exists e\in B'_0 \cap \langle M,B \rangle$
finite projection with $\|e-e_B\|_{2,Tr} \leq \varepsilon(\delta)$ and
$(B_0e)'\cap
e\langle M, B \rangle e =\Cal Z(B_0)e$. But by ([Po6]) the
latter condition implies there exists $A_1\subset B_0$ abelian such that
$A_1e$ is maximal abelian in $e\langle M, B \rangle e$. By ([K2]) there
exists a projection $P\in A_1=A_1e$ such that $P$ is equivalent to
the support projection of $ee_Be \in e\langle M, B \rangle e$.
In particular, $P$ is majorized by $e_B$. Also,
$P$, $e$ and $e_B$ are $\varepsilon(\delta)$-close one to another. By (1.2 in [C2]),
there exists a partial isometry $V \in \langle M,B \rangle$
such that $V$ is $\varepsilon(\delta)$-close to $e_B$, $V^*V=P \in A_1\subset B_0'$
and $VV^* \leq e_B$.
Like in ([Chr]) and in
the proof of A.1, if $p_0\in B_0$ and $p\in B$ denote the
support projections of $V^*V$ in $B_0$ and respectively
$VV^*$ in $B$ then
there exists
a unital isomorphism $\theta$ of $p_0B_0p_0$ into $pBp$ such that
$Vb_0 = \theta(b_0)V, \forall b_0\in p_0B_0p_0$. If we now take
the partial isomtery $v=\Phi(V)|\Phi(V)|^{-1}\in M$,
then we still have $vb_0 = \theta(b_0)v, \forall b_0\in p_0B_0p_0$
and $v$ is $\varepsilon(\delta)$-close to $1$ (using $\|\Phi(V)-1\|_1
\leq \|V-e_B\|_{1,Tr}$ and applying 2.1 in [C2]). Since $v^*v
\in (p_0B_0p_0)'\cap p_0Mp_0=\Cal Z(B_0)p_0$ and $vv^*\in
\theta(p_0B_0p_0)'\cap pMp$, letting $q=vv^*$, we are done.

Step 2. If $p_0, p, q, v, \theta$ are as in Step 1,
then $vB_0v^* = \theta(p_0B_0p_0)q$, so by
spatiality we have:
$$
q(\theta(p_0B_0p_0)'\cap pMp)q= (vB_0v^*)'\cap qMq
$$
$$
=v(p_0B_0p_0' \cap p_0Mp_0)v^* = v\Cal Z(B_0)v^* = \Cal
Z(\theta(p_0B_0p_0))q.
$$
In particular, $q(\theta(p_0B_0p_0)'\cap pBp)q
=\Cal
Z(\theta(p_0B_0p_0))q$.
Since $Z(\theta(p_0B_0p_0)) \subset \theta(p_0B_0p_0)'\cap pBp$
this implies that there exists a normal
conditional expectation $E$ of $\theta(p_0B_0p_0)'\cap pBp$
onto $Z(\theta(p_0B_0p_0))$ such that
$qxq=E(x)q, \forall x\in \theta(p_0B_0p_0)'\cap pBp$.

Let $p'\in \theta(p_0B_0p_0)'\cap pBp$ be the minimal projection
such that $qp'=q$. By replacing
if necessary $\theta$ by $\theta(\cdot)q'$,
(while leaving $v$ unchanged) we may thus assume
$p'=p$. Thus, if $a
\in \theta(p_0B_0p_0)'\cap pBp$ satisfies $aq=0$ then
the support of $a^*a$ follows majorized
by $p-p'=0$, implying that $a=0$ and showing that
$E$ is faithful. Since
$q$ implements the normal faithful conditional expectation $E$
of $\theta(p_0B_0p_0)'\cap pBp$ onto $Z(\theta(p_0B_0p_0))$,
it follows that the weak closure of sp$\{x q y\mid x, y
\in  \theta(p_0B_0p_0)'\cap pBp\}$ is a finite von Neumann
subalgebra $Q$ of $pMp$ with $qQq \simeq Z(\theta(p_0B_0p_0))$.
Since $q$ has support $1$ in $Q$, this shows that $Q$ is type I$_{fin}$.
But $Q$ contains $(\theta(p_0B_0p_0)'\cap pBp)1_Q$,
which is isomorphic to $\theta(p_0B_0p_0)'\cap pBp$.
Thus, the latter follows type I$_{fin}$ as well.

Let $q'\in Z(\theta(p_0B_0p_0)) (\subset
\Cal Z(\theta(p_0B_0p_0)'\cap pBp))$
be the maximal projection with
$$
q'Z(\theta(p_0B_0p_0))=q'(\theta(p_0B_0p_0)'\cap pBp).
$$
It follows that
there exists $b\in L^2(\theta(p_0B_0p_0)'\cap pBp)(p-q')$ with
$E(b)=0$ and $E(b^*b)=p-q'$ (see e.g.,
[Po2]). This shows that $bqb^*$ is a projection
orthogonal to $q(p-q')$ and equivalent to $q(p-q')$,
while still under $p-q'$. Thus
$$
\tau(q(p-q')) = \tau(bq(p-q')b^*)
$$
$$
\leq \tau((1-q)(p-q')) \leq \tau(1-q) \leq \varepsilon.
$$
Thus, $1-\varepsilon - \tau(q') \leq \tau(p-q') \leq 2\varepsilon$,
implying that $\tau(q') \geq 1-3\varepsilon$. This shows
that by
``cutting everything''  by $q'$ we may assume
$\theta(p_0B_0p_0)'\cap pBp = \Cal Z(\theta(p_0B_0p_0))$.

Since $B_0$ is regular in $M$, $p_0B_0p_0$ is regular in $p_0Mp_0$
(see e.g. [JPo])
and thus, by spatiality,  $\theta(p_0B_0p_0)q$ is regular
in $qMq$. Since $\theta(p_0B_0p_0)
\ni b \rightarrow bq \in \theta(p_0B_0p_0)q$
is an isomorphism, for each $u \in \Cal N_{qMq}(\theta(p_0B_0p_0)q)$
there exists an automorphism $\sigma_u$ of $\theta(p_0B_0p_0)$
such that $ubqu^*= \sigma_u(b)q, \forall b\in \theta(p_0B_0p_0)$.
Thus, $ub= \sigma_u(b)u, \forall b\in \theta(p_0B_0p_0)$.

By applying $E_B$ to both sides of this equality,
it follows that $E_B(u)b = \sigma_u(b)E_B(u),
\forall b\in \theta(p_0B_0p_0)$. By also taking into
account that $\theta(p_0B_0p_0)'\cap pBp \subset \theta(p_0B_0p_0)$,
this shows that if $B_1 \subset pBp$ denotes the von Neumann
algebra generated by the normalizer of $\theta(p_0B_0p_0)$
in $pBp$ then  $E_B(\Cal N_{qMq}(\theta(p_0B_0p_0)q)\subset B_1$.
By the regularity of $\theta(p_0B_0p_0)q$ in $qMq$,
this entails $E_B(qMq) \subset B_1$ as well. Since $q\leq p$ and
$\tau(q) \geq 1-\varepsilon$, we thus have
$pBp \subset_\varepsilon B_1 \subset pBp$. Taking into
account that $pBp$ is a factor, this implies there exists
a projection $p''\in \Cal Z(B_1)$ with $\tau(p'') \geq 1-2\varepsilon$
such that $B_1p'' = p''Bp''$.

By cutting with $p''$ we may thus also assume $\theta(p_0B_0p_0)$
is regular in $pBp$. Since $pBp'\cap pMp = \Bbb Cp$,
this implies $\Cal N_1 =
\Cal N_{pBp}(\theta(p_0B_0p_0))$ satisfies $\Cal N_1'\cap pMp = \Bbb Cp$.
Since $\Cal N_1$ also normalizes the algebras $\Cal Z(\theta(p_0B_0p_0))
= \theta(p_0B_0p_0)'\cap pBp$ and $\theta(p_0B_0p_0)'\cap pMp$,
it follows that it acts ergodically on both. By
ergodicity, $\theta(p_0B_0p_0)'\cap pMp$ follows either
homogeneous of type I$_{fin}$ or of type II$_1$. Since $q(\theta(p_0B_0p_0)'
\cap pMp)q=\Cal Z(\theta(p_0B_0p_0))q$ is abelian and $\tau(q) > 1/2$
(for $\varepsilon$ chosen sufficiently small),
$\theta(p_0B_0p_0)'\cap pMp$ follows abelian.

Denote $A_0=\Cal Z(\theta(p_0B_0p_0))$, $A_1=\theta(p_0B_0p_0)'\cap pMp$,
$N_0=pBp$ and $Q_0$ the factor generated by $\Cal N_1$
and $A_1$ in $pMp$. Thus, we have $N_0'\cap Q_0 = \Bbb C$ and
the non-degenerate commuting square
$$
\matrix
N_0 &\subset& Q_0&\\
\cup && \cup && \\
A_0 & \subset& A_1&
\endmatrix
$$
(Recall that we also have $q\in A_1, A_1q=A_0q$ and
$\tau(q) \geq 1-\varepsilon$.)

Thus, if $e=e^{Q_0}_{N_0}$ denotes the Jones
projection corresponding to the inclusion $N_0\subset Q_0$
then $A_0\subset A_1 \subset \langle A_1, e \rangle$ is
the basic construction
for $A_0\subset A_1$. Since $\Cal Z(\langle A_1, e \rangle)=A_0$
and since $\Cal N_1$ acts on $A_0 \subset A_1$
with the action on $A_0$ being ergodic, it
follows that $\langle A_1, e \rangle$ is
homogeneous of type I. But
$q(A_1eA_1)q= A_0(qeq)A_0$, and since $[A_0, qeq]=0$ this implies
$q \langle A_1, e \rangle q = A_0qeq$.
Thus, $q\langle A_1, e \rangle q$ is abelian. Equivalently,
$q$ is an abelian projection in $\langle A_1, e \rangle$.
But then $q$ is majorised by $e$ in $\langle A_1, e \rangle$.
Thus $q$ is majorised by $e$ in $\langle Q_0, e \rangle$ as well,
showing that $q$ is finite in $\langle Q_0, e \rangle$.

But $q$ enters finitely many times in $1_Q$, in
the factor $Q_0$, which is a subalgebra of $
\langle Q_0, e \rangle$. Thus
$\langle Q_0, e \rangle$ is a finite factor and $\tau(e) \geq \tau(q)
\geq 1-\varepsilon > 1/2$.
By Jones Theorem, $e=1$ and $N_0=Q_0$. In particular,
$q \in \theta(p_0B_0p_0)$, so $q=p$. Thus, $v^*v = p_0 \in B_0,
vv^* = p \in B$ and $v(p_0B_0p_0)v^* \subset pBp$. Since
the normaliser of $B_0$ acts ergodically
on the center of $B_0$ and $B$ is a factor, there exists a unitary
element $u \in M$ such that $up_0 = v$ and $uB_0u^* \subset B$.
But then $\|1-u\|_2 \leq \|1-v\|_2+\|v-u\|_2 \leq 2 \varepsilon^{1/2}
=\varepsilon_0$.
\hfill Q.E.D.

Our last conjugacy result, somewhat technical, is needed in
the proof of 4.3.2$^\circ$.

\proclaim{A.3. Theorem} Let $M$ be a type ${\text{\rm II}}_1$
factor and $P, Q \subset M$ von Neumann subalgebras. Assume there
exists a group of unitary elements $\Cal U_0 \subset P$ that
normalizes $Q$ and satisfies $N_0'\cap M = Q'\cap \Cal Z(N_0)$,
where $N_0=\Cal U_0''$. If $Q\subset_{\varepsilon_0} P$, for some
$\varepsilon_0 < 1/2$, then there exists a non-zero projection $p
\in Q'\cap \Cal Z(N_0)$ such that $Qp \subset P$.
\endproclaim
\vskip .1in \noindent {\it Proof}. Let $M \subset^{e_P} \langle M,
e_P \rangle$ be the basic construction for $P\subset M$, with $Tr$
and $\Phi$ the canonical trace and weight, respectively, as in
1.3.1. The statement is equivalent to proving that there exists $p
\in Q'\cap \Cal Z(N_0)$, $p\neq 0$, such that $[Qp,e_P]=0$.

Let $k$ be the unique element of minimal norm $\| \quad \|_{2,Tr}$
in $K=\overline{\text{\rm co}}^w \{ue_Pu^* \mid u \in \Cal
U(Q)\}$. Note that $0 \leq k \leq 1, Tr(k) \leq 1$. Also, since
for $u \in \Cal U(Q)$ we have
$$
\|e_P-ue_Pu^*\|_{2,Tr}^2= 2 - 2 \|E_P(u)\|_2^2 = 2\|u-E_P(u)\|_2^2
\leq 2 \varepsilon_0^2,
$$
by taking convex combinations and weak limits it follows that $\|k
- e_P\|^2_{2,Tr} \leq 2\varepsilon_0^2 < 1/2$.

Since $u K u^* = K$ and $\|uku^*\|_{2, Tr} = \|k\|_{2, Tr},
\forall u\in \Cal U(Q)$, by the uniqueness of $k$ as the element
of minimal norm $\|\quad\|_{2,Tr}$ in $K$, it follows that
$uku^*=k, \forall u\in \Cal U(Q)$. Thus $[k, Q]=0$. Moreover, if
$v \in \Cal U_0 \subset P$ then $[v, e_P]=0$ and $vQ v^* = Q$,
implying that $v(ue_Pu^*)v^* = (vuv^*)e_P(vu^*v^*) \subset K$,
$\forall u \in \Cal U(Q)$. Thus, $vKv^*=K$ and so, by the
uniqueness of $k$, $[k, v]=0$. Since $\Cal U_0$ generates $N_0$, it
follows that $k$ and all its spectral projections commute with
both $Q$ and $N_0=\Cal U_0''$.

Together with $[e_P, N_0]=0$ this yields $[ke_P, N_0]=0$ and
further on, by applying the operator valued weight $\Phi$ of
$\langle M, e_P \rangle$ on $M$ (which is $M$-bimodular, thus
$N_0$-bimodular as well) and letting $a=\Phi(ke_P)$, gives $[a,
N_0]=0$. Equivalently, $a \in N_0'\cap M = Q'\cap \Cal Z(N_0)$.
Since $\Cal Z(N_0)\subset N_0 \subset P$, $a\in P$ and so
$[a,e_P]=0$. Together with $ae_P=ke_P$, this entails
$ae_P=e_Pae_P=e_Pke_P \geq 0$, and so $a \geq 0$. In particular,
$a=a^*$. Thus, $ke_P=ae_P=(ae_P)^*=(ke_P)^*=e_Pk$, showing that
$[k, e_P]=0$.

Let now $e_1$ be the spectral projection of $k$ corresponding to
the set $\{1\}$. Thus $e_1=e_1k \in \overline{\text{\rm co}}^w
\{u(e_1e_P)u^* \mid u \in \Cal U_0\}$, showing that $e_1 \leq
e_P$. Thus, if $p=\Phi(e_1)$ then $p$ is a projection in $P$ with
$e_1 = pe_P$, $[p, Q\vee N_0]=0$ and $[e_P, Qp]=0$. Thus, we are
done, provided we can show that $p\neq 0$.

Assume by contradiction $e_1 = 0$. We show that this implies that
for any spectral projection $e$ of $k$, $ee_P$ is majorized by
$e(1-e_P)$ in $\langle M, e_P\rangle$. Indeed, for if this is not
the case then there exists a projection $z$ in $\Cal Z(\langle M,
e_P \rangle)$ and a partial isometry $w\in \langle M, e_P \rangle$
such that $w^*w \lneqq zee_P$, $ww^* = ze(1-e_P)$. If we denote
$b=\Phi(w)$, then $be_P=w$ and so
$$
bb^*=\Phi(ww^*)=\Phi(ze(1-e_P)) \in N_0'\cap M = Q'\cap \Cal
Z(N_0).
$$
Similarly, $q= \Phi(eze_P)$ is a projection in $P$ which commutes
with $N_0$, thus lying in $Q'\cap \Cal Z(N_0)\subset P$. Noticing
that $bb^* \geq be_Pb^* = z e(1-e_P)$ and that the morphism
$Q'\cap \Cal Z(N_0) \ni x \mapsto xze(1-e_P)$ has support $q$
(because $e_1=0$), it follows that $bb^* \geq q$. Thus
$$
\tau(q)=Tr(zee_P) \gneqq Tr(w^*w) = Tr(ww^*)=\tau(bb^*) \geq
\tau(q),
$$
a contradiction.

In particular, since $ee_P \prec e(1-e_P)$ for any spectral
projection $e$ of $k$, we have $\|k(1-e_P)\|_{2,Tr} \geq
\|ke_P\|_{2,Tr}$. By Pythagora, this gives
$$
\tau((1-k)^2) + \tau(k^2) \leq \|ke_P-e_P\|_{2,Tr}^2 +
\|k(1-e_P)\|_{2,Tr}^2  = \|k-e_P\|_{2,Tr}^2 < 1/2
$$
Thus $0 > \tau(2(1-k)^2 + 2k^2 - 1)=\tau(1-4k +
4k^2)=\tau((1-2k)^2)$. This final contradiction ends the proof of
the Theorem. \hfill Q.E.D.

\head References\endhead

\item{[A-De]} C. Anantharam-Delaroche: {\it On
Connes' property} (T) {\it for von Neumann algebras}, Math. Japon.
{\bf 32} (1987), 337-355.

\item{[At]} M. Atiyah: {\it Elliptic operators, discrete groups and
von Neumann algebras}, In Colloque ``Analyse et Topologie'',
in honor of Henri Cartan, Ast\'erisque, {\bf 32} (1976), 43-72.

\item{[Bi]} D. Bisch: {\it A note on intermediate subfactors},
Pac. J. Math., {\bf 163} (1994), 201-215.

\item{[Bo]} F. Boca: {\it On the method for constructng irreducible
finite index subfactors of Popa},
Pac. J. Math., {\bf 161} (1993), 201-231.

\item{[B]} A. Borel: {\it The $L^2$-cohomology
of negatively curved Riemannian symmetric spaces},
Acad. Sci. Fenn. Ser. A, Math., {\bf 10} (1985), 95-105.

\item{[Bu]} M. Burger: {\it Kazhdan constants for}
SL$(3,\Bbb Z)$, J. Reine Angew. Math. {\bf 413} (1991), 36-67.

\item{[dCaH]} J. de Canni\'ere, U. Haagerup: {\it
Multipliers of the Fourier algebra of some simple Lie groups
and their discrete subgroups}, Amer. J. Math. {\bf 107} (1984), 455-500.

\item{[ChGr]} J. Cheeger, M. Gromov: {\it $L^2$-cohomology and group
cohomology}, Topology {\bf 25} (1986), 189-215.

\item{[CCJJV]} P.-A. Cherix, M. Cowling,
P. Jolissaint, P. Julg, A. Valette: ``Groups with the Haagerup property
(Gromov's a-T-menability)'', book.

\item{[Cho]} M. Choda: {\it Group factors of the
Haagerup type}, Proc. Japan Acad.,
{\bf 59} (1983), 174-177.

\item{[Chr]}
E. Christensen: {\it Subalgebras of a finite algebra},
Math. Ann. {\bf 243} (1979), 17-29.

\item{[C1]} A. Connes: {\it A type II$_1$
factor with countable fundamental group}, J. Operator
Theory {\bf 4} (1980), 151-153.

\item{[C2]} A. Connes: {\it Classification of injective factors},
Ann. of Math.,
{\bf 104} (1976), 73-115.

\item{[C3]} A. Connes: {\it Classification des facteurs},
Proc. Symp. Pure Math.
{\bf 38}
(Amer. Math. Soc., 1982), 43-109.

\item{[C4]} A. Connes: {\it Sur la th\'eory non-commutative de
l'int\'egration}, in ``Alg\`ebres d'op\'erateurs, S\'eminaire
Les Plans-sur-Bex, 1978'', pp. 19-143, Lecture Notes 725,
Springer, Berlin, 1979.

\item{[C5]} A. Connes: {\it Sur la classification des facteurs de type} II,  C. R. Acad. Sci. Paris S\'er. I Math.,
{\bf 281} (1975), A13-A15.

\item{[C6]} A. Connes: {\it Factors and geometry}, lecture at MSRI,
May 1'st, 2001.

\item{[C7]} A. Connes: {\it Correspondences}, hand-written
notes, 1980.

\item{[CJ]} A. Connes, V.F.R. Jones: {\it Property} (T)
{\it for von Neumann algebras}, Bull. London Math. Soc. {\bf 17} (1985),
57-62.

\item{[CW]} A. Connes, B. Weiss: {\it Property $\Gamma$ and
asymptotically
invariant sequences}, Israel. J. Math. {\bf 37} (1980), 209-210.

\item{[CowH]} M. Cowling, U. Haagerup: {\it Completely bounded
multipliers and the Fourier algebra of a simple Lie group of
real rank one}, Invent. Math. {\bf 96} (1989), 507-549.

\item{[DeKi]} C. Delaroche, Kirilov: {\it Sur les relations entre
l'espace dual d'un groupe et la structure de ses sous-groupes fermes},
Se. Bourbaki, 20'e ann\'ee, 1967-1968, no. 343, juin 1968.

\item{[D]}
J. Dixmier: {\it Sous-anneaux ab\'{e}liens maximaux
dans les facteurs de type fini}, Ann. Math. {\bf 59} (1954), 279-286.

\item{[Dy]} H. Dye: {\it On groups of measure preserving transformations},
I, II, Amer. J. Math. {\bf 81} (1959), 119-159, and {\bf 85} (1963),
551-576.

\item{[Dyk]} K. Dykema: {\it
Interpolated free group factors}, Duke Math J. {\bf 69} (1993), 97-119.

\item{[Ey]} P. Eymard: ``Moyennes invariantes et repr\'esentations
unitaires'', Lecture Notes in Math, {\bf 300}, Springer-Verlag,
Berlin, 1972.

\item{[FM]} J. Feldman, C.C. Moore: {\it Ergodic equivalence relations,
cohomology, and von Neumann algebras I, II}, Trans. Amer. Math.
Soc. {\bf 234} (1977), 289-324, 325-359.

\item{[Fu]} A. Furman: {\it Gromov's measure equivalence and rigidity of
higher rank lattices}, Ann. of Math. {\bf 150} (1999), 1059-1081.

\item{[G1]} D. Gaboriau: {\it Cout des r\'elations d'\'equivalence
et des groupes}, Invent. Math. {\bf 139} (2000), 41-98.

\item{[G2]} D. Gaboriau: {\it Invariants $\ell^2$ de relations
d'\'equivalence et de groupes}, preprint May 2001.

\item{[G3]} D. Gaboriau: {\it Approximate dimension for
equivalence relations and groups}, preprint.

\item{[Ge]} L. Ge: {\it Prime factors},
Proc. Natl. Acad. Sci. USA, {\bf
93} (1996), 12762-12763.

\item{[GoNe]} V. Y. Golodets, N. I. Nesonov: T{\it -property and
nonisomorphic factors of type} II {\it and} III, J. Funct. Analysis
{\bf 70} (1987), 80-89.

\item{[Gr]} M. Gromov: {\it Rigid transformation groups},
``G\'eometrie diff\'erentielle'' (Paris 1986), 65-139, Travaux en Cours, 33,
Hermann, Paris, 1988.

\item{[H]} U. Haagerup: {\it An example of non-nuclear C$^*$-algebra
which has the metric approximation property}, Invent. Math.
{\bf 50} (1979), 279-293.

\item{[dHVa]} P. de la Harpe, A. Valette: ``La propri\'et\'e T
de Kazhdan pour les
groupes localement compacts'', Ast\'erisque {\bf 175} (1989).

\item{[Hj]} G. Hjorth: {\it A lemma for cost attained},
UCLA preprint 2002.

\item{[HkS]} R. Hoegh-Krohn, T. Skjelbred: {\it Classification of
$C^*$-algebras admitting ergodic actions of the two-dimensional torus},
J. reine engev. Math., {\bf 328} (1981), 1-8.

\item{[ILP]} M. Izumi, R. Longo and S. Popa: {\it A Galois correspondence
for compact groups of automorphisms of von Neumann algebras with
a generalization to Kac algebras}, J. Funct. Analysis, 155 (1998), 25-63.

\item{[Jo1]} P. Jolissait: {\it Haagerup approximation property
for von Neumann algebras},
preprint 2001.

\item{[Jo2]} P. Jolissaint: {\it On the relative property} (T),
preprint 2001.

\item{[J1]} V.F.R. Jones : {\it Index for subfactors}, Invent. Math.
{\bf 72} (1983), 1-25.

\item{[J2]} V.F.R. Jones : {\it A converse to Ocneanu's theorem},
J. Operator Theory {\bf 4} (1982), 21-23.

\item{[JPo]} V.F.R. Jones, S. Popa: {\it Some properties of
MASAs in factors}, in ``Invariant subspaces and other topics'', pp. 89-102, Operator Theory: Adv. Appl. {\bf 6}, Birkhäuser, 1982.

\item{[K1]} R.V. Kadison: {\it Problems on von Neumann algebras},
Baton Rouge Conference 1967.

\item{[K2]} R.V. Kadison: {\it Diagonalizing matrices}, Amer. Math.
Journal (1984), 1451-1468.

\item{[KaWe]} V. Kaftal, G. Weiss: {\it Compact
derivations relative to semifinite von Neumann algebras},
J. Funct. Anal. {\bf 62} (1985), 202-220.

\item{[Kaz]} D. Kazhdan: {\it Connection of the dual space of a group
with the structure of its closed subgroups}, Funct. Anal. and its Appl.
{\bf1} (1967), 63-65.

\item{[KoY]} H. Kosaki, S. Yamagami: {\it Irreducible bimodules associated
with crossed product algebras}, preprint.

\item{[L\"u]} W. L\"uck: {\it Dimension theory of arbitrary
modules over finite von Neumann algebras and
$L^2$-Betti numbers} II. {\it Applications to Grothendieck groups,
$L^2$-Euler characteristics and Burnside groups},
J. Reine Angew. Math., {\bf 496} (1998), 213-236.

\item{[Ma]} G. Margulis: {\it Finitely-additive invariant measures
on Euclidian spaces}, Ergodic. Th. and Dynam. Sys. {\bf 2} (1982),
383-396.

\item{[McD]} D. McDuff: {\it Central sequences
and the hyperfinite factor}, Proc. London Math. Soc. {\bf 21} (1970),
443-461.

\item{[MvN]} F. Murray, J. von Neumann: {\it Rings of operators IV},
Ann. Math. {\bf44} (1943), 716-808.

\item{[PiPo]} M. Pimsner, S. Popa: {\it Entropy and index for subfactors},
Ann. Scient.
Ec. Norm. Sup. {\bf 19} (1986), 57-106.

\item{[Po1]} S. Popa: {\it Correspondences}, INCREST preprint 1986,
unpublished.

\item{[Po2]} S. Popa: ``Classification of subfactors and of their
endomorphisms'', CBMS Lecture Notes Series, Vol. {\bf 86}, 1994.

\item{[Po3]} S. Popa: {\it Classification
of subfactors of type} II, Acta Math.
{\bf 172} (1994), 163-255.

\item{[Po4]} S. Popa: {\it Free independent
sequences in type II$_1$ factors and
related problems}, Asterisque {\bf 232} (1995), 187-202.

\item{[Po5]} S. Popa: {\it Some properties of the
symmetric enveloping algebras with applications to amenability and
property} T, Documenta Math. {\bf 4} (1999), 665-744.

\item{[Po6]} S. Popa: {\it On a problem of
R. V. Kadison on maximal abelian subalgebras},
Invent. Math. {\bf 65} (1981), 269-281.

\item{[Po7]} S. Popa: {\it The relative Dixmier property for
inclusions of von Neumann algebras of finite index}, Ann. Ec. Norm. Sup.,
{\bf 32} (1999), 743-767.

\item{[Po8]} S. Popa: {\it Notes on Cartan subalgebras
in type} II$_1$ {\it factors}, Math. Scand. {\bf 57} (1985), 171-188.

\item{[Po9]} S. Popa: {\it On the distance between MASA's
in type} II$_1$ {\it factors}, in ``Mathematical Physics in
Mathematics and Physics'' pp. 321-324, Fields Inst. Comm. {\bf 30}.

\item{[PoRa]} S. Popa, F. Radulescu: {\it Derivations of von Neumann
algebras into the compact ideal space of a semifinite algebra},
Duke Math. J. {\bf 57} (1988), 485-518.

\item{[PoSh]} S. Popa, D. Shlyakhtenko: {\it Cartan subalgebras and bimodule
decomposition of} II$_1$ {\it factors}, Math. Scand. (2003).

\item{[Ra]} F. Radulescu: {\it Random matrices,
amalgamated free products and
subfactors of the von Neumann algebra of a free group},
Invent. Math. {\bf 115} (1994), 347--389.

\item{[Ri]} M. Rieffel: {\it $C^*$-algebras associated
with irrational rotations}, Pacific. J. Math. {\bf 93} (1981), 415-429.

\item{[S]} Sakai: ``C$^*$-algebras and W$^*$-algebras'', Springer-Verlag,
Berlin-Heidelberg-New York, 1971.

\item{[Sa]} Sauvageot: {\it Sur le produit tensoriel relatif
d'espaces de Hilbert}, J. Oprator Theory {\bf 9} (1983), 237-252.

\item{[Sc]} K. Schmidt: {\it Asymptotically invariant
sequences and an action of $SL(2, \Bbb Z)$ on the
$2$-sphere}, Israel. J. Math. {\bf 37} (1980), 193-208.

\item{[Sha]} Y. Shalom: {\it Bounded generation
and Kazhdan's property} (T), Publ. Math. I.H.E.S. (2001).

\item{[Sh]} D. Shlyakhtenko: {\it Free quasi-free states},
Pacific J. Math. {\bf 177} (1997),
329-368.

\item{[Si]} I.M. Singer: {\it
Automorphisms of finite factors}, Amer. J. Math. {\bf 77}, (1955). 117-133.

\item{[Va]} A. Valette: {\it Semi-direct products with
the property (T)}, preprint, 2001.

\item{[V1]} D. Voiculescu: {\it
Circular and semicircular systems and free product
factors}, Prog. in Math. {\bf 92}, Birkhauser, Boston,
1990, pp. 45-60.

\item{[V2]} D. Voiculescu: {\it The analogues of entropy and
of Fisher's information theory in free probability} II,
Invent. Math. {\bf 118} (1994), 411-440 .

\item{[W]} A. Wassermann: {\it Coactions and Yang-Baxter equations
for ergodic actions and subfactors}, 
London Math. Soc. Lect. Notes Ser.,
{\bf 136}, Cambridge Univ. Press 1988, pp. 202-236.

\item{[Zi]} R. Zimmer: ``Ergodic theory and semisimple groups'',
Birkha\"user-Verlag, Boston 1984.

\enddocument